\documentclass[11pt]{article}

\usepackage[english]{babel}

\usepackage[T1]{fontenc}
\usepackage{lmodern}

\let\Sho\TH
\usepackage[utf8]{inputenc}
\usepackage[a4paper,margin=2.5cm]{geometry}
\usepackage[dvipsnames]{xcolor}
\usepackage{ltablex}  
\keepXColumns
\usepackage{eufrak}
\usepackage{chemist}
\usepackage{amsmath,amsfonts,amsthm,amssymb}

\DeclareMathOperator{\sgn}{sgn}

\usepackage{blindtext}
\usepackage{mathtools}
\usepackage{mathrsfs}
\newcommand{\defeq}{\mathrel{\mathop:}=}

\usepackage{booktabs,tabularx}
\usepackage{apxproof}
\usepackage{environ}
\usepackage{dsfont}
\usepackage{etoolbox}
\usepackage{fancyvrb}
\usepackage{ifthen}
\usepackage{kvoptions}
\usepackage{bibunits}
\usepackage{comment}
\usepackage{chemformula}
\usepackage{mathrsfs}
\usepackage{float}
\usepackage{xcolor}
\usepackage{enumerate}
\usepackage[toc,page]{appendix}
\usepackage{graphicx,color}
\usepackage{hyperref}
\usepackage{mathtools}
\usepackage{setspace}
\usepackage{cite}

\makeatletter
\renewcommand*\env@matrix[1][\arraystretch]{%
  \edef\arraystretch{#1}%
  \hskip -\arraycolsep
  \let\@ifnextchar\new@ifnextchar
  \array{*\c@MaxMatrixCols c}}
\makeatother

\newtheorem{theorem}{Theorem}[section]
\newtheorem{lemma}[theorem]{Lemma}
\newtheorem{corollary}[theorem]{Corollary}
\newtheorem{assumption}{Assumption}
\newtheorem{proposition}[theorem]{Proposition}

\numberwithin{equation}{section}

\newtheorem*{theorem*}{Theorem}
\newtheorem*{proposition*}{Proposition}
\newtheorem*{corollary*}{Corollary}
\newtheorem*{lemma*}{Lemma}

\let\OLDthebibliography\thebibliography
\renewcommand\thebibliography[1]{
  \OLDthebibliography{#1}
  \setlength{\parskip}{0pt}
  \setlength{\itemsep}{0pt plus 0.3ex}
}

\theoremstyle{definition}
\newtheorem{definition}[theorem]{Definition}
\newtheorem{remark}[theorem]{Remark}

\newcommand{\ignore}[1]{}

\usepackage{amsmath}
\usepackage{graphicx}
\title{Existence, uniqueness and moment bounds \\
for a spatial model of Muller's ratchet}

\author{João Luiz de Oliveira Madeira\thanks{Department of Statistics, University of Oxford, UK} \and Marcel Ortgiese\footnote{Department of Mathematical Sciences, University of Bath, UK} \and Sarah Penington\footnotemark[2]}

\begin{document}


\maketitle

\begin{abstract}
In this article, we consider a generalisation of the spatial Muller's ratchet introduced by Foutel-Rodier and Etheridge~\cite{foutel2020spatial}. This particle system is a spatial model of an asexual population, with birth and death rates that depend on the local population density. Particles live in discrete demes and migrate to neighbouring demes. 
Each particle carries some number of mutations (its `type'), and additional mutations can occur during birth events.
Mutations are assumed to be deleterious, i.e.~carrying a higher number of mutations results in a lower birth rate.
Our main result shows that this interacting particle system can be constructed even when the total initial number of particles is infinite. We also prove moment bounds on the local density of particles;
these bounds are a crucial ingredient of the proof of a law of large numbers result for the particle system in the companion article~\cite{madeiraPreprintmuller}.
The construction of the particle system uses a sequence of approximating processes. 
Proving weak convergence of this sequence of processes is non-trivial because the particle system is non-monotone and interactions are non-local in type space.
The uniqueness of the limit relies on a delicate coupling argument.
\end{abstract}


\section{Introduction}
\label{Paper01_introduction}

In this article, we rigorously construct a reaction-diffusion system of interacting particles with infinitely many types which is a generalisation of a model introduced by Foutel-Rodier and Etheridge in~\cite{foutel2020spatial}, called the \emph{spatial Muller's ratchet}. Muller's ratchet is a mechanism first proposed in~\cite{muller1964relation} to explain the evolution of sexual reproduction. Consider a population of constant size in which, at each reproduction event, a mutation -- i.e.~a change in the genome transmitted from parent to offspring -- can occur with positive probability. Most such mutations are deleterious, meaning that they typically reduce fitness. Here, fitness refers to an individual’s ability to survive and reproduce in a given environment~\cite{loewe2010population}. In asexual reproduction, chromosomes are inherited as indivisible units, so the number of deleterious mutations along an ancestral lineage can only increase. This leads to a gradual decline in overall population fitness~\cite{etheridge2009often}. When all individuals in the most-adapted class acquire at least one additional deleterious mutation, the minimum mutational load in the population increases -- this is known as a click of the ratchet. In contrast, sexual reproduction enables recombination, which can reduce mutation load along lineages, thereby preventing further clicks. From a mathematical perspective, models of Muller's ratchet have been extensively studied in populations without spatial structure~\cite{casanova2022quasi,pfaffelhuber2012muller,haigh1978accumulation,etheridge2009often,mariani2020metastability}.

In spatially structured populations, however, the probability of survival of an individual depends not only on its fitness, but also on the local population density around this individual~\cite{foutel2020spatial}. In this setting, the question of whether or not deleterious mutations could accumulate is complex (related to a phenomenon known as gene surfing), and biological data is inconclusive~\cite{peischl2013accumulation,peischl2015expansion,peischl2016genetic,do2015no,simons2014deleterious}. To investigate this question, Foutel-Rodier and Etheridge proposed a spatial version of Muller's ratchet in~\cite{foutel2020spatial}. In this model, the population is subdivided into demes indexed by~$\mathbb{Z}$. For $x \in \mathbb{Z}$, $t \geq 0$ and $k \in \mathbb{N}_{0}$, let $\eta_{k}(t,x)$ indicate the number of particles carrying exactly~$k$ mutations and living in deme $x$ at time~$t$, and let $\| \eta(t,x) \|_{\ell_{1}} \defeq \sum_{k = 0}^{\infty} \eta_{k}(t,x)$. In the model proposed in~\cite{foutel2020spatial}, for some fixed $m>0$, each particle (representing an individual in the population) independently migrates at rate~$m$, moving to one of the two neighbouring demes with equal probability. Moreover, a particle located in deme $x \in \mathbb{Z}$ at time $t \geq 0$ and carrying $k \in \mathbb{N}_{0}$ mutations gives birth to a new particle at rate
\begin{equation*}
    (1-s)^{k}\Big(B\vert \vert \eta(t-,x) \vert \vert_{\ell_{1}} + 1\Big),
\end{equation*}
and dies at rate
\begin{equation*}
   \vert \vert \eta(t-,x) \vert \vert_{\ell_{1}}\Big(B\vert \vert \eta(t-,x) \vert \vert_{\ell_{1}} + 1\Big).
\end{equation*}
Here, the parameter $s \in (0,1)$ represents the impact on fitness of each deleterious mutation, and the parameter $B \in [0, \infty)$ captures the balance between cooperation and competition in the local population dynamics: larger values of $B$ indicate a greater role of cooperation. Take $\mu \in (0,1)$. After a birth event, an offspring particle is added to the same deme as its parent; with probability $1 - \mu$, the offspring inherits the same number of mutations as its parent, and with probability~$\mu$, it accumulates an additional mutation. In~\cite{foutel2020spatial}, Foutel-Rodier and Etheridge used this model to make theoretical predictions about how the spatial expansion of the population impacts the propagation of deleterious mutations. To make these predictions, they (non-rigorously) derive scaling limits and study the system's long-term behaviour.

Our first goal in this paper is to rigorously construct (a generalisation) of the spatial Muller's ratchet, in particular for starting configurations consisting of infinitely many particles.
For an initial configuration consisting of finitely many particles, it suffices to show that the total number of particles remains almost surely finite over fixed time intervals in order to establish that the model introduced in~\cite{foutel2020spatial} can be formally constructed as a finite system of stochastic differential equations driven by Poisson processes~\cite{garcia2006spatial}. Nevertheless, for modelling purposes, it is often more natural to consider infinite particle systems when studying large real-world populations -- a classical approach in the theory of interacting particle systems~\cite{durrett1999stochastic, durrett1994importance, liggett1972existence}. However, constructing the spatial Muller's ratchet as an infinite particle system poses significant challenges. In particular, there are no \emph{a priori} bounds on the number of particles per site, and the birth and death rates are unbounded and depend on the local population density. Moreover, the system is non-monotone: since the death rate at a deme depends on the total population size in the deme,
fitter particles (particles carrying few mutations) can be eliminated through local competition with less fit particles, i.e.~the non-monotonicity arises from non-local interactions in the type space. We will explain this issue in more detail in Section~\ref{Paper01_ssec:non-monotonicity} below. When interpreting numerical simulations -- such as those in~\cite{foutel2020spatial} -- it is generally assumed that the interacting particle system of interest, starting with an infinite number of particles, can be well approximated by a system with a finite initial number of particles. While this assumption can often be justified for monotone interacting particle systems using standard techniques, its validity is not immediately clear in the case of spatial birth-death processes with non-monotone, non-local interactions.

To overcome these challenges, we derive a general recipe for constructing semigroups of Markov processes taking values in Polish spaces which are not necessarily locally compact as the limit of an approximating sequence of semigroups. Unlike the standard theory of Feller processes, our approach extends to complete, separable metric spaces that are not locally compact. This allows us to construct spatial interacting particle systems without imposing any \emph{a priori} bounds on particle density. We will then apply this recipe with an approximating sequence of particle systems in which particles outside a large spatial box are frozen and particles carrying too many mutations cannot give birth.
A crucial ingredient will be a moment bound 
that gives us local control over particle numbers in the approximating particle systems.
To show uniqueness, we use
a careful coupling argument that encodes two realisations, starting from different initial conditions, of one of the approximating particle systems in terms of susceptible, infected, and partially recovered particles. 
The moment bounds on local particle densities for the approximating particle systems will then imply moment bounds on the limiting process, which is our second main result.
To the best of our knowledge, there are no prior results constructing infinite interacting particle systems with non-monotone, non-local interactions, unbounded numbers of particles per site, and unbounded transition rates. We hope that the ideas developed here will be applicable to other spatial birth-death processes with non-monotone interactions, which often exhibit similar challenges to those encountered in the spatial Muller's ratchet~\cite{penington2018spreading,etheridge2019genealogical}.

In a companion paper~\cite{madeiraPreprintmuller}, 
we derive a functional law of large numbers for the spatial Muller's ratchet, confirming the non-rigorous calculations in~\cite{foutel2020spatial}. Both of our main results here, existence and the moment bounds, will be essential in our proofs in the companion paper. 
There, we also discuss the implications for various biologically relevant examples. We also analyse the spreading speed of a population invading an empty habitat and use tracer dynamics to investigate the question of gene surfing.

\medskip

\noindent \textbf{Structure of the article.} In Section~\ref{Paper01_model_description}, we describe the model and state our main results. We discuss the lack of monotonicity in Section~\ref{Paper01_ssec:non-monotonicity}. Section~\ref{Paper01_related_work} briefly reviews some techniques commonly used in the construction of interacting particle systems and explains why they are not directly applicable in our setting. In Section~\ref{Paper01_subsection_heuristics}, we informally present the key ideas behind the proofs of our results.
In Section~\ref{Paper01_subsection_notation}, we give a glossary of frequently used notation.

The proofs of our results are contained in Sections~\ref{Paper01_section_functional_analysis_feller_non_locally_compact}-\ref{Paper01_section_spread_infection}.
In Section~\ref{Paper01_section_functional_analysis_feller_non_locally_compact}, we derive, in Theorem~\ref{Paper01_general_thm_convergence_feller_non_locally_compact}, a general recipe for constructing weakly Feller càdlàg Markov processes in non-locally compact state spaces, without assuming that the infinitesimal generator or its associated semigroup satisfies a global Lipschitz condition. In the remaining sections, we apply this recipe to the spatial Muller's ratchet. In Section~\ref{Paper01_section_correlation_functions}, we derive moment estimates for a sequence of stochastic processes that approximate the spatial Muller's ratchet. These estimates are used in Section~\ref{Paper01_formal_construction_generator_section} to prove that this sequence of processes is tight and converges weakly to a unique in law Markov process. We will show in~Section~\ref{Paper01_formal_construction_generator_section} that this Markov process satisfies a certain martingale problem. A key step in proving uniqueness of the weak limit is a coupling argument, comparing approximating processes started from different configurations using an encoding in terms of susceptible, infected, and partially recovered particles. To our knowledge, this coupling argument -- presented in detail in Section~\ref{Paper01_section_spread_infection} -- is new and potentially applicable to other non-monotone interacting particle systems.
Some technical details and proofs are deferred to the appendix.

\medskip

\noindent \textbf{Notation.}
In the remainder of the article, we will use the following general notation. For a complete and separable metric space $(\mathcal{S},d_{\mathcal S})$, we denote the space of $\mathcal{S}$-valued càdlàg processes by $\mathscr{D}\left([0, \infty), \mathcal{S}\right)$. 
Recall that since $(\mathcal{S}, d_{\mathcal{S}})$ is a Polish space, the set $\mathscr{D}([0, \infty), \mathcal{S})$ equipped with the Skorokhod $J_{1}$-topology is also a Polish space (see~\cite[Theorem~3.5.6]{ethier2009markov}).
We write $\mathscr{C}(\mathcal{S}, \mathbb{R})$ for the space of continuous real-valued functions defined on $\mathcal{S}$, and we write $\mathscr{C}_{b}(\mathcal{S}, \mathbb{R})$ for the space of bounded continuous real-valued functions. We write $\mathscr{B}_{b}(\mathcal{S}, \mathbb{R})$ for the set of real-valued bounded Borel measurable functions on $\mathcal{S}$. Then $\mathscr{B}_{b}(\mathcal{S}, \mathbb{R})$, equipped with the norm $\vert \vert \cdot \vert \vert_{L_{\infty}(\mathcal{S}; \mathbb{R})}$, where $ \vert \vert \phi \vert \vert_{L_{\infty}(\mathcal{S}; \mathbb{R})} = \sup_{\boldsymbol{\xi} \in \mathcal{S}}\vert  \phi(\boldsymbol{\xi}) \vert$, is a Banach space. We also let $\ell_{1}$ denote the space of summable real sequences, i.e. 
\begin{equation*}
    \ell_{1} = \ell_{1}\left(\mathbb{N}_{0}\right) \defeq \left\{z = \left(z_{k}\right)_{k \in \mathbb{N}_{0}} \in \mathbb{R}^{\mathbb{N}_{0}}: \, \vert \vert z \vert \vert_{\ell_1} \defeq \sum_{k = 0}^{\infty}\vert z_{k} \vert < \infty \right\},
\end{equation*}
and we let $\ell_1^+ \subset \ell_1$ denote the subset of summable sequences with non-negative entries. Recall that $(\ell_1,\vert \vert \cdot \vert \vert_{\ell_1})$ is a complete and separable Banach space. For a stochastic process $(\eta(t))_{t \geq 0}$ taking values in $\mathcal{S}$, for $t \geq 0$, $\boldsymbol{\eta} \in \mathcal{S}$ and a Borel subset $A \subseteq \mathcal{S}$, we write
\begin{equation*}
    \mathbb{P}_{\boldsymbol{\eta}}(\eta(t) \in A) \defeq \mathbb{P}(\eta(t) \in A \vert \, \eta(0) = \boldsymbol{\eta}),
\end{equation*}
whenever the right-hand side is well defined. We let $\mathbb{E}_{\boldsymbol{\eta}}$ be the expectation operator corresponding to the measure $\mathbb{P}_{\boldsymbol{\eta}}$. For a sequence of $\mathcal{S}$-valued random variables $({\eta}^{n})_{n \in \mathbb{N}}$ and for a random variable $\eta \in \mathcal{S}$, we use the notation $\eta^n \xrightarrow[]{\mathcal{D}} \eta$ as $n\to \infty$ to indicate that $({\eta}^{n})_{n \in \mathbb{N}}$ converges weakly on $(\mathcal{S},d_{\mathcal S})$ to $\eta$ as $n \rightarrow \infty$. 

We  let $\sgn: \mathbb{R} \rightarrow \{-1,0,1\}$ denote the sign function. For a finite set $A$, we denote the number of elements of $A$ by~$\# A$. For $K \in \mathbb{N}_0$, we let $[K]_0 \defeq \{0,1, \ldots,K\}$ and $[K] \defeq [K]_0 \setminus \{0\}$.
For $x\in \mathbb{R}$, we let $x^+ \defeq x \vee 0$.

It will also be convenient to introduce some notation regarding estimates for expressions. For a set $\mathcal{S}$, let $f,g: \mathcal{S} \rightarrow [0,\infty)$ be non-negative real-valued functions on $\mathcal{S}$. We write $f \lesssim g$ to indicate that there exists a constant $C > 0$ not depending on any parameters such that for every $x \in \mathcal{S}$ we have $f(x) \leq Cg(x)$. 
We write $f \lesssim_{p_1,\ldots , p_k} g$ to indicate that there exists $C>0$ that depends only on the parameters $p_1,\ldots , p_k$ such that for every $x \in \mathcal{S}$ we have $f(x) \leq Cg(x)$. 

\section{Model description and main results} \label{Paper01_model_description}

We begin by giving an informal description of the (generalised) spatial Muller's ratchet.
Our model consists of particles living on the rescaled one-dimensional lattice $L^{-1}\mathbb{Z}$, where $L > 0$ is a space renormalisation parameter, and we refer to each point $x \in L^{-1}\mathbb{Z}$ as a deme. Each particle (representing an individual in the population) carries a unique chromosome and is characterised by two features: its number of deleterious mutations and its spatial location. 
For each $t \ge 0$, $x \in L^{-1}\mathbb{Z}$ and $k\in \mathbb N_0$, let $\eta_{k}(t,x)$ denote the number of particles carrying exactly $k$ mutations and living in deme $x$ at time $t$, so that 
$\eta(t,x) := \left(\eta_{k}(t,x)\right)_{k \in \mathbb{N}_{0}}$ characterises the particles in this deme. Hence, the total number of particles living at time $t$ in deme $x$ is $\vert \vert \eta(t,x) \vert \vert_{\ell_1}
= \sum_{k = 0}^{\infty} \eta_{k}(t,x)$.

Let $m > 0$ be the migration rate, and for each $k\in \mathbb N_0$ we let $s_{k} \geq 0$ denote the fitness parameter of a particle carrying $k$ mutations. Let $q_{+},q_{-}: [0,\infty) \rightarrow [0,\infty)$ be non-negative functions, $N \in \mathbb{N}$ a scaling parameter for the carrying capacity of the environment, and $\mu \in [0,1]$ the mutation probability. We aim to construct a Markov process $(\eta(t))_{t \geq 0} = (\eta_{k}(t,x): \, k \in \mathbb{N}_{0}, \, x \in L^{-1}\mathbb{Z})_{t \geq 0}$ evolving (informally) as follows:
\begin{itemize}
    \item \textit{Migration events}: For each $t \geq 0$ and $x \in L^{-1}\mathbb{Z}$, each particle living in deme $x$ at time $t$ makes a jump independently at rate $m$, to a uniformly chosen deme from $\{x - L^{-1}, \, x + L^{-1}\}$.
    \item \textit{Reproduction events}: For each $t \geq 0$, $x \in L^{-1}\mathbb{Z}$ and $k \in \mathbb{N}_{0}$, each particle carrying $k$ mutations living in deme $x$ at time $t$ reproduces independently at rate $\displaystyle s_{k}q_{+}\left(\vert \vert \eta(t-,x) \vert \vert_{\ell_{1}} / N\right)$. After such an event, a new offspring particle will be added in deme $x$. With probability $1 - \mu$, the offspring  will carry $k$ mutations, and with probability $\mu$ it will carry $k + 1$ mutations. Note that in both cases, the parent particle {remains}  alive and keeps the same number of mutations.
    \item \textit{Death events}: For each $t \geq 0$ and $x \in L^{-1}\mathbb{Z}$, each particle living in deme $x$ at time $t$ dies independently at rate $\displaystyle q_{-}\left(\vert \vert \eta(t-,x) \vert \vert_{\ell_{1}} / N\right)$, regardless of its number of mutations. When a particle dies, it is removed from the process.
\end{itemize}

We now state some conditions on
the parameters of our model that we will assume throughout. 
Our first assumption is that mutations are deleterious, i.e.~particles that carry more mutations have lower fitness.

\begin{assumption} [Fitness parameters]\label{Paper01_assumption_fitness_sequence}
The sequence of fitness parameters $\left(s_{k}\right)_{k \in \mathbb{N}_{0}}$ satisfies the following conditions:
\begin{enumerate}[(i)]
\item $s_{0} = 1$;
\item $s_{k} \geq 0$ for all $k \in \mathbb{N}_{0}$;
\item $\left(s_{k}\right)_{k \in \mathbb{N}_{0}}$ is monotonically non-increasing, i.e.~$s_{k} \geq s_{k+1}$ for all $k \in \mathbb{N}_{0}$;
\item $\displaystyle \lim_{k \rightarrow \infty} s_{k} = 0$.
\end{enumerate}
\end{assumption}

Since in our model we do not impose \emph{a priori} bounds on the number of particles per deme, we need to make assumptions on the birth and death rates which will guarantee that the population 
remains locally controlled.

\begin{assumption} [Birth and death polynomial rates]\label{Paper01_assumption_polynomials}
The functions $q_{+}, q_{-}: [0,\infty) \rightarrow [0,\infty)$ are polynomials satisfying $0 \leq \deg q_{+} < \deg q_{-}$.
\end{assumption}

Note also that the parameter $N$, used together with $q_+$ and $q_-$ to determine the birth and death rates, can be thought of as being proportional to the local carrying capacity of the population, as observed in~\cite{foutel2020spatial}.

We think of $\xi=(\xi_j(x))_{j\in \mathbb N_0,\, x\in L^{-1}\mathbb Z}\in (\mathbb N_0^{\mathbb N_0})^{L^{-1}\mathbb Z}$ as representing a particle configuration in which $\xi_j(x)$ is the number of particles in deme $x$ carrying exactly $j$ mutations, for each $x\in L^{-1}\mathbb Z$ and $j\in \mathbb N_0$.
In order to define an appropriate state space for the interacting particle system, 
set
\begin{equation} \label{Paper01_definition_state_space_formal}
    \mathcal{S} \defeq \Bigg\{\boldsymbol{\xi} \in (\mathbb{N}_{0}^{\mathbb{N}_{0}})^{L^{-1}\mathbb{Z}}: \; \vert \vert \vert \boldsymbol{\xi} \vert \vert \vert_{\mathcal{S}} \defeq \sum_{x \in L^{-1}\mathbb{Z}} \frac{\vert \vert \xi(x) \vert \vert_{\ell_{1}}}{(1 + \vert x \vert)^{2}} < \infty\Bigg\}.
\end{equation}
Then $\vert \vert \vert \cdot \vert \vert \vert_{\mathcal{S}}$ induces a metric $d_{\mathcal{S}}: \mathcal{S} \times \mathcal{S} \rightarrow [0, \infty)$ on $\mathcal{S}$ given by
\begin{equation} \label{Paper01_definition_semi_metric_state_space}
\begin{aligned}
   d_{\mathcal{S}}: \mathcal{S} \times \mathcal{S} & \rightarrow [0, \infty) \\
    (\boldsymbol{\xi}, \boldsymbol{\zeta}) & \mapsto  d_{\mathcal{S}}(\boldsymbol{\xi}, \boldsymbol{\zeta}) \defeq \vert \vert \vert \boldsymbol{\xi} - \boldsymbol{\zeta} \vert \vert \vert_{\mathcal{S}} = \sum_{x \in L^{-1}\mathbb{Z}} \frac{\vert \vert (\xi - \zeta)(x) \vert \vert_{\ell_{1}}}{(1 + \vert x \vert)^{2}},
\end{aligned}
\end{equation}
where for each $x \in L^{-1}\mathbb{Z}$, the difference $(\xi - \zeta)(x)$ is defined as the usual difference between elements of~$\ell_{1}$. In Section~\ref{Paper01_subsection_topological_properties_state_space} in the appendix, we will verify that the metric space $(\mathcal{S},d_{\mathcal{S}})$ is complete and separable. 

\begin{remark}
   Our results would extend to other choices of state space~$\mathcal{S}$. In particular, the power $2$ in the normalisation in~\eqref{Paper01_definition_state_space_formal} and~\eqref{Paper01_definition_semi_metric_state_space} could be replaced by any power $p > 1$. In fact, as long as $\mathcal S$ was chosen in such a way that the number of particles per deme could grow at most subexponentially with distance from the origin for configurations in $\mathcal S$, our methods would apply without modification. However, from a modelling perspective, the state space~$\mathcal{S}$ as defined in~\eqref{Paper01_definition_state_space_formal} is already sufficiently general and includes all biologically relevant configurations.
\end{remark}

For technical reasons that we will discuss further in Section~\ref{Paper01_subsection_heuristics} below, we can only construct a càdlàg version of our process for initial conditions in a particular subspace of $\mathcal{S}$.

\begin{assumption} [Initial condition] \label{Paper01_assumption_initial_condition}
The set of initial configurations will be given by the following set $\mathcal S_0\subset \mathcal S$. Let
\begin{equation} \label{Paper01_set_initial_configurations}
    \mathcal{S}_{0} \defeq \Bigg\{{\eta}= (\eta_j(x))_{j\in \mathbb N_0,\, x\in L^{-1}\mathbb Z}\in \mathcal{S}: \; \sup_{x \in L^{-1}\mathbb{Z}} \vert \vert \eta(x) \vert \vert_{\ell_{1}} < \infty \quad \textrm{and} \quad \lim_{k \rightarrow \infty} \, \sup_{x \in L^{-1}\mathbb{Z}} \, \sum_{j = k}^{\infty} \, \eta_{j}(x) = 0\Bigg\}.
\end{equation}
\end{assumption}

We aim to construct the particle system described informally at the start of this section as a càdlàg $\mathcal S$-valued Markov process associated to an infinitesimal generator~$\mathcal{L}$ {that we will now define}.
For~$x \in L^{-1}\mathbb{Z}$ and~$k \in \mathbb{N}_{0}$, let $\boldsymbol{e}_{k}^{(x)} \in \mathcal{S}$ denote the configuration consisting of one particle carrying exactly $k$ mutations in deme $x$.  Moreover, for $N \in \mathbb{N}$, we let $q_{+}^{N}, q_{-}^{N}: [0,\infty) \rightarrow [0,\infty)$ be given by, for all $u \in [0,\infty)$,
\begin{equation} \label{Paper01_scaled_polynomials_carrying_capacity}
    q_{+}^{N}(u) \defeq q_{+}\left(\frac{u}{N}\right) \quad \textrm{and} \quad q_{-}^{N}(u) \defeq q_{-}\left(\frac{u}{N}\right).
\end{equation}
Also, let the sequences of maps $(F^{b}_{k})_{k \in \mathbb{N}_{0}}$ and $(F^{d}_{k})_{k \in \mathbb{N}_{0}}$ be given by, for $k \in \mathbb{N}_0$ and $u =(u_j)_{j\in \mathbb N_0}\in \ell_1^+$,
    \begin{equation} \label{Paper01_discrete_birth_rate}
        F^{b}_{k}(u) \defeq \left(s_{k}(1-\mu)u_{k} + \mathds{1}_{\{k \geq 1\}} s_{k-1}\mu u_{k-1}\right)q^{N}_{+}\left({\vert \vert u \vert \vert_{\ell_{1}}}\right),
    \end{equation}
and
\begin{equation} \label{Paper01_discrete_death_rate}
        F^{d}_{k}(u) \defeq u_{k}q^{N}_{-}\left({\vert \vert u \vert \vert_{\ell_{1}}}\right).
    \end{equation}
Then, for a suitable subset $\mathscr{C}_*({\mathcal{S}, \mathbb{R}}) \subset \mathscr{C}(\mathcal{S}, \mathbb{R})$ that we will define in~\eqref{Paper01_definition_Lipschitz_function}, we define the operator $\mathcal L:\mathscr{C}_*({\mathcal{S}, \mathbb{R}}) \to \mathscr{C}(\mathcal{S}, \mathbb{R})$ as follows.
For any $\phi \in \mathscr{C}_*({\mathcal{S}, \mathbb{R}})$ and any $\boldsymbol{\xi} =(\xi_k(x))_{k\in \mathbb N_0,\, x\in L^{-1}\mathbb Z}\in \mathcal{S}$, we let
\begin{equation} \label{Paper01_generator_foutel_etheridge_model}
    \mathcal{L}\phi(\boldsymbol{\xi}) = \frac{m}{2}\mathcal{L}_{m}\phi(\boldsymbol{\xi}) + \mathcal{L}_{r}\phi(\boldsymbol{\xi}),
\end{equation}
where $\mathcal{L}_{m}$ can be thought of as corresponding to the migration process and $\mathcal{L}_{r}$ as corresponding to the birth-death process, and $\mathcal{L}_{m}$ and $\mathcal{L}_{r}$ are defined as follows:
\begin{equation} \label{Paper01_infinitesimal_generator}
\begin{aligned}
    (\mathcal{L}_{m}\phi)(\boldsymbol{\xi}) & \defeq \sum_{x \in L^{-1}\mathbb{Z}} \; \sum_{k = 0}^{\infty} \; \sum_{z \in \{-L^{-1}, L^{-1}\}}\xi_{k}(x) \left(\phi\Big(\boldsymbol{\xi} + \boldsymbol{e}_{k}^{(x+z)} - \boldsymbol{e}_{k}^{(x)}\Big) - \phi(\boldsymbol{\xi}) \right), \\
    (\mathcal{L}_{r}\phi)(\boldsymbol{\xi}) & \defeq \sum_{x \in L^{-1}\mathbb{Z}} \; \sum_{k = 0}^{\infty} \Big( F^{b}_k(\xi(x))\left(\phi\Big(\boldsymbol{\xi} + \boldsymbol{e}_{k}^{(x)}\Big) - \phi(\boldsymbol{\xi})\right) + F^{d}_k(\xi(x))\left(\phi\Big(\boldsymbol{\xi}- \boldsymbol{e}_{k}^{(x)}\Big) - \phi(\boldsymbol{\xi})\right)\Big).
\end{aligned}    
\end{equation}
We will prove in Lemma~\ref{Paper01_convergence_infinitesimal_generator_lipischitz_semi_metric} that the right-hand sides of the identities in~\eqref{Paper01_infinitesimal_generator} are well defined for any $\phi \in \mathscr{C}_*(\mathcal{S}, \mathbb{R})$. 

As we will discuss in Section~\ref{Paper01_ssec:non-monotonicity} below, the particle system that we want to construct is non-monotone. 
Together with the lack of \emph{a priori} bounds on both the number of particles per deme and the per capita birth and death rates, as well as the non-local interactions in the type space, this makes the construction of the process highly non-standard.

Our strategy will be to construct the desired Markov process as the weak limit of a sequence of approximating Feller processes~$(\eta^{n})_{n \in \mathbb{N}}$. More precisely, let $(\lambda_{n})_{n \in \mathbb{N}} \subset (0, \infty)$ and $(K_{n})_{n \in \mathbb{N}} \subseteq \mathbb{N}$ denote increasing sequences such that
\begin{equation} \label{Paper01_scaling_parameters_for_discrete_approximation}
    \lim_{n \rightarrow \infty} \, \lambda_{n} = \infty \quad \textrm{and} \quad  \lim_{n \rightarrow \infty} \, K_{n} = \infty.
\end{equation}
For each $n \in \mathbb{N}$, let $\Lambda_{n} \defeq L^{-1}\mathbb{Z} \, \cap \, [-\lambda_n, \lambda_n] $ and let $(\eta^{n}(t))_{t \geq 0}$ be a Markov process on $(\mathcal{S},d_{\mathcal{S}})$ with associated infinitesimal generator~$\mathcal{L}^{n}$ given as follows:
for every~$n \in \mathbb{N}$, any~$\phi \in \mathscr{C}(\mathcal{S}, \mathbb{R})$ and any~$\boldsymbol{\xi} = (\xi_{k}(x))_{k \in \mathbb{N}_{0}, \, x \in L^{-1}\mathbb{Z}}\in \mathcal{S}$, let
\begin{equation} \label{Paper01_generator_foutel_etheridge_model_restriction_n}
    \mathcal{L}^{n}\phi(\boldsymbol{\xi}) = \frac{m}{2}\mathcal{L}^{n}_{m}\phi(\boldsymbol{\xi}) + \mathcal{L}^{n}_{r}\phi(\boldsymbol{\xi}),
\end{equation}
where
\begin{equation} \label{Paper01_infinitesimal_generator_restriction_n}
\begin{aligned}
(\mathcal{L}^{n}_{m}\phi)(\boldsymbol{\xi}) & \defeq \sum_{x \in \Lambda_n} \; \sum_{k = 0}^{\infty} \; \sum_{z \in \{-L^{-1}, L^{-1}\}}  \mathds{1}_{\{x + z \in \Lambda_n\}}\xi_{k}(x) \left(\phi\Big(\boldsymbol{\xi} + \boldsymbol{e}_{k}^{(x+z)} - \boldsymbol{e}_{k}^{(x)}\Big) - \phi(\boldsymbol{\xi}) \right),
\\ 
(\mathcal{L}^{n}_{r}\phi)(\boldsymbol{\xi}) & \defeq \sum_{x \in \Lambda_n} \Bigg(\sum_{k = 0}^{K_n} F^b_k(\xi(x))\left(\phi\Big(\boldsymbol{\xi} + \boldsymbol{e}_{k}^{(x)}\Big) - \phi(\boldsymbol{\xi})\right) + \sum_{k = 0}^{\infty} F^d_k(\xi(x))\left(\phi\Big(\boldsymbol{\xi}- \boldsymbol{e}_{k}^{(x)}\Big) - \phi(\boldsymbol{\xi})\right)\Bigg),
\end{aligned}    
\end{equation}
and where~$(F^b_k)_{k \in \mathbb{N}_0}$ and $(F^d_k)_{k \in \mathbb{N}_0}$ are given by~\eqref{Paper01_discrete_birth_rate} and~\eqref{Paper01_discrete_death_rate} respectively.
Let $\{P^n_t\}_{t\ge 0}$ denote the semigroup of probability measures associated to $(\eta^n(t))_{t\ge 0}$.

Informally, the approximating process $\eta^n$ can be described as follows: particles outside the box $\Lambda_n$ are frozen. All particles living inside $\Lambda_n$
can migrate and die, but particles on the boundary of $\Lambda_n$ can only migrate to a deme inside $\Lambda_n$. Finally, only particles that carry $K_n$ or fewer mutations are able to reproduce. Since the number of demes in $\Lambda_n$ is finite, and only particles carrying at most $K_n$ mutations can reproduce, the process $(\eta^{n}(t))_{t \geq 0}$ can be thought of as a Markov process with finitely many types (ignoring the frozen particles outside $\Lambda_n$).

By using the fact that by Assumption~\ref{Paper01_assumption_polynomials}, the birth and death polynomials satisfy $0 \leq \deg q_{+} < \deg q_{-}$, and by  adapting the technique of correlation functions introduced by Boldrighini, De~Masi, Pellegrinotti and Presutti in~\cite{boldrighini1987collective}, we will verify in Section~\ref{Paper01_section_correlation_functions} that for any $n \in \mathbb{N}$ and any $\boldsymbol{\xi} \in \mathcal{S}$, conditioned on $\eta^{n}(0) = \boldsymbol{\xi}$, the Markov process $(\eta^{n}(t))_{t \geq 0}$ is well-defined as an $\mathcal{S}$-valued càdlàg process and does not explode in finite time (see Proposition~\ref{Paper01_bound_total_mass}).

To state the main results of this article, we need to introduce some extra notation. Define the set of bounded continuous cylindrical functions as follows:
\begin{align} \label{Paper01_general_definition_cylindrical_function}
    \mathscr{C}^{\textrm{cyl}}_{b,*}(\mathcal{S}, \mathbb{R}) & \defeq \Big\{\phi :\mathcal{S} \rightarrow \mathbb{R} \text{ such that }\; \exists J \in \mathbb{N}, \, (x_{1}, \ldots, x_{J}) \in (L^{-1}\mathbb{Z})^{J}, \, (k_{1}, \ldots, k_{J}) \in \left(\mathbb{N}_{0}\right)^{J}, \notag \\ &  \quad \quad \quad A\subset (\mathbb N_0 )^J \text{ finite, and } \psi\in \mathscr{C}_b((\mathbb N_0)^J,\mathbb R) \text{ with }\psi(x)=\psi(y) \; \forall x,y\notin A, \notag \\ 
    & \quad \quad \quad \textrm{with~}  \phi(\boldsymbol{\xi}) = \psi(\xi_{k_{1}}(x_{1}), \ldots, \xi_{k_{J}}(x_{J})) \; \forall \boldsymbol{\xi} = (\xi_{k}(x))_{k \in \mathbb{N}_{0}, \, x \in L^{-1}\mathbb{Z}} \in \mathcal{S}\Big\}.
\end{align}
Moreover, let $\mathscr{C}_{*}(\mathcal{S}, \mathbb{R})$ denote the following subset of $\mathscr{C}(\mathcal{S}, \mathbb{R})$ that can be thought of as a set of functions whose discrete derivatives decrease polynomially fast in distance from the origin:
\begin{equation} \label{Paper01_definition_Lipschitz_function}
\begin{aligned}
    & \mathscr{C}_*(\mathcal{S}, \mathbb{R}) \\ \, & \defeq \Bigg\{\phi \in \mathscr{C}(\mathcal{S}, \mathbb{R}): \\
    &\qquad \sup_{\substack{x \in L^{-1}\mathbb{Z}, \, k \in \mathbb{N}_{0}}} \;\sup_{\boldsymbol{\zeta}=(\zeta_j(y))_{j\in \mathbb N_0,\, y\in L^{-1}\mathbb Z} \in \mathcal{S}} \Bigg( \Big\vert \phi\Big(\boldsymbol{\zeta} + \boldsymbol{e}^{(x)}_{k}\Big) - \phi(\boldsymbol{\zeta}) \Big\vert +  \mathds{1}_{\{\zeta_{k}(x)>0\}}\Big\vert \phi\Big(\boldsymbol{\zeta} - \boldsymbol{e}^{(x)}_{k}\Big) - \phi(\boldsymbol{\zeta}) \Big\vert \\ 
    & \hspace{2.5cm} + \mathds{1}_{\{\zeta_{k}(x)>0\}} \sum_{a = 1}^{2} \Big\vert \phi\Big(\boldsymbol{\zeta} + \boldsymbol{e}^{(x+ (-1)^{a}L^{-1})}_{k} - \boldsymbol{e}^{(x)}_{k}\Big) - \phi(\boldsymbol{\zeta}) \Big\vert \Bigg) (1 + \vert x \vert)^{2(1 + \deg q_{-})} < \infty \Bigg\}.
\end{aligned}
\end{equation}
We will prove in Lemma~\ref{Paper01_characterisation_lipschitz_functions_semi_metric}(v) in the appendix that~$\mathscr{C}_{b,*}^{\textrm{cyl}}(\mathcal{S}, \mathbb{R}) \subseteq \mathscr{C}_{*}(\mathcal{S}, \mathbb{R})$. 

For a complete and separable metric space $(\mathcal A,d_{\mathcal A})$, and an $\mathcal A$-valued càdlàg Markov process $(\eta(t))_{t \geq 0}$ constructed on a probability space $(\Omega, \mathcal{F}, \mathbb{P})$, we define its natural filtration $\{\mathcal{F}_{t}^{\eta}\}_{t \geq 0}$ by $\mathcal{F}_{t}^{\eta} \defeq \sigma(\eta(u): u \leq t)$ for $t\ge 0$, and define its right-continuous filtration $\{\mathcal{F}_{t+}^{\eta}\}_{t \geq 0}$, for any $t \geq 0$, as
\begin{equation} \label{eq:rightctsfiltration}
    \mathcal{F}^{\eta}_{t+} \defeq  \bigcap_{T > t} \mathcal{F}^{\eta}_{T}.
\end{equation}
Let $\{P_t\}_{t \geq 0}$ be the semigroup associated to $(\eta(t))_{t \geq 0}$. Following the exposition in~\cite[Chapter~4]{ethier2009markov}, we say that $(\eta(t))_{t \geq 0}$ is a strong Markov process with respect to the filtration $\{\mathcal{F}_{t+}^{\eta}\}_{t \geq 0}$ if for any almost surely finite $\{\mathcal{F}_{t+}^{\eta}\}_{t \geq 0}$-stopping time $\tau$, any $t \geq 0$ and any $\phi \in \mathscr{B}_{b}(\mathcal{A}, \mathbb{R})$, 
\begin{equation} \label{Paper01_strong_markov_property}
    \mathbb{E}\Big[\phi(\eta(t + \tau)) \, \Big\vert \, \mathcal{F}_{\tau+}^{\eta} \Big] = (P_t\phi)(\eta(\tau)) \, \textrm{ almost surely}.
\end{equation}
Moreover, in Section~\ref{Paper01_subsection_appendix_feller_non_locally_compact_polish_space} of the appendix we recall the definition of (weakly) Feller semigroups on non-locally compact Polish spaces (see Definition~\ref{Paper01_definition_feller_non_locally_compact}), as well as the definition of a Feller process (at the end of Section~\ref{Paper01_subsection_appendix_feller_non_locally_compact_polish_space}). We can now state our first main result, which shows that there exists a unique limit of the sequence of processes $\left((\eta^n(t))_{t \geq 0}\right)_{n \in \mathbb{N}}$.

\begin{theorem} \label{Paper01_thm_existence_uniqueness_IPS_spatial_muller}
Take $L,m>0$, $N\in \mathbb N$, $\mu\in [0,1]$ and $(s_k)_{k\in \mathbb N_0}$, $q_+$, $q_-$ satisfying Assumptions~\ref{Paper01_assumption_fitness_sequence} and~\ref{Paper01_assumption_polynomials}, and recall the definitions of $\mathcal S$, $\mathcal S_0$ and $\mathcal L$ in~\eqref{Paper01_definition_state_space_formal},~\eqref{Paper01_set_initial_configurations} and~\eqref{Paper01_generator_foutel_etheridge_model} respectively.
Then for any choice of increasing sequences $(\lambda_n)_{n\in \mathbb N}$ and $(K_n)_{n\in \mathbb N}$ satisfying~\eqref{Paper01_scaling_parameters_for_discrete_approximation}, 
    for any $T \geq 0$, the linear operator $P_T : \mathscr{C}_b(\mathcal{S}, \mathbb{R}) \rightarrow \mathscr{C}_b(\mathcal{S}, \mathbb{R})$ given by, for $\phi \in \mathscr{C}_b(\mathcal{S}, \mathbb{R})$ and $\boldsymbol{\eta} \in \mathcal{S}$,
    \begin{equation*}
        (P_T\phi)(\boldsymbol{\eta}) \defeq \lim_{n \rightarrow \infty} \,  (P^n_T\phi)(\boldsymbol{\eta}),
    \end{equation*}
    where $P^n_T$ is defined after~\eqref{Paper01_infinitesimal_generator_restriction_n},
    is a well-defined bounded linear operator which does not depend on the choice of $(\lambda_n)_{n\in \mathbb N}$ and $(K_n)_{n\in \mathbb N}$, and the family $\{P_t\}_{t \geq 0}$ is a (weakly) Feller semigroup on $\mathcal{S}$. Moreover, for any $\boldsymbol{\eta} \in \mathcal{S}_{0}$, conditioning on $\eta^{n}(0) = \boldsymbol{\eta}$ for every $n \in \mathbb{N}$, the sequence of processes $(\eta^{n}(t))_{t \geq 0}$ converges weakly with respect to the $J_{1}$-topology on $\mathscr{D}([0, \infty), \mathcal{S})$ as $n \rightarrow \infty$ to an $\mathcal S$-valued càdlàg Markov process $(\eta(t))_{t \geq 0}$ with $\eta(0) = \boldsymbol{\eta}$ almost surely and with finite dimensional distributions given by the semigroup $\{P_t\}_{t \geq 0}$. Also, $(\eta(t))_{t \geq 0}$ is a strong Markov process with respect to its right-continuous filtration $\{\mathcal{F}^{\eta}_{t+}\}_{t \geq 0}$, and for any $\phi \in \mathscr{C}_*(\mathcal{S}, \mathbb{R})$, the process $(M^{\phi}(t))_{t \geq 0}$ given by
    \begin{equation*}
       M^{\phi}(T) \defeq \phi(\eta(T)) -\phi(\boldsymbol{\eta}) - \int_{0}^{T} (\mathcal{L}\phi)(\eta(t-)) \, dt \quad \text{for }T \geq 0
    \end{equation*}
    is a càdlàg $\{\mathcal{F}^{\eta}_{t+}\}_{t\ge 0}$-square integrable martingale.
\end{theorem}

For any $\boldsymbol{\eta} \in \mathcal S_0$, we will refer to 
the $\mathcal{S}$-valued càdlàg Markov process $(\eta(t))_{t \geq 0}$ with $\eta(0) = \boldsymbol{\eta}$ defined in Theorem~\ref{Paper01_thm_existence_uniqueness_IPS_spatial_muller} as the \emph{spatial Muller's ratchet}.

Our second main result shows that the spatial Muller's ratchet satisfies certain moment estimates on the local density of particles, uniformly in the spatial rescaling, the migration rate and the scaling parameter for the carrying capacity. As part of the proof of this result, we will prove a similar result for the approximating process $(\eta^n_t)_{t \geq 0}$ which will be essential for proving the tightness required in the proof of Theorem~\ref{Paper01_thm_existence_uniqueness_IPS_spatial_muller}.

\begin{theorem} \label{Paper01_thm_moments_estimates_spatial_muller_ratchet}
Take $q_+$ and $q_-$ satisfying Assumption~\ref{Paper01_assumption_polynomials}; then
    for any $p \geq 1$, there is a non-decreasing function $C_{p}: [0,\infty) \rightarrow [1,\infty)$ such that the following holds.
 For any $L,m>0$, $N\in \mathbb N$, $\mu\in [0,1]$ and $(s_k)_{k\in \mathbb N_0}$ satisfying Assumption~\ref{Paper01_assumption_fitness_sequence}, recalling the definition of $\mathcal S_0$ in~\eqref{Paper01_set_initial_configurations},
    for any $\boldsymbol{\eta} =(\eta(x))_{x\in L^{-1}\mathbb Z} \in \mathcal{S}_{0}$ and $T \geq 0$,
\begin{equation*}
  \sup_{x \in L^{-1}\mathbb{Z}} \mathbb{E}_{\boldsymbol{\eta}}\left[{\left\vert\left\vert \eta(T,x) \right\vert\right\vert_{\ell_{1}}^{p}}\right] \leq C_{p}(T) \Big(\sup_{x \in L^{-1}\mathbb{Z}} \, {\vert \vert \eta(x) \vert \vert_{\ell_{1}}^{p}} + N^p\Big),
\end{equation*}
where $\mathbb P_{\boldsymbol{\eta}}$ is the probability measure under which $(\eta(t))_{t \geq 0}=(\eta(t,x))_{t\ge 0,\, x\in L^{-1}\mathbb Z}$ is the $\mathcal S$-valued càdlàg Markov process defined in Theorem~\ref{Paper01_thm_existence_uniqueness_IPS_spatial_muller} with $\eta(0) = \boldsymbol{\eta}$ almost surely, and $\mathbb E_{\boldsymbol{\eta}}$ is the corresponding expectation.
\end{theorem}

\begin{remark}
\begin{enumerate}[(i)]
\item In our companion article~\cite{madeiraPreprintmuller}, we establish a law of large numbers for the spatial Muller's ratchet when scaling the parameters $L,m$ and the initial number of particles with $N$. Under this scaling, Theorem~\ref{Paper01_thm_moments_estimates_spatial_muller_ratchet} implies that the rescaled number of particles remains (locally) controlled, which is essential for our proof of the partial differential equation limit in~\cite{madeiraPreprintmuller}.

\item We note that the arguments used in the proofs of Theorems~\ref{Paper01_thm_existence_uniqueness_IPS_spatial_muller} and~\ref{Paper01_thm_moments_estimates_spatial_muller_ratchet} remain valid under more general assumptions. In particular, Assumption~\ref{Paper01_assumption_fitness_sequence}(iii) is not essential for our main results. Furthermore, instead of requiring the per-capita birth and death rates $q_+, q_- : [0,\infty)\to[0,\infty)$ to be polynomials satisfying Assumption~\ref{Paper01_assumption_polynomials}, it suffices that $q_+$ and $q_-$ are non-negative locally Lipschitz functions such that
\begin{enumerate}[(1)]
    \item $ \liminf_{u\to\infty} \displaystyle \frac{q_-(u)}{u\,q_+(u)+1} > 0$.
    \item  $\lim_{u\to\infty} q_-(u) = \infty$.
    \item There exists $u^*\in(0,\infty)$ for which $q_-(u^{(1)}) > q_-(u^{(2)})$ whenever $u^{(1)} > u^{(2)} > u^*$.
    \item There exists a non-negative polynomial $p: [0,\infty) \rightarrow [0, \infty)$ such that $q_-(u) \leq p(u)$, for all $u \in [0, \infty)$.
\end{enumerate}
The same arguments also yield the construction of the spatial Muller's ratchet on the lattice $L^{-1}\mathbb{Z}^d$ for any spatial dimension $d\in\mathbb{N}$.

\item    We highlight that Theorems~\ref{Paper01_thm_existence_uniqueness_IPS_spatial_muller} and~\ref{Paper01_thm_moments_estimates_spatial_muller_ratchet} allow us to represent the spatial Muller's ratchet in terms of an infinite system of stochastic differential equations driven by Poisson point processes through an application of the Markov Mapping Theorem (see for instance~\cite[Theorem~A.2]{etheridge2019genealogical} and the references therein). Since this construction is not needed for the purposes of this article, we omit the details.
\end{enumerate}
\end{remark}

\section{Discussion}\label{Paper01_sec:discussion}

\subsection{Non-monotonicity}\label{Paper01_ssec:non-monotonicity}

One of the main issues that makes the approximation of the spatial Muller's ratchet by a particle system with finitely many particles difficult is that our process is \emph{non-monotone}
due to the competition between particles carrying different numbers of mutations.

We start by recalling the definition of a monotone interacting particle system on the rescaled lattice $L^{-1}\mathbb Z$. Let $(\xi(t))_{t \geq 0}$ be a particle system taking values in $\mathcal{X}^{L^{-1}\mathbb{Z}}$, where $\mathcal{X}$ is some complete, separable metric space equipped with a partial order relation~$\preccurlyeq$. We say that $(\xi(t))_{t \geq 0}$ is \textit{monotone} if for any two configurations $\boldsymbol{\xi}^{(1)}=(\xi^{(1)}(x))_{x\in L^{-1}\mathbb Z}, \boldsymbol{\xi}^{(2)}=(\xi^{(2)}(x))_{x\in L^{-1}\mathbb Z} \in \mathcal{X}^{L^{-1}\mathbb{Z}}$ such that $\xi^{(1)}(x) \preccurlyeq \xi^{(2)}(x)$ for all $x \in L^{-1}\mathbb{Z}$, it is possible to construct realisations $(\xi^{(1)}(t))_{t \geq 0}=(\xi^{(1)}(t,x))_{t \geq 0,\, x\in L^{-1}\mathbb Z}$ and $(\xi^{(2)}(t))_{t \geq 0}=(\xi^{(2)}(t,x))_{t \geq 0,\, x\in L^{-1}\mathbb Z}$ of the particle system on the same probability space in such a way that $\xi^{(1)}(0) = \boldsymbol{\xi}^{(1)}$, $\xi^{(2)}(0) = \boldsymbol{\xi}^{(2)}$ and for all $t \geq 0$,
\begin{equation*}
    \xi^{(1)}(t,x) \preccurlyeq  \xi^{(2)}(t,x) \; \mbox{for all} \, x \in L^{-1}\mathbb{Z}.
\end{equation*}
Take a partial order relation~$\preccurlyeq$ on $\ell_1$ satisfying the following properties:
\begin{itemize}
    \item for any $x\in \ell_1^+$ with $x\neq 0$ we have $x\not\preccurlyeq 0.$
    \item for any $x^{(1)}=(x^{(1)}_k)_{k\in \mathbb N_0},$ $x^{(2)}=(x^{(2)}_k)_{k\in \mathbb N_0}\in \ell_1 \cap \mathbb N_0 ^{\mathbb N_0}$ with $x^{(1)}_k\le x^{(2)}_k$ $\forall k\in \mathbb N_0$ we have $x^{(1)}\preccurlyeq x^{(2)}$.
    \end{itemize}
 To see heuristically that the spatial Muller's ratchet, as informally defined at the start of Section~\ref{Paper01_model_description}, is not monotone in general, 
 consider the classical Fisher-KPP dynamics: that is, suppose the polynomials $q_{+}$ and $q_{-}$ are given, for all $U \in [0, \infty)$, by
\[
    q_{+}(U) \equiv 1 
    \quad \text{and} \quad 
    q_{-}(U) \equiv U.
\]
(Note that Assumption~\ref{Paper01_assumption_polynomials} holds for this choice of $q_{+}$ and $q_{-}$.)
Also, for simplicity, suppose that $\mu$ and $s_1$ are very small.
We call a particle carrying $k$ mutations a `type-$k$' particle. Using the notation defined before~\eqref{Paper01_scaled_polynomials_carrying_capacity}, let $\boldsymbol{\eta}^{(1)}=\boldsymbol{e}^{(0)}_0$ and $\boldsymbol{\eta}^{(2)}=\boldsymbol{e}^{(0)}_0+\boldsymbol{e}^{(0)}_1$, and take an arbitrary coupling between realisations $(\eta^{(1)}(t))_{t\ge 0}$ and $(\eta^{(2)}(t))_{t\ge 0}$ with initial conditions $\boldsymbol{\eta}^{(1)}$ and $\boldsymbol{\eta}^{(2)}$.

{For $i \in \{1,2\}$, let $E_i$ denote the event that the type-$0$ particle in $\eta^{(i)}$ dies before it gives birth to another particle, and let $\tau^{(i)}$ denote the first time the type-$0$ particle in $\eta^{(i)}$ reproduces or dies. Note that for each $k\in \mathbb N_0$, the birth rate for every type-$k$ particle is $s_k$, and the death rate for a particle in a deme containing $n$ particles is $n/N$. Therefore, the probability that the initial type-$0$ particle dies before reproducing is strictly higher in  $(\eta^{(2)}(t))_{t\ge 0}$ than in $(\eta^{(1)}(t))_{t\ge 0}$, i.e.~$\mathbb{P}(E_2) > \mathbb{P}(E_1)$. Hence, under any coupling, the event $E_2 \cap (E_1)^c$ must have strictly positive probability.}

\begin{figure}[t]
\centering
\includegraphics[width=0.75\linewidth,trim=0cm 0cm 0cm 0cm,clip=true]{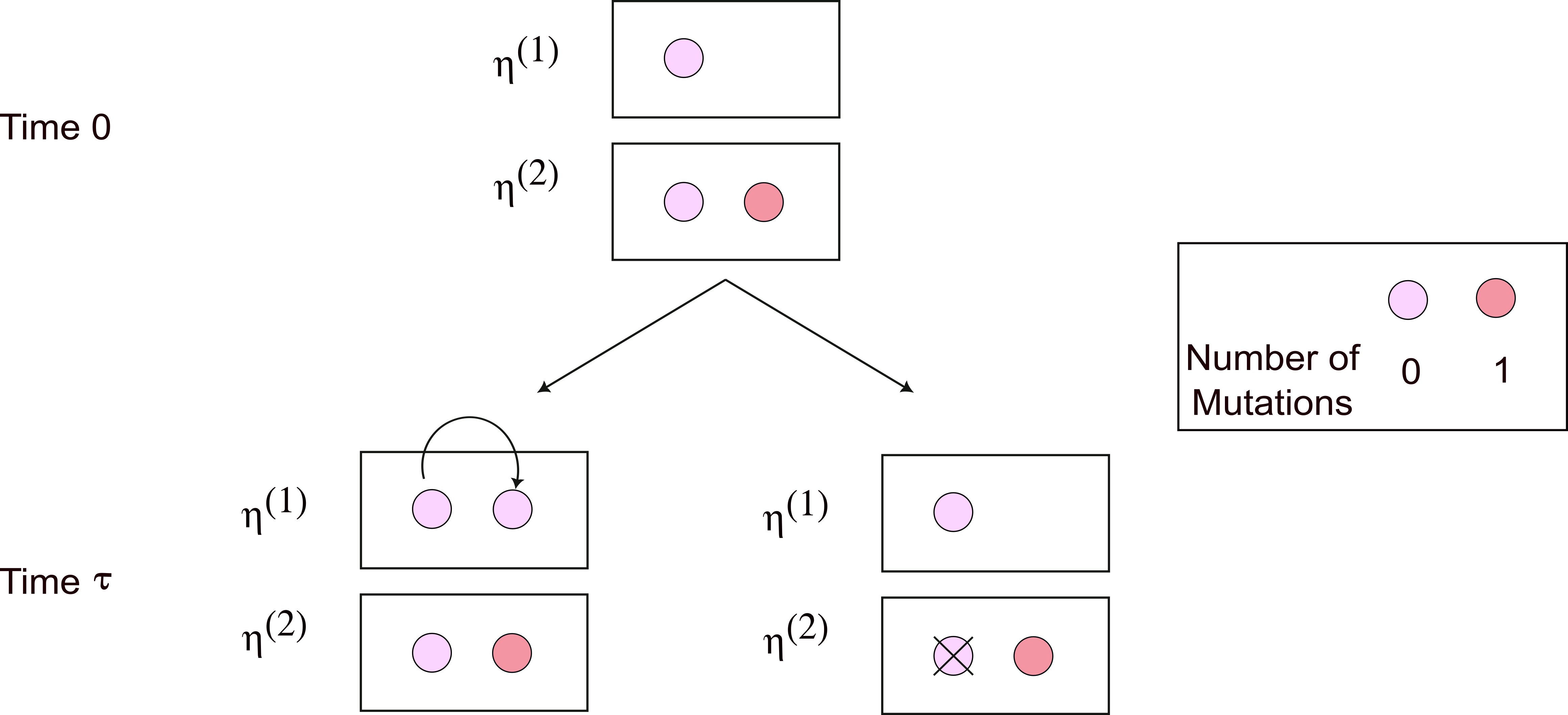}
\caption{\label{Paper01_non_monotone_particle_system_realisation}~Coupling between two realisations of the spatial Muller's ratchet $(\eta^{(1)}(t))_{t \geq 0}$ and $(\eta^{(2)}(t))_{t \geq 0}$, started from initial configurations $\eta^{(1)}(0) = \boldsymbol{e}^{(0)}_0$ and $\eta^{(2)}(0) = \boldsymbol{e}^{(0)}_0 + \boldsymbol{e}^{(0)}_1$. Then, by letting $\tau$ denote the first time at which the type-$0$ particle dies or reproduces in either process, with positive probability either the type-$0$ particle in $\eta^{(1)}$ reproduces at time $\tau$, or the type-$0$ particle in $\eta^{(2)}$ dies at time $\tau$.}
\end{figure}

{By assuming that~$\mu$ and~$s_1$ are sufficiently small (as mentioned above), it follows that  with positive probability, letting $\tau = \tau^{(1)} \wedge \tau^{(2)}$ denote the first time at which the type-$0$ particle dies or reproduces in either process, either (see Figure~\ref{Paper01_non_monotone_particle_system_realisation})
\begin{itemize}
    \item $\eta^{(2)}(\tau)$ has no type-$0$ particles and at most one other particle, and $\eta^{(1)}(\tau)$ consists of at least one type-$0$ particle and no other particles, or 
    \item $\eta^{(1)}(\tau)$ consists of two type-$0$ particles and no other particles, and $\eta^{(2)}(\tau)$ consists of at most one type-$0$ particle and at most one other particle.
\end{itemize}}
Then by assuming that $N \in \mathbb{N}$ is large, and $\mu, s_1 \ll N^{-1}$, and so particles that are not type-$0$ reproduce at a very low rate and new mutations are rare, and type-$0$ particles in low density regions have birth rate higher than their death rate, one can convince oneself that, conditioning on the event $E_2 \cap (E_1)^c$, with positive probability, there exists $t>0$ such that
\[
\eta^{(2)}(t)=0 \quad \text{and}\quad \eta^{(1)}(t)\neq 0,
\]
and therefore the system is non-monotone.
We emphasise that this is intended purely as an illustrative example and not a formal proof of non-monotonicity.

\subsection{Comparison with the literature} \label{Paper01_related_work}

Interacting particle systems consisting of finitely many particles can be constructed directly as Markov processes taking values in $\mathbb{Z}^{d}$; see for example~\cite[Section~4.2]{ethier2009markov}. On the other hand, the standard methods for the construction of infinite particle systems can be roughly divided into the following categories:
\begin{enumerate}[(A)]
    \item Methods based on the theory of Feller processes in locally compact Polish spaces~\cite{liggett1972existence,iscoe1986weighted};
    \item Methods based on monotonicity~\cite[Chapter~IV]{demasi1991mathematical};
    \item Methods based on the Lipschitz property of the infinitesimal generator and its associated semigroup as in~\cite[Chapter~IX]{liggett1985interacting} and~\cite[Chapter~13]{chen2004markov};
    \item Methods based on stochastic differential equations driven by Poisson measures on Polish spaces~\cite{kurtz1999particle,kurtz1996weak}.
\end{enumerate}

The first group of methods (A) includes the standard machinery used for spatial interacting particle systems with finitely many types and a bounded number of particles per site (see, for instance,~\cite{liggett1972existence} and~\cite[Chapter~I]{liggett1985interacting}). In this case, the state space in which the particle system evolves is $(\mathbb{Z}^d)^\mathcal{K}$, for some compact set $\mathcal{K}$, and hence the state space is 
compact with respect to the product topology (see~\cite[Section~I.3]{liggett1985interacting}). Therefore, one can use the standard theory of Feller processes, including the application of the Hille--Yosida theorem, to prove existence, uniqueness and regularity  properties of the interacting particle system. It can then be established that if the transition rates between states with different numbers of particles at each spatial location are uniformly bounded, the process is well-defined (see~\cite[Theorem~I.3.9]{liggett1985interacting}). Since in our scenario the number of particles per deme is not uniformly bounded, we cannot apply this approach directly.

An alternative method, also based on the theory of Feller processes in locally compact Polish spaces, was developed in~\cite{iscoe1986weighted}. In that work, Iscoe constructs measure-valued branching processes starting from configurations consisting of infinitely many particles. His idea was to construct the branching process as a weighted measure-valued process on~$\dot{\mathbb{R}}^d$, the one-point compactification of~$\mathbb{R}^d$, for some $d \in \mathbb{N}$. To follow his approach, it is crucial to characterise the behaviour of the process `at infinity'. To illustrate why this technique cannot be applied to construct the spatial Muller's ratchet, let $\mathcal{Z} \defeq L^{-1}\mathbb{Z} \times \mathbb{N}_0$, and denote its one-point compactification by $\dot{\mathcal{Z}}$. The state space~$\mathcal{S}$ defined in~\eqref{Paper01_definition_state_space_formal} can be viewed as the space of weighted, positive, integer-valued measures on~$\mathcal{Z}$. By extending~$\mathcal{S}$ to the space~$\mathcal{S}'$ of such measures on~$\dot{\mathcal{Z}}$, one can define a metric~$d_{\mathcal{S}'}$ such that~$(\mathcal{S}', d_{\mathcal{S}'})$ is a locally compact, complete, and separable metric space. To apply the standard theory of Feller processes in this setting, it is then necessary to characterise the behaviour of the spatial Muller's ratchet `at the infinity' in~$\dot{\mathcal{Z}}$. The difficulty with applying this strategy in our setting is that there are in fact two qualitatively different types of ‘infinity' to consider, at which the behaviour of the particle system should be different: one corresponds to the accumulation of particles far from the origin in the spatial direction (i.e.~in~$L^{-1}\mathbb{Z}$), and the other corresponds to the accumulation of particles carrying a high number of mutations. Since the reproduction rate converges to $0$ as the number of mutations increases (see Assumption~\ref{Paper01_assumption_fitness_sequence}(iv)), the behaviour of particles in these two directions toward infinity diverges significantly. Therefore, the arguments in~\cite{iscoe1986weighted} cannot be applied directly to the spatial Muller's ratchet.

Another possible strategy for constructing infinite particle systems relies on a monotonicity property and the use of Carathéodory's extension theorem. This method corresponds to the category~(B) listed at the beginning of this section, and was used for example by De Masi and Presutti in~\cite[Chapter~IV]{demasi1991mathematical} to construct one-type interacting particle systems with polynomial birth and death rates as the weak limit of a sequence of particle systems with truncated birth rates. However, because of the non-monotonicity of our particle system, we cannot use this technique directly.

Category~(C) includes methods that rely on a global Lipschitz condition for the infinitesimal generator of the Markov process -- a property that does not hold in general, especially for particle systems with non-monotone interactions. This approach is used by Liggett in~\cite[Section~IX.1]{liggett1985interacting} and by Chen in~\cite[Chapter~13]{chen2004markov} to construct infinite particle systems. Interestingly, this construction can be applied to some interacting particle systems with multiple (but finitely many) types of particles (see~\cite[Example~13.39]{chen2004markov}), but the topology used by Chen cannot be extended to particle systems with infinitely many types of particles. Moreover, the existence and uniqueness results stated by Liggett in~\cite[Theorem~IX.1.14]{liggett1985interacting} and by Chen in~\cite[Theorems~13.17 and~13.18]{chen2004markov} do not imply the existence of càdlàg versions of the Markov process or the strong Markov property.

Finally, the methods in category~(D) rely on the characterisation of infinite particle systems as infinite-dimensional stochastic differential equations. Although this is a powerful method to prove existence and uniqueness of reaction-diffusion Markov processes, the proof of uniqueness also requires a global Lipschitz property of the reaction part of the infinitesimal generator~\cite{kurtz1999particle,kurtz1996weak}, which is not satisfied in general by the infinitesimal generator $\mathcal{L}$ defined in~\eqref{Paper01_generator_foutel_etheridge_model} and~\eqref{Paper01_infinitesimal_generator}.

\subsection{Overview of the proof of Theorems~\ref{Paper01_thm_existence_uniqueness_IPS_spatial_muller} and~\ref{Paper01_thm_moments_estimates_spatial_muller_ratchet}} \label{Paper01_subsection_heuristics}

In this section, we give an overview of the main ideas in the proofs of Theorems~\ref{Paper01_thm_existence_uniqueness_IPS_spatial_muller} and~\ref{Paper01_thm_moments_estimates_spatial_muller_ratchet}.

In Section~\ref{Paper01_section_functional_analysis_feller_non_locally_compact}, we will establish, in Theorem~\ref{Paper01_general_thm_convergence_feller_non_locally_compact}, a general recipe for constructing a semigroup corresponding to a càdlàg strong Markov process taking values in a non-locally compact Polish space
$(\mathcal{S}, d_{\mathcal{S}})$. As is standard in the construction of infinite particle systems, our approach is based on approximating the target Markov process $(\eta(t))_{t \geq 0}$ by a sequence of $\mathcal{S}$-valued càdlàg Feller processes $((\eta^n(t))_{t \geq 0})_{n \in \mathbb{N}}$. Although the recipe in Section~\ref{Paper01_section_functional_analysis_feller_non_locally_compact} is rather general, for the purposes of our proof outline in this section we will restrict our attention to the construction of the spatial Muller's ratchet and therefore use the notation introduced in Section~\ref{Paper01_model_description}. Roughly speaking, there are two main properties that we must establish in order to apply our recipe in Section~\ref{Paper01_section_functional_analysis_feller_non_locally_compact} to the sequences of processes $(\eta^n)_{n\in \mathbb N}$:
\begin{enumerate}[(i)]
    \item For any $\boldsymbol{\eta} \in \mathcal{S}$, 
    conditioning on $\eta^n(0) = \boldsymbol{\eta}$ $\forall n \in \mathbb{N}$,
    the one-dimensional distributions of $((\eta^n(t))_{t \geq 0})_{n \in \mathbb{N}}$ converge weakly.
    \item For any $\boldsymbol{\eta} \in \mathcal{S}_{0}$, conditioning on $\eta^{n}(0) = \boldsymbol{\eta}$ $\forall n \in \mathbb{N}$, the sequence of processes $(\eta^{n})_{n \in \mathbb{N}}$ is tight in $\mathscr{D}([0, \infty), \mathcal{S})$.
\end{enumerate}
These properties are shown in Section~\ref{Paper01_section_functional_analysis_feller_non_locally_compact} to imply that there is a limiting process $(\eta(t))_{t\ge 0}$.
Indeed, we use property~(i) to construct a family of probability measures $\{P_t\}_{t \geq 0}$ on $\mathcal{S}$ from the sequence of semigroups $(\{P^n_t\}_{t \geq 0})_{n \in \mathbb{N}}$ associated to $(\eta^n)_{n \in \mathbb{N}_0}$. By requiring the existence of a (sufficiently large) subset $\mathfrak{A} \subset \mathscr{C}_b(\mathcal{S}, \mathbb{R})$ such that for every $\phi \in \mathfrak{A}$ and $t\ge 0$, $(P^n_t\phi)_{n \in \mathbb{N}}$ converges uniformly on compact subsets of $\mathcal{S}$, we are able to establish the semigroup and other regularity properties for $\{P_t\}_{t \geq 0}$. By standard arguments in stochastic analysis (see~\cite[Theorem~4.1.1]{ethier2009markov}), these properties are enough to guarantee the existence of a Markov process taking values in~$\mathcal{S}$ whose finite-dimensional distributions are determined by the semigroup $\{P_t\}_{t \geq 0}$. If we were then to assume that $\mathcal{S}$ were locally compact, we could construct càdlàg modifications of this Markov process directly~\cite[Theorem~4.2.5]{ethier2009markov}. However, since our state space $\mathcal S$ is not locally compact, we cannot use this argument. Instead, we rely on condition~(ii), namely the tightness of the sequence of processes $(\eta^n)_{n \in \mathbb{N}}$ in $\mathscr{D}([0, \infty),\mathcal{S})$ for initial conditions in $\mathcal{S}_0$.

To establish properties (i) and~(ii) for the spatial Muller's ratchet, we begin by proving moment bounds for $(\eta^{n}(t))_{t \geq 0}$ in Section~\ref{Paper01_section_correlation_functions}. We will analyse an auxiliary process, in which all particles carry no mutations. For $n \in \mathbb{N}$, denote this auxiliary process by $(\zeta^{n}(t))_{t\geq 0} = (\zeta^{n}(t,x), x \in \Lambda_n)_{t \geq 0}$ taking values in $\mathbb{N}_{0}^{\Lambda_n}$. 
Recall from Assumption~\ref{Paper01_assumption_fitness_sequence},~\eqref{Paper01_infinitesimal_generator_restriction_n} and~\eqref{Paper01_discrete_birth_rate} that particles in $(\eta^n(t))_{t\ge 0}$ carrying larger numbers of mutations have lower birth rates.
We will define the process $\left(\zeta^{n}(t)\right)_{t \geq 0}$  in Section~\ref{Paper01_section_correlation_functions} in such a way that the death and migration rates for particles in $\zeta^{n}$ are the same as those for particles in $\eta^{n}$, and the birth rates for particles in $\zeta^{n}$ are at least as large as the birth rates for particles in $\eta^{n}$. We will then be able to couple the two processes in such a way that $\zeta^{n}(t,x) \geq \vert \vert \eta^{n}(t,x) \vert \vert_{\ell_{1}}$ for all $t \geq 0$ and $x \in \Lambda_n$ almost surely. Hence, in order to prove moment bounds on $\eta^{n}$, it will suffice to establish the corresponding bounds for $\zeta^{n}$.

The core technique that we will use to prove moment bounds on $\zeta^n$ is the method of correlation functions developed by Boldrighini, De Masi, Pellegrinotti and Presutti in \cite{boldrighini1987collective, demasi1991mathematical}. This method was applied in~\cite{boldrighini1987collective, demasi1991mathematical} to analyse the local behaviour of reaction-diffusion systems in which both birth and death rates are polynomial functions of the local number of particles. Heuristically, since, by Assumption~\ref{Paper01_assumption_polynomials}, the birth and death rate polynomials $q_+$ and $q_-$ are such that $0 \leq \deg q_{+} < \deg q_{-}$, whenever the local number of particles is very large, the rate of death events will be much higher than the rate of reproduction events. This will allow us to control the number of particles locally, even though the number of particles per deme is not uniformly bounded \emph{a priori}. The key difference between the arguments in~\cite[Chapter~IV]{demasi1991mathematical} and our results in Section~\ref{Paper01_section_correlation_functions} is that we keep an explicit dependence in our estimates on the parameter $N$, the scaling parameter for the carrying capacity. While tracking this dependence in the moment bounds is not needed for the construction of the spatial Muller's ratchet, 
it will allow us to prove the uniform estimate in Theorem~\ref{Paper01_thm_moments_estimates_spatial_muller_ratchet}, which in turn
will be essential in our companion article~\cite{madeiraPreprintmuller}. 

More precisely, our moment bounds on $\zeta^n$ will first imply (see Proposition~\ref{Paper01_bound_total_mass}) that for every $n \in \mathbb{N}$, any initial condition $\boldsymbol{\eta} \in \mathcal{S}$, and any time $T \geq 0$,
\begin{equation*}
    \mathbb{E}_{\boldsymbol{\eta}}\Bigg[\sup_{t \leq T} \sum_{x \in \Lambda_n} \vert \vert \eta^{n}(t,x) \vert \vert_{\ell_1}\Bigg] < \infty.
\end{equation*}
In Proposition~\ref{Paper01_bound_total_mass}, we will also establish that for every $p,n \in \mathbb{N}$ and any $\boldsymbol{\eta}=(\eta(x))_{x\in L^{-1}\mathbb Z} \in \mathcal{S}$ and $T\ge 0$,
\begin{equation} \label{Paper01_heuristics_correlation_functions_ii}
     \sup_{t \leq T} \; \sup_{x \in \Lambda_n} \mathbb{E}_{\boldsymbol{\eta}}\left[\left\vert\left\vert \eta^{n}(t,x) \right\vert\right\vert_{\ell_{1}}^{p}\right] \lesssim_p \sup_{x \in \Lambda_n} \, \vert \vert \eta(x) \vert \vert_{\ell_{1}}^{p}+N^{p}.
\end{equation}
For a general initial configuration $\boldsymbol{\eta} \in \mathcal{S}$, the right-hand side of~\eqref{Paper01_heuristics_correlation_functions_ii} may diverge as $n \to \infty$. For instance, recalling~\eqref{Paper01_definition_state_space_formal}, consider the configuration $\boldsymbol{\eta} = (\eta_{k}(x))_{k \in \mathbb{N}_{0}, \, x \in L^{-1}\mathbb{Z}} \in \mathcal{S}$ given by, for every $x \in L^{-1}\mathbb{Z}$,
\begin{equation*}
    \eta_{0}(x) = \lfloor (1 + \vert x \vert )^{1/2} \rfloor \; \textrm{and} \; \eta_{k}(x) = 0 \; \forall  k \in \mathbb{N}.
\end{equation*}
Then, by the definition of $\Lambda_n$ below~\eqref{Paper01_scaling_parameters_for_discrete_approximation}, we have
\begin{equation*}
    \lim_{n \rightarrow \infty} \sup_{x \in \Lambda_n} \, \vert \vert \eta(x) \vert \vert_{\ell_{1}} = \infty.
\end{equation*}
This possibility of divergence leads to technical challenges in constructing the semigroup associated with the spatial Muller's ratchet in Section~\ref{Paper01_formal_construction_generator_section}. Note, however, that by~\eqref{Paper01_set_initial_configurations}, for $\boldsymbol{\eta} \in \mathcal{S}_{0}$, the number of particles per deme in $\boldsymbol{\eta}$ is uniformly bounded, and~\eqref{Paper01_heuristics_correlation_functions_ii} yields the stronger estimate
\begin{equation} \label{Paper01_heuristics_correlation_functions_iii}
     \sup_{n \in \mathbb{N}} \; \sup_{t \leq T} \; \sup_{x \in L^{-1}\mathbb{Z}} \mathbb{E}_{\boldsymbol{\eta}}\left[\left\vert\left\vert \eta^{n}(t,x) \right\vert\right\vert_{\ell_{1}}^{p}\right] \lesssim_p \sup_{x \in L^{-1}\mathbb{Z}} \, \vert \vert \eta(x) \vert \vert_{\ell_{1}}^{p}+N^{p}.
\end{equation}
We highlight that~\eqref{Paper01_heuristics_correlation_functions_iii} is still a ‘local' estimate, in the sense that the supremum over space is outside the expectation. Nevertheless, by using~\eqref{Paper01_heuristics_correlation_functions_iii} and standard stochastic analysis techniques, including the Jakubowski criterion~\cite[Theorem~3.1]{jakubowski1986skorokhod}, in Section~\ref{Paper01_formal_construction_generator_section} we will be able to establish, for any $\boldsymbol{\eta} \in \mathcal{S}_0$,
tightness of the sequence $(\eta^n)_{n \in \mathbb{N}}$ in $\mathscr{D}([0, \infty), \mathcal{S})$ when conditioning on $\eta^n(0) = \boldsymbol{\eta} $ for every $n$, i.e.~we establish property~(ii) stated at the beginning of this subsection. It is important to note that our argument relies on the fact that the right-hand side of~\eqref{Paper01_heuristics_correlation_functions_iii} is finite, which holds only when the initial condition $\boldsymbol{\eta}$ belongs to $\mathcal{S}_0$. As a result, we are unable to establish tightness in $\mathscr{D}([0, \infty), \mathcal{S})$ for general initial data in $\mathcal{S}$.

For the same reason, we must use an alternative argument to prove tightness of one-dimensional distributions for general initial data in $\mathcal S$ (which we require in order to establish the convergence in property~(i) above). Recall the definition of $\vert \vert \vert \cdot \vert \vert \vert_{\mathcal{S}}$ in~\eqref{Paper01_definition_state_space_formal}. Since by Assumption~\ref{Paper01_assumption_polynomials} we have $0 \leq \deg q_+ < \deg q_-$, for any initial condition in $\mathcal{S}$, the expected increment in $\vert \vert \vert \cdot \vert \vert \vert_{\mathcal{S}}$ after a finite period of time $T \in [0, \infty)$ is uniformly bounded (see the proof of Proposition~\ref{Paper01_estimates_holding_uniformly_compact_subsets_S_N}). We will also establish in Proposition~\ref{Paper01_topological_properties_state_space} in the appendix that compact subsets of $(\mathcal{S},d_{\mathcal S})$ are $\vert \vert \vert \cdot \vert \vert \vert_{\mathcal{S}}$-bounded. Combining these ideas,~we will establish in Proposition~\ref{Paper01_estimates_holding_uniformly_compact_subsets_S_N} that for every $p \in \mathbb{N}$, any compact set $\mathcal{K} \subset \mathcal{S}$, and any $T \geq 0$,
    \begin{equation*}
        \sup_{n \in \mathbb{N}}\; \sup_{\boldsymbol{\eta} \in \mathcal{K}} \; \sup_{t \leq T} \; \mathbb{E}_{\boldsymbol{\eta}}\Big[\vert \vert \vert \eta^{n}(t) \vert \vert \vert_{\mathcal{S}}^{p}\Big] < \infty.
    \end{equation*}
    By the definition of $\vert \vert \vert \cdot \vert \vert \vert_{\mathcal{S}}$ in~\eqref{Paper01_definition_state_space_formal}, we will then conclude that there exists $C_{p,\mathcal{K},T} > 0$, which does not depend on $N$, such that for every $x \in L^{-1}\mathbb{Z}$,
    \begin{equation} \label{Paper01_control_moments_improvement_general_compact_ii}
          \sup_{n \in \mathbb{N}}\; \sup_{\boldsymbol{\eta} \in \mathcal{K}} \; \sup_{t \leq T} \; \mathbb{E}_{\boldsymbol{\eta}}\Big[\vert \vert \eta^{n}(t,x) \vert \vert_{\ell_1}^{p}\Big] \leq C_{p,\mathcal{K},T} (1 + \vert x \vert)^{2p}.
    \end{equation}
    Unlike estimate~\eqref{Paper01_heuristics_correlation_functions_ii}, estimate~\eqref{Paper01_control_moments_improvement_general_compact_ii} gives a bound for general $x\in L^{-1}\mathbb Z$ that is uniform in $n$. This will allow us to establish tightness of the one-dimensional distributions of $\eta^n$ for any initial condition in $\mathcal{S}$.
    
    After establishing tightness, to complete the proof of property~(i), i.e.~to show that the one-dimensional distributions of $(\eta^n)_{n \in \mathbb{N}}$ converge for any initial condition in $\mathcal S$, we must show uniqueness of the limit. The standard argument for proving uniqueness consists of establishing that perturbations in the initial condition that occur far away from the origin in space do not immediately propagate to the origin, or in other words, that ‘infinity does not influence the origin'. This property is usually required for the construction of infinite particle systems with unbounded numbers of particles per site (see for instance~\cite[Equation~IX.1.17]{liggett1985interacting} and~\cite[Equation~13.17]{chen2004markov}) or with unbounded transition rates (see e.g.~\cite[Proof of Proposition~3.5]{zbMATH07976502}).

    The main difficulty in establishing this property for the spatial Muller's ratchet lies in the lack of monotonicity explained in Section~\ref{Paper01_ssec:non-monotonicity}, as well as the fact that there are infinitely many types of particles. Because particles of different types, living in the same deme at the same time, interact with one another in the sense that they influence each other's birth and death rates, for configurations $\boldsymbol{\eta}^{(1)} \neq \boldsymbol{\eta}^{(2)}$, it is difficult to couple two realisations $(\eta^{n,(1)}(t))_{t \geq 0}$ and $(\eta^{n,(2)}(t))_{t \geq 0}$ of the process $\eta^n$, conditioned on $\eta^{n,(1)}(0) = \boldsymbol{{\eta}}^{(1)}$ and $\eta^{n,(2)}(0) = \boldsymbol{{\eta}}^{(2)}$, in such a way that after a period of time $T > 0$, we can establish a useful bound on 
    \[\mathbb{E}_{(\boldsymbol{\eta}^{(1)}, \boldsymbol{\eta}^{(2)})}\Big[d_{\mathcal{S}}(\eta^{n, (1)}(T), \eta^{n, (2)}(T))\Big]
    \]
    as a function of $d_{\mathcal{S}}(\boldsymbol{\eta}^{(1)}, \boldsymbol{\eta}^{(2)})$ uniformly in $n \in \mathbb{N}$. This difficulty 
    is also observed in other interacting particle systems with non-monotone, non-local interactions~\cite{maillard2024branching, etheridge2019genealogical}.

     To overcome this challenge, we construct a coupling in Section~\ref{Paper01_section_spread_infection} by encoding the difference between~$(\eta^{n,(1)}(t))_{t \geq 0}$ and~$(\eta^{n,(2)}(t))_{t \geq 0}$ in terms of \emph{infected} and \emph{partially recovered} particles. Infected and partially recovered particles will represent particles that are alive in only one of the processes, while \emph{susceptible} particles are shared by both $(\eta^{n,(1)}(t))_{t \geq 0}$ and~$(\eta^{n,(2)}(t))_{t \geq 0}$. Importantly, the coupling will be constructed in such a way that the number of partially recovered particles in each deme will be almost surely the same for both processes for all times, which means that only the infected particles will cause differences in the reproduction and death rates of susceptible particles. Under this representation, to prove that ‘infinity does not influence the origin', it suffices to show that the propagation speed in space of infected and partially recovered particles is bounded.

    This will be achieved in Proposition~\ref{Paper01_proposition_bound_spread_infection}. The main difficulty arises when an infected particle enters a region with many susceptible particles. Although Assumption~\ref{Paper01_assumption_polynomials} guarantees that, in such regions, the death rate dominates the birth rate (so infected particles typically die before producing many descendants), a single infected particle may still cause the death of many susceptibles in one process but not in the other.
    
    To control this asymmetry, the coupling is designed so that whenever an infected particle causes the death of a susceptible particle in a high-density region, it is immediately converted into a \emph{partially recovered} particle, which can no longer trigger further one-sided deaths; see Figure~\ref{Paper01_coupling_infection}(b). Moreover, the rate at which partially recovered particles can become reinfected is uniformly bounded. These two features ensure the desired bound on the propagation speed.
    
    Since this coupling is non-standard and, to the best of our knowledge, new, we now briefly outline it heuristically and relate it to existing work on the spread of infections in interacting particle systems. The following discussion is not required for the rigorous proofs in Section~\ref{Paper01_section_spread_infection}.

    \subsubsection{Coupling in terms of susceptible, infected and partially recovered particles}

    To couple two realisations of the process~$(\eta^{n,(1)}(t))_{t \geq 0}$ and~$(\eta^{n,(2)}(t))_{t \geq 0}$ started from different configurations, we give each particle a ‘class' in addition to its type (number of mutations) and spatial location. There will be five different classes of particles:
    \begin{itemize}
    \item Class 0 (\emph{susceptible}) particles;
    \item Class 1 and class 2 (\emph{infected}) particles;
    \item Class $1*$ and class $2*$ (\emph{partially recovered}) particles.
    \end{itemize}
    The class 0 particles are shared by both processes $(\eta^{n,(1)}(t))_{t\ge 0}$ and $(\eta^{n,(2)}(t))_{t\ge 0}$, while the class 1 and class $1*$ particles are exclusive to the process $(\eta^{n,(1)}(t))_{t\ge 0}$, and the class 2 and class $2*$ particles are exclusive to the process $(\eta^{n,(2)}(t))_{t\ge 0}$. Additionally, each partially recovered class $1*$ (respectively, class $2*$) particle is uniquely associated with a class $2*$ (respectively, class $1*$) particle, referred to as its \emph{dual particle}, with the same spatial location (but not necessarily with the same number of mutations). Using this construction, the fact that ‘infinity does not influence the origin' will be proved by establishing an upper bound on the propagation speed of infected and partially recovered particles.

\begin{figure}[p]
\centering
\includegraphics[width=\linewidth,trim=0cm 0cm 0cm 0cm,clip=true]{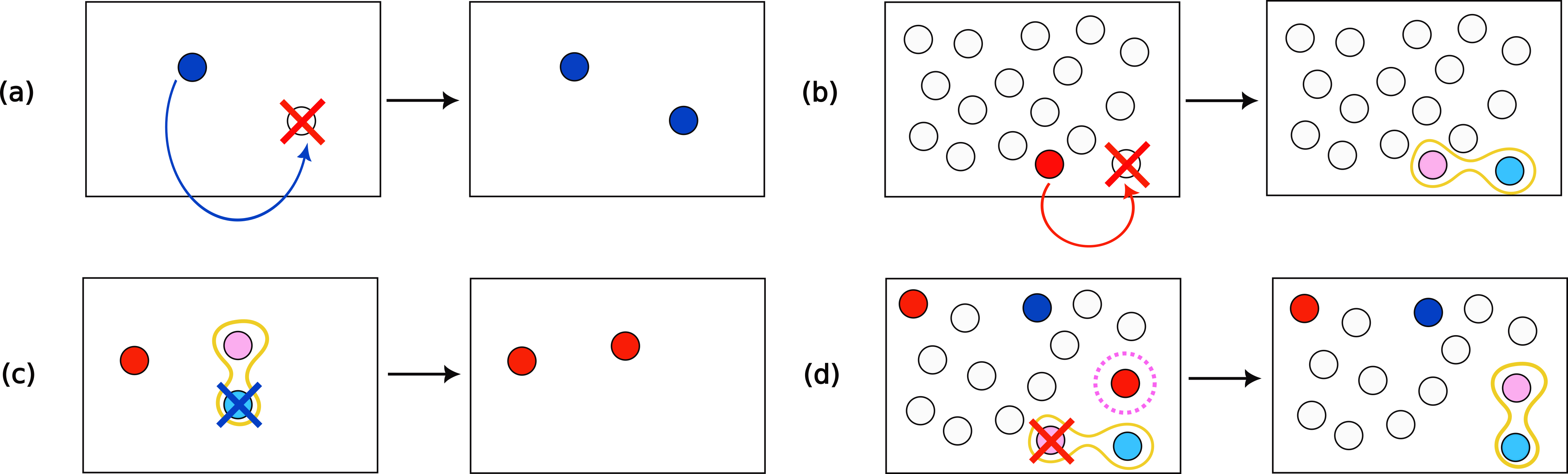}
\caption{\label{Paper01_coupling_infection}~Representation of some of the birth and death events in the coupling described in this section. Different colours represent different particle classes but do not distinguish between types (i.e.~the number of mutations). Class~$0$ particles are shown in white; class~$1$ and class~$1*$ particles are shown in dark and light blue, respectively; and class~$2$ and class~$2*$ particles are shown in red and pink, respectively. Dual pairs (a class~$1*$ and a class~$2*$ particle) are circled in yellow. Death events occurring exclusively in $(\eta^{n,(1)}(t))_{t \geq 0}$ are marked with dark blue crosses, while those occurring exclusively in $(\eta^{n,(2)}(t))_{t \geq 0}$ are marked with red crosses. Arrows indicate which infected particle caused a transmission event, coloured according to the class of the transmitting particle. The figure illustrates four types of events:\\
\textit{(a)~Transmission without recovery}. The presence of an extra class~$1$ particle increases the death rate in $\eta^{n,(2)}$ relative to $\eta^{n,(1)}$, causing a susceptible particle to die in~$\eta^{n,(2)}$, but not in~$\eta^{n,(1)}$. We say that the class~$1$ particle `transmits' the infection to a susceptible particle, which then becomes class~$1$. \\
\textit{(b)~Transmission with partial recovery}. An extra class~$2$ particle increases the death rate in~$\eta^{n,(2)}$ relative to~$\eta^{n,(1)}$. As a result, the class~$2$ particle transmits the infection to a susceptible particle, which becomes class~$1$. However, a `partial recovery' then occurs: the newly infected particle becomes class~$1*$, forming a dual pair with the transmitting particle, which becomes class~$2*$.\\
\textit{(c)~Death of a partially recovered particle with reinfection}. The presence of an extra class~$2$ particle increases the death rate in $\eta^{n,(1)}$ relative to $\eta^{n,(2)}$, causing a class~$1*$ particle to die, but not its class~$2*$ dual pair particle. Since there are no class~$1$ particles at the deme when the death occurs, `reinfection' takes place: the class~$2*$ dual pair particle becomes class~$2$. \\
\textit{(d)~Death of a partially recovered particle with replacement}. Here, an extra class~$2$ particle causes the death of a class~$2*$ particle, but not its class~$1*$ dual pair particle. One of the class~$2$ particles is then chosen uniformly at random (indicated by the pink circle) to become a new class~$2*$ particle, replacing the particle that died in its dual pair.\\
We will show that transmission without recovery and reinfection events only occur at demes where the number of particles is sufficiently small (see Lemma~\ref{Paper01_lem:cutoffconseq}).}
\end{figure}

    We construct the coupling in such a way that particles of classes 0, 1 and 2 migrate independently from each other, while each dual pair of particles of classes $1*$ and $2*$ always migrates together, and independently from other particles and other pairs. To define the particle system in such a way that $\eta^{n,(1)}$ and $\eta^{n,(2)}$ are both realisations of $\eta^n$, the following birth and death events are allowed to occur (see Figure~\ref{Paper01_coupling_infection}):
    \begin{itemize}
        \item Reproduction of partially recovered or infected particles, producing infected particles.
        \item `Uninfected' reproduction of susceptible particles, producing susceptible particles (corresponding to reproduction events that happen in both processes $\eta^{n,(1)}$ and $\eta^{n,(2)}$).
        \item `Infected' reproduction of susceptible particles, producing infected particles (corresponding to reproduction events only happening in one of the processes $\eta^{n,(1)}$ and $\eta^{n,(2)}$ due to a difference in the number of class~$1$ and class~$2$ particles).
        \item Death of infected particles.
        \item Death of susceptible particles that occur in both processes.
        \item `Transmission' events: these events correspond to death events of susceptible particles that happen in only one of the processes, due to a difference in the number of class~$1$ and class~$2$ particles. They are classified into events with or without `partial recovery' (see Figure~\ref{Paper01_coupling_infection}(a) and~(b)).
        \item {Death of dual pairs of particles}: These events correspond to the simultaneous death of both the class~$1*$ particle and the class~$2*$ particle in a dual pair.
        \item {Death of partially recovered particles with `reinfection' or `replacement'}: These events correspond to death events that occur to only one of the partially recovered particles of a dual pair, due to a difference in the number of class~$1$ and class~$2$ particles (see Figure~\ref{Paper01_coupling_infection}(c) and~(d)). We construct the coupling in such a way that for each $i \in \{1,2\}$, when a death event occurs to a  class~$i*$ particle due to a higher number of class $i$ particles than class $3-i$ particles, then we choose one of the class~$i$ particles to replace the $i*$ particle in the dual pair. 
    \end{itemize}

To bound the propagation speed of infected and partially recovered particles, we will adapt an argument used by Kesten and Sidoravicius to bound the spread of infection in a moving population~\cite[Theorem~1]{kesten2005spread}. The main difference between the model studied by Kesten and Sidoravicius and our interacting particle system is that in~\cite{kesten2005spread} there is no reproduction of particles, while in our model particles reproduce and die at rates that depend on the local number of particles. Hence, their result is not directly applicable to the coupling introduced in this section. In fact, although there are many works on the spread of infection in interacting particle systems (see~\cite{stauffer2022} and the references therein), we are not aware of any results concerning the spread of infection in spatial birth-death systems with unbounded reproduction  rates and without uniform bounds on the number of particles per site.

In~\cite[Theorem~1]{kesten2005spread}, Kesten and Sidoravicius associate to each infected particle a Poisson clock, in such a way that at each ‘ring' of the clock, the particle either jumps to another deme or transmits the infection. This construction allows them to apply a first moment method and bound the propagation speed of infected particles. Baldasso and Stauffer use a similar argument in~\cite[Proposition~3.1]{baldasso2022local} to bound the spread of infection in a population moving according to a biased random walk. Assuming the original infected particle is located at the origin at time~$0$, in both~\cite[Theorem~1]{kesten2005spread} and~\cite[Proposition~3.1]{baldasso2022local}, the authors bound the spreading speed of infection by observing that, for $t \geq 0$, every site that is infected at time $t$ can be reached by a concatenation of trajectories of random walks. Since, by using the strong Markov property, the concatenation of random walks is also a random walk, the upper bound on the speed follows by bounding the number of jumps a random walk performs over finite periods of time.

To adapt this idea to our setting, we must overcome the fact that, in our case, the rate at which an infected or partially recovered particle generates another infected particle, either by reproducing, by inducing an `infected' reproduction event for a susceptible particle, or by `transmitting' the infection, is not deterministically bounded. 
We are still able to control the number of infected and partially recovered particles generated, for the following reasons.
In demes where the number of particles is very large, infected or partially recovered particles die at a much higher rate than the rate at which they cause reproduction events, due to Assumption~\ref{Paper01_assumption_polynomials}. Transmission events can also happen with high rate in such demes; however, when an infected particle transmits the infection at a deme with a high number of particles, the particle becomes partially recovered and can no longer transmit the infection (see Figure~\ref{Paper01_coupling_infection}). 
Moreover, partially recovered particles can only be reinfected when at demes with a low number of particles (see Figure~\ref{Paper01_coupling_infection}).
At demes with a low number of particles, the total rate of birth and death events is bounded, which implies a bound on the rate at which new infected particles are generated or partially recovered particles are reinfected at such demes.
By combining these ideas, in Lemma~\ref{Paper01_lem:boundPE} and Corollary~\ref{Paper01_corollary_almost_final_bound_spread_infection} we will be able to control the number of infected or partially recovered particles that are generated over finite periods of time. We will then use similar arguments to those in~\cite[Theorem~1]{kesten2005spread} and~\cite[Proposition~3.1]{baldasso2022local}, to prove in Proposition~\ref{Paper01_proposition_bound_spread_infection} that `the origin is not influenced by infinity'. This will complete the proof of uniqueness of the limit of the one-dimensional distributions of $\eta^n$ and establish property~(i).

\subsection{Glossary} \label{Paper01_subsection_notation}

Here we list frequently used notation. 
Recall that we also introduced some standard general notation at the end of Section~\ref{Paper01_introduction}.
In the second column of the table we give a brief heuristic description, and in the third column we refer to the section or equation where the notation is introduced.

\setlength{\extrarowheight}{1pt}   
\renewcommand{\arraystretch}{1.1}  

\begin{tabularx}{\textwidth}{@{}lXr@{}}
\toprule
\textbf{Notation} & \textbf{Meaning} & \textbf{Defn./Sect.} \\
\midrule
\endfirsthead

\toprule
\textbf{Notation} & \textbf{Meaning} & \textbf{Defn./Sect.} \\
\midrule
\endhead

\midrule
\multicolumn{3}{r}{\emph{(continued on next page)}} \\
\midrule
\endfoot

\bottomrule
\endlastfoot

$L$ & space renormalisation parameter & Section~\ref{Paper01_model_description} \\
$m$ & migration rate & Section~\ref{Paper01_model_description} \\
$(s_k)_{k \in \mathbb{N}_{0}}$ & sequence of fitness parameters & Section~\ref{Paper01_model_description} \\
$N$ & scaling parameter for carrying capacity & Section~\ref{Paper01_model_description} \\
$\mu$ & mutation probability & Section~\ref{Paper01_model_description} \\
$q_+$ & birth rate polynomial & Section~\ref{Paper01_model_description} \\
$q_+^N$ & scaled birth rate polynomial & \eqref{Paper01_scaled_polynomials_carrying_capacity} \\
$q_-$ & death rate polynomial & Section~\ref{Paper01_model_description} \\
$q_-^N$ & scaled death rate polynomial & \eqref{Paper01_scaled_polynomials_carrying_capacity} \\ $  \bar{q}^{N}$ & upper bound on the total birth and death rates for the population in a given deme & \eqref{Paper01_bar_polynomial} \\
$F^b_k$ & total rate of birth of particles carrying $k$ mutations & \eqref{Paper01_discrete_birth_rate} \\
$F^d_k$ & total rate of death of particles carrying $k$ mutations & \eqref{Paper01_discrete_death_rate} \\
$\mathcal{S}$ & space of configurations; state space for spatial Muller's ratchet & \eqref{Paper01_definition_state_space_formal} \\
$\vert \vert \vert \cdot \vert \vert \vert_{\mathcal{S}}$ & total mass function on configurations in $\mathcal S$ & \eqref{Paper01_definition_state_space_formal} \\ $d_{\mathcal{S}}$ & metric on $\mathcal{S}$ & \eqref{Paper01_definition_semi_metric_state_space} \\
$\mathcal{S}_{0}$ & set of possible initial configurations & \eqref{Paper01_set_initial_configurations} \\
$\eta_k(x)$ & number of particles in deme $x$ carrying $k$ mutations in configuration $\boldsymbol{\eta}=(\eta_j(y))_{j\in \mathbb N_0,\, y\in L^{-1}\mathbb Z}\in \mathcal S$ & Section~\ref{Paper01_model_description} \\
$\Big((\eta^{n}(t))_{t \geq 0}\Big)_{n \in \mathbb{N}}$ & sequence of Feller processes that approximates the spatial Muller's ratchet & before \eqref{Paper01_generator_foutel_etheridge_model_restriction_n} \\
$\lambda_n$ & radius of the spatial region outside of which particles in $(\eta^{n}(t))_{t\ge 0}$ are frozen & \eqref{Paper01_scaling_parameters_for_discrete_approximation} \\
$\Lambda_n$ & spatial region outside of which particles in $(\eta^{n}(t))_{t\ge 0}$ are frozen & after~\eqref{Paper01_scaling_parameters_for_discrete_approximation} \\
$K_n$ & only particles carrying at most $K_n$ mutations can reproduce in $(\eta^{n}(t))_{t \geq 0}$ & \eqref{Paper01_scaling_parameters_for_discrete_approximation} \\
$\mathcal{L}^{n}$ & infinitesimal generator of $(\eta^{n}(t))_{t \geq 0}$ & \eqref{Paper01_generator_foutel_etheridge_model_restriction_n} \\
$\mathcal{L}^{n}_{m}$ & migration part of $\mathcal{L}^{n}$ & \eqref{Paper01_infinitesimal_generator_restriction_n} \\
$\mathcal{L}^{n}_{r}$ & reaction part of $\mathcal{L}^{n}$ & \eqref{Paper01_infinitesimal_generator_restriction_n} \\ $\mathcal{L}$ & infinitesimal generator of the spatial Muller's ratchet & \eqref{Paper01_generator_foutel_etheridge_model} \\ 
$\mathcal{L}_{m}$ & migration part of $\mathcal{L}$ & \eqref{Paper01_infinitesimal_generator} \\
$\mathcal{L}_{r}$ & reaction part of $\mathcal{L}$ & \eqref{Paper01_infinitesimal_generator} \\
$\boldsymbol{e}^{(x)}_{k}$ & particle configuration consisting of only one particle carrying exactly $k$ mutations in deme $x$ & before \eqref{Paper01_scaled_polynomials_carrying_capacity} \\ 
$\{P^n_t\}_{t \geq 0}$ & semigroup associated with $(\eta^{n}(t))_{t \geq 0}$
& after~\eqref{Paper01_infinitesimal_generator_restriction_n} \\
$\mathscr{C}_b^{\textrm{cyl}}(\mathcal{S}, \mathbb{R})$ & set of bounded cylindrical functions & \eqref{Paper01_general_definition_cylindrical_function} \\ $\mathscr{C}_*(\mathcal{S}, \mathbb{R})$ & set of continuous functions whose discrete derivatives decay polynomially fast in space & \eqref{Paper01_definition_Lipschitz_function} \\
$\Big\{\mathcal{F}^{\eta}_{t}\Big\}_{t \geq 0}$ & natural filtration of the stochastic process $(\eta(t))_{t \geq 0}$ & before~\eqref{eq:rightctsfiltration} \\ 
$\Big\{\mathcal{F}^{\eta}_{t+}\Big\}_{t \geq 0}$ & right-continuous filtration of the stochastic process $(\eta(t))_{t \geq 0}$ & before~\eqref{eq:rightctsfiltration} \\ 
\end{tabularx}

\section{Construction and characterisation of Feller semigroups on non-locally compact Polish spaces}
\label{Paper01_section_functional_analysis_feller_non_locally_compact}

In this section, we will state and prove some general results regarding the construction of Feller processes that have the strong Markov property taking values in Polish spaces that are not locally compact. Although, as discussed in Section~\ref{Paper01_subsection_heuristics}, we will subsequently use these results in the specific context of the spatial Muller's ratchet, in this section we will state them in a more general context that should also be applicable to other non-monotone interacting particle systems. As explained in Section~\ref{Paper01_subsection_heuristics}, the main focus in this section will be the construction of the semigroup generated by the transition probabilities on appropriate function spaces.

Let $(\mathcal{S}, d_{\mathcal{S}})$ be a complete and separable metric space, and let $\mathcal{S}_{0} \subseteq \mathcal{S}$. Let $\mathcal{D}(\mathcal{L}) \subseteq \mathscr{C}_{b}(\mathcal{S}, \mathbb{R})$ be a vector subspace of the space of continuous bounded functions from $\mathcal{S}$ to $\mathbb{R}$, and let $\mathcal{L}: \mathcal{D}(\mathcal{L}) \rightarrow \mathscr{C}(\mathcal{S}, \mathbb{R})$ be an unbounded linear operator. We will refer to $\mathcal{D}(\mathcal{L})$ as the domain of $\mathcal{L}$. Heuristically, $\mathcal{L}$ should be thought of as the infinitesimal generator of a non-monotone reaction-diffusion particle system defined on some unbounded spatial domain. We will show that under certain assumptions, $\mathcal{L}$ ‘generates' a Feller semigroup associated to this particle system (see Section~\ref{Paper01_subsection_appendix_feller_non_locally_compact_polish_space} for a brief overview of Feller semigroups on non-locally compact Polish spaces), and the particle system is a càdlàg (strong) Markov process. Our first assumption is that we can characterise the action of $\mathcal{L}$ on a sufficiently large (and regular) subset of $\mathscr{C}_{b}(\mathcal{S}, \mathbb{R})$.

\begin{assumption}[Set of test functions] \label{Paper01_assumption_candidate_set_functions_domain_generator}
    There exists $\mathfrak{A} \subseteq \mathcal{D}(\mathcal{L})$ such that
    \begin{enumerate}[(i)]
        \item For any $\phi \in \mathfrak{A}$ and $a > 0$, the map $\phi_{a} \in \mathscr{C}_{b}(\mathcal{S}, \mathbb{R})$ given by, for any $\boldsymbol{\eta} \in \mathcal{S}$,
        \begin{equation*}
            \phi_{a}(\boldsymbol{\eta}) \defeq \sgn(\phi(\boldsymbol{\eta})) (\vert \phi(\boldsymbol{\eta})\vert \wedge a)
        \end{equation*}
        is also an element of $\mathfrak{A}$;
        \item The constant function $\boldsymbol{1}: \mathcal{S} \rightarrow \{1\}$ is an element of $\mathfrak{A}$;
        \item $\mathfrak{A}$ is an algebra over $\mathbb{R}$, i.e.~for any $\phi, \psi \in \mathfrak{A}$ and $a,b \in \mathbb{R}$, $a\phi + b\psi$ and $\phi \psi$ are also elements of $\mathfrak{A}$;
        \item $\mathfrak{A}$ separates the points of $\mathcal{S}$, i.e.~for any $\boldsymbol{\eta}, \boldsymbol{\xi} \in \mathcal{S}$ such that $\boldsymbol{\eta} \neq \boldsymbol{\xi}$, there exists $\phi \in \mathfrak{A}$ such that $\phi(\boldsymbol{\eta}) \neq \phi(\boldsymbol{\xi})$.
    \end{enumerate}
\end{assumption}

By an application of the Stone-Weierstrass theorem (see for e.g.~\cite[Section~4.7]{folland1999real}) and conditions (ii) to (iv) in Assumption~\ref{Paper01_assumption_candidate_set_functions_domain_generator}, for any compact set $\mathcal{K} \subset \mathcal{S}$, $\mathfrak{A}$ is a dense subset of $\mathscr{C}(\mathcal{K}, \mathbb{R})$ in the topology of uniform convergence on $\mathcal{K}$. For ease of reference, we state this result below without proof.

\begin{theorem}[The Stone-Weierstrass theorem] \label{Paper01_corollary_stone_weierstrass_theorem}
    Suppose $\mathfrak{A}$ satisfies Assumption~\ref{Paper01_assumption_candidate_set_functions_domain_generator}.
    For any $\Phi \in \mathscr{C}_{b}(\mathcal{S}, \mathbb{R})$, any compact subset $\mathcal{K} \subset \mathcal{S}$ and any $\varepsilon > 0$, there exists $\phi = \phi(\Phi, \mathscr{K}, \varepsilon) \in \mathfrak{A}$ such that
    \begin{equation*}
        \sup_{\boldsymbol{\eta} \in \mathcal{K}} \; \vert \Phi(\boldsymbol{\eta}) - \phi(\boldsymbol{\eta}) \vert \leq \varepsilon.
    \end{equation*}
\end{theorem}

Test functions commonly used in the construction of interacting particle systems -- such as bounded globally Lipschitz functions, cylindrical functions (in the case where the particle system lives on the lattice $\mathbb{Z}^d$), and integrals against smooth, compactly supported test functions (when the system is defined on $\mathbb{R}^d$) -- satisfy Assumption~\ref{Paper01_assumption_candidate_set_functions_domain_generator}~\cite{liggett1985interacting,chen2004markov,etheridge2023looking}.

 To allow us to prove the semigroup property, we will need a sequence of Markov processes that approximates the particle system that we wish to construct. This will be our next assumption. Recall that $\mathcal{S}_{0} \subseteq \mathcal{S}$ is a fixed subset of $\mathcal{S}$. Recall the definition of $\mathbb{P}_{\boldsymbol{\eta}}$ and $\mathbb{E}_{\boldsymbol{\eta}}$ introduced at the end of Section~\ref{Paper01_introduction}.

\begin{assumption}[Tight sequence of processes] \label{Paper01_assumption_tightness_sequence_processes}
    For each $n \in \mathbb{N}$, let $\mathcal{L}^{n}:\mathcal D(\mathcal{L}^{n})\to \mathscr{C}(\mathcal S,\mathbb R)$ be the infinitesimal generator of a càdlàg strong Markov process $(\eta^{n}(t))_{t \geq 0}$ taking values in $\mathcal{S}$, with corresponding Feller semigroup $\{P^{n}_t\}_{t \geq 0}$. The sequence $(\eta^{n})_{n \in \mathbb{N}}$ is assumed to satisfy the following conditions:
    \begin{enumerate}[(i)]
        \item For every $n \in \mathbb{N}$, $\mathfrak{A} \subseteq \mathcal{D}(\mathcal{L}^{n})$. Moreover, for any $\phi \in \mathfrak{A}$ and $\boldsymbol{\eta} \in \mathcal{S}$,
        \begin{equation*}
            \lim_{t \rightarrow 0^{+}} \; \frac{(P^{n}_{t}\phi)(\boldsymbol{\eta}) - \phi(\boldsymbol{\eta})}{t} = (\mathcal{L}^{n}\phi)(\boldsymbol{\eta}) \quad \textrm{and} \quad \lim_{n \rightarrow \infty} \; (\mathcal{L}^{n}\phi)(\boldsymbol{\eta}) = (\mathcal{L}\phi)(\boldsymbol{\eta}) ;
        \end{equation*}
        \item For any $\boldsymbol{\eta} \in \mathcal{S}_{0}$, conditioning on $\eta^{n}(0) = \boldsymbol{\eta}$ for every $n \in \mathbb{N}$, the sequence of processes $\left((\eta^{n}(t))_{t \geq 0}\right)_{n \in \mathbb{N}}$ is tight in the space of càdlàg functions $\mathscr{D}([0,\infty), \mathcal{S})$;
        \item For any $T \geq 0$, any compact subset $\mathcal{K} \subset \mathcal{S}$ and any $\varepsilon > 0$, there exists a compact subset $\mathscr{K}' \defeq \mathscr{K}'(T, \mathcal{K}, \varepsilon) \subset \mathcal{S}$ such that
        \begin{equation*}
            \inf_{\boldsymbol{\eta} \in \mathcal{K}} \; \inf_{n \in \mathbb{N}} \; \inf_{t \in [0,T]} \; \mathbb{P}_{\boldsymbol{\eta}}\left(\eta^{n}(t) \in \mathscr{K}'\right) \geq 1 - \varepsilon;
        \end{equation*}
        \item For any $\phi \in \mathfrak{A}$, $\boldsymbol{\eta} \in \mathcal{S}$ and $T \geq 0$, the limit  $\lim_{n \rightarrow \infty} (P^{n}_{T}\phi)(\boldsymbol{\eta})$ is well defined. Moreover, the convergence is uniform on compact time intervals and compact subsets of $\mathcal{S}$, i.e.~for any $T \geq 0$, any $\phi \in \mathfrak{A}$, any compact subset $\mathcal{K} \subset \mathcal{S}$ and any $\varepsilon > 0$, there exists $n(T, \phi, \mathcal{K}, \varepsilon) \in \mathbb{N}$ such that
        \begin{equation*}
            \sup_{n_{1}, n_{2} \geq n(T, \phi, \mathcal{K}, \varepsilon)} \; \sup_{\boldsymbol{\eta} \in \mathcal{K}} \; \sup_{t \in [0,T]} \; \left\vert (P^{n_{1}}_{t}\phi)(\boldsymbol{\eta}) - (P^{n_{2}}_{t}\phi)(\boldsymbol{\eta}) \right\vert \leq \varepsilon.
        \end{equation*}
        \ignore{
        \item For any $\phi \in \mathfrak{A}$, any $T \geq 0$, any compact subset $\mathcal{K} \subset \mathcal{S}$ and any $\varepsilon > 0$, there exists $\delta \defeq \delta(\varepsilon, \phi, \mathcal{K}, T) > 0$ such that if $\boldsymbol{\eta}, \boldsymbol{\xi} \in \mathcal{K}$ satisfy $d_{\mathcal{S}}(\boldsymbol{\eta}, \boldsymbol{\xi}) \leq \delta$, then
        \begin{equation*}
            \sup_{n \in \mathbb{N}} \, \sup_{t \in [0,T]} \, \vert (P^{n}_{t}\phi)(\boldsymbol{\eta}) - (P^{n}_{t}\phi)(\boldsymbol{\xi})\vert \leq \varepsilon. 
        \end{equation*}}
    \end{enumerate}
\end{assumption}

We now show that under Assumptions~\ref{Paper01_assumption_candidate_set_functions_domain_generator} and~\ref{Paper01_assumption_tightness_sequence_processes}, we can define a semigroup of bounded linear operators on $\mathscr{C}_{b}(\mathcal{S}, \mathbb{R})$ corresponding to a càdlàg strong Markov process with initial configuration in $\mathcal{S}_{0}$. Recall the definition of the strong Markov property from~\eqref{Paper01_strong_markov_property}, and observe that we do not require the filtration to be complete to define the strong Markov property.

\begin{theorem} \label{Paper01_general_thm_convergence_feller_non_locally_compact}
Suppose $\mathcal L:\mathcal D(\mathcal L)\to \mathscr{C}(\mathcal S,\mathbb R)$, $\mathfrak{A} \subseteq \mathcal D(\mathcal L)$, $(\mathcal L^n)_{n\in \mathbb N}$, $(\eta^n)_{n\in \mathbb N}$ and $(\{P^n_t\}_{t\ge 0})_{n\in \mathbb N}$ satisfy
    Assumptions~\ref{Paper01_assumption_candidate_set_functions_domain_generator} and~\ref{Paper01_assumption_tightness_sequence_processes}. Then
    \begin{enumerate}[(i)]
        \item For any $\boldsymbol{\eta} \in \mathcal{S}$ and $T \geq 0$, there is a unique $\mathcal{S}$-valued random variable $\eta(T)$ such that, conditioning on $\eta^{n}(0) = \boldsymbol{\eta}$ for every $n \in \mathbb{N}$, $\eta^{n}(T) \xrightarrow[]{\mathcal{D}} \eta(T)$ in $\mathcal{S}$ as $n \rightarrow \infty$. We denote the probability measure associated to the law of $\eta(T)$ by $\mathbb{P}_{\boldsymbol{\eta}}$ and the corresponding expectation operator by $\mathbb{E}_{\boldsymbol{\eta}}$.
        \item For any $\phi \in \mathscr{C}_{b}(\mathcal{S}, \mathbb{R})$ and $T \geq 0$, define $P_{T}\phi: \mathcal{S} \rightarrow \mathbb{R}$ by setting, for any $\boldsymbol{\eta} \in \mathcal{S}$,
    \begin{equation*}
        (P_{T}\phi)(\boldsymbol{\eta}) \defeq \mathbb{E}_{\boldsymbol{\eta}}\left[\phi\left(\eta(T)\right)\right] = \lim_{n \rightarrow \infty} \; (P^{n}_{T}\phi)(\boldsymbol{\eta}).
    \end{equation*}
     Then $P_{T}\phi \in \mathscr{C}_{b}(\mathcal{S}, \mathbb{R})$ for any $\phi \in \mathfrak{A}$ and $T \geq 0$.
      \item For any non-empty compact subset $\mathcal{K} \subset \mathcal{S}$, $\phi \in \mathscr{C}_{b}(\mathcal{S}, \mathbb{R})$ and $T \geq 0$,
     \begin{equation*}
         \lim_{n \rightarrow \infty} \; \sup_{\boldsymbol{\eta} \in \mathcal{K}} \; \sup_{t \in [0,T]} \; \Big\vert (P_{t}\phi)(\boldsymbol{\eta}) - (P^{n}_{t}\phi)(\boldsymbol{\eta})\Big\vert = 0.
     \end{equation*}
     \item For any $\phi \in \mathscr{C}_{b}(\mathcal{S}, \mathbb{R})$, the map $(T, \boldsymbol{\eta}) \mapsto (P_{T}\phi)(\boldsymbol{\eta})$ for $(T, \boldsymbol{\eta}) \in [0, \infty) \times \mathcal{S}$ is continuous with respect to the product topology on $[0, \infty) \times \mathcal{S}$. In particular, for any $\phi \in \mathscr{C}_{b}(\mathcal{S}, \mathbb{R})$ and $T \geq 0$, $P_{T}\phi \in \mathscr{C}_{b}(\mathcal{S}, \mathbb{R})$.
     \item For any $T \geq 0$, $P_{T}: \mathscr{C}_{b}(\mathcal{S}, \mathbb{R}) \rightarrow \mathscr{C}_{b}(\mathcal{S}, \mathbb{R})$ defined as in assertion~$(iv)$ is a bounded linear operator with $\vert \vert P_{T} \vert \vert_{\textrm{op}} \leq 1$, where $\vert \vert \cdot \vert \vert_{\textrm{op}}$ indicates the operator norm. Moreover, $\{P_{t}\}_{t \geq 0}$ is a semigroup, i.e.~$P_{T}P_{t} = P_{T+t}$ for any $T,t \geq 0$.
     \item For any $\boldsymbol{\eta} \in \mathcal{S}_{0}$, conditioning on $\eta^{n}(0) = \boldsymbol{\eta}$ for every $n \in \mathbb{N}$, the sequence of processes $(\eta^{n}(t))_{t \geq 0}$ converges weakly in $\mathscr{D}([0, \infty), \mathcal{S})$ to a càdlàg Markov process $(\eta(t))_{t \geq 0}$ with $\eta(0) = \boldsymbol{\eta}$ almost surely, and whose transition probabilities correspond to the Feller semigroup $\{P_{t}\}_{t \geq 0}$, i.e.~such that for any $I \in \mathbb{N}$, $(\phi_{i})_{i = 1}^{I} \subset \mathscr{C}_{b}(\mathcal{S}, \mathbb{R})$ and $0 \leq t_{1} \leq t_{2} \leq \cdots \leq t_{I}$,
    \begin{equation*}
    \begin{aligned}
        & \mathbb{E}_{\boldsymbol{\eta}}\Big[\phi_{1}(\eta(t_{1}))\phi_{2}(\eta(t_{2})) \cdots \phi_{I}(\eta(t_{I}))\Big] \\ & \quad = P_{t_{1}}\Big(\phi_{1}P_{(t_2-t_1)}\Big(\phi_{2}P_{(t_3-t_2)}\Big(\cdots \Big(\phi_{I-1}P_{(t_{I} - t_{I-1})}\phi_{I}\Big) \cdots\Big)\Big)\Big)(\boldsymbol{\eta}).
    \end{aligned}
    \end{equation*}

     \item For any $\phi \in \mathfrak{A}$ and $\boldsymbol{\eta} \in \mathcal{S}$,
        \begin{equation*}
            \lim_{t \rightarrow 0^{+}} \; \frac{(P_{t}\phi)(\boldsymbol{\eta}) - \phi(\boldsymbol{\eta})}{t} = (\mathcal{L}\phi)(\boldsymbol{\eta}).
        \end{equation*}
     
     \item For $\boldsymbol{\eta} \in \mathcal{S}_{0}$, let $(\eta(t))_{t \geq 0}$ be the càdlàg Markov process with $\eta(0) = \boldsymbol{\eta}$ defined in~$(vi)$. Then $(\eta(t))_{t \geq 0}$ is a strong Markov process with respect to the filtration~$\{\mathcal{F}_{t+}^{\eta}\}_{t \geq 0}$.
    \end{enumerate}
\end{theorem}

We highlight that assertions~(i) - (vii) can be understood as a generalisation of existence theorems for some classes of infinite interacting particle systems as in~\cite[Theorem~IX.1.14]{liggett1985interacting} and~\cite[Theorem~13.8]{chen2004markov}. Moreover, assertion~(viii) is a generalisation of the strong Markov property of Feller processes taking values in locally compact state spaces (see for instance~\cite[Theorem~4.2.5]{ethier2009markov}).

\begin{proof}
    We will divide the proof into steps corresponding to each of the assertions (i) - (viii).

    \medskip

    \noindent \underline{Proof of assertion~(i):}

    \medskip

    By Assumption~\ref{Paper01_assumption_tightness_sequence_processes}(iii), for any $T \geq 0$, the sequence of $\mathcal{S}$-valued random variables $(\eta^{n}(T))_{n \in \mathbb{N}}$ is tight. Let $\eta^{(1)}(T)$ and $\eta^{(2)}(T)$ be any two weak limits of subsequences of $(\eta^{n}(T))_{n \in \mathbb{N}}$. By Assumption~\ref{Paper01_assumption_tightness_sequence_processes}(iv), for any $\phi \in \mathfrak{A}$, we have
    \begin{equation} \label{Paper01_separating_class_functions_identity}
        \mathbb{E}_{\boldsymbol{\eta}}\left[\phi(\eta^{(1)}(T))\right] = \mathbb{E}_{\boldsymbol{\eta}}\left[\phi(\eta^{(2)}(T))\right]. 
    \end{equation}
    Since by Assumption~\ref{Paper01_assumption_candidate_set_functions_domain_generator}(iv), $\mathfrak{A}$ separates points of $\mathcal{S}$, we conclude that $\mathfrak{A}$ is a separating class of continuous bounded functions (see~\cite[Theorem~3.4.5]{ethier2009markov}), and therefore~\eqref{Paper01_separating_class_functions_identity} implies that $\eta^{(1)}(T)$ and $\eta^{(2)}(T)$ have the same distribution, as desired.

    \medskip

    \noindent \underline{Proof of assertion~(ii):}

    \medskip
    
    We first observe that by assertion~(i), we have $\lim_{n \rightarrow \infty} \, (P^{n}_{T}\phi)(\boldsymbol{\eta}) = \mathbb{E}_{\boldsymbol{\eta}}[\phi(\eta(T))]$ for any $\phi \in \mathscr{C}_{b}(\mathcal{S}, \mathbb{R})$ and $T \geq 0$. We now prove that for any $T \geq 0$ and any $\phi \in \mathfrak{A}$, $P_{T}\phi \in \mathscr{C}_{b}(\mathcal{S}, \mathbb{R})$. Let $T \geq 0$ and $\phi \in \mathfrak{A}$ be fixed. By Assumption~\ref{Paper01_assumption_tightness_sequence_processes}(iv), the sequence $(P^{n}_{T}\phi)_{n \in \mathbb{N}}$ converges to $P_{T}\phi$ in the topology of uniform convergence on compact sets. Since $(\mathcal{S}, d_{\mathcal{S}})$ is a metric space, by standard results in topology (see for instance~\cite[Theorem~43.11]{willard2012general}), to conclude that $P_{T}\phi$ is a continuous function from $\mathcal{S}$ to $\mathbb{R}$, it suffices to establish that for every $n \in \mathbb{N}$,
    \begin{equation} \label{Paper01_continuity_criteria_compact_open_topology}
        P_{T}^{n}\phi \in \mathscr{C}(\mathcal{S}, \mathbb{R}).
    \end{equation}
    By Assumption~\ref{Paper01_assumption_tightness_sequence_processes}, for every $n \in \mathbb{N}$, $\{P^{n}_{t}\}_{t \geq 0}$ is a Feller semigroup (see Definition~\ref{Paper01_definition_feller_non_locally_compact}), and therefore~\eqref{Paper01_continuity_criteria_compact_open_topology} holds. This shows that $P_{T}\phi \in \mathscr{C}(\mathcal{S}, \mathbb{R})$. Finally, since $\phi \in \mathfrak{A}$ implies that $\phi$ is bounded, and hence that $P^{n}_{T}\phi$ is uniformly bounded in $n \in \mathbb{N}$ and $T \geq 0$, we conclude that $P_{T}\phi \in \mathscr{C}_{b}(\mathcal{S}, \mathbb{R})$, which completes the proof.

    \medskip

     \noindent \underline{Proof of assertion~(iii):}
    
    \medskip

     Observe that for $\phi \in \mathfrak{A}$, the desired limit follows directly from  Assumption~\ref{Paper01_assumption_tightness_sequence_processes}(iv). To tackle the case of general $\phi \in \mathscr{C}_{b}(\mathcal{S}, \mathbb{R})$, let $(\varepsilon_{j})_{j \in \mathbb{N}}$ be a sequence of strictly positive real numbers such that $\lim_{j \rightarrow \infty} \, \varepsilon_{j} = 0$, and take a compact subset $\mathcal{K} \subset \mathcal{S}$ and $T \geq 0$. By Assumption~\ref{Paper01_assumption_tightness_sequence_processes}(iii), for any $j \in \mathbb{N}$, there exists a compact set $\mathscr{K}_{j} \defeq \mathscr{K}_{j}(T, \mathcal{K}, j) \subset \mathcal{S}$ such that
     \begin{equation} \label{Paper01_intermediate_step_high_probability_sequence_compacts}
         \inf_{n \in \mathbb{N}} \; \inf_{\boldsymbol{\eta} \in \mathcal{K}} \; \inf_{t \in [0,T]} \; \mathbb{P}_{\boldsymbol{\eta}} \Big(\eta^{n}(t) \in \mathscr{K}_{j}\Big) \geq 1 - \varepsilon_{j},
     \end{equation}
     which implies by assertion~(i) of this theorem that we also have
     \begin{equation}
     \label{Paper01_compact_containment_condition_general_argument}
     \inf_{\boldsymbol{\eta} \in \mathcal{K}} \; \inf_{t \in [0,T]} \; \mathbb{P}_{\boldsymbol{\eta}} \Big(\eta(t) \in \mathscr{K}_{j}\Big) \geq 1 - \varepsilon_{j}.
     \end{equation}
     By Theorem~\ref{Paper01_corollary_stone_weierstrass_theorem}, there exists a sequence $(\phi^{(j)})_{j \in \mathbb{N}}$ in~$\mathfrak{A}$ such that for each $j \in \mathbb{N}$,
     \begin{equation}
     \label{Paper01_application_stone_weierstrass}
         \sup_{\boldsymbol{\xi} \in \mathscr{K}_{j}} \; \left\vert \phi(\boldsymbol{\xi}) - \phi^{(j)}(\boldsymbol{\xi})\right\vert \leq \varepsilon_{j}.
     \end{equation}
      Since $\phi$ is uniformly bounded, by using Assumption~\ref{Paper01_assumption_candidate_set_functions_domain_generator}(i) and by replacing $\phi^{(j)}$ with $\phi^{(j)}_{\vert \vert \phi \vert \vert_{L_{\infty}(\mathcal{S}; \mathbb{R})}}$ if necessary, we can assume that
    \begin{equation} \label{Paper01_uniform_boundedness_approximating_sequence_maps}
        \sup_{j \in \mathbb{N}} \;  \|\phi^{(j)} \|_{L_{\infty}(\mathcal{S}; \mathbb{R})} \leq  \|\phi \|_{L_{\infty}(\mathcal{S}; \mathbb{R})}.
    \end{equation}
     Thus, for any $\boldsymbol{\eta} \in \mathcal{K} $, $n \in \mathbb{N}$, $t \in [0,T]$ and $j \in \mathbb{N}$, we have by the triangle inequality
     \begin{equation} \label{Paper01_intermediate_estimate_uniform_convergence_semigroups_in_n_over_compacts}
     \begin{aligned}
         & \left\vert (P_{t}\phi)\left(\boldsymbol{\eta}\right) - \left(P^{n}_{t}\phi\right)\left(\boldsymbol{\eta}\right)\right\vert \\ & \quad \leq \left\vert (P_{t}\phi)\left(\boldsymbol{\eta}\right) - (P_{t}\phi^{(j)})\left(\boldsymbol{\eta}\right)\right\vert + \left\vert (P^{n}_{t}\phi^{(j)})\left(\boldsymbol{\eta}\right) - \left(P^{n}_{t}\phi\right)\left(\boldsymbol{\eta}\right)\right\vert + \left\vert (P_{t}\phi^{(j)})\left(\boldsymbol{\eta}\right) - (P_{t}^{n}\phi^{(j)})\left(\boldsymbol{\eta}\right)\right\vert.
     \end{aligned}
     \end{equation}
     We will bound each term on the right-hand side of~\eqref{Paper01_intermediate_estimate_uniform_convergence_semigroups_in_n_over_compacts} separately. For the first term, note that for any $\boldsymbol{\eta} \in \mathcal{K}$, $t \in [0,T]$ and $j \in \mathbb{N}$, we have
    \begin{equation} \label{Paper01_uniform_estimate_compact_open_topology}
    \begin{aligned}
        & \left\vert (P_{t}\phi)\left(\boldsymbol{\eta}\right) - (P_{t}\phi^{(j)})\left(\boldsymbol{\eta}\right)\right\vert \\ & \quad  \leq \mathbb{E}_{\boldsymbol{\eta}}\left[\left\vert \phi(\eta(t)) - \phi^{(j)}(\eta(t))\right\vert \right] \\ & \quad = \mathbb{E}_{\boldsymbol{\eta}}\left[\left\vert \phi(\eta(t)) - \phi^{(j)}(\eta(t))\right\vert \cdot \mathds{1}_{\{\eta(t) \in \mathscr{K}_{j}\}} \right]  + \mathbb{E}_{\boldsymbol{\eta}}\left[\left\vert \phi(\eta(t)) - \phi^{(j)}(\eta(t))\right\vert \cdot \mathds{1}_{\{\eta(t) \not\in \mathscr{K}_{j}\}} \right] \\ & \quad \leq \varepsilon_{j}\Big(1 + 2 \vert \vert \phi \vert\vert_{L_{\infty}(\mathcal{S}; \mathbb{R})}\Big),
    \end{aligned}
    \end{equation}
    where for the last inequality we applied~\eqref{Paper01_application_stone_weierstrass},~\eqref{Paper01_compact_containment_condition_general_argument} and the uniform bound on $(\phi^{(j)})_{j \in \mathbb{N}}$ in~\eqref{Paper01_uniform_boundedness_approximating_sequence_maps}. By the same argument, using~\eqref{Paper01_intermediate_step_high_probability_sequence_compacts} in place of~\eqref{Paper01_compact_containment_condition_general_argument}, we also have that for $j\in \mathbb N$,
    \begin{equation}
    \label{Paper01_uniform_estimate_compact_open_topology_ii}
        \sup_{\boldsymbol{\eta} \in \mathcal{K}} \; \sup_{n \in \mathbb{N}} \; \sup_{t \in [0,T]} \; \left\vert (P^{n}_{t}\phi^{(j)})\left(\boldsymbol{\eta}\right) - \left(P^{n}_{t}\phi\right)\left(\boldsymbol{\eta}\right)\right\vert  \leq \varepsilon_{j}\Big(1 + 2 \vert \vert \phi \vert\vert_{L_{\infty}(\mathcal{S}; \mathbb{R})}\Big).
    \end{equation}
    For the third term on the right-hand side of~\eqref{Paper01_intermediate_estimate_uniform_convergence_semigroups_in_n_over_compacts}, by Assumption~\ref{Paper01_assumption_tightness_sequence_processes}(iv), there exists a sequence of positive integers $(n_{j})_{j \in \mathbb{N}}$ such that for any $j \in \mathbb{N}$ and $n \geq n_{j}$, we have
     \begin{equation} \label{Paper01_almost_finishing_proof_assertion_iii}
         \sup_{\boldsymbol{\eta} \in \mathcal{K}} \; \sup_{t \in [0,T]} \; \left\vert \left(P_{t}\phi^{(j)}\right)\left(\boldsymbol{\eta}\right) - \left(P_{t}^{n}\phi^{(j)}\right)\left(\boldsymbol{\eta}\right)\right\vert \leq \varepsilon_{j}.
     \end{equation}
     Therefore, applying~\eqref{Paper01_uniform_estimate_compact_open_topology},~\eqref{Paper01_uniform_estimate_compact_open_topology_ii} and~\eqref{Paper01_almost_finishing_proof_assertion_iii} to~\eqref{Paper01_intermediate_estimate_uniform_convergence_semigroups_in_n_over_compacts}, and recalling that by construction $\varepsilon_{j} \rightarrow 0$ as $j \rightarrow \infty$, we conclude that assertion~(iii) holds.
    
    \medskip

    \noindent \underline{Proof of assertion~(iv):}

    \medskip

   Since $\mathcal{S}$ is a complete and separable metric space, the product topology on $[0, \infty) \times \mathcal{S}$ is also metrisable, and therefore, to prove our claim, it is enough to verify that for any $T \geq 0$, any compact subset $\mathcal{K} \subset \mathcal{S}$ and any $\phi \in \mathscr{C}_{b}(\mathcal{S}, \mathbb{R})$, the restriction of the map $(t, \boldsymbol{\eta}) \mapsto (P_{t}\phi)(\boldsymbol{\eta})$ to $[0,T] \times \mathcal{K}$ is continuous (see~\cite[Theorem~43.9 and Lemma~43.10]{willard2012general}). Since, by Assumption~\ref{Paper01_assumption_tightness_sequence_processes}, for each $n \in \mathbb{N}$, 
   $\{P^n_t\}_{t\ge 0}$ is a Feller semigroup,
   for any $n \in \mathbb{N}$, the map $(t, \boldsymbol{\eta}) \mapsto (P^{n}_{t}\phi)(\boldsymbol{\eta})$ is continuous (see Section~\ref{Paper01_subsection_appendix_feller_non_locally_compact_polish_space} for a brief overview of Feller semigroups on non-locally compact metric spaces). Therefore, by assertion~(iii) of this theorem, the map $(t, \boldsymbol{\eta}) \mapsto (P_{t}\phi)(\boldsymbol{\eta})$ restricted to $[0,T] \times \mathcal{K}$ is the limit of a sequence of continuous functions in the topology of uniform convergence, and therefore it is also continuous, as desired.
   The fact that $P_T \phi \in \mathscr C_b(\mathcal S,\mathbb R)$ follows directly since $\phi$ is bounded and so $P_T \phi$ is bounded by assertion~(ii).
     
    \medskip

    \noindent \underline{Proof of assertion~(v):}

    \medskip

    The fact that for any $T \geq 0$, $P_{T}: \mathscr{C}_{b}(\mathcal{S}, \mathbb{R}) \rightarrow \mathscr{C}_{b}(\mathcal{S}, \mathbb{R})$ is a bounded linear operator and that $\vert \vert P_{T} \vert \vert_{\textrm{op}} \leq 1$ follows directly from assertion~(ii), i.e.~follows from $(P_{T}\phi)(\boldsymbol{\eta}) = \mathbb{E}_{\boldsymbol{\eta}}[\phi(\eta(T))]$. It remains to prove that $\{P_{t}\}_{t \geq 0}$ is a semigroup. Since by Assumption~\ref{Paper01_assumption_tightness_sequence_processes}, for each $n \in \mathbb{N}$, $\{P^{n}_{t}\}_{t \geq 0}$ is a semigroup of bounded linear operators, by assertion~(ii) of this theorem, it will suffice to prove that for any $\phi \in \mathscr{C}_{b}(\mathcal{S}, \mathbb{R})$, $T,t\geq 0$ and $\boldsymbol{\eta} \in \mathcal{S}$, we have
    \begin{equation} \label{Paper01_equivalence_characterisation_semigroup}
        \lim_{n \rightarrow \infty} \left\vert (P_{T}P_{t}\phi)(\boldsymbol{\eta}) - (P^{n}_{T}P^{n}_{t}\phi)(\boldsymbol{\eta})\right\vert = 0.
    \end{equation}
    By the triangle inequality, we can write
    \begin{equation} \label{Paper01_simple_triangle_inequality_semigroup_proof}
        \left\vert (P_{T}P_{t}\phi)(\boldsymbol{\eta}) - (P^{n}_{T}P^{n}_{t}\phi)(\boldsymbol{\eta})\right\vert \leq \left\vert (P_{T}- P^{n}_{T})(P_{t}\phi)(\boldsymbol{\eta})\right\vert + \left\vert P^{n}_{T}(P_{t} - P^{n}_{t})(\phi)(\boldsymbol{\eta}) \right\vert.
    \end{equation}
    We will tackle each of the terms on the right-hand side of~\eqref{Paper01_simple_triangle_inequality_semigroup_proof} separately. Starting with the first term, observe that by assertion~(iv) of this theorem, we have $P_{t}\phi \in \mathscr{C}_{b}(\mathcal{S}, \mathbb{R})$. Hence, the weak convergence stated in assertion~(i) of this theorem yields
    \begin{equation*}
        \lim_{n \rightarrow \infty} \; \left\vert (P_{T}- P^{n}_{T})(P_{t}\phi)(\boldsymbol{\eta})\right\vert = 0.
    \end{equation*}
    To bound the second term on the right-hand side of~\eqref{Paper01_simple_triangle_inequality_semigroup_proof}, as in the proof of assertion~(iii), let $(\varepsilon_{j})_{j \in \mathbb{N}}$ denote a sequence of strictly positive real numbers converging to $0$, and, using Assumption~\ref{Paper01_assumption_tightness_sequence_processes}(iii), let $(\mathscr{K}_{j})_{j \in \mathbb{N}}$ denote a sequence of compact subsets of $\mathcal{S}$ satisfying, for each $j \in \mathbb{N}$,
     \begin{equation} \label{Paper01_intermediate_sequence_compact_sets_containing_high_probability}
         \inf_{n \in \mathbb{N}} \; \mathbb{P}_{\boldsymbol{\eta}} \Big(\eta^{n}(T) \in \mathscr{K}_{j}\Big) \geq 1 - \varepsilon_{j}.
     \end{equation}
    By assertion~(iii) of this theorem, there exists a strictly increasing sequence of positive integers $(n_{j})_{j \in \mathbb{N}} \defeq (n_{j}(t,\phi,\boldsymbol{\eta}))_{j \in \mathbb{N}}$ such that for any $j \in \mathbb{N}$ and $n \geq n_{j}$,
    \begin{equation} \label{Paper01_converging_intermediate_step_semigroups_on_compacts}
        \sup_{n \geq n_{j}} \; \sup_{\boldsymbol{\xi} \in \mathscr{K}_{j}} \; \vert (P_{t}\phi)(\boldsymbol{\xi}) - (P_{t}^{n}\phi)(\boldsymbol{\xi}) \vert \leq \varepsilon_{j}.
    \end{equation}
    Therefore, for $n \geq n_{j}$, we have
    \begin{equation} \label{Paper01_almost_finastep_to_pove_semigroup_property}
    \begin{aligned}
        \left\vert P^{n}_{T}(P_{t} - P^{n}_{t})(\phi)(\boldsymbol{\eta}) \right\vert & = \left\vert \mathbb{E}_{\boldsymbol{\eta}}\Big[P_{t}\phi(\eta^{n}(T)) - P^{n}_{t}\phi(\eta^{n}(T))\Big] \right\vert \\ & \leq \mathbb{E}_{\boldsymbol{\eta}}\Big[\Big\vert P_{t}\phi(\eta^{n}(T)) - P^{n}_{t}\phi(\eta^{n}(T)) \Big\vert \cdot \mathds{1}_{\{\eta^{n}(T) \in \mathscr{K}_{j}\}}\Big] \\ & \quad \quad + \mathbb{E}_{\boldsymbol{\eta}}\Big[\Big\vert P_{t}\phi(\eta^{n}(T)) - P^{n}_{t}\phi(\eta^{n}(T)) \Big\vert \cdot \mathds{1}_{\{\eta^{n}(T) \not\in \mathscr{K}_{j}\}}\Big] \\ & \leq \varepsilon_{j}\Big(1 + 2 \vert \vert \phi \vert \vert_{L_{\infty}(\mathcal{S}; \mathbb{R})}\Big),
    \end{aligned}
    \end{equation}
    where for the last inequality we used~\eqref{Paper01_converging_intermediate_step_semigroups_on_compacts} to bound the first term, and~\eqref{Paper01_intermediate_sequence_compact_sets_containing_high_probability} and the fact that $\vert \vert P_{T} \vert \vert_{\textrm{op}} \leq 1$ and $\vert \vert P^{n}_{T} \vert \vert_{\textrm{op}} \leq 1$  to bound the second term. Thus, by taking $j \rightarrow \infty$ in~\eqref{Paper01_almost_finastep_to_pove_semigroup_property}, it follows from~\eqref{Paper01_simple_triangle_inequality_semigroup_proof} that~\eqref{Paper01_equivalence_characterisation_semigroup} holds.

    \medskip

    \noindent \underline{Proof of assertion~(vi):}

    \medskip
    
    Fix $\boldsymbol{\eta} \in \mathcal{S}_{0}$ and, for each $n \in \mathbb{N}$, let $(\eta^{n}(t))_{t \geq 0}$ be the Markov process with sample paths in $\mathscr{D}([0, \infty), \mathcal{S})$ and associated semigroup $\{P^{n}_{t}\}_{t \geq 0}$ with $\eta^{n}(0) = \boldsymbol{\eta}$ almost surely. Since by Assumption~\ref{Paper01_assumption_tightness_sequence_processes}(ii) this sequence is tight in $\mathscr{D}([0, \infty), \mathcal{S})$, in order to prove uniqueness of the limit it will suffice to establish convergence of the finite-dimensional distributions~\cite[Theorem~3.7.8]{ethier2009markov}, i.e.~that for any $I \in \mathbb{N}$, $(\phi_{i})_{i = 1}^{I} \subset \mathscr{C}_{b}(\mathcal{S}, \mathbb{R})$ and $0 \leq t_{1} \leq t_{2} \leq \cdots \leq t_{I}$, the sequence
    \begin{equation*}
        \Big(\mathbb{E}_{\boldsymbol{\eta}}\left[\phi_{1}(\eta^{n}(t_{1}))\phi_{2}(\eta^{n}(t_{2})) \cdots \phi_{I}(\eta^{n}(t_{I}))\right]\Big)_{n \in \mathbb{N}}
    \end{equation*}
    converges as $n \rightarrow \infty$. Thus, the convergence of the finite-dimensional distributions will be proved after establishing that for any $I \in \mathbb{N}$, $(\phi_{i})_{i = 1}^{I} \subset \mathscr{C}_{b}(\mathcal{S}, \mathbb{R})$ and $0 \leq t_{1} \leq t_{2} \leq \cdots \leq t_{I}$, we have
    \begin{equation} \label{Paper01_convergence_fdd_markovian_process}
    \begin{aligned}
        & \lim_{n \rightarrow \infty} \; \mathbb{E}_{\boldsymbol{\eta}}\Big[\phi_{1}(\eta^{n}(t_{1}))\phi_{2}(\eta^{n}(t_{2})) \cdots \phi_{I}(\eta^{n}(t_{I}))\Big] \\ & \quad = P_{t_{1}}\Big(\phi_{1}P_{(t_2-t_1)}\Big(\phi_{2}P_{(t_3-t_2)}\Big(\cdots \Big(\phi_{I-1}P_{(t_{I} - t_{I-1})}\phi_{I}\Big) \cdots\Big)\Big)\Big)(\boldsymbol{\eta}).
    \end{aligned}
    \end{equation}

    We will prove~\eqref{Paper01_convergence_fdd_markovian_process} using an induction argument on $I$. Observe that assertion~(ii) of this theorem guarantees that for any $t \geq 0$ and any $\phi \in \mathscr{C}_{b}(\mathcal{S}, \mathbb{R})$, we have
    \begin{equation*}
        \lim_{n \rightarrow \infty} \; \mathbb{E}_{\boldsymbol{\eta}}\Big[\phi(\eta^{n}(t))\Big] = (P_{t}\phi)(\boldsymbol{\eta}).
    \end{equation*}
    Suppose now that for some $I \in \mathbb{N}$, identity~\eqref{Paper01_convergence_fdd_markovian_process} holds. Then, for any $(\phi_{i})_{i = 1}^{I+ 1} \subset \mathscr{C}_{b}(\mathcal{S}, \mathbb{R})$, any $0 \leq t_{1} \leq t_{2} \leq \cdots \leq t_{I} \leq t_{I+1}$, and any $n \in \mathbb{N}$, by the fact that, by Assumption~\ref{Paper01_assumption_tightness_sequence_processes}, $\eta^{n}$ is a Markov process, applying the Markov property at time $t_{I}$ and using the tower property of conditional expectation, we have
    \begin{equation} \label{Paper01_intermediate_equation_proof_convergence_fdds_markovian}
    \begin{aligned}
        & \mathbb{E}_{\boldsymbol{\eta}}\Big[\phi_{1}(\eta^{n}(t_{1}))\phi_{2}(\eta^{n}(t_{2})) \cdots \phi_{I}(\eta^{n}(t_{I}))\phi_{I+1}(\eta^{n}(t_{I+1}))\Big] \\ & \quad = \mathbb{E}_{\boldsymbol{\eta}}\Big[\phi_{1}(\eta^{n}(t_{1}))\phi_{2}(\eta^{n}(t_{2})) \cdots \phi_{I}(\eta^{n}(t_{I}))P^{n}_{(t_{I+1} - t_{I})}\phi_{I+1}(\eta^{n}(t_{I}))\Big] \\ & \quad = \mathbb{E}_{\boldsymbol{\eta}}\Big[\phi_{1}(\eta^{n}(t_{1}))\phi_{2}(\eta^{n}(t_{2})) \cdots \phi_{I}(\eta^{n}(t_{I}))(P^{n}_{(t_{I+1} - t_{I})} - P_{(t_{I+1} - t_{I})})\phi_{I+1}(\eta^{n}(t_{I}))\Big] \\ & \quad \quad \quad  + \mathbb{E}_{\boldsymbol{\eta}}\Big[\phi_{1}(\eta^{n}(t_{1}))\phi_{2}(\eta^{n}(t_{2})) \cdots \Big(\phi_{I}P_{(t_{I+1} - t_{I})}\phi_{I+1}\Big)(\eta^{n}(t_{I}))\Big].
    \end{aligned}
    \end{equation}
    We will tackle the terms on the right-hand side of~\eqref{Paper01_intermediate_equation_proof_convergence_fdds_markovian} separately. Starting with the last term, note that by assertion~(iv) of this theorem, $P_{(t_{I+1} - t_{I})}\phi_{I+1} \in \mathscr{C}_{b}(\mathcal{S}, \mathbb{R})$. Therefore, by the induction hypothesis, we conclude that
    \begin{equation*}
    \begin{aligned}
        & \lim_{n \rightarrow \infty} \; \mathbb{E}_{\boldsymbol{\eta}}\Big[\phi_{1}(\eta^{n}(t_{1}))\phi_{2}(\eta^{n}(t_{2})) \cdots \Big(\phi_{I}P_{(t_{I+1} - t_{I})}\phi_{I+1}\Big)(\eta^{n}(t_{I}))\Big] \\ & \quad = P_{t_{1}}\Big(\phi_{1}P_{(t_2-t_1)}\Big(\phi_{2}P_{(t_3-t_2)}\Big(\cdots \Big(\phi_{I}P_{(t_{I+1} - t_{I})}\phi_{I+1}\Big) \cdots\Big)\Big)\Big)(\boldsymbol{\eta}).
    \end{aligned}
    \end{equation*}
    Hence, to show that~\eqref{Paper01_convergence_fdd_markovian_process} holds with $I$ replaced with $I+ 1$, and so conclude the induction argument, it will be enough to verify that the first term on the right-hand side of~\eqref{Paper01_intermediate_equation_proof_convergence_fdds_markovian} converges to $0$ as $n \rightarrow \infty$.
    
    As in the proof of assertion~(iii), let $(\varepsilon_{j})_{j \in \mathbb{N}}$ be a sequence of strictly positive real numbers converging to $0$, and, using Assumption~\ref{Paper01_assumption_tightness_sequence_processes}(iii), let $(\mathscr{K}_{j})_{j \in \mathbb{N}} = (\mathscr{K}_{j}(\boldsymbol{\eta}, t_{I}))_{j \in \mathbb{N}}$ be a sequence of compact subsets of $\mathcal{S}$ such that for each $j \in \mathbb{N}$,
    \begin{equation} \label{Paper01_sequence_compacts_containing_high_probability_proof_markov_property}
        \inf_{n \in \mathbb{N}} \, \mathbb{P}_{\boldsymbol{\eta}}\Big(\eta^{n}(t_{I}) \in \mathscr{K}_{j}\Big) \geq 1 - \varepsilon_{j}.
    \end{equation}
    By repeating the argument used for~\eqref{Paper01_converging_intermediate_step_semigroups_on_compacts} in the proof of assertion~(v), there exists a strictly increasing sequence of positive integers $(n_{j})_{j \in \mathbb{N}} \defeq (n_{j}(t_{I+1} - t_{I},\phi,\boldsymbol{\eta}))_{j \in \mathbb{N}}$ such that for every $j \in \mathbb{N}$,
    \begin{equation} \label{Paper01_intermediate_step_markov_property}
        \sup_{n \geq n_{j}} \; \sup_{\boldsymbol{\xi} \in \mathscr{K}_{j}} \; \vert (P_{(t_{I+1} - t_{I})}\phi_{I+1})(\boldsymbol{\xi}) - (P_{(t_{I+1} - t_{I})}^{n}\phi_{I+1})(\boldsymbol{\xi}) \vert \leq \varepsilon_{j}.
    \end{equation}
    Therefore, for every $j \in \mathbb{N}$, for any $n \geq n_{j}$, using that $\phi_{i} \in \mathscr{C}_{b}(\mathcal{S}, \mathbb{R})$ for all $i \leq I$ in the first inequality,
    \begin{align} \label{Paper01_almost_almost_final_step_markov_property}
        & \Big\vert \mathbb{E}_{\boldsymbol{\eta}}\Big[\phi_{1}(\eta^{n}(t_{1}))\phi_{2}(\eta^{n}(t_{2})) \cdots \phi_{I}(\eta^{n}(t_{I}))(P^{n}_{(t_{I+1} - t_{I})} - P_{(t_{I+1} - t_{I})})\phi_{I+1}(\eta^{n}(t_{I}))\Big] \Big\vert \notag \\ 
        & \quad \leq \Bigg(\prod_{i = 1}^{I} \vert \vert \phi_{i} \vert \vert_{L_{\infty}(\mathcal{S}; \mathbb{R})}\Bigg) \mathbb{E}_{\boldsymbol{\eta}}\Big[\Big\vert(P_{(t_{I+1} - t_{I})}\phi_{I+1})(\eta^{n}(t_{I})) - (P_{(t_{I+1} - t_{I})}^{n}\phi_{I+1})(\eta^{n}(t_{I}))\Big\vert\Big]
        \notag \\ 
        & \quad = \Bigg(\prod_{i = 1}^{I} \vert \vert \phi_{i} \vert \vert_{L_{\infty}(\mathcal{S}; \mathbb{R})}\Bigg) \mathbb{E}_{\boldsymbol{\eta}}\Big[\Big\vert(P_{(t_{I+1} - t_{I})}\phi_{I+1})(\eta^{n}(t_{I})) - (P_{(t_{I+1} - t_{I})}^{n}\phi_{I+1})(\eta^{n}(t_{I}))\Big\vert \cdot \mathds{1}_{\{\eta^{n}(t_{I}) \in \mathscr{K}_{j}\}}\Big] \notag \\
        & \quad \quad \, + \Bigg(\prod_{i = 1}^{I} \vert \vert \phi_{i} \vert \vert_{L_{\infty}(\mathcal{S}; \mathbb{R})}\Bigg) \mathbb{E}_{\boldsymbol{\eta}}\Big[\Big\vert(P_{(t_{I+1} - t_{I})}\phi_{I+1})(\eta^{n}(t_{I})) - (P_{(t_{I+1} - t_{I})}^{n}\phi_{I+1})(\eta^{n}(t_{I}))\Big\vert \cdot \mathds{1}_{\{\eta^{n}(t_{I}) \not\in \mathscr{K}_{j}\}}\Big] \notag \\
        & \quad \leq \varepsilon_{j} \Bigg(\prod_{i = 1}^{I} \vert \vert \phi_{i} \vert \vert_{L_{\infty}(\mathcal{S}; \mathbb{R})}\Bigg)\Big(1 + 2\vert \vert \phi_{I+1} \vert \vert_{L_{\infty}(\mathcal{S}; \mathbb{R})}\Big),
    \end{align}
    where for the last inequality we used~\eqref{Paper01_intermediate_step_markov_property} to bound the first term, and the bounds on the operator norms $\vert \vert P_{t_{I+1} - t_{I}} \vert \vert_{\textrm{op}} \leq 1$ (by assertion (v)) and $\vert \vert P^{n}_{t_{I+1} - t_{I}} \vert \vert_{\textrm{op}} \leq 1$ (by Assumption~\ref{Paper01_assumption_tightness_sequence_processes}) for every $n \in \mathbb{N}$, together with~\eqref{Paper01_sequence_compacts_containing_high_probability_proof_markov_property} for the second term.
    Then, by taking $j \rightarrow \infty$ in~\eqref{Paper01_almost_almost_final_step_markov_property}, we conclude that
    \begin{equation*}
        \lim_{n \rightarrow \infty} \, \mathbb{E}_{\boldsymbol{\eta}}\Big[\phi_{1}(\eta^{n}(t_{1}))\phi_{2}(\eta^{n}(t_{2})) \cdots \phi_{I}(\eta^{n}(t_{I}))(P^{n}_{(t_{I+1} - t_{I})} - P_{(t_{I+1} - t_{I})})\phi_{I+1}(\eta^{n}(t_{I}))\Big] = 0.
    \end{equation*}
    This completes the proof of~\eqref{Paper01_convergence_fdd_markovian_process}, i.e.~that the finite-dimensional distributions of the sequence $(\eta^{n})_{n \in \mathbb{N}}$ converge.
    
    As discussed before~\eqref{Paper01_convergence_fdd_markovian_process}, we now have convergence of the sequence of processes $(\eta^{n})_{n \in \mathbb{N}}$ in $\mathscr{D}([0, \infty), \mathcal{S})$ to a unique process $(\eta(t))_{t \geq 0}$.
    Let $(\eta(t))_{t \geq 0}$ be the limiting càdlàg process. Observe also that $\{P_{t}\}_{t \geq 0}$ is a Feller semigroup. Indeed, assertion~(v) of this theorem implies conditions~(ii) and~(iii) of Definition~\ref{Paper01_definition_feller_non_locally_compact}, assertion~(ii) implies condition~(iv), and assertion~(iv) of this theorem implies conditions~(i),(v) and~(vi). By identity~\eqref{Paper01_convergence_fdd_markovian_process}, the transition probabilities of $(\eta(t))_{t \geq 0}$ are given by~$\{P_{t}\}_{t \geq 0}$, and therefore, by Lemma~\ref{Paper01_markov_property_feller_semigroup}, $(\eta(t))_{t \geq 0}$ is a Markov process. This completes the proof of assertion~(vi) of this theorem.

    \medskip

    \noindent \underline{Proof of assertion~(vii):}

    \medskip

    Let $\phi \in \mathfrak{A}$ and $\boldsymbol{\eta} \in \mathcal{S}$ be fixed. Then, observe that by assertion (ii),
    \begin{equation} \label{Paper01_double_limit}
         \lim_{t \rightarrow 0^{+}} \; \frac{(P_{t}\phi)(\boldsymbol{\eta}) - \phi(\boldsymbol{\eta})}{t} =  \lim_{t \rightarrow 0^{+}} \;  \lim_{n \rightarrow \infty} \; \frac{(P^{n}_{t}\phi)(\boldsymbol{\eta}) - \phi(\boldsymbol{\eta})}{t}.
    \end{equation}
    To compute the term on right-hand side of~\eqref{Paper01_double_limit}, we will prove that we can interchange the order of the limits. To follow this approach, by the Moore-Osgood theorem (see~\cite[Theorem~7.11]{rudin1976principles}), we need to establish the following conditions:
    \begin{enumerate}[(1)]
        \item The following limit holds for any $n \in \mathbb{N}$:
        \begin{equation*}
            \lim_{t \rightarrow 0^{+}} \; \frac{(P^{n}_{t}\phi)(\boldsymbol{\eta}) - \phi(\boldsymbol{\eta})}{t} = (\mathcal{L}^{n}\phi)(\boldsymbol{\eta});
    \end{equation*}
        \item There exists $T >0$ such that the following limit holds uniformly in $t \in (0,T]$:
        \begin{equation*}
        \lim_{n \rightarrow \infty} \; (P^{n}_{t}\phi)(\boldsymbol{\eta}) =  (P_{t}\phi)(\boldsymbol{\eta}).
    \end{equation*}
    
    \end{enumerate}
    Condition~(1) follows directly from Assumption~\ref{Paper01_assumption_tightness_sequence_processes}(i), while condition~(2) follows from assertion~(iii) of this theorem. Therefore, we can interchange the order of the limits on the right-hand side of~\eqref{Paper01_double_limit}, obtaining
    \begin{equation*}
    \begin{aligned}
         \lim_{t \rightarrow 0^{+}} \; \frac{(P_{t}\phi)(\boldsymbol{\eta}) - \phi(\boldsymbol{\eta})}{t} & = \lim_{n \rightarrow \infty} \; \lim_{t \rightarrow 0^{+}} \;   \frac{(P^{n}_{t}\phi)(\boldsymbol{\eta}) - \phi(\boldsymbol{\eta})}{t} \\ & = \lim_{n \rightarrow \infty} \;(\mathcal{L}^{n}\phi)(\boldsymbol{\eta}) \\ & = (\mathcal{L}\phi)(\boldsymbol{\eta}),
    \end{aligned}
    \end{equation*}
    where the last line follows by Assumption~\ref{Paper01_assumption_tightness_sequence_processes}(i). Hence, assertion~(vii) holds.
    
    \medskip
    
    \noindent \underline{Proof of assertion~(viii):}

    \medskip 
    
    We start by mentioning that the main differences between our proof of assertion~(viii) and the proof of the strong Markov property for Feller processes in locally compact metric spaces (see~\cite[Theorem~4.2.5]{ethier2009markov}) are that we need the tightness condition stated in Assumption~\ref{Paper01_assumption_tightness_sequence_processes}(ii), and also that, since we do not assume that the semigroup~$\{P_{t}\}_{t \geq 0}$ is strongly continuous in the norm of uniform convergence, we also use Assumption~\ref{Paper01_assumption_tightness_sequence_processes}(iii) and assertion~(iii) of this theorem.
    
    We now prove assertion~(viii). Fix $\boldsymbol{\eta} \in \mathcal{S}_{0}$. Recall that $\mathscr{C}_b(\mathcal{S}, \mathbb{R})$ is boundedly and pointwise dense (bp-dense) in $\mathscr{B}_b(\mathcal{S},\mathbb{R})$, i.e.~any bounded Borel–measurable function can be obtained as the pointwise limit of a uniformly bounded sequence $(\phi_n)_{n \in \mathbb{N}} \subset \mathscr{C}_b(\mathcal{S},\mathbb{R})$ (see~\cite[Proposition~3.4.2]{ethier2009markov}). Hence, by~\eqref{Paper01_strong_markov_property}, it will suffice to establish that for any almost surely finite $\{\mathcal{F}_{t+}^{\eta}\}_{t\ge 0}$-stopping time $\tau$, any $T \geq 0$ and any $\phi \in \mathscr{C}_{b}(\mathcal{S}, \mathbb{R})$, we have
    \begin{equation} \label{eq:strongMarkovinproof}
        \mathbb{E}_{\boldsymbol{\eta}}\left[\phi(\eta(\tau + T)) \,  \Big\vert \, \mathcal{F}^{\eta}_{\tau+}\right] = (P_{T}\phi)(\eta(\tau)) \, \textrm{ almost surely}.
    \end{equation}
    
    We start by considering the case in which the almost surely finite $\{\mathcal{F}_{t+}^{\eta}\}_{t\ge 0}$-stopping time $\tau$ is discrete, i.e.~suppose $\tau$ is concentrated on $\{T_{i}: \, i \in \mathbb{N}\} \subset [0, \infty)$. Fix $T > 0$ and $\phi \in \mathscr{C}_{b}(\mathcal{S}, \mathbb{R})$. Observe that for any $J > 0$, the random variable $\tau \wedge J$ is also an $\{\mathcal{F}_{t+}^{\eta}\}_{t\ge 0}$-stopping time. Let $(\varepsilon_{j})_{j \in \mathbb{N}}$ be a sequence of strictly positive real numbers converging to $0$. Since by assertion~(vi) of this theorem, $(\eta(t))_{t \geq 0}$ is the weak limit of a tight sequence of processes in $\mathscr{D}([0, \infty), \mathcal{S})$, by Remark~3.7.3 in~\cite{ethier2009markov}, for any $j \in \mathbb{N}$, there exists a compact subset $\mathcal{K}_{j} = \mathcal{K}_{j}(J,T,\boldsymbol{\eta}) \subset \mathcal{S}$ such that
    \begin{equation*}
        \mathbb{P}_{\boldsymbol{\eta}}\Big(\eta(t) \in \mathcal{K}_{j} \quad \forall \, t \in [0, J+T]\Big) \geq 1 - \varepsilon_{j}.
    \end{equation*}
    In particular, this implies that 
    \begin{equation} \label{Paper01_bound_probability_compact_stopping_time}
        \inf_{\delta \in (0,T)} \; \mathbb{P}_{\boldsymbol{\eta}}\Big(\eta((\tau \wedge J) + \delta) \in \mathcal{K}_{j} \Big) \geq 1 - \varepsilon_{j}.
    \end{equation}
    Observe that for every $i \in \mathbb{N}$ and any $\delta > 0$, 
    \begin{equation*}
        \{\tau \wedge J = T_{i}\} \in \mathcal{F}^\eta_{T_{i}+} \subseteq \mathcal{F}^\eta_{T_{i}+\delta}.
    \end{equation*}
   Let $J > T_i$. For any $\mathscr{A} \in \mathcal{F}_{\tau+}^{\eta}$, any $\delta > 0$ and every $i \in \mathbb{N}$,
    \begin{equation} \label{Paper01_measurability_stopping_times}
        \mathscr{A} \cap \{\tau \wedge J = T_{i}\} \in \mathcal{F}^{\eta}_{T_{i}+} \subseteq \mathcal{F}^{\eta}_{T_{i} + \delta}.
    \end{equation}
    By dominated convergence, we conclude that for every $i \in \mathbb{N}$,
    \begin{equation} \label{Paper01_first_step_strong_markov}
    \begin{aligned}
        \mathbb{E}_{\boldsymbol{\eta}}\left[\phi(\eta(\tau + T)) \cdot \mathds{1}_{\mathscr{A} \cap \{\tau = T_{i}\}}\right] & = \lim_{J \rightarrow \infty} \; \mathbb{E}_{\boldsymbol{\eta}}\left[\phi(\eta((\tau \wedge J) + T)) \cdot \mathds{1}_{ \mathscr{A} \cap \{\tau \wedge J = T_{i}\}}\right] \\ & = \lim_{J \rightarrow \infty} \; \mathbb{E}_{\boldsymbol{\eta}}\left[\phi(\eta(T_{i} + T)) \cdot \mathds{1}_{ \mathscr{A} \cap \{\tau \wedge J = T_{i}\}}\right] \\ & = \lim_{J \rightarrow \infty} \; \mathbb{E}_{\boldsymbol{\eta}}\left[(P_{T - \delta} \phi)(\eta(T_{i} + \delta)) \cdot  \mathds{1}_{ \mathscr{A} \cap \{\tau \wedge J = T_{i}\}}\right],
    \end{aligned}
    \end{equation}
    where the last identity holds for any $\delta \in (0, T)$ by~\eqref{Paper01_measurability_stopping_times}, the Markov property of $(\eta(t))_{t \geq 0}$ and the tower property of conditional expectation. 
    
    Now take $J > T_{i}$ fixed. By assertion~(iv) of this theorem, there exists a strictly decreasing sequence $(\delta_{j})_{j \in \mathbb{N}}$ with $\delta_{j} \in (0,T)$ for every $j \in \mathbb{N}$ and $\lim_{j \rightarrow \infty} \delta_{j} = 0$ such that for any $j \in \mathbb{N}$ and $\delta \in (0, \delta_{j})$, we have
    \begin{equation} \label{Paper01_bound_uniform_compact_subsets_distance_semigroup}
        \sup_{\boldsymbol{\xi} \in \mathcal{K}_{j}} \; \left\vert (P_{T - \delta}\phi)(\boldsymbol{\xi}) - (P_{T}\phi)(\boldsymbol{\xi})\right\vert \leq \varepsilon_{j}.
    \end{equation}
    We now rewrite the term on the right-hand side of~\eqref{Paper01_first_step_strong_markov}, obtaining, for $i,j \in \mathbb{N}$ and $\delta \in (0,T \wedge \delta_{j})$, using that $J>T_i$ and so $\{\tau\wedge J=T_i\}=\{\tau=T_i\}$ in the first term,
    \begin{equation} \label{Paper01_intermediate_step_strong_markov_property}
    \begin{aligned}
        & \mathbb{E}_{\boldsymbol{\eta}}\left[(P_{T - \delta} \phi)(\eta(T_{i} + \delta)) \cdot \mathds{1}_{\mathscr{A} \cap \{\tau \wedge J = T_{i}\}}\right] \\ 
        & \quad = \mathbb{E}_{\boldsymbol{\eta}}\left[(P_{T} \phi)(\eta(T_{i} + \delta)) \cdot \mathds{1}_{\mathscr{A} \cap \{\tau  = T_{i}\}}\right] \\ 
        & \quad \quad \quad + \mathbb{E}_{\boldsymbol{\eta}}\left[(P_{T-\delta} \phi - P_{T} \phi)(\eta(T_{i} + \delta)) \cdot \mathds{1}_{\mathscr{A} \cap \{\tau \wedge J = T_{i}\} \cap  \{\eta((\tau \wedge J) + \delta) \in \mathcal{K}_{j}\}}\right] \\ & \quad \quad \quad + \mathbb{E}_{\boldsymbol{\eta}}\left[(P_{T-\delta} \phi - P_{T} \phi)(\eta(T_{i} + \delta)) \cdot \mathds{1}_{\mathscr{A} \cap \{\tau \wedge J = T_{i}\} \cap  \{\eta((\tau \wedge J) + \delta) \not\in \mathcal{K}_{j}\}}\right].
    \end{aligned}
    \end{equation}
    We will tackle each term on the right-hand side of~\eqref{Paper01_intermediate_step_strong_markov_property} separately. Starting with the first term, since by assertion~(iv) of this theorem we have $P_{T}\phi \in \mathscr{C}_{b}(\mathcal{S}, \mathbb{R})$, and since $(\eta(t))_{t \geq 0}$ is right-continuous by assertion~(vi), by dominated convergence we have that
    \begin{equation} \label{Paper01_intermediate_limit_strong_markov_property}
    \begin{aligned}
        \lim_{\delta \rightarrow 0^{+}} \; \mathbb{E}_{\boldsymbol{\eta}}\left[(P_{T} \phi)(\eta(T_{i} + \delta)) \cdot \mathds{1}_{\mathscr{A} \cap \{\tau = T_{i}\}}\right] = \mathbb{E}_{\boldsymbol{\eta}}\left[(P_{T} \phi)(\eta(T_{i})) \cdot \mathds{1}_{\mathscr{A} \cap \{\tau = T_{i}\}}\right].
    \end{aligned}
    \end{equation}
    To bound the second term on the right-hand side of~\eqref{Paper01_intermediate_step_strong_markov_property}, we observe that by~\eqref{Paper01_bound_uniform_compact_subsets_distance_semigroup} and since $\delta \in (0,\delta_j)$,
    \begin{equation} \label{Paper01_intermediate_intermediate_second_term_strong_markov}
    \begin{aligned}
        \left\vert \mathbb{E}_{\boldsymbol{\eta}}\left[(P_{T-\delta} \phi - P_{T} \phi)(\eta(T_{i} + \delta)) \cdot \mathds{1}_{\mathscr{A} \cap \{\tau \wedge J = T_{i}\} \cap  \{\eta((\tau \wedge J) + \delta) \in \mathcal{K}_{j}\}}\right] \right\vert  \leq \varepsilon_{j}.
    \end{aligned}
    \end{equation}
    For the third term on the right-hand side of~\eqref{Paper01_intermediate_step_strong_markov_property}, we use the fact that $\phi \in \mathscr{C}_{b}(\mathcal{S}, \mathbb{R})$, and that by assertion~(v) of this theorem $\vert \vert P_{T} \vert \vert_{\textrm{op}} \leq 1$ and $\vert \vert P_{T-\delta} \vert \vert_{\textrm{op}} \leq 1$ in the second inequality, and then~\eqref{Paper01_bound_probability_compact_stopping_time} in the third inequality, to obtain
    \begin{equation} \label{Paper01_intermediate_intermediate_third_term_strong_markov}
    \begin{aligned}
         & \left\vert \mathbb{E}_{\boldsymbol{\eta}}\left[(P_{T-\delta} \phi - P_{T} \phi)(\eta(T_{i} + \delta)) \cdot \mathds{1}_{\mathscr{A} \cap \{\tau \wedge J = T_{i}\} \cap  \{\eta((\tau \wedge J) + \delta) \not\in \mathcal{K}_{j}\}}\right] \right\vert \\ & \quad \leq \mathbb{E}_{\boldsymbol{\eta}}\left[\Big\vert (P_{T-\delta} \phi - P_{T} \phi)(\eta(T_{i} + \delta)) \big\vert \cdot \mathds{1}_{\mathscr{A} \cap \{\tau \wedge J = T_{i}\} \cap  \{\eta((\tau \wedge J) + \delta) \not\in \mathcal{K}_{j}\}}\right] \\ & \quad \leq 2 \vert \vert \phi \vert \vert_{L_{\infty}(\mathcal{S}; \mathbb{R})} \mathbb{P}_{\boldsymbol{\eta}} \Big(\eta((\tau \wedge J) + \delta) \not\in \mathcal{K}_{j} \Big) \\ & \quad \leq 2 \varepsilon_{j} \vert \vert \phi \vert \vert_{L_{\infty}(\mathcal{S}; \mathbb{R})}.
    \end{aligned}
    \end{equation}
    Therefore, by applying~\eqref{Paper01_intermediate_step_strong_markov_property} to~\eqref{Paper01_first_step_strong_markov}, then using~\eqref{Paper01_intermediate_limit_strong_markov_property},~\eqref{Paper01_intermediate_intermediate_second_term_strong_markov} and~\eqref{Paper01_intermediate_intermediate_third_term_strong_markov}, taking the limit $\delta \rightarrow 0^+$ and using that $\varepsilon_j\to 0$ as $j\to \infty$, we conclude that for every $i \in \mathbb{N}$,
    \begin{equation} \label{eq:towardsSMP}
        \mathbb{E}_{\boldsymbol{\eta}}\left[\phi(\eta(\tau + T)) \cdot \mathds{1}_{\mathscr{A} \cap \{\tau = T_{i}\}}\right] = \mathbb{E}_{\boldsymbol{\eta}}\left[(P_{T} \phi)(\eta(T_{i})) \cdot \mathds{1}_{\mathscr{A} \cap \{\tau = T_{i}\}}\right].
    \end{equation}
    Since we assumed $\tau$ to be a $\{\mathcal{F}^{\eta}_{t+}\}_{t\ge 0}$-stopping time concentrated on $\{T_i:i\in \mathbb N\}$, we conclude that for any $\mathscr{A} \in \mathcal{F}^{\eta}_{\tau+}$, any $T > 0$ and any $\phi \in \mathscr{C}_{b}(\mathcal{S}, \mathbb{R})$, using~\eqref{eq:towardsSMP} in the second line,
     \begin{equation} \label{Paper01_strong_markov_property_discrete_stopping_time}
     \begin{aligned}
     \mathbb{E}_{\boldsymbol{\eta}}\left[\phi(\eta(\tau + T)) \cdot \mathds{1}_{\mathscr{A}}\right] & =  \sum_{i \in \mathbb{N}} \mathbb{E}_{\boldsymbol{\eta}}\left[\phi(\eta(\tau + T)) \cdot \mathds{1}_{\mathscr{A}\cap\{\tau = T_{i}\}}\right] \\ &  = \sum_{i \in \mathbb{N}} \mathbb{E}_{\boldsymbol{\eta}}\left[(P_{T} \phi)(\eta(T_{i})) \cdot \mathds{1}_{\mathscr{A}\cap \{\tau = T_{i}\}}\right] \\ & = \mathbb{E}_{\boldsymbol{\eta}}\left[(P_{T} \phi)(\eta(\tau)) \cdot \mathds{1}_{\mathscr{A}}\right],
     \end{aligned}
    \end{equation}
    which implies~\eqref{eq:strongMarkovinproof}. This proves the strong Markov property with respect to almost surely finite discrete stopping times.

    For a general almost surely finite $\{\mathcal{F}^{\eta}_{t+}\}_{t\ge 0}$-stopping time $\tau$, we observe that by standard results in stochastic processes (see~\cite[Proposition~2.1.3]{ethier2009markov}), there exists a non-increasing sequence of almost surely finite discrete $\{\mathcal{F}^{\eta}_{t+}\}_{t\ge 0}$-stopping times $(\tau_{i})_{i \in \mathbb{N}}$ such that $\tau_{i} \rightarrow \tau$ as $i \rightarrow \infty$ almost surely. Let $\mathscr{A} \in \mathcal{F}^{\eta}_{\tau+}$ be fixed. Since $\tau_{i} \geq \tau$ for every $i \in \mathbb{N}$, we conclude (see~\cite[Proposition~2.1.4]{ethier2009markov}) that $\mathscr{A} \in \mathcal{F}^{\eta}_{\tau_{i}+}$ for every $i \in \mathbb{N}$. For any $\phi \in \mathscr{C}_{b}(\mathcal{S}, \mathbb{R})$ and $T \geq 0$, since $(\eta(t))_{t \geq 0}$ is right-continuous by assertion~(vi) of this theorem, and then by dominated convergence we have
    \begin{equation*}
    \begin{aligned}
        \mathbb{E}_{\boldsymbol{\eta}}\left[\phi(\eta(\tau + T)) \cdot \mathds{1}_{\mathscr{A}}\right] & = \mathbb{E}_{\boldsymbol{\eta}}\left[\lim_{i \rightarrow \infty} \; \phi(\eta(\tau_{i} + T)) \cdot \mathds{1}_{\mathscr{A}}\right] \\ & = \lim_{i \rightarrow \infty} \; \mathbb{E}_{\boldsymbol{\eta}}\left[ \phi(\eta(\tau_{i} + T)) \cdot \mathds{1}_{\mathscr{A}}\right] \\ & = \lim_{i \rightarrow \infty} \; \mathbb{E}_{\boldsymbol{\eta}}\left[(P_{T} \phi)(\eta(\tau_{i})) \cdot \mathds{1}_{\mathscr{A}}\right] \\ 
        & = \mathbb{E}_{\boldsymbol{\eta}}\left[ (P_{T} \phi)(\eta(\tau)) \cdot \mathds{1}_{\mathscr{A}}\right],
    \end{aligned}
    \end{equation*}
    where we used~\eqref{Paper01_strong_markov_property_discrete_stopping_time} and that $\mathcal A\in \mathcal F^\eta_{\tau_i+}$ for each $i\in \mathbb N$ on the third identity, and dominated convergence together with the fact that by assertion~(iv) of this theorem $P_{T}\phi$ is continuous, and by assertion~(vi) the process $(\eta(t))_{t \geq 0}$ is right-continuous to derive the last identity. We then conclude that~\eqref{eq:strongMarkovinproof} holds for any almost surely finite $\{\mathcal{F}^{\eta}_{t+}\}_{t\ge 0}$-stopping time, which means $(\eta(t))_{t \geq 0}$ is a strong Markov process with respect to the filtration $\{\mathcal{F}^{\eta}_{t+}\}_{t\ge 0}$, as desired.
\end{proof}

\section{Moment estimates} 
\label{Paper01_section_correlation_functions}

Our main goal in this section will be to prove moment estimates on the particle system $\eta^n$ that will be used in the proofs of Theorems~\ref{Paper01_thm_existence_uniqueness_IPS_spatial_muller} and~\ref{Paper01_thm_moments_estimates_spatial_muller_ratchet}.
Recall that as introduced in Section~\ref{Paper01_model_description}, 
$q_+$ and $q_-$ are the birth and death rate polynomials,
$N \in \mathbb{N}$ is a scaling parameter for the carrying capacity, $L > 0$ is the space renormalisation parameter, $m > 0$ is the migration rate, $\mu\in [0,1]$ is the mutation probability, and $(s_k)_{k\in \mathbb N_0}$ is the sequence of fitness parameters for the process $(\eta^{n}(t))_{t \geq 0}$, for each $n \in \mathbb{N}$. Also recall the definitions of $\mathcal S$ and $\mathcal S_0$ in~\eqref{Paper01_definition_state_space_formal} and~\eqref{Paper01_set_initial_configurations} respectively, and
for each $n \in \mathbb{N}$, the definition of $\Lambda_n$ after~\eqref{Paper01_scaling_parameters_for_discrete_approximation}, and that the infinitesimal generator of $\eta^{n}$ is given by $\mathcal{L}^{n}$, defined in~\eqref{Paper01_generator_foutel_etheridge_model_restriction_n} and~\eqref{Paper01_infinitesimal_generator_restriction_n}.

\begin{proposition} \label{Paper01_bound_total_mass}
Take $q_+$ and $q_-$ satisfying Assumption~\ref{Paper01_assumption_polynomials}; then for any $p\in \mathbb N$, there is a non-decreasing function $C_p:[0,\infty)\to [1,\infty)$ such that the following holds.
For any $L>0$, $m>0$, $N\in \mathbb N$, $\mu\in [0,1]$ and $(s_k)_{k\in \mathbb N_0}$ satisfying Assumption~\ref{Paper01_assumption_fitness_sequence}, for any $n\in \mathbb N$ and $\boldsymbol{\eta}\in \mathcal S$, there exists a well-defined $\mathcal S$-valued c\`adl\`ag (strongly) Feller process $(\eta^n(t))_{t\ge 0}=(\eta^n(t,x))_{t\ge 0,\, x\in L^{-1}\mathbb Z}$ with generator $\mathcal L^n$ and with $\eta^n(0)=\boldsymbol{\eta}$ almost surely.
For any~$n \in \mathbb{N}$,~$\boldsymbol{\eta} \in \mathcal{S}$ and $T \geq 0$,
\begin{equation} \label{Paper01_discrete_approximation_does_not_explode}
    \mathbb{E}_{\boldsymbol{\eta}}\Bigg[\sup_{t \leq T} \sum_{x \in \Lambda_n} \vert \vert \eta^{n}(t,x) \vert \vert_{\ell_{1}}\Bigg] < \infty.
\end{equation}
Moreover, for any $n\in \mathbb N$, $\boldsymbol{\eta}=(\eta(x))_{x\in L^{-1}\mathbb Z} \in \mathcal{S}$ and $T \geq 0$,
\begin{equation} \label{Paper01_crude_estimate_moments_general_configuration}
    \sup_{x \in \Lambda_n} \mathbb{E}_{\boldsymbol{\eta}}\left[{\left\vert\left\vert \eta^{n}(T,x) \right\vert\right\vert_{\ell_{1}}^{p}}\right] \leq C_{p}(T) \Bigg(\sup_{x \in \Lambda_n} \, {\vert \vert \eta(x) \vert \vert_{\ell_{1}}^{p}} + N^p\Bigg).
\end{equation}
In particular, for $\boldsymbol{\eta}=(\eta(x))_{x\in L^{-1}\mathbb Z} \in \mathcal{S}_{0}$ and $T\ge 0$,
\begin{equation} \label{Paper01_moments_estimates_uniform_space}
    \sup_{n \in \mathbb{N}} \; \sup_{x \in L^{-1}\mathbb{Z}} \mathbb{E}_{\boldsymbol{\eta}}\left[{\left\vert\left\vert \eta^{n}(T,x) \right\vert\right\vert_{\ell_{1}}^{p}}\right] \leq C_{p}(T) \Bigg(\sup_{x \in L^{-1}\mathbb{Z}} \,{\vert \vert \eta(x) \vert \vert_{\ell_{1}}^{p}} + N^p\Bigg).
\end{equation}
\end{proposition}

\ignore{\begin{remark} \label{Paper01_remark_parameter_N_estimates}
    Observe that we could have stated the result by saying the non-decreasing function $C_{p}: \mathbb{R}_{+} \rightarrow \mathbb{R}_{+}$ does not depend on the parameter $N \in \mathbb{N}$, and simply multiplying both sides of the estimate of Proposition~\ref{Paper01_bound_total_mass} by $N^{p}$. We kept the division by $N^{p}$ on both sides of the estimate since this formulation is used in our companion paper~\cite{madeiraPreprintmuller}, where we derive a functional law of large numbers.
\end{remark}}
As mentioned in Section~\ref{Paper01_subsection_heuristics},
the core technique that we will use in the proof of Proposition~\ref{Paper01_bound_total_mass} is the method of correlation functions developed by Boldrighini, De Masi, Pellegrinotti and Presutti in \cite{boldrighini1987collective, demasi1991mathematical}, where it was used to characterise the local behaviour of reaction-diffusion systems in which both birth and death rates are polynomial functions of the local number of particles. The main differences between the model in~\cite{boldrighini1987collective, demasi1991mathematical} and our particle system are that we scale the birth and death rates using the parameter $N \in \mathbb{N}$, and that we have to deal with infinitely many types of particles.

To prove Proposition~\ref{Paper01_bound_total_mass}, our strategy will be to study an auxiliary process, in which all particles have no mutations, which will be a stochastic upper bound for $\eta^n$. For $n \in \mathbb{N}$, we write this auxiliary process as 
$$(\zeta^{n}(t))_{t\geq 0} = (\zeta^{n}(t,x), x \in \Lambda_n)_{t \geq 0},$$
taking values in $\left(\mathbb{N}_{0}\right)^{\Lambda_n}$. 
To simplify notation, for $n\in \mathbb N$, let
\begin{equation} \label{Paper01_modified_state_space_simplified_version_without_types}
    \mathcal{X}^{n} \defeq \left(\mathbb{N}_{0}\right)^{\Lambda_n}.
\end{equation}
For an initial particle configuration
$\boldsymbol{\eta}=(\eta(x))_{x\in L^{-1}\mathbb Z} \in \mathcal{S}$, we take the initial state of $(\zeta^{n}(t))_{t\geq 0}$ as
$\zeta^n(0)=(\zeta(x))_{x\in L^{-1}\mathbb Z}\in \mathcal{X}^n$, where
\begin{equation} \label{Paper01_deriving_initial_condition_auxiliary_process}
    \zeta(x) = \vert \vert \eta(x) \vert \vert_{\ell_1} \; \forall \, x \in \Lambda_n.
\end{equation}
Recall the definitions of~$q^N_+$ and~$q^N_-$ from~\eqref{Paper01_scaled_polynomials_carrying_capacity}. For~$n \in \mathbb{N}$, the process $\left(\zeta^{n}(t)\right)_{t \geq 0}$ evolves as follows:
\begin{itemize}
    \item \textit{Migration events}: For each $t \geq 0$ and $x \in \Lambda_n$, each particle living in deme $x$ at time $t$ independently attempts to make a jump at rate $m$ to a uniformly chosen deme from $\{x - L^{-1}, \, x + L^{-1}\}$. If the chosen deme lies outside $\Lambda_n$, the particle immediately returns to its original position. Otherwise, the jump is carried out.
    \item \textit{Reproduction events}: For each $t \geq 0$ and $x \in \Lambda_n $, each particle living in deme $x$ at time $t$ reproduces independently at rate $ q^N_{+}\left({\zeta^{n}(t-,x)}\right).$ After such an event, a new particle is added in deme $x$.
    \item \textit{Death events}:  For each $t \geq 0$ and $x \in \Lambda_n$, each particle living in deme $x$ at time $t$ dies independently at rate $ q_{-}^N\left({\zeta^{n}(t-,x)}\right).$ When a particle dies, it is simply removed from the process.
\end{itemize}

Note that, in contrast to the process $(\eta^{n}(t))_{t \geq 0}$, the process $(\zeta^{n}(t))_{t \geq 0}$ is monotone. Moreover, observe that on the box~$\Lambda_{n}$, i.e.~on the spatial region where particles in $\eta^{n}$ are not frozen, the death and migration rates for particles in $\zeta^{n}$ are the same as those for particles in $\eta^{n}$. Since the birth rates for particles in $(\zeta^{n}(t))_{t \geq 0}$ are bounded below by the rates for particles in $(\eta^{n}(t))_{t \geq 0}$ due to the mutations in $\eta^n$ being deleterious (more precisely because $s_{k} \leq 1$ for all $k \in \mathbb{N}$), for $\boldsymbol{\eta}\in \mathcal S$ and $\boldsymbol{\zeta}\in \mathcal X^n$ as defined in~\eqref{Paper01_deriving_initial_condition_auxiliary_process},  for each $n \in \mathbb{N}$, we can couple the processes $(\zeta^{n}(t))_{t \geq 0}$ and $(\eta^{n}(t))_{t \geq 0}$ in such a way that $\zeta^n(0)=\boldsymbol{\zeta}$, $\eta^n(0)=\boldsymbol{\eta}$ and
\begin{equation} \label{eq:zetaaboveeta}
    \zeta^{n}(t,x) \geq \vert \vert \eta^{n}(t,x) \vert \vert_{\ell_{1}}
\end{equation}
for all $t \geq 0$ and $x \in \Lambda_{n}$ almost surely. Then, by establishing moment bounds for $\zeta^n$, we will have corresponding bounds for $\eta^n$, which imply that $\eta^n$ can be constructed by standard arguments as a well-posed càdlag $\mathbb{N}_0^{d}$-valued Markov process for some $d\in \mathbb N$ (see~\cite[Theorem~6.4.1]{ethier2009markov} for this construction). Hence, to prove Proposition~\ref{Paper01_bound_total_mass}, it will suffice to prove the corresponding moment bounds for the process $(\zeta^{n}(t))_{t \geq 0}$, for every~$n \in \mathbb{N}$.

Since there is only one type of particle in $(\zeta^{n}(t))_{t \geq 0}$, we can adapt the ideas of De Masi and Presutti in~\cite[Chapter~IV]{demasi1991mathematical}. Following their approach, we start by noting that it is not clear a priori that the number of particles in $(\zeta^{n}(t))_{t \geq 0}$ does not explode in finite time. To avoid this problem, we construct a sequence of processes $(\zeta^{n, \kappa})_{\kappa \in \mathbb{N}}$ that will converge weakly to $\zeta^{n}$ as $\kappa \rightarrow \infty$. The parameter $\kappa \in \mathbb{N}$ will indicate the maximum local population density at which reproduction events will be allowed to occur. More precisely, for each $n,\kappa \in \mathbb{N}$ we let
\begin{equation} \label{eq:zetankappadefn}
    (\zeta^{n, \kappa}(t))_{t\geq 0} = (\zeta^{n,\kappa}(t,x), x \in \Lambda_n)_{t \geq 0}
\end{equation}
be a Markov process taking values in $\mathcal X^n$ with $\zeta^{n,\kappa}(0)=\zeta^n(0)$ defined as follows.
Migration and death events occur with the same rates as in the definition of $(\zeta^{n}(t))_{t \geq 0}$ above. 
For reproduction events, the rates differ from the definition of $\zeta^n$.
For each $t \geq 0$ and $x \in \Lambda_n$, each particle living in deme $x$ at time $t$ reproduces independently at rate $\displaystyle q^{N,\kappa}_{+}\left({\zeta^{n, \kappa}(t-,x)}\right)$, where $q^{N,\kappa}_{+}:[0,\infty)\to [0,\infty)$ is given by
\begin{equation*}
        q^{N,\kappa}_{+}(i) = \left\{\begin{array}{ccc} q_{+}(i/N) & \textrm{if} & i \leq N\kappa, \\ 0 & \textrm{if} & i > N\kappa. \end{array}\right.
\end{equation*}
After such an event, a new particle is added to deme $x$.

For $\boldsymbol{\eta}\in \mathcal S$ and $n,\kappa \in \mathbb N$, write $\mathbb P_{\boldsymbol{\eta}}$ for the probability measure under which $\zeta^{n,\kappa}(0)=\boldsymbol{\zeta}$, where $\boldsymbol{\zeta}=(\zeta(x))_{x\in \Lambda_n}$ is given by~\eqref{Paper01_deriving_initial_condition_auxiliary_process}, and write $\mathbb E_{\boldsymbol{\eta}}$ for the corresponding expectation.
Note that for each $n, \kappa \in \mathbb{N}$ and $\boldsymbol{\zeta}\in \mathcal X^n$, the process $\left(\zeta^{n, \kappa}(t)\right)_{t \geq 0}$ conditioned on $\zeta^{n,\kappa}(0)=\boldsymbol{\zeta}$ is well defined, since the birth rate of each particle is bounded above by $\max_{i \in [0, \kappa]}q_{+}(i)$, and the initial number of particles living in $\Lambda_n$ is finite (see~\cite[Theorem~6.4.1]{ethier2009markov}). For the same reason, for any $\boldsymbol{\eta} \in \mathcal{S}$, any $n, \kappa \in \mathbb{N}$, any $p \geq 1$ and any $T \geq 0$ we have
\begin{equation} \label{Paper01_well_posedness_truncated_process}
    \mathbb{E}_{\boldsymbol{\eta}}\Bigg[\sup_{0 \leq t \leq T} \Bigg(\sum_{x \in \Lambda_n} \zeta^{n, \kappa}(t,x)\Bigg)^{p}\Bigg] < \infty.
\end{equation}
By~\eqref{Paper01_well_posedness_truncated_process}, $(\zeta^{n, \kappa}(t))_{t \geq 0}$ is well defined as a jump process taking values in $\mathcal X^n$. For $n,\kappa\in \mathbb N$, let $\hat{\mathcal{L}}^{n}$ and $\hat{\mathcal{L}}^{n,\kappa}$ be the infinitesimal generators of $\zeta^{n}$ and of $\zeta^{n, \kappa}$ respectively. For any $\phi \in \mathscr{C}(\mathcal X^n, \mathbb{R})$ and any $\boldsymbol{\zeta} \in \mathcal X^n$, we can write
\begin{equation} \label{Paper01_infinitesimal_generator_auxiliary_process}
    \hat{\mathcal{L}}^{n} \phi(\boldsymbol{\zeta}) = \frac{m}{2} \hat{\mathcal{L}}^{n}_{m} \phi(\boldsymbol{\zeta}) + \hat{\mathcal{L}}^{n}_{r} \phi(\boldsymbol{\zeta}) \quad \textrm{ and } \quad \hat{\mathcal{L}}^{n,\kappa} \phi(\boldsymbol{\zeta}) = \frac{m}{2} \hat{\mathcal{L}}_{m}^{n} \phi(\boldsymbol{\zeta}) + \hat{\mathcal{L}}^{n,\kappa}_{r} \phi(\boldsymbol{\zeta}),
\end{equation}
where $\hat{\mathcal{L}}_{m}^{n}$ indicates the migration part (which is identical for both generators), and $\hat{\mathcal{L}}^{n}_{r}$ and $\hat{\mathcal{L}}^{n, \kappa}_{r}$ indicate the reaction parts of the generators of $\zeta^{n}$ and $\zeta^{n,\kappa}$, respectively. For conciseness, we do not write the definitions of $\hat{\mathcal{L}}_{m}^{n}$, $\hat{\mathcal{L}}_{r}^{n}$ and $\hat{\mathcal{L}}^{n, \kappa}_{r}$ explicitly.

We will establish estimates for $\left(\zeta^{n, \kappa}(t)\right)_{t \geq 0}$ that hold uniformly in $\kappa \in \mathbb{N}$. Bearing this aim in mind, we start by modifying the definition of Poisson polynomials given in \cite{demasi1991mathematical, boldrighini1987collective, perrut2000hydrodynamic} in the following way.

\begin{definition} \label{Paper01_poisson_polynomial}
For $j \in \mathbb Z$ and $N \in \mathbb{N}$, we define the \emph{scaled Poisson polynomial} of order $j$ and scaling parameter $N$ as the function $D^{N}_{j}: \mathbb Z \rightarrow \mathbb{Q}$ given by, for all $i \in \mathbb Z$,
\begin{equation*}
    D^{N}_{j}(i) = \left\{\begin{array}{lc}
    1 & \textrm{if } j \le 0 \text{ or }i<0, \\
    \displaystyle \frac{i  (i-1) \cdots  (i-j+1)}{N^{j}} & \textrm{otherwise.} \end{array}\right.
\end{equation*}
\end{definition}

Note that the denominator $N^{j}$ does not appear in the definition of Poisson polynomials in~\cite[Equation~2.53]{demasi1991mathematical}. For each $N,n \in \mathbb{N}$, we define the function
\begin{equation} \label{Paper01_definition_duality_function}
\begin{aligned}
    \mathfrak{D}^{N,n}: \mathcal{X}^{n}  \times \mathcal{X}^{n} & \rightarrow \mathbb{R}_{+} \\ (\boldsymbol{\xi}, \boldsymbol{\zeta})=((\xi(x))_{x\in \Lambda_n},(\zeta(x))_{x\in \Lambda_n}) & \mapsto \mathfrak{D}^{N,n}(\boldsymbol{\xi}, \boldsymbol{\zeta}) \defeq \prod_{x \in \Lambda_n} D^{N}_{\xi(x)}(\zeta(x)),
\end{aligned}
\end{equation}
where $D^{N}_j$ is given in Definition~\ref{Paper01_poisson_polynomial}. 

The map $ \mathfrak{D}^{N,n}$ will be used as a duality function between the process $\zeta^{n,\kappa}$ and a migration process $(\xi^{n}(t))_{t \geq 0}$ taking values in $\mathcal{X}^{n}$. The infinitesimal generator of $(\xi^{n}(t))_{t \geq 0}$ is given by $\frac{m}{2}\hat{\mathcal{L}}_{m}^{n}$; recall from~\eqref{Paper01_infinitesimal_generator_auxiliary_process} that this is the migration part of the infinitesimal generator of $(\zeta^{n,\kappa}(t))_{t \geq 0}$. In other words, particles in $(\xi^{n}(t))_{t \geq 0}$ can only undergo migration events, which are described as follows. For any $t \geq 0$ and $x \in \Lambda_n$, each particle living in deme $x$ at time $t$ independently attempts to make a jump at rate $m$, to a uniformly chosen deme from $\{x - L^{-1}, \, x + L^{-1}\}$. If the chosen deme lies outside $\Lambda_n$, the particle immediately returns to its original position; otherwise, the jump is carried out. As there are no birth events, the total number of particles in $(\xi^{n}(t))_{t \geq 0}$ remains finite for all $t \geq 0$. In order to properly introduce the duality relation between $(\xi^{n}(t))_{t \geq 0}$ and $(\zeta^{n,\kappa}(t))_{t \geq 0}$, we start by observing that the action of the generator $\hat{\mathcal{L}}^n_{m}$ on the function $\mathfrak{D}^{N,n}(\boldsymbol{\xi}, \boldsymbol{\zeta})$ is symmetric with respect to the two variables.

\begin{lemma} \label{Paper01_duality_xi_Lambda_N}
For every $n \in \mathbb{N}$ and all $\boldsymbol{\xi}, \boldsymbol{\zeta} \in \mathcal{X}^{n}$,
\begin{equation*}
    \left(\frac{m}{2}\hat{\mathcal{L}}_{m}^{n}\mathfrak{D}^{N,n}\left(\boldsymbol{\xi},\cdot \right)\right)\left(\boldsymbol{\zeta}\right) = \left(\frac{m}{2}\hat{\mathcal{L}}^{n}_{m}\mathfrak{D}^{N,n}\left(\cdot,\boldsymbol{\zeta} \right)\right)\left(\boldsymbol{\xi}\right).
\end{equation*}
\end{lemma}

We highlight that this result (for unscaled Poisson polynomials) is stated in the proof of Proposition~2.9.4 in~\cite{demasi1991mathematical}. Since the proof of this fact is not provided in~\cite{demasi1991mathematical}, we prove Lemma~\ref{Paper01_duality_xi_Lambda_N} below.

\begin{proof}
Take $\boldsymbol{\zeta}=(\zeta(x))_{x\in \Lambda_n}$, $\boldsymbol{\xi}=(\xi(x))_{x\in \Lambda_n}\in \mathcal X^n$.
We will directly compute the action of $\hat{\mathcal{L}}^{n}_{m}$ on each of the variables of $\mathfrak{D}^{N,n}(\boldsymbol{\xi}, \boldsymbol{\zeta})$. Recall from after~\eqref{Paper01_scaling_parameters_for_discrete_approximation} that $\Lambda_n = L^{-1}\mathbb{Z} \cap [-\lambda_n, \lambda_n]$. By~\eqref{Paper01_definition_duality_function}, after rearranging terms, note that
\begin{equation} \label{eq:generatoronDNn}
\left(\hat{\mathcal{L}}^{n}_{m}\mathfrak{D}^{N,n}(\boldsymbol{\xi}, \cdot)\right)(\boldsymbol{\zeta}) = \sum_{x \in \Lambda_n\setminus \{L^{-1} \lfloor L \lambda_n \rfloor\}} \Bigg(\prod_{ x^{*} \in \Lambda_n\setminus  \{x,  \, x + L^{-1}\}} D^{N}_{\xi(x^{*})}(\zeta({x^{*}}))\Bigg) \cdot \Upsilon^{N}(\boldsymbol{\xi}, \boldsymbol{\zeta}, x),
\end{equation}
where, for $x\in \Lambda_n\setminus \{L^{-1} \lfloor L \lambda_n \rfloor\}$, $\Upsilon^{N}(\boldsymbol{\xi}, \boldsymbol{\zeta}, x)$ is given by
\begin{align} \label{Paper01_definition_upsilon}
    & \Upsilon^{N}(\boldsymbol{\xi}, \boldsymbol{\zeta}, x)  \notag \\
    & \quad \defeq \zeta(x) \Big(D^{N}_{\xi(x + L^{-1})}(\zeta(x + L^{-1}) + 1)D^{N}_{\xi(x)}(\zeta(x) - 1)  - D^{N}_{\xi(x + L^{-1})}(\zeta(x + L^{-1}))D^{N}_{\xi(x)}(\zeta({x}))\Big) \notag \\
    & \quad \quad  \quad + \zeta(x + L^{-1})\Big(D^{N}_{\xi(x+L^{-1})}(\zeta(x+L^{-1}) - 1)D^{N}_{\xi({x})}(\zeta({x}) + 1)  \notag \\
    & \quad \quad \quad \quad \quad \quad \quad \quad\quad \quad \quad \quad \quad \quad - D^{N}_{\xi(x+L^{-1})}(\zeta(x + L^{-1}))D^{N}_{\xi(x)}(\zeta(x))\Big).
\end{align}
By Definition~\ref{Paper01_poisson_polynomial}, in the case $\zeta(x)\neq 0$ and $\xi(x),$ $\xi(x+L^{-1})\notin \{0,1\}$, we can rewrite the term in brackets on the first line of the definition of $\Upsilon^{N}(\boldsymbol{\xi}, \boldsymbol{\zeta}, x)$ as
\begin{align*}
    & D^{N}_{\xi(x + L^{-1})}(\zeta(x + L^{-1}) + 1)D^{N}_{\xi(x)}(\zeta(x) - 1)  - D^{N}_{\xi(x + L^{-1})}(\zeta(x + L^{-1}))D^{N}_{\xi(x)}(\zeta({x})) \\
    & \quad = \frac{(\zeta(x + L^{-1}) + 1) \cdots (\zeta(x+L^{-1}) - \xi(x + L^{-1}) + 2) \cdot (\zeta(x) - 1) \cdots (\zeta(x) - \xi(x))}{N^{\xi(x + L^{-1})} \cdot N^{\xi(x)}} \\ 
    & \quad \quad \quad - \frac{\zeta(x + L^{-1}) \cdots (\zeta(x+L^{-1}) - \xi(x + L^{-1}) + 1) \cdot \zeta(x) \cdots (\zeta(x) - \xi(x) + 1)}{N^{\xi(x + L^{-1})} \cdot N^{\xi(x)}} \\ 
    & \quad = \left(\frac{\zeta(x + L^{-1}) \cdots (\zeta(x+L^{-1}) - \xi(x + L^{-1}) + 2) \cdot  (\zeta(x) - 1) \cdots  (\zeta(x) - \xi(x) + 1)}{N^{2 + (\xi(x + L^{-1})-1) + (\xi(x) - 1)}}\right) \\ & \quad \quad \quad \quad \quad \quad \quad \cdot\Big((\zeta(x + L^{-1}) + 1)(\zeta(x) - \xi(x)) - \zeta(x)(\zeta(x + L^{-1})-\xi(x + L^{-1}) + 1)\Big).
\end{align*}
Therefore, applying Definition~\ref{Paper01_poisson_polynomial} again to the first term in brackets on the right-hand side, in the case $\zeta(x)\neq 0$ and $\xi(x),$ $\xi(x+L^{-1})\notin \{0,1\}$,
\begin{align} \label{eq:DNxizetacalc}
     & D^{N}_{\xi(x + L^{-1})}(\zeta(x + L^{-1}) + 1)D^{N}_{\xi(x)}(\zeta(x) - 1)  - D^{N}_{\xi(x + L^{-1})}(\zeta(x + L^{-1}))D^{N}_{\xi(x)}(\zeta({x})) \notag \\ 
     & \quad = \frac{D^{N}_{\xi(x + L^{-1})-1}(\zeta(x + L^{-1}))D^{N}_{\xi(x)-1}(\zeta(x) - 1) }{N^{2}} \notag \\ 
     & \quad \quad \quad \cdot \Big((\zeta(x + L^{-1}) + 1)(\zeta(x) - \xi(x)) - \zeta(x)(\zeta(x + L^{-1})-\xi(x + L^{-1}) + 1)\Big) \notag \\ 
     & \quad = \frac{D^{N}_{\xi(x + L^{-1})-1}(\zeta(x + L^{-1}))D^{N}_{\xi(x)-1}(\zeta(x) - 1) }{N^{2}} \Big(\zeta(x)\xi(x + L^{-1}) - \xi(x)(\zeta(x + L^{-1}) + 1)\Big).
\end{align}
By repeating the same calculation in the cases $\xi(x)=1$ and $\xi(x+L^{-1})=1$, we see that~\eqref{eq:DNxizetacalc} holds as long as
$\zeta(x)\neq 0$ and $\xi(x),$ $\xi(x+L^{-1})\neq 0$.
Multiplying both sides of~\eqref{eq:DNxizetacalc} by $\zeta(x)$, and again using Definition~\ref{Paper01_poisson_polynomial}
(considering the case $\xi(x)=1$ separately), we obtain that in the case $\zeta(x)\neq 0$ and $\xi(x),$ $\xi(x+L^{-1})\neq 0$,
\begin{align} \label{Paper01_first_step_duality_argument_correlation}
     & \zeta(x)\Big(D^{N}_{\xi(x + L^{-1})}(\zeta(x + L^{-1}) + 1)D^{N}_{\xi(x)}(\zeta(x) - 1)  - D^{N}_{\xi(x + L^{-1})}(\zeta(x + L^{-1}))D^{N}_{\xi(x)}(\zeta({x}))\Big) \notag \\ 
     & \quad = \frac{D^{N}_{\xi(x + L^{-1})-1}(\zeta(x + L^{-1}))D^{N}_{\xi(x)-1}(\zeta(x)) }{N^{2}} \notag \\ 
     & \quad \quad \quad \cdot \Big(\zeta(x) - \xi(x) + 1\Big)\Big(\zeta(x)\xi(x + L^{-1}) - \xi(x)(\zeta(x + L^{-1}) + 1)\Big).
\end{align}
Note that in fact~\eqref{Paper01_first_step_duality_argument_correlation} also holds in the case $\zeta(x)= 0$ and $\xi(x),$ $\xi(x+L^{-1})\neq 0$, because both sides are equal to $0$.
By an analogous argument, in the case $\xi(x),$ $\xi(x+L^{-1})\neq 0$,
\begin{align} \label{Paper01_second_step_duality_argument_correlation}
    & \zeta(x + L^{-1})\Big(D^{N}_{\xi(x+L^{-1})}(\zeta(x+L^{-1}) - 1)D^{N}_{\xi({x})}(\zeta({x}) + 1) - D^{N}_{\xi(x+L^{-1})}(\zeta(x + L^{-1}))D^{N}_{\xi(x)}(\zeta(x)) \Big) \notag \\ 
    & \quad = \frac{D^{N}_{\xi(x + L^{-1})-1}(\zeta(x + L^{-1}))D^{N}_{\xi(x)-1}(\zeta(x)) }{N^{2}} \notag \\ 
    & \quad \quad \quad \cdot \Big(\zeta(x+L^{-1}) - \xi(x+L^{-1}) + 1\Big)\Big(\zeta(x+L^{-1}) \xi(x) - \xi(x+L^{-1})(\zeta(x) +1)\Big).
\end{align}
Hence, in the case $\xi(x),$ $\xi(x + L^{-1}) \neq 0$, by applying~\eqref{Paper01_first_step_duality_argument_correlation} and~\eqref{Paper01_second_step_duality_argument_correlation} to~\eqref{Paper01_definition_upsilon}, $\Upsilon^{N}(\boldsymbol{\xi}, \boldsymbol{\zeta}, x)$ can be rewritten as
\begin{equation} \label{Paper01_proof_duality_first_with_lambda_ii}
    \Upsilon^{N}(\boldsymbol{\xi}, \boldsymbol{\zeta}, x) = \frac{D^{N}_{\xi(x)-1}(\zeta(x))D^{N}_{\xi(x + L^{-1})-1}(\zeta(x+L^{-1}))}{N^{2}} \cdot Y^{N}(\boldsymbol{\xi}, \boldsymbol{\zeta}, x),
\end{equation}
where $Y^{N}(\boldsymbol{\xi}, \boldsymbol{\zeta}, x)$ is given by
\begin{equation*}
\begin{aligned}
    & Y^{N}(\boldsymbol{\xi}, \boldsymbol{\zeta}, x) \\ & \quad = \Big(\zeta(x) - \xi(x) + 1\Big)\Big(\zeta(x) \xi(x+L^{-1}) - \zeta(x+L^{-1})\xi(x) - \xi(x)\Big) \\ & \quad \quad \quad + \Big(\zeta(x+L^{-1}) - \xi(x+L^{-1}) + 1\Big)\Big(\zeta(x+L^{-1}) \xi(x) - \zeta(x)\xi(x+L^{-1}) - \xi(x+L^{-1})\Big).
\end{aligned}
\end{equation*}
On the other hand, by~\eqref{Paper01_definition_duality_function} again,
\begin{equation} \label{eq:genDNndag}
\left(\hat{\mathcal{L}}^{n}_{m}\mathfrak{D}^{N,n}(\cdot, \boldsymbol{\zeta})\right)(\boldsymbol{\xi}) =  \sum_{x \in \Lambda_n\setminus \{L^{-1} \lfloor L \lambda_n \rfloor\}} \Bigg(\prod_{x^{*} \in \Lambda_n\setminus \{x, \, x + L^{-1}\}} D^{N}_{\xi(x^{*})}(\zeta(x^{*}))\Bigg) \cdot \widetilde{\Upsilon}^{N}(\boldsymbol{\xi}, \boldsymbol{\zeta}, x),
\end{equation}
where, for $x \in \Lambda_n\setminus \{L^{-1} \lfloor L \lambda_n \rfloor\}$,  $\widetilde{\Upsilon}^{N}(\boldsymbol{\xi}, \boldsymbol{\zeta}, x)$ is given by
\begin{align} \label{eq:tildeUpsdefn}
   & \widetilde{\Upsilon}^{N}(\boldsymbol{\xi}, \boldsymbol{\zeta}, x) \notag \\ 
   & \quad \defeq \xi(x) \Big(D^{N}_{\xi(x + L^{-1})+1}(\zeta(x + L^{-1}))D^{N}_{\xi(x) - 1}(\zeta(x)) - D^{N}_{\xi(x + L^{-1})}(\zeta(x+L^{-1}))D^{N}_{\xi(x)}(\zeta(x))\Big) \notag \\ 
   & \quad \quad \quad + \xi(x + L^{-1})\Big(D^{N}_{\xi(x + L^{-1})-1}(\zeta(x + L^{-1}))D^{N}_{\xi(x)+1}(\zeta(x)) \notag \\ 
   & \quad \quad \quad \quad \quad \quad \quad \quad \quad \quad \quad - D^{N}_{\xi(x + L^{-1})}(\zeta(x + L^{-1}))D^{N}_{\xi(x)}(\zeta(x))\Big).
\end{align}

By~\eqref{eq:generatoronDNn} and~\eqref{eq:genDNndag}, to complete the proof it will suffice to verify that $ {\Upsilon}^{N}(\boldsymbol{\xi}, \boldsymbol{\zeta}, x) =  \widetilde{\Upsilon}^{N}(\boldsymbol{\xi}, \boldsymbol{\zeta}, x)$ for each
$x \in \Lambda_n\setminus \{L^{-1} \lfloor L \lambda_n \rfloor\}$.
By an argument analogous to the one we used to derive~\eqref{Paper01_proof_duality_first_with_lambda_ii}, in the case $\xi(x),\xi(x + L^{-1}) \neq 0$, we can rewrite $\widetilde{\Upsilon}^{N}(\boldsymbol{\xi}, \boldsymbol{\zeta}, x)$, obtaining
\begin{equation*}
\begin{aligned}
    \widetilde{\Upsilon}^{N}(\boldsymbol{\xi}, \boldsymbol{\zeta}, x) = \frac{D^{N}_{\xi(x)-1}(\zeta(x))D^{N}_{\xi(x + L^{-1})-1}(\zeta(x+L^{-1}))}{N^{2}} \cdot {Y}^{N}(\boldsymbol{\xi}, \boldsymbol{\zeta}, x) = {\Upsilon}^{N}(\boldsymbol{\xi}, \boldsymbol{\eta}, x).
\end{aligned}
\end{equation*}
Hence, we conclude that $\Upsilon^{N}(\boldsymbol{\xi}, \boldsymbol{\zeta}, x) = \widetilde{\Upsilon}^{N}(\boldsymbol{\xi}, \boldsymbol{\zeta}, x)$ for any $x \in \Lambda_n\setminus \{L^{-1} \lfloor L \lambda_n \rfloor\}$
such that $\xi(x), \xi(x + L^{-1}) \neq 0$. In the case $\xi(x) \wedge \xi(x + L^{-1}) =0$, 
the expressions for $\Upsilon^{N}$ and $\widetilde{\Upsilon}^{N}$ in~\eqref{Paper01_definition_upsilon} and~\eqref{eq:tildeUpsdefn} simplify, and by an easy calculation
we obtain $\Upsilon^{N}(\boldsymbol{\xi}, \boldsymbol{\zeta}, x) = \widetilde{\Upsilon}^{N}(\boldsymbol{\xi}, \boldsymbol{\zeta}, x)$ for every $x \in \Lambda_n\setminus \{L^{-1} \lfloor L \lambda_n \rfloor\}$, which completes the proof.
\end{proof}

We now define a family of correlation functions; recall the definition of $\mathfrak D^{N,n}$ in~\eqref{Paper01_definition_duality_function}.

\begin{definition} \label{Paper01_correlation_functions}
For $N,n,\kappa \in \mathbb{N}$, we define the \emph{scaled correlation function} $\rho^{N,n, \kappa}$ by
\begin{equation*}
\begin{aligned}
    \rho^{N,n, \kappa}: \mathcal{S} \times \mathcal{X}^{n} \times [0, \infty) & \rightarrow [0,\infty) \\
    (\boldsymbol{\eta}, \boldsymbol{\xi}, T) & \mapsto \rho^{N,n, \kappa}(\boldsymbol{\eta}, \boldsymbol{\xi}, T) = \mathbb{E}_{\boldsymbol{\eta}}\left[\mathfrak{D}^{N,n}\left(\boldsymbol{\xi}, \zeta^{n,\kappa}(T)\right)\right].
\end{aligned}
\end{equation*}
\end{definition}

Recall from above~\eqref{Paper01_well_posedness_truncated_process} that we use the subscript $\boldsymbol{\eta}$ to indicate that the initial condition of $(\zeta^{n,\kappa}(t))_{t \geq 0}$ is given by $\boldsymbol{\zeta}$ as defined in~\eqref{Paper01_deriving_initial_condition_auxiliary_process}. We also highlight that Definition~\ref{Paper01_correlation_functions} is analogous to Equation~(4.12) in~\cite{demasi1991mathematical}. Moreover, note that due to~\eqref{Paper01_well_posedness_truncated_process}, for any $N,n,\kappa \in \mathbb N$ we have $\rho^{N,n,\kappa}(\boldsymbol{\eta},\boldsymbol{\xi},T)<\infty$ $\forall \boldsymbol{\eta}\in \mathcal S$, $\boldsymbol{\xi} \in \mathcal{X}^{n}$ and $T\ge 0$. We now state and prove a result analogous to Lemma~4.2.2 in~\cite{demasi1991mathematical}.
Recall the definition of the migration process $\xi^n$ before the statement of Lemma~\ref{Paper01_duality_xi_Lambda_N}. 
\begin{lemma} \label{Paper01_rewriting_definition_correlation_function}
For $N,n,\kappa \in \mathbb N$,
for any $\boldsymbol{\eta} \in \mathcal{S}$, any $\boldsymbol{\xi} \in \mathcal{X}^{n}$, and any $T > 0$, 
\begin{equation*}
\begin{aligned}
    \rho^{N, n, \kappa}(\boldsymbol{\eta}, \boldsymbol{\xi}, T) = & \sum_{\bar{\boldsymbol{\xi}} \in \mathcal{X}^{n}} \mathbb{P}_{\boldsymbol{\xi} }(\xi^n(T) = \bar{\boldsymbol{\xi}}) \rho^{N,n, \kappa}(\boldsymbol{\eta}, \bar{\boldsymbol{\xi} }, 0) \\ & + \int_{0}^{T} \Bigg(\sum_{\bar{\boldsymbol{\xi} } \in \mathcal{X}^{n}} \mathbb{P}_{\boldsymbol{\xi} }(\xi^n(T - t) = \bar{\boldsymbol{\xi} })\mathbb{E}_{\boldsymbol{\eta}}\left[\Big(\hat{\mathcal{L}}^{n,\kappa}_{r}\mathfrak{D}^{N,n}(\bar{\boldsymbol{\xi} }, \cdot)\Big)\left(\zeta^{n, \kappa}(t)\right)\right]\Bigg) dt,
\end{aligned}
\end{equation*}
where $\hat{\mathcal{L}}^{n,\kappa}_{r}$ is the reaction part of the generator of $\zeta^{n,\kappa}$ as defined in~\eqref{Paper01_infinitesimal_generator_auxiliary_process}.
\end{lemma}

\begin{proof}
Recall that $\hat{\mathcal{L}}^{n,\kappa}$ is the generator of $\zeta^{n,\kappa}$, and recall
the definition of $\hat{\mathcal{L}}_{m}^{n}$ in~\eqref{Paper01_infinitesimal_generator_auxiliary_process}. Our strategy will be to prove our claim by using a standard stochastic chain rule for time-dependent test functions (see Lemma~\ref{Paper01_general_integration_by_parts_formula}). For $T \ge 0$, $n \in \mathbb{N}$ and $\boldsymbol{\xi} \in \mathcal{X}^{n}$, we define the map $g^{n,T}(\cdot, \boldsymbol{\xi}, \cdot): [0,T] \times \mathcal{X}^{n} \rightarrow [0,\infty)$ by
\begin{equation} \label{Paper01_definition_auxiliary_map_chain_rule_formula_correlation_function}
    (t,\boldsymbol{\zeta}) \mapsto g^{n,T}(\boldsymbol{\zeta}, \boldsymbol{\xi}, t) \defeq \sum_{\bar{\boldsymbol{\xi}} \in \mathcal{X}^{n}} \mathbb{P}_{\boldsymbol{\xi}}(\xi^n(T-t) = \bar{\boldsymbol{\xi}}) \mathfrak{D}^{N,n}(\bar{\boldsymbol{\xi}}, \boldsymbol{\zeta}).
\end{equation}
Observe that by Definition~\ref{Paper01_correlation_functions}, for $N,n,\kappa\in \mathbb N$, $\boldsymbol{\eta}\in \mathcal S$, $\boldsymbol{\xi}\in \mathcal X^n$ and $T\ge 0$ we can write
\begin{equation} \label{Paper01_equivalence_test_function_correlation_funciton}
    \rho^{N,n,\kappa}(\boldsymbol{\eta}, \boldsymbol{\xi},T) = \mathbb{E}_{\boldsymbol{\eta}}\Big[g^{n,T}(\zeta^{n,\kappa}(T), \boldsymbol{\xi}, T)\Big].
\end{equation}
Now fix $N,n,\kappa\in \mathbb N$, $\boldsymbol{\eta}\in \mathcal S$, $\boldsymbol{\xi}\in \mathcal X^n$ and $T> 0$.
We claim that the following statements hold:
\begin{enumerate}[(i)]
    \item $\displaystyle \sup_{t_{1}, t_{2} \in [0,T]} \, \mathbb{E}_{\boldsymbol{\eta}}\Big[g^{n,T}(\zeta^{n,\kappa}(t_1), \boldsymbol{\xi}, t_2)\Big] < \infty$.
    \item The map $(t, \boldsymbol{\zeta}) \mapsto \frac{\partial}{\partial t} (g^{n,T}(\boldsymbol{\zeta}, \boldsymbol{\xi}, \cdot))(t)$ is continuous with respect to the product topology on $(0,T) \times \mathcal{X}^{n}$, and satisfies the following identity, for $t \in (0,T)$ and $\boldsymbol{\zeta} \in \mathcal{X}^{n}$:
    \begin{equation*}
        \frac{\partial}{\partial t} (g^{n,T}(\boldsymbol{\zeta}, \boldsymbol{\xi}, \cdot))(t) = - \frac{m}{2} \sum_{\bar{\boldsymbol{\xi}} \in \mathcal{X}^{n}} \mathbb{P}_{\boldsymbol{\xi}}(\xi^n(T-t) = \bar{\boldsymbol{\xi}}) (\hat{\mathcal{L}}^{n}_{m}\mathfrak{D}^{N,n}(\cdot, \boldsymbol{\zeta}))(\bar{\boldsymbol{\xi}}).
    \end{equation*}
    \item For any $t_{1} \in [0,T]$, the process $(M^{n,\kappa,T,t_1}_{\boldsymbol{\xi}}(t))_{t \geq 0}$ given by
    \begin{equation*}
        M^{n,\kappa,T,t_1}_{\boldsymbol{\xi}}(t) \defeq g^{n,T}(\zeta^{n,\kappa}(t), \boldsymbol{\xi},t_1) - g^{n,T}(\zeta^{n,\kappa}(0), \boldsymbol{\xi},t_1) - \int_{0}^{t} \Big(\hat{\mathcal{L}}^{n,\kappa}g^{n,T}(\cdot, \boldsymbol{\xi},t_1)\Big)(\zeta^{n,\kappa}(t'-)) dt'
    \end{equation*}
    is a càdlàg martingale with respect to the filtration $\Big\{\mathcal{F}^{\zeta^{n,\kappa}}_{t+}\Big\}_{t \geq 0}$.
    \item The map $[0,T] \times \mathcal{X}^n \ni (t,\boldsymbol{\zeta}) \mapsto \Big(\hat{\mathcal{L}}^{n,\kappa}g^{n,T}(\cdot, \boldsymbol{\xi},t)\Big)(\boldsymbol{\zeta})$ is in~$\mathscr{C}([0,T] \times \mathcal{X}^n, \mathbb{R})$, and is given by
    \begin{equation*}
    \begin{aligned}
    & \Big(\hat{\mathcal{L}}^{n,\kappa}g^{n,T}(\cdot, \boldsymbol{\xi},t)\Big)(\boldsymbol{\zeta}) \\ & \quad =  \sum_{\bar{\boldsymbol{\xi}} \in \mathcal{X}^n} \mathbb{P}_{\boldsymbol{\xi}}(\xi^n (T-t) = \bar{\boldsymbol{\xi}})  \left(\frac{m}{2}(\hat{\mathcal{L}}^{n}_{m}\mathfrak{D}^{N,n}(\cdot, \boldsymbol{\zeta}))(\bar{\boldsymbol{\xi}}) + (\hat{\mathcal{L}}^{n,\kappa}_{r}\mathfrak{D}^{N,n}(\bar{\boldsymbol{\xi}}, \cdot))(\boldsymbol{\zeta}) \right).
    \end{aligned}
    \end{equation*}
    \item
    $
    \!
    \begin{aligned}[t]
           \sup_{t_{1},t_{2} \in [0,T]} \; \mathbb{E}_{\boldsymbol{\eta}}\Bigg[\left(\frac{\partial}{\partial t} g^{n,T}(\zeta^{n,\kappa}(t_{1}), \boldsymbol{\xi}, \cdot)\right)^{2}(t_{2}) + \Big(\hat{\mathcal{L}}^{n,\kappa}g^{n,T}(\cdot, \boldsymbol{\xi}, t_{2})\Big)^{2} (\zeta^{n,\kappa}(t_{1})) \Bigg] &< \infty.
        \end{aligned}
    $
\end{enumerate}
Since the arguments for the proof of claims~(i)-(v) are standard (see for instance~\cite[Section~2.6]{demasi1991mathematical}), we will only sketch them and omit the details. Indeed, claim~(i) follows from applying~\eqref{Paper01_well_posedness_truncated_process} to~\eqref{Paper01_definition_auxiliary_map_chain_rule_formula_correlation_function}; claim~(ii) follows from the fact that the process $(\xi^n(t))_{t \geq 0}$ is a strongly continuous Feller process taking values in the locally compact space $\mathcal{X}^{n}$ defined in~\eqref{Paper01_modified_state_space_simplified_version_without_types}; claim~(iii) follows from~\eqref{Paper01_well_posedness_truncated_process}, the standard theory of Feller processes in finite dimensional state spaces (see for instance~\cite[Theorem~8.4]{darling2008differential}), and from the fact that the birth and death rates of $\zeta^{n,\kappa}$ are bounded above by polynomials of the number of particles at a deme. Claim~(iv) follows from the expression for $\hat{\mathcal{L}}^{n,\kappa}$ in~\eqref{Paper01_infinitesimal_generator_auxiliary_process} and from Lemma~\ref{Paper01_duality_xi_Lambda_N}. Finally, claim~(v) follows from claims~(ii) and~(iv), and from~\eqref{Paper01_well_posedness_truncated_process}. 

Claims (i)-(v) allow us to apply a general stochastic chain rule formula, see Lemma~\ref{Paper01_general_integration_by_parts_formula}, to the map $g^{n,T}$, yielding that under $\mathbb P_{\boldsymbol{\eta}}$, the stochastic process $(M^{n,\kappa,T}_{\boldsymbol{\xi}}(t))_{t \in [0,T]}$ given by, for $t \in [0,T]$,
\begin{equation} \label{Paper01_local_martingale_definition_to_rewrite_correlation_function}
\begin{aligned}
M^{n,\kappa,T}_{\boldsymbol{\xi}}(t) & \defeq g^{n,T}(\zeta^{n,\kappa}(t), \boldsymbol{\xi}, t) - g^{n,T}(\zeta^{n,\kappa}(0), \boldsymbol{\xi}, 0) \\ & \quad \quad - \int_{0}^{t} \left(\left(\frac{\partial}{\partial t'} g(\zeta^{n,\kappa}(t'-),\boldsymbol{\xi},\cdot)\right)(t') + \Big(\hat{\mathcal{L}}^{n,\kappa}g^{n,T}(\cdot, \boldsymbol{\xi},t')\Big) (\zeta^{n,\kappa}(t'-))\right) \, dt',
\end{aligned}
\end{equation}
is a càdlàg martingale with respect to the filtration~$\Big\{\mathcal{F}^{\zeta^{n,\kappa}}_{t+}\Big\}_{t \geq 0}$. By taking $t = T$ in~\eqref{Paper01_local_martingale_definition_to_rewrite_correlation_function}, and then by taking expectations on both sides of~\eqref{Paper01_local_martingale_definition_to_rewrite_correlation_function}, using claims~(ii) and (iv), and identity~\eqref{Paper01_equivalence_test_function_correlation_funciton}, and rearranging terms, we get
\begin{align*}
 &\rho^{N,n,\kappa}(\boldsymbol{\eta}, \boldsymbol{\xi},T) \\
 &\quad = \mathbb{E}_{\boldsymbol{\eta}}\left[ g^{n,T}(\zeta^{n,\kappa}(0), \boldsymbol{\xi}, 0)\right] + \int_{0}^{T} \mathbb{E}_{\boldsymbol{\eta}}\Bigg[\sum_{\bar{\boldsymbol{\xi}} \in \mathcal{X}^n} \mathbb{P}_{\boldsymbol{\xi}}(\xi^n (T-t) = \bar{\boldsymbol{\xi}}) (\hat{\mathcal{L}}^{n,\kappa}_{r}\mathfrak{D}^{N,n}(\bar{\boldsymbol{\xi}}, \cdot))(\zeta^{n,\kappa}(t-))\Bigg] \, dt.
\end{align*}
The proof is then completed by applying Fubini's theorem together with the definition of $g^{n,T}$ in~\eqref{Paper01_definition_auxiliary_map_chain_rule_formula_correlation_function}.
\end{proof}

We will need another lemma to allow us to bound the correlation functions; this lemma is analogous to Lemma~4.2.3 in~\cite{demasi1991mathematical}. Before introducing this result, we define, for any $n\in \mathbb N$ and $x \in \Lambda_n$, the configuration $\boldsymbol{E}^{(x)} \in \mathcal{X}^{n}$ to denote the configuration consisting of exactly one particle at deme $x$ and no particles elsewhere.
Moreover, for $n \in \mathbb{N}$ and $\boldsymbol{\xi}=(\xi(x))_{x\in \Lambda_n} \in \mathcal{X}^{n}$, denote by $\boldsymbol{\xi}^{(x)} = \xi(x)\boldsymbol{E}^{(x)}$ the element of $\mathcal{X}^{n}$ which has exactly $\xi(x)$ particles in deme $x \in \Lambda_n$ and none elsewhere. 
\begin{lemma} \label{Paper01_first_bound_correlation_functions}
For every $N,n,\kappa \in \mathbb{N}$ and all $\boldsymbol{\xi}, \boldsymbol{\zeta} \in \mathcal{X}^{n}$,
\begin{equation*}
    \hat{\mathcal{L}}^{ n,\kappa}_{r}\mathfrak{D}^{N,n}(\boldsymbol{\xi}, \cdot)(\boldsymbol{\zeta}) = \sum_{x \in \Lambda_n} \mathfrak{D}^{N,n}(\boldsymbol{\xi} - \boldsymbol{\xi}^{(x)}, \boldsymbol{\zeta})  \left(\hat{\mathcal{L}}^{n,\kappa}_{r}\mathfrak{D}^{N,n}(\boldsymbol{\xi}^{(x)}, \cdot) (\boldsymbol{\zeta})\right).
\end{equation*}
Furthermore, there is a non-negative non-decreasing function $\vartheta: \mathbb{N}_{0} \rightarrow [0,\infty)$ with $\vartheta(0) = 0$, which does not depend on $N$, $n$ or $\kappa$, such that for every $N, n,\kappa \in \mathbb{N}$, every $\boldsymbol{\zeta}, \boldsymbol{\xi} \in \mathcal{X}^{n}$ with $\boldsymbol{\xi}=(\xi(y))_{y\in \Lambda_n}$ and every $x \in \Lambda_n$,
\begin{equation} \label{Paper01_estimates_coefficients_reaction_part_generator}
    \hat{\mathcal{L}}^{n,\kappa}_{r}\mathfrak{D}^{N,n}(\boldsymbol{\xi}^{(x)}, \cdot) (\boldsymbol{\zeta}) \leq \left\{\arraycolsep=1.4pt\def\arraystretch{1.5}\begin{array}{ll} 0 & \textrm{if } \xi(x) = 0, \\ \vartheta\left(\xi(x)\right) \cdot \mathfrak{D}^{N,n}(\boldsymbol{\xi}^{(x)}, \boldsymbol{\zeta}) & \textrm{if } \xi(x) = 1, \\ \vartheta\left(\xi(x)\right) \cdot \mathfrak{D}^{N,n}(\boldsymbol{\xi}^{(x)} - \boldsymbol{E}^{({x})}, \boldsymbol{\zeta}) & \textrm{if } \xi(x) \geq 1. \end{array}\right.
\end{equation}
\end{lemma}

\begin{remark} \label{Paper01_when_xi_has_one_particle_on_x_more_inequalities_satisfied}
    Note that by the definition of $\mathfrak{D}^{N,n}$ in~\eqref{Paper01_definition_duality_function}, the third inequality of~\eqref{Paper01_estimates_coefficients_reaction_part_generator} implies that if the deme $x \in \Lambda_n$ and the configuration $\boldsymbol{\xi} =(\xi(y))_{y\in \Lambda_n}\in \mathcal{X}^{n}$ are such that $\xi{(x)} = 1$, then for any $\boldsymbol{\zeta} \in \mathcal{X}^{n}$,
    \begin{equation*}
         \hat{\mathcal{L}}^{n,\kappa}_{r}\mathfrak{D}^{N,n}(\boldsymbol{\xi}^{(x)}, \cdot) (\boldsymbol{\zeta}) \leq \vartheta(1).
    \end{equation*}
\end{remark}

We highlight that the main difference between the result proved in~\cite{demasi1991mathematical} and our Lemma~\ref{Paper01_first_bound_correlation_functions} is that our bounds are uniform in $N\in \mathbb N$, the scaling parameter for
the carrying capacity that appears in the per capita birth and death rates in our particle system. This, however, does not lead to any technical difficulties  in the proof of Lemma~\ref{Paper01_first_bound_correlation_functions}, since our definition of scaled Poisson polynomials given by Definition~\ref{Paper01_poisson_polynomial} already includes the carrying capacity parameter $N$. Since the proof of Lemma~4.2.3 in~\cite{demasi1991mathematical} is provided in~\cite{demasi1991mathematical}, we omit the proof of Lemma~\ref{Paper01_first_bound_correlation_functions} here.

Our next lemma provides a bound on the correlation functions that holds uniformly in $N$.
Recall the definition of $\rho^{N,n,\kappa}$ in Definition~\ref{Paper01_correlation_functions}.

\begin{lemma} \label{Paper01_uniform_bound_correlation_functions}
For any $n\in \mathbb N$ and $\boldsymbol{\xi}=(\xi(x))_{x\in \Lambda_n} \in \mathcal{X}^{n}$, let $$\vert \vert \vert \boldsymbol{\xi} \vert \vert \vert_{\mathcal{X}^{n}} \defeq \sum_{x \in \Lambda_n} \xi(x)$$ denote the total number of particles in $\boldsymbol{\xi}$. Moreover, for each $N,n,\kappa,p \in \mathbb{N}$, for any $\boldsymbol{\eta} \in \mathcal{S}$ and $T \geq 0$, let
\begin{equation*}
    \Psi^{N,n,\kappa}(\boldsymbol{\eta},p,T) \defeq \sup_{t \leq T} \; \sup_{\substack{\{\boldsymbol{\xi} \in \mathcal{X}^{n}: \; \vert \vert \vert \boldsymbol{\xi} \vert \vert \vert_{\mathcal{X}^{n}} \leq p\}}} \; \rho^{N,n,\kappa}(\boldsymbol{\eta}, \boldsymbol{\xi}, t).
\end{equation*}
Then there exists a positive function $\Psi: \mathbb{N} \times [0,\infty) \rightarrow (0,\infty)$, which does not depend on $N,n,\kappa,L$ or $m$, such that for each $N, n, p \in \mathbb{N}$ and any $\boldsymbol{\eta}=(\eta(x))_{x\in L^{-1}\mathbb Z}\in \mathcal{S}$ and $T \geq 0$,
\begin{equation*}
    \sup_{\kappa \in  \mathbb{N}} \, \Psi^{N,n,\kappa}(\boldsymbol{\eta}, p,T) \leq \Psi(p,T)\Bigg(\sup_{x \in \Lambda_n} \, \frac{\vert \vert \eta(x) \vert \vert_{\ell_{1}}^{p}}{N^{p}} + 1 \Bigg).
\end{equation*}
\end{lemma}

\begin{proof}
Let $N,n,\kappa \in \mathbb{N}$, $\boldsymbol{\eta}=(\eta(x))_{x\in L^{-1}\mathbb Z} \in \mathcal{S}$, $p \in \mathbb{N}$ and $T \geq 0$ be fixed. Notice that by~\eqref{Paper01_well_posedness_truncated_process}, $\Psi^{N,n,\kappa}(\boldsymbol{\eta},p,T)$ is well defined and finite. 
Recall from above~\eqref{Paper01_well_posedness_truncated_process} that under $\mathbb P_{\boldsymbol{\eta}}$,
the initial condition of $(\zeta^{n,\kappa}(t))_{t \geq 0}$ is given by~\eqref{Paper01_deriving_initial_condition_auxiliary_process}; hence by Definition~\ref{Paper01_poisson_polynomial} and~\eqref{Paper01_definition_duality_function}, for any $\boldsymbol{\xi} \in \mathcal{X}^{n}$ such that $\vert \vert \vert \boldsymbol{\xi} \vert \vert \vert_{\mathcal{X}^{n}} \leq p$,
\begin{equation} \label{Paper01_initial_uniform_bound_correlation_function}
    \rho^{N,n,\kappa}(\boldsymbol{\eta},\boldsymbol{\xi},0) \leq \sup_{x \in \Lambda_n} \, \frac{\vert \vert \eta(x) \vert \vert_{\ell_{1}}^{p}}{N^{p}} + 1.
\end{equation}
Note that for any $t \in [0,T]$ and $\bar{\boldsymbol{\xi}}=(\bar{\xi}(x))_{x\in \Lambda_n} \in \mathcal{X}^{n}$ such that $\vert \vert \vert \bar{\boldsymbol{\xi}} \vert \vert \vert_{\mathcal{X}^{n}} \leq p$, we have by Lemma~\ref{Paper01_first_bound_correlation_functions}, using that $\vartheta(0)=0$ and so terms with $\bar{\xi}(x)=0$ do not contribute to the sum, that
\begin{equation} \label{Paper01_auxiliary_estimate_correlation_functions}
\begin{aligned}
    & \mathbb{E}_{\boldsymbol{\eta}}\left[\hat{\mathcal{L}}^{n,\kappa}_{r}\mathfrak{D}^{N,n}(\bar{\boldsymbol{\xi}}, \cdot )\left(\zeta^{n, \kappa}(t)\right)\right] \\ 
    & \quad \leq \mathbb{E}_{\boldsymbol{\eta}}\Bigg[\sum_{\{x \in \Lambda_n: \, \bar{\xi}(x) \neq 0\}} \mathfrak{D}^{N,n}\left(\bar{\boldsymbol{\xi}} - \bar{\boldsymbol{\xi}}^{(x)}, \zeta^{n, \kappa}(t)\right) \cdot \vartheta(\bar{\xi}(x)) \cdot \mathfrak{D}^{N,n}\left(\bar{\boldsymbol{\xi}}^{(x)} - \boldsymbol{E}^{(x)}, \zeta^{n, \kappa}(t)\right)\Bigg] \\ 
    & \quad \leq \vartheta(p) \sum_{\{x \in \Lambda_n: \, \bar{\xi}(x) \neq 0\}} \mathbb{E}_{\boldsymbol{\eta}}\left[\mathfrak{D}^{N,n}\left(\bar{\boldsymbol{\xi}} - \bar{\boldsymbol{\xi}}^{(x)}, \zeta^{n,\kappa}(t)\right) \cdot \mathfrak{D}^{N,n}\left(\bar{\boldsymbol{\xi}}^{(x)} - \boldsymbol{E}^{(x)}, \zeta^{n,\kappa}(t)\right)\right] \\ 
    & \quad \leq \vartheta(p) \sum_{\{x \in \Lambda_n : \, \bar{\xi}(x) \neq 0\}} \mathbb{E}_{\boldsymbol{\eta}}\left[\mathfrak{D}^{N,n}\left(\bar{\boldsymbol{\xi}} - \boldsymbol{E}^{(x)}, \zeta^{n,\kappa}(t)\right)\right] \\ 
    & \quad \leq p\vartheta(p) \Psi^{N,n,\kappa}({\boldsymbol{\eta}},p,t),
\end{aligned}
\end{equation}
where in the second inequality above we used the fact that the function $\vartheta: \mathbb{N}_{0} \rightarrow [0,\infty)$ defined in Lemma~\ref{Paper01_first_bound_correlation_functions} is non-decreasing, in the third inequality we used~\eqref{Paper01_definition_duality_function} and the definitions of $\bar{\xi}^{(x)}$ and $E^{(x)}$ before Lemma~\ref{Paper01_first_bound_correlation_functions}, and in the last inequality we used the fact that $\vert \vert \vert \bar{\boldsymbol{
\xi}} \vert \vert \vert_{\mathcal{X}^{n}} \leq p$. Thus by applying Lemma~\ref{Paper01_rewriting_definition_correlation_function}, estimate~\eqref{Paper01_auxiliary_estimate_correlation_functions}, and identity~\eqref{Paper01_initial_uniform_bound_correlation_function}, we conclude that for all $\boldsymbol{\eta}=(\eta(x))_{x\in L^{-1}\mathbb Z} \in \mathcal{S}$, $N, n, \kappa, p \in \mathbb{N}$ and $T \geq 0$,
\begin{equation} \label{Paper01_almost_final_bound_pre_gronwall_correlation}
    \Psi^{N,n,\kappa}(\boldsymbol{\eta}, p,T) \leq \Bigg(\sup_{x \in \Lambda_n} \, \frac{\vert \vert \eta(x) \vert \vert_{\ell_{1}}^{p}}{N^{p}} + 1\Bigg) + \int_{0}^{T} p \vartheta(p)  \Psi^{N,n,\kappa}(\boldsymbol{\eta}, p,t) \; dt.
\end{equation}
Since~\eqref{Paper01_well_posedness_truncated_process} implies that $\Psi^{N,n,\kappa}(\boldsymbol{\eta}, p,T) < \infty$, we can apply Gr\"{o}nwall's lemma to~\eqref{Paper01_almost_final_bound_pre_gronwall_correlation}, which completes the proof.
\end{proof}

We are now ready to prove Proposition~\ref{Paper01_bound_total_mass}. Since the proof follows the proofs of Propositions~4.2.1 and~4.2.5 in~\cite{demasi1991mathematical} closely, we will omit some of the details.

\begin{proof}[Proof of Proposition~\ref{Paper01_bound_total_mass}]
Fix $q_+$ and $q_-$ satisfying Assumption~\ref{Paper01_assumption_polynomials}.
    By the fact that for each $n \in \mathbb{N}$, the process without mutations $(\zeta^{n}(t))_{t \geq 0}$ stochastically dominates the process $(\eta^{n}(t))_{t \geq 0}$ in the sense of the coupling in~\eqref{eq:zetaaboveeta}, the bound~\eqref{Paper01_discrete_approximation_does_not_explode} will be proved after establishing that for any $L>0$, $m>0$ and $N\in \mathbb N$,
    for any $n\in \mathbb N$, $\boldsymbol{\eta} \in \mathcal{S}$ and $T \geq 0$,
    \begin{equation} \label{Paper01_supremum_time_total_number_particles_zeta_n}
        \mathbb{E}_{\boldsymbol{\eta}}\Bigg[\sup_{t \le T} \; \sum_{x \in \Lambda_n} \zeta^{n}(t,x) \Bigg] < \infty.
    \end{equation}
    Observe that since $(\zeta^{n}(t))_{t \geq 0}$ stochastically dominates $(\eta^{n}(t))_{t \geq 0}$, estimate~\eqref{Paper01_supremum_time_total_number_particles_zeta_n} also implies that~$(\eta^{n}(t))_{t \geq 0}$ is well defined as a (strongly) Feller process (see~\cite[Theorem~6.4.1]{ethier2009markov}). Moreover, the bound~\eqref{Paper01_crude_estimate_moments_general_configuration} will be proved after establishing that for every $p \in \mathbb{N}$, there exists a non-decreasing function $C_{p}: [0,\infty) \rightarrow [1,\infty)$ such that for any $L>0$, $m>0$ and $N\in \mathbb N$,
    for any $n\in \mathbb N$, $\boldsymbol{\eta}=(\eta(x))_{x\in L^{-1}\mathbb Z} \in \mathcal{S}$ and $T \geq 0$,
    \begin{equation} \label{Paper01_moments_bounds_auxiliary_process}
        \sup_{x \in \Lambda_n} \mathbb{E}_{\boldsymbol{\eta}}\left[\left(\frac{\zeta^{n}(T,x)}{N}\right)^{p}\right] \leq C_{p}(T) \Bigg(\sup_{x \in \Lambda_n} \, \frac{\vert \vert \eta(x) \vert \vert_{\ell_{1}}^{p}}{N^{p}} + 1\Bigg).
    \end{equation}
    Estimates~\eqref{Paper01_supremum_time_total_number_particles_zeta_n} and~\eqref{Paper01_moments_bounds_auxiliary_process} will be proved by following the steps below:
    \begin{enumerate}[(1)]
        \item First, we prove that for any $\boldsymbol{\eta}\in \mathcal S$, taking $\boldsymbol{\zeta}$ as in~\eqref{Paper01_deriving_initial_condition_auxiliary_process} and conditioning on $\zeta^{n,\kappa}(0)=\boldsymbol{\zeta}$ for each $\kappa \in \mathbb N$,
        $(\zeta^{n,\kappa}(t))_{t \geq 0}$ converges in distribution on $\mathscr{D}([0, \infty), \mathcal X^n)$ as $\kappa \rightarrow \infty$ to $(\zeta^{n}(t))_{t \geq 0}$ with $\zeta^n(0)=\boldsymbol{\zeta}$;
        \item By the convergence in distribution and Lemma~\ref{Paper01_uniform_bound_correlation_functions}, it is immediate that we can define the correlation function $\rho^{N,n}$ for the process $\zeta^{n}$ analogously to Definition~\ref{Paper01_correlation_functions}, and that, by Fatou's lemma, $\rho^{N,n}$ also satisfies the bounds provided by Lemma~\ref{Paper01_uniform_bound_correlation_functions}.
        \end{enumerate}
        
        By the definition of $\zeta^{n, \kappa}$ below~\eqref{eq:zetankappadefn}, the convergence in  distribution stated in step $(1)$ will follow (see~\cite[Proposition~4.2.1]{demasi1991mathematical}) after establishing that for every $N,n \in \mathbb{N}$, $L,m>0$ and any $\boldsymbol{\eta}\in \mathcal S$ and $T > 0$,
         \begin{equation} \label{Paper01_supremum_time_total_number_particles_zeta_n_kappa}
            \sup_{\kappa\in \mathbb N} \mathbb{E}_{\boldsymbol{\eta}}\Bigg[\sup_{t \leq T} \; \sum_{x \in \Lambda_n} {\zeta^{n,\kappa}(t,x)}\Bigg] < \infty.
        \end{equation}
        Estimate~\eqref{Paper01_supremum_time_total_number_particles_zeta_n_kappa} follows from the representation of the total number of particles of $(\zeta^{n,\kappa}(t))_{t \geq 0}$ in terms of a martingale problem, the fact that the expectation of the birth and death rates of $(\zeta^{n}(t))_{t \geq 0}$ can be written in terms of sums of correlation functions, and then from an application of Lemma~\ref{Paper01_uniform_bound_correlation_functions}. Since the proof of an analogous result is shown carefully in~\cite[Proposition~4.2.1]{demasi1991mathematical}, we omit the details. By Skorokhod's representation theorem, we can construct the sequence of processes $(\zeta^{n,\kappa})_{\kappa \in \mathbb{N}}$ and $\zeta^{n}$ on the same probability space so that convergence occurs almost surely. By using this construction, and then by combining Fatou's lemma and~\eqref{Paper01_supremum_time_total_number_particles_zeta_n_kappa}, estimate~\eqref{Paper01_supremum_time_total_number_particles_zeta_n} is proved.
        
        Moving now to step $(2)$, we observe that by Lemma~\ref{Paper01_uniform_bound_correlation_functions} and by the definition of the correlation functions in Definition~\ref{Paper01_correlation_functions}, for every~$p \in \mathbb{N}$, there exists a non-decreasing function $C_{p}: [0, \infty) \rightarrow [1, \infty)$ such that for any $n \in \mathbb{N}$, $\boldsymbol{\eta}=(\eta(x))_{x\in L^{-1}\mathbb Z} \in \mathcal{S}$ and $T \geq 0$,
        \begin{equation} \label{Paper01_moments_bound_auxilliary_process_cut_off}
            \sup_{\kappa \in \mathbb{N}} \; \sup_{x \in \Lambda_n} \mathbb{E}_{\boldsymbol{\eta}}\left[\left(\frac{\zeta^{n,\kappa}(T,x)}{N}\right)^{p}\right] \leq C_{p}(T) \Bigg(\sup_{x \in \Lambda_n} \, \frac{\vert \vert \eta(x) \vert \vert_{\ell_{1}}^{p}}{N^{p}} + 1\Bigg).
        \end{equation}
        As in the proof of~\eqref{Paper01_supremum_time_total_number_particles_zeta_n}, since $(\zeta^{n,\kappa})_{\kappa \in \mathbb{N}}$ converges in distribution to~$\zeta^{n}$ as $\kappa \rightarrow \infty$, by Skorokhod's representation theorem, and then by combining Fatou's lemma and~\eqref{Paper01_moments_bound_auxilliary_process_cut_off}, estimate~\eqref{Paper01_moments_bounds_auxiliary_process} is proved.

        As explained at the beginning of this proof, estimates~\eqref{Paper01_discrete_approximation_does_not_explode} and~\eqref{Paper01_crude_estimate_moments_general_configuration} follow from~\eqref{Paper01_supremum_time_total_number_particles_zeta_n} and~\eqref{Paper01_moments_bounds_auxiliary_process}. Observing that by the definition of $\mathcal{S}_{0}$ in~\eqref{Paper01_set_initial_configurations}, for every $\boldsymbol{\eta}=(\eta(x))_{x\in L^{-1}\mathbb Z} \in \mathcal{S}_{0}$,
        \begin{equation*}
            \sup_{n \in \mathbb{N}} \; \sup_{x \in \Lambda_n} \frac{\vert \vert \eta(x) \vert \vert_{\ell_{1}}^{p}}{N^{p}} = \sup_{\substack{x \in L^{-1}\mathbb{Z}}} \frac{\vert \vert \eta(x) \vert \vert_{\ell_{1}}^{p}}{N^{p}} < \infty,
        \end{equation*}
        we conclude that~\eqref{Paper01_moments_estimates_uniform_space} follows from~\eqref{Paper01_crude_estimate_moments_general_configuration}, which completes the proof.
\end{proof}

\section{Construction of the discrete spatial Muller's ratchet} \label{Paper01_formal_construction_generator_section}

In this section, we will prove Theorems~\ref{Paper01_thm_existence_uniqueness_IPS_spatial_muller} and~\ref{Paper01_thm_moments_estimates_spatial_muller_ratchet}, i.e.~in particular, we will prove that the sequence of $\mathcal{S}$-valued Feller processes $((\eta^{n}(t))_{t \geq 0})_{n \in \mathbb{N}}$ introduced in Section~\ref{Paper01_model_description}  converges weakly on $\mathscr{D}([0,\infty), \mathcal{S})$ as $n \rightarrow \infty$ to the spatial Muller's ratchet $(\eta(t))_{t \geq 0}$. Recall that the space of configurations $\mathcal{S}$ is defined in~\eqref{Paper01_definition_state_space_formal}, and the set of initial configurations $\mathcal{S}_{0} \subset \mathcal{S}$ is defined in~\eqref{Paper01_set_initial_configurations}. Also, 
take $L>0$, $m>0$, $\mu\in [0,1]$ and a
sequence of fitness parameters $(s_{k})_{k \in \mathbb{N}_{0}}$ satisfying Assumption~\ref{Paper01_assumption_fitness_sequence}, and take birth and death rate polynomials $q_{+}$ and $q_{-}$ satisfying Assumption~\ref{Paper01_assumption_polynomials}. Throughout this section, we fix 
$N \in \mathbb{N}$ and, to simplify notation, we omit explicit dependence on $N$ whenever it is clear from context.

Recall from Proposition~\ref{Paper01_bound_total_mass} that $(\eta^n(t))_{t\ge 0}$ is a well-defined $\mathcal S$-valued c\`adl\`ag Feller process with generator $\mathcal L^n$ (defined in \eqref{Paper01_generator_foutel_etheridge_model_restriction_n} and~\eqref{Paper01_infinitesimal_generator_restriction_n}).
For $n \in \mathbb{N}$ and  $\boldsymbol{\xi}=(\xi_k(x))_{x\in L^{-1}\mathbb Z, \, k\in \mathbb N_0} \in \mathcal{S}$, we can think of $(\eta^{n}(t))_{t \geq 0}=(\eta^{n}_k(t,x))_{t \geq 0,\, x\in L^{-1}\mathbb Z,\, k\in \mathbb N_0}$ conditioned on $\eta^{n}(0) = \boldsymbol{\xi}$ as an $\mathbb{N}_{0}^{d}$-valued process, for some $d \in \mathbb{N}$, by ignoring particles living outside $\Lambda_n$ which do not move, reproduce or die, and tracking $\eta^{n}_{k}(t,x)$ for $x \in \Lambda_n$ and 
$$k \leq K_n \vee \max\Big\{k' \in \mathbb{N}_0 : \sum_{y\in \Lambda_n}\xi_{k'}(y)>0\Big\}.$$
As defined after~\eqref{Paper01_infinitesimal_generator_restriction_n},
for $\phi \in \mathscr{C}_{b}(\mathcal{S},\mathbb{R})$, $T \geq 0$ and $\boldsymbol{\xi} \in \mathcal{S}$, we let
\begin{equation*}
    (P^{n}_{T}\phi)(\boldsymbol{\xi}) \defeq \mathbb{E}_{\boldsymbol{\xi}}[\phi(\eta^{n}(T))].
\end{equation*}
By Proposition~\ref{Paper01_bound_total_mass}, $\{P^{n}_{t}\}_{t \geq 0}$ is a Feller semigroup on $\mathcal{S}$ (see Definition~\ref{Paper01_definition_feller_non_locally_compact} of Feller semigroups on non-locally compact Polish spaces) with infinitesimal generator $\mathcal{L}^{n}$ defined in~\eqref{Paper01_infinitesimal_generator_restriction_n} and~\eqref{Paper01_generator_foutel_etheridge_model_restriction_n}.

It will sometimes be convenient to use a construction of $(\eta^{n}(t))_{t \geq 0}$ in terms of families of independent Poisson processes. Since for any $\boldsymbol{\eta} \in \mathcal{S}$ and $n \in \mathbb{N}$, the number of types of particles in the process $(\eta^{n}(t))_{t \geq 0}$ conditioned on $\eta^n(0)=\boldsymbol{\eta}$ is deterministically bounded, and since by~\eqref{Paper01_discrete_approximation_does_not_explode} in Proposition~\ref{Paper01_bound_total_mass}, the number of particles in $(\eta^{n}(t))_{t \geq 0}$ living in the box $\Lambda_n$ remains almost surely finite during finite periods of time, this construction is standard (see for instance~\cite[Theorem~2.8]{garcia2006spatial}). Take three families of i.i.d.~Poisson random measures on $[0, \infty) \times [0, \infty)$ with intensity measure given by the Lebesgue measure $\lambda$. These families are denoted $\{\mathcal{Q}^{x,x+z}_{k}: \, x \in L^{-1}\mathbb{Z}, \, z \in \{-L^{-1},L^{-1}\}, \, k \in \mathbb{N}_{0}\}$, $\left\{\mathcal{R}^{x}_{k}: \, x \in L^{-1}\mathbb{Z}, \, k \in \mathbb{N}_{0}\right\}$ and $\left\{\mathcal{D}^{x}_{k}: \, x \in L^{-1}\mathbb{Z}, \, k \in \mathbb{N}_{0}\right\}$, and will be used in the construction of migration, reproduction, and death events, respectively. We recall the definition of the discrete birth and death rates $F^{b}_{k}, F^{d}_{k}: \ell_{1}^{+} \rightarrow [0, \infty)$ given by~\eqref{Paper01_discrete_birth_rate} and~\eqref{Paper01_discrete_death_rate} respectively. We also recall from before~\eqref{Paper01_scaled_polynomials_carrying_capacity} that for~$x \in L^{-1}\mathbb{Z}$ and~$k \in \mathbb{N}_{0}$, $\boldsymbol{e}_{k}^{(x)} \in \mathcal{S}$ denotes the configuration consisting of one particle carrying exactly $k$ mutations in deme $x$.
Then, using e.g.~\cite[Theorem~2.8]{garcia2006spatial}, for $n\in \mathbb N$ and $\boldsymbol{\eta}\in \mathcal S$, we can construct $(\eta^n(t))_{t\ge 0}=(\eta^n_k(t,x))_{t\ge 0, \, x\in L^{-1}\mathbb Z, \, k\in \mathbb N_0}$ conditioned on $\eta^n(0)=\boldsymbol{\eta}$
 in such a way that for any $\phi \in \mathscr{C}\left(\mathcal{S}, \mathbb{R}\right)$ and $T \geq 0$, we have almost surely
\begin{equation} \label{Paper01_Poisson_process_representation_discrete_approximation}
\begin{aligned}
    & \phi\left(\eta^{n}(T)\right) - \phi\left(\boldsymbol{\eta}\right) \\ 
    & \quad = \sum_{x \in \Lambda_n} \, \sum_{k = 0}^{\infty} \, \sum_{z \in \{-L^{-1}, L^{-1}\}} \mathds{1}_{\{x + z \in \Lambda_n\}}\int_{[0,T]\times [0,\infty)} \Big(\phi\Big(\eta^{n}(t-) + \boldsymbol{e}^{(x + z)}_{k} - \boldsymbol{e}^{(x)}_{k}\Big) - \phi(\eta^{n}(t-))\Big) \\[-6mm] & \quad \quad \quad \quad \quad \quad \quad \quad \quad \quad \quad \quad \quad \quad \quad \quad \quad \quad \quad \quad \quad \quad \quad \quad \quad \quad \quad \cdot \mathds{1}_{\left\{y \leq \frac{m}{2}\eta^{n}_{k}(t-,x)\right\}} \mathcal{Q}^{x,x+z}_{k}(dt \times dy) \\ 
    & \quad \quad + \sum_{x \in \Lambda_n} \, \sum_{k = 0}^{K_n} \, \int_{[0,T]\times [0,\infty)} \left(\phi\left(\eta^{n}(t-) + \boldsymbol{e}^{(x)}_{k}\right) - \phi(\eta^{n}(t-))\right) \cdot \mathds{1}_{\left\{y \leq F^{b}_{k}(\eta^{n}(t-,x))\right\}} \mathcal{R}^{x}_{k}(dt \times dy) \\ 
    & \quad \quad + \sum_{x \in \Lambda_n} \, \sum_{k = 0}^{\infty} \, \int_{[0,T]\times [0,\infty)} \left(\phi\left(\eta^{n}(t-) - \boldsymbol{e}^{(x)}_{k}\right) - \phi(\eta^{n}(t-))\right) \cdot \mathds{1}_{\left\{y \leq F^{d}_{k}(\eta^{n}(t-,x))\right\}} \mathcal{D}^{x}_{k}(dt \times dy).
\end{aligned}
\end{equation}
For $n\in \mathbb N$, $\phi \in \mathscr{C}(\mathcal{S}, \mathbb{R})$ and $\boldsymbol{\eta}\in \mathcal S$, let $(M^{n,\phi,\boldsymbol{\eta}}(t))_{t \geq 0}$ be the stochastic process obtained  by subtracting the intensity measures from the Poisson measures in~\eqref{Paper01_Poisson_process_representation_discrete_approximation}, i.e.~for any $T \geq 0$, we let
\begin{equation} \label{Paper01_local_martingale_problem_discrete_approximation}
    M^{n,\phi,\boldsymbol{\eta}}(T) \defeq \phi(\eta^{n}(T)) - \phi(\boldsymbol{\eta}) - \int_{0}^{T} (\mathcal{L}^{n}\phi)(\eta^{n}(t-)) \, dt.
\end{equation}
We will show in the proof of Lemma~\ref{Paper01_auxilliary_lemma_easy_martingale_representation} below that
the following estimate holds for suitable $\phi:\mathcal S\to \mathbb R$, and for $\boldsymbol{\xi}\in \mathcal S$, $n\in \mathbb N$ and $T \geq 0$:
\begin{align} \label{Paper01_easy_criterion_martingale_feller_process_N_d}
&E^n(\phi,\boldsymbol{\xi},T) :=\notag \\
    &\sum_{x \in \Lambda_n} \sum_{k = 0}^{\infty} \, \sum_{z \in \{-L^{-1}, L^{-1}\}} \mathds{1}_{\{x + z \in \Lambda_n\}} \frac{m}{2} \int_{0}^{T} \mathbb{E}_{\boldsymbol{\xi}}\Big[  \Big\vert \phi\Big(\eta^{n}(t-) + \boldsymbol{e}^{(x + z)}_{k} - \boldsymbol{e}^{(x)}_{k}\Big) - \phi(\eta^{n}(t-))\Big\vert^2 \eta^{n}_{k}(t-,x) \Big] dt \notag \\ 
    & \qquad + \sum_{x \in \Lambda_n} \, \sum_{k = 0}^{K_n} \, \int_{0}^{T} \mathbb{E}_{\boldsymbol{\xi}} \Big[\Big\vert\phi\left(\eta^{n}(t-) + \boldsymbol{e}^{(x)}_{k}\right) - \phi(\eta^{n}(t-))\Big\vert^2 F^{b}_{k}(\eta^{n}(t-,x)) \Big] \, dt \notag  \\ 
    & \qquad + \sum_{x \in \Lambda_n} \, \sum_{k = 0}^{\infty} \, \int_{0}^{T} \mathbb{E}_{\boldsymbol{\xi}} \Big[\Big\vert \phi\left(\eta^{n}(t-) - \boldsymbol{e}^{(x)}_{k}\right) - \phi(\eta^{n}(t-))\Big\vert^2 F^{d}_{k}(\eta^{n}(t-,x))\Big] \, dt \notag \\ & \quad < \infty.
\end{align}
By standard results in the theory of Markov processes (see for instance~\cite[Proposition~8.7]{darling2008differential}),
if~\eqref{Paper01_easy_criterion_martingale_feller_process_N_d} holds then it follows that
$(M^{n,\phi,\boldsymbol{\xi}}(t))_{t \in [0,T]}$ is a square-integrable càdlàg martingale under $\mathbb P_{\boldsymbol{\xi}}$. Our first result of this section is an application of~\eqref{Paper01_easy_criterion_martingale_feller_process_N_d}. Recall the definition of $\mathscr{C}_{*}(\mathcal{S}, \mathbb{R})$ from~\eqref{Paper01_definition_Lipschitz_function}, and the definition of $(F^{b}_{k})_{k \in \mathbb{N}_{0}}$ and $(F^{d}_{k})_{k \in \mathbb{N}_{0}}$ from~\eqref{Paper01_discrete_birth_rate} and~\eqref{Paper01_discrete_death_rate}, respectively. For $u \in [0,\infty)$, let
     \begin{equation} \label{Paper01_bar_polynomial}
         \bar{q}^{N}(u) \defeq u(q^{N}_{+}(u) + q^{N}_{-}(u)),
     \end{equation}
     where $q^{N}_{+},q^{N}_{-}: [0, \infty) \rightarrow [0, \infty)$ are the scaled polynomials defined in~\eqref{Paper01_scaled_polynomials_carrying_capacity}.
     Note that since by Assumption~\ref{Paper01_assumption_fitness_sequence}, $s_{k} \leq 1$ for every $k \in \mathbb{N}_{0}$, we have that for any $u = (u_{k})_{k \in \mathbb{N}_{0}} \in \ell_{1}^{+}$,
     \begin{equation} \label{Paper01_simple_bound_total_birth_death_rates_deterministic}
         \sum_{k = 0}^{\infty} (F^{b}_{k}(u) + F^{d}_{k}(u)) \leq \bar{q}^{N}(\vert \vert u \vert \vert_{\ell_{1}}).
     \end{equation}
For a càdlàg local martingale $(M(t))_{t \geq 0}$, we will denote its quadratic variation by $([M](t))_{t \geq 0}$ and its predictable bracket process by $(\langle M \rangle(t))_{t \geq 0}$ (see Section~\ref{Paper01_subsection_integral_poisson}, and see~\eqref{Paper01_quadratic_variation_stochastic_integral_poisson_point_process} for the definition of these processes in the context of integrals with respect to Poisson point processes).

\begin{lemma} \label{Paper01_auxilliary_lemma_easy_martingale_representation}
    For $\phi \in \mathscr{C}_*(\mathcal{S}, \mathbb{R})$, $\boldsymbol{\xi} \in \mathcal{S}$ and $n \in \mathbb{N}$, conditioning on $\eta^{n}(0) = \boldsymbol{\xi}$, the process $(M^{n,\phi,\boldsymbol{\xi}}(t))_{t \geq 0}$ defined in~\eqref{Paper01_local_martingale_problem_discrete_approximation} is a square-integrable càdlàg martingale. Moreover, its predictable bracket process is given by, for $T \geq 0$,
   \begin{align} \label{Paper01_predictable_bracket_process_martingale_lipschitz_functions_discrete_IPS}
        & \left\langle M^{n,\phi,\boldsymbol{\xi}}\right\rangle(T) \notag \\ 
        & \;  = \sum_{\substack{x \in \Lambda_n}} \sum_{k = 0}^{\infty} \, \sum_{z \in \{-L^{-1},L^{-1}\}} \,\int_{0}^{T} \mathds{1}_{\{x + z \in \Lambda_n\}} \frac{m}{2}\eta^{n}_{k}(t-,x) \Big(\phi\Big(\eta^{n}(t-) + e^{(x+ z)}_{k} - e^{(x)}_{k}\Big) - \phi(\eta^{n}(t-))\Big)^{2} \, dt \notag \\ 
        &  \quad \; + \sum_{\substack{x \in \Lambda_n}} \, \sum_{k = 0}^{K_n}  \, \int_{0}^{T} F^{b}_{k}(\eta^{n}(t-,x))\left(\phi\left(\eta^{n}(t-) + e^{(x)}_{k}\right) - \phi(\eta^{n}(t-))\right)^{2} \, dt \notag \\ 
        & \quad \; + \sum_{\substack{x \in \Lambda_n}} \, \sum_{k = 0}^{\infty} \, \int_{0}^{T} F^{d}_{k}(\eta^{n}(t-,x))\left(\phi\left(\eta^{n}(t-) - e^{(x)}_{k}\right) - \phi(\eta^{n}(t-))\right)^{2} \, dt.
    \end{align}
    Furthermore, the quadratic variation of $(M^{n,\phi,\boldsymbol{\xi}}(t))_{t \geq 0}$ is given by, for $T \geq 0$,
    \begin{align} \label{Paper01_definition_quadratic_variation_martingale_problem}
         & \left[M^{n,\phi,\boldsymbol{\xi}}\right](T) \notag \\ 
         & \quad = \sum_{x \in \Lambda_n} \, \sum_{k = 0}^{\infty} \, \sum_{z \in \{-L^{-1}, L^{-1}\}} \, \int_{[0,T]\times [0,\infty)} \mathds{1}_{\{x + z \in \Lambda_n\}}\Big(\phi\Big(\eta^{n}(t-) + e^{(x + z)}_{k} - e^{(x)}_{k}\Big) - \phi(\eta^{n}(t-))\Big)^{2} \notag \\[-6mm] 
         & \quad \quad \quad \quad \quad \quad \quad \quad \quad \quad \quad \quad \quad \quad \quad \quad \quad \quad \quad \quad \cdot \mathds{1}_{\left\{y \leq \frac{m}{2}\eta^{n}_{k}(t-,x)\right\}} \mathcal{Q}^{x,x+z}_{k}(dt \times dy) \notag \\ 
         & \quad \quad + \sum_{x \in \Lambda_n} \, \sum_{k = 0}^{K_n} \, \int_{[0,T]\times [0,\infty)} \left(\phi\left(\eta^{n}(t-) + e^{(x)}_{k}\right) - \phi(\eta^{n}(t-))\right)^{2} \cdot \mathds{1}_{\left\{y \leq F^{b}_{k}(\eta^{n}(t-,x))\right\}} \mathcal{R}^{x}_{k}(dt \times dy) \notag \\ 
         & \quad \quad + \sum_{x \in \Lambda_n} \, \sum_{k = 0}^{\infty} \, \int_{[0,T]\times [0,\infty)} \left(\phi\left(\eta^{n}(t-) - e^{(x)}_{k}\right) - \phi(\eta^{n}(t-))\right)^{2} \cdot \mathds{1}_{\left\{y \leq F^{d}_{k}(\eta^{n}(t-,x))\right\}} \mathcal{D}^{x}_{k}(dt \times dy).
    \end{align}
\end{lemma}
\begin{proof}
    We start by establishing that~\eqref{Paper01_easy_criterion_martingale_feller_process_N_d} holds for any $\phi \in \mathscr{C}_*(\mathcal{S}, \mathbb{R})$, $\boldsymbol{\xi} \in \mathcal{S}$, $n \in \mathbb{N}
    $ and $T \geq 0$. Take $\boldsymbol{\xi} \in \mathcal{S}$, $n \in \mathbb{N}
    $ and $T \geq 0$, and recall the definition of $\boldsymbol{e}^{(x)}_{k}$ from before~\eqref{Paper01_scaled_polynomials_carrying_capacity}. 
     By estimate~\eqref{Paper01_crude_estimate_moments_general_configuration} in Proposition~\ref{Paper01_bound_total_mass} and the fact that $\bar{q}^{N}$, defined in~\eqref{Paper01_bar_polynomial}, is a non-negative polynomial, we have
     \begin{equation} \label{Paper01_trivial_estimate_polynomial}
         \sup_{t \leq T} \, \sup_{x \in \Lambda_n} \mathbb{E}_{\boldsymbol{\xi}}\Big[\bar{q}^{N}(\vert \vert \eta^{n}(t, x) \vert \vert_{\ell_{1}})\Big] < \infty.
     \end{equation}
     By the definition of $\mathscr{C}_*(\mathcal{S}, \mathbb{R})$ in~\eqref{Paper01_definition_Lipschitz_function}, then by applying~\eqref{Paper01_simple_bound_total_birth_death_rates_deterministic}, and then finally~\eqref{Paper01_trivial_estimate_polynomial} and estimate~\eqref{Paper01_crude_estimate_moments_general_configuration} in Proposition~\ref{Paper01_bound_total_mass}, to the expression in~\eqref{Paper01_easy_criterion_martingale_feller_process_N_d}, we conclude that for any $\phi \in \mathscr{C}_*(\mathcal{S}, \mathbb{R})$, there exists $C_{\phi} > 0$ such that
     \begin{align*}
         & E^n(\phi,\boldsymbol{\xi},T) \\ & \quad \leq \sum_{x \in \Lambda_n} \frac{C_\phi(m \vee 1)}{(1 + \vert x \vert)^{4(1 + \deg q_{-})}}\int_{0}^{T} \mathbb{E}_{\boldsymbol{\xi}}\Bigg[\vert \vert \eta^{n}(t-,x) \vert \vert_{\ell_{1}} + \sum_{k = 0}^{K_n} F^{b}_{k}(\eta^{n}(t-, x)) + \sum_{k = 0}^{\infty} F^{d}_{k}(\eta^{n}(t-, x)) \Bigg] \, dt \\ & \quad \leq \sum_{x \in \Lambda_n} \frac{C_\phi(m \vee 1)}{(1 + \vert x \vert)^{4(1 + \deg q_{-})}}\int_{0}^{T} \mathbb{E}_{\boldsymbol{\xi}}\Big[\vert \vert \eta^{n}(t-,x) \vert \vert_{\ell_{1}} + \bar{q}^{N}(\vert \vert \eta^{n}(t-,x) \vert \vert_{\ell_{1}})\Big] \, dt \\ & \quad < \infty.
     \end{align*}
     Hence, as explained after~\eqref{Paper01_easy_criterion_martingale_feller_process_N_d}, we conclude that for $\phi \in \mathscr{C}_*(\mathcal{S}, \mathbb{R})$, conditioning on $\eta^n(0)=\boldsymbol{\xi}$, the process $(M^{n,\phi,\boldsymbol{\xi}}(t))_{t \geq 0}$ defined in~\eqref{Paper01_local_martingale_problem_discrete_approximation} is a square-integrable martingale. Identities~\eqref{Paper01_predictable_bracket_process_martingale_lipschitz_functions_discrete_IPS} and~\eqref{Paper01_definition_quadratic_variation_martingale_problem} then follow by using~\eqref{Paper01_Poisson_process_representation_discrete_approximation} and~\eqref{Paper01_quadratic_variation_stochastic_integral_poisson_point_process} in Appendix~\ref{Paper01_subsection_integral_poisson}.
\end{proof}

We will also need in this section the following version of the Burkholder-Davis-Gundy (BDG) inequality for càdlàg martingales; the inequality in this form is stated on page~527 of~\cite{mueller1995stochastic}, and in~(53) in~\cite{durrett2016genealogies}. For a square-integrable càdlàg martingale $(M(t))_{t \geq 0}$ with $M(0) = 0$, for $T \geq 0$ and $p \geq 2$,
    \begin{equation} \label{Paper01_BDG_inequality_cadlag_version}
        \mathbb{E}\left[\sup_{t \leq T} \left\vert M(t) \right\vert^{p} \right] \lesssim_{p} \mathbb{E}\left[\left(\left\langle M\right\rangle(T)\right)^{p/2} + \sup_{t \leq T} \left\vert M(t) - M(t-) \right\vert^{p}\right].
    \end{equation}

For $x \in L^{-1}\mathbb{Z}$, $\boldsymbol{\xi}=(\xi_k(y))_{y\in L^{-1}\mathbb Z,\, k\in \mathbb N_0} \in \mathcal{S}$ and $\mathcal{I} \subseteq \mathbb{N}_{0}$, let
\begin{equation} \label{Paper01_number_particles_position_x_number_mutations_in_I}
    \phi^{(x)}_{\mathcal{I}}(\boldsymbol{\xi}) \defeq \sum_{k \in \mathcal{I}} \xi_{k}(x).
\end{equation}
Observe from~\eqref{Paper01_definition_Lipschitz_function} that $\phi^{(x)}_{\mathcal{I}} \in \mathscr{C}_*(\mathcal{S}, \mathbb{R})$ for every $x \in L^{-1}\mathbb{Z}$ and~$\mathcal{I} \subseteq \mathbb{N}_{0}$. In what follows, for $n\in \mathbb N$, let $(X^{n}(t))_{t \geq 0}$ denote a simple symmetric random walk with reflecting boundaries on $\Lambda_n$, with total jump rate $m$, i.e.~a particle at $x\in \Lambda_n$ attempts to jump at rate $m$ to a uniformly chosen deme from $\{x-L^{-1},x+L^{-1}\}$, and if the chosen deme is outside $\Lambda_n$ then no jump occurs. For $t \geq 0$ and $x,y \in \Lambda_n$, we let
\begin{equation}  \label{eq:RWinitcond}
    \mathbb{P}_{x}(X^{n}(t) = y) \defeq  \mathbb{P}(X^{n}(t) = y \vert \,  X^{n}(0) = x).
\end{equation}
For notational convenience, we define, for $x,y \in L^{-1}\mathbb{Z}$ such that $\vert x \vert \vee \vert y \vert > \lambda_n$ and $t \geq 0$,
\begin{equation*}
    \mathbb{P}_{x}(X^{n}(t) = y) \defeq \left\{\begin{array}{ll}
        0 & \textrm{if } x \neq y,  \\
        1 & \textrm{otherwise.} 
    \end{array} \right.
\end{equation*}
We will also need the following result, which can be thought of as a Green's function representation for $\phi^{(x)}_{\mathcal{I}}$.
Recall the definitions of $\mathcal L^n$ in~\eqref{Paper01_generator_foutel_etheridge_model_restriction_n} and $\mathcal L^n_r$ and $\mathcal L^n_m$ in~\eqref{Paper01_infinitesimal_generator_restriction_n}.
\begin{lemma} \label{Paper01_greens_function_representation}
    For $\boldsymbol{\eta} \in \mathcal{S}$, $n \in \mathbb{N}$, $\mathcal{I} \subseteq \mathbb{N}_{0}$, $x \in \Lambda_n$ and $T \geq 0$, 
    \begin{align*}
        \mathbb{E}_{\boldsymbol{\eta}}\Bigg[\sum_{k \in \mathcal{I}} \eta^{n}_{k}(T,x)\Bigg] & = \sum_{y \in \Lambda_n} \mathbb{P}_{x}(X^{n}(T) = y) \phi^{(y)}_{\mathcal{I}}(\boldsymbol{\eta}) \\ & \quad \quad + \sum_{y \in \Lambda_n} \, \int_{0}^{T} \mathbb{P}_{x}(X^{n}(T-t) = y) \mathbb{E}_{\boldsymbol{\eta}}\Big[(\mathcal{L}^{n}_{r}\phi^{(y)}_{\mathcal{I}})(\eta^{n}(t-))\Big] \, dt.
    \end{align*}
\end{lemma}

\begin{proof}
    Our argument will be similar to the one used in the proof of Lemma~\ref{Paper01_rewriting_definition_correlation_function}. For $n \in \mathbb{N}$, $T \geq 0$, $\mathcal{I} \subseteq \mathbb{N}_{0}$ and $x \in \Lambda_n$, we define the map $g^{n,T,x}_{\mathcal{I}}: [0,T] \times \mathcal{S} \rightarrow [0,\infty)$ by letting
    \begin{equation} \label{Paper01_auxiliary_function_green_representation}
        g^{n,T,x}_{\mathcal{I}}(t,\boldsymbol{\zeta}) \defeq \sum_{y \in \Lambda_n} \mathbb{P}_{x}(X^n (T-t) = y) \phi^{(y)}_{\mathcal{I}}(\boldsymbol{\zeta}) \quad \forall (t,\boldsymbol{\zeta}) \in [0,T] \times \mathcal{S}.
    \end{equation}
    As in the proof of Lemma~\ref{Paper01_rewriting_definition_correlation_function}, our strategy will be to apply a general stochastic chain rule formula for $g^{n,T,x}_{\mathcal{I}}$ (see Lemma~\ref{Paper01_general_integration_by_parts_formula}). 
    Fix $\boldsymbol{\eta}\in \mathcal S$, $n\in \mathbb N$, $\mathcal I\subseteq \mathbb N_0$, $x\in \Lambda_n$ and $T\ge 0$.
    We start by claiming that the following statements hold:
    \begin{enumerate}[(i)]
    \item $\displaystyle \sup_{t_{1}, t_{2} \in [0,T]} \, \mathbb{E}_{\boldsymbol{\eta}}\Big[g^{n,T,x}_{\mathcal{I}}(t_1, \eta^{n}(t_2))\Big] < \infty$.
    \item The map $(t, \boldsymbol{\zeta}) \mapsto \frac{\partial}{\partial t} (g^{n,T,x}_{\mathcal{I}}(\cdot,\boldsymbol{\zeta}))(t)$ is continuous with respect to the product topology on $(0,T) \times \mathcal{S}$, and satisfies the following identity, for $t \in (0,T)$ and $\boldsymbol{\zeta} \in \mathcal{S}$:
    \begin{equation*}
        \Bigg(\frac{\partial}{\partial t} g^{n,T,x}_{\mathcal{I}}(\cdot,\boldsymbol{\zeta})\Bigg)(t) = - \frac{m}{2} \sum_{y \in \Lambda_n} \mathbb{P}_{x}(X^n(T-t) = y) (\mathcal{L}^{n}_{m}\phi^{(y)}_{\mathcal{I}})(\boldsymbol{\zeta}).
    \end{equation*}
    \item For any $t \in [0,T]$, $g^{n,T,x}_{\mathcal I}(t, \cdot) \in \mathscr{C}_*(\mathcal{S}, \mathbb{R})$.
    \item The map $(t,\boldsymbol{\zeta}) \mapsto \Big({\mathcal{L}}^{n}g_{\mathcal I}^{n,T,x}(t, \cdot)\Big)(\boldsymbol{\zeta})$ is continuous with respect to the product topology on $[0,T] \times \mathcal{S}$, and satisfies the following identity for $t \in [0,T]$ and $\boldsymbol{\zeta} \in \mathcal{S}$:
    \begin{equation*}
    \begin{aligned}
    \Big({\mathcal{L}}^{n}g_{\mathcal I}^{n,T,x}(t, \cdot)\Big)(\boldsymbol{\zeta}) =  \sum_{y \in \Lambda_n} \mathbb{P}_{x}(X^n(T-t) = y)  \left(\frac{m}{2}(\mathcal{L}^{n}_{m}\phi^{(y)}_{\mathcal{I}})(\boldsymbol{\zeta}) + ({\mathcal{L}}^{n}_{r}\phi^{(y)}_{\mathcal{I}})(\boldsymbol{\zeta}) \right).
    \end{aligned}
    \end{equation*}
    \item
$
    \!
    \begin{aligned}[t]
 \sup_{t_{1},t_{2} \in [0,T]} \; \mathbb{E}_{\boldsymbol{\eta}}\Bigg[\left(\frac{\partial}{\partial t} g_{\mathcal I}^{n,T,x}(\cdot, \eta^{n}(t_{2}))\right)^{2}(t_{1}) + \Big({\mathcal{L}}^{n}g_{\mathcal I}^{n,T,x}(t_1, \cdot)\Big)^{2} (\eta^{n}(t_{2})) \Bigg] < \infty.
        \end{aligned}
    $
\end{enumerate}
    As the proofs of claims (i)–(v) are standard, we will sketch the main ideas and omit the details. Claim~(i) follows from the definition of $g_{\mathcal I}^{n,T,x}$ in~\eqref{Paper01_auxiliary_function_green_representation} and estimate~\eqref{Paper01_crude_estimate_moments_general_configuration} in Proposition~\ref{Paper01_bound_total_mass}. To prove claim~(ii), we first notice that letting the infinitesimal generator of $X^n$ be denoted by $\frac{m}{2}\mathcal L^n_X$,
    and recalling the definition of $\mathcal L^n_m$ in~\eqref{Paper01_infinitesimal_generator_restriction_n},
   for any $\boldsymbol{\zeta} =(\zeta_k(z))_{z\in L^{-1}\mathbb Z,\, k\in \mathbb N_0}\in \mathcal{S}$ and $y \in \Lambda_n$,
     \begin{equation} \label{Paper01_duality_migration_relation_trivial_case}
        (\mathcal{L}^{n}_{X} \phi_{\mathcal{I}}^{(\cdot)}(\boldsymbol{\zeta}))(y) = \sum_{k \in \mathcal{I}} \, \sum_{z \in \{-L^{-1},L^{-1}\}} \mathds{1}_{\{y + z \in \Lambda_n\}}(\zeta_{k}(y + z) - \zeta_{k}(y)) =  (\mathcal{L}^{n}_{m} \phi_{\mathcal{I}}^{(y)})(\boldsymbol{\zeta}).
    \end{equation}
Then, claim~(ii) follows by differentiating the semigroup associated to $(X^{n}(t))_{t \geq 0}$, and then by applying~\eqref{Paper01_duality_migration_relation_trivial_case}. Claim~(iii) follows from the definition of $g^{n,T,x}_{\mathcal{I}}$ in~\eqref{Paper01_auxiliary_function_green_representation}, and the definition of $\mathscr{C}_*(\mathcal{S}, \mathbb{R})$ in~\eqref{Paper01_definition_Lipschitz_function}. Claim~(iv) follows by applying the definition of $\mathcal{L}^{n}$ in~\eqref{Paper01_infinitesimal_generator_restriction_n} to~\eqref{Paper01_auxiliary_function_green_representation}. To establish claim~(v), we recall the definitions of $(F^{b}_{k})_{k \in \mathbb{N}_{0}}$ and $(F^{d}_{k})_{k \in \mathbb{N}_{0}}$ in~\eqref{Paper01_discrete_birth_rate} and~\eqref{Paper01_discrete_death_rate}, respectively, and of $\bar{q}^{N}$ in~\eqref{Paper01_bar_polynomial}. Observe that by~\eqref{Paper01_duality_migration_relation_trivial_case}, for $\boldsymbol{\zeta} \in \mathcal{S}$ and $y \in \Lambda_n$ we have
    \begin{equation} \label{Paper01_simple_computation_migration_term_greens}
    \vert (\mathcal{L}^{n}_{m}\phi^{(y)}_{\mathcal{I}})(\boldsymbol{\zeta}) \vert \leq \sum_{k \in \mathcal{I}} \sum_{z \in \{-L^{-1}, L^{-1}\}} \, \mathds{1}_{\{y + z \in \Lambda_n\}}(\zeta_{k}(y) + \zeta_{k}(y + z)),
    \end{equation}
    and by~\eqref{Paper01_infinitesimal_generator_restriction_n}, for $\boldsymbol{\zeta} \in \mathcal{S}$ and $y \in \Lambda_n$ we have
    \begin{equation} \label{Paper01_simple_computation_reaction_term_greens}
           \vert (\mathcal{L}^{n}_{r}\phi^{(y)}_{\mathcal{I}})(\boldsymbol{\zeta}) \vert \leq \sum_{k \in \mathcal{I}} \, \left( F_{k}^{b}(\zeta(y)) +  F_{k}^{d}(\zeta(y))\right) \leq \sum_{k = 0}^{\infty} \, \left(F_{k}^{b}(\zeta(y)) +  F_{k}^{d}(\zeta(y)) \right) \leq \bar{q}^{N}(\vert \vert \zeta(y) \vert \vert_{\ell_1}),
    \end{equation}
   where the third inequality follows from~\eqref{Paper01_simple_bound_total_birth_death_rates_deterministic}. Then, claim~(v) follows from applying claims~(ii) and~(iv), and then from combining estimates~\eqref{Paper01_simple_computation_migration_term_greens} and~\eqref{Paper01_simple_computation_reaction_term_greens} with estimate~\eqref{Paper01_crude_estimate_moments_general_configuration} in Proposition~\ref{Paper01_bound_total_mass}.

   Claims (i)-(v) and Lemma~\ref{Paper01_auxilliary_lemma_easy_martingale_representation} allow us to apply a general stochastic chain rule formula, see Lemma~\ref{Paper01_general_integration_by_parts_formula}, to the map $g^{n,T,x}_{\mathcal{I}}$, obtaining that conditioning on $\eta^n(0)=\boldsymbol{\eta}$, the stochastic process $(M^{n,T,x,\boldsymbol{\eta}}_{\mathcal{I}}(t))_{t \in [0,T]}$ given by, for $t \in [0,T]$,
\begin{equation} \label{Paper01_local_martingale_definition_to_rewrite_greens_function}
\begin{aligned}
M^{n,T,x,\boldsymbol{\eta}}_{\mathcal{I}}(t) & \defeq g^{n,T,x}_{\mathcal{I}}(t, \eta^{n}(t)) - g^{n,T,x}_{\mathcal{I}}(0,\boldsymbol{\eta}) \\ & \quad \quad - \int_{0}^{t} \left(\left(\frac{\partial}{\partial t'} g^{n,T,x}_{\mathcal I}(\cdot, \eta^{n}(t'-))\right)(t') + \Big({\mathcal{L}}^{n}g^{n,T,x}_{\mathcal{I}}(t',\cdot)\Big) (\eta^{n}(t'-))\right) \, dt'
\end{aligned}
\end{equation}
is a càdlàg martingale with respect to the filtration~$\Big\{\mathcal{F}^{\eta^{n}}_{t+}\Big\}_{t \geq 0}$. By taking $t = T$ in~\eqref{Paper01_local_martingale_definition_to_rewrite_greens_function}, taking expectations on both sides of~\eqref{Paper01_local_martingale_definition_to_rewrite_greens_function}, using claims~(ii) and (iv), and then by rearranging terms, we get
\begin{equation*}
 \mathbb{E}_{\boldsymbol{\eta}}\Big[\phi^{(x)}_{\mathcal I}(\eta^n(T))\Big] = g^{n,T,x}_{\mathcal{I}}(0,\boldsymbol{\eta}) + \int_{0}^{T} \mathbb{E}_{\boldsymbol{\eta}}\Bigg[\sum_{y \in \Lambda_n} \mathbb{P}_{x}(X^n (T-t) = y) ({\mathcal{L}}^{n}_{r}\phi^{(y)}_{\mathcal{I}})(\eta^{n}(t-))\Bigg] \, dt.
\end{equation*}
The proof is completed by applying Fubini's theorem together with the definitions of $g^{n,T,x}_{\mathcal{I}}$ in~\eqref{Paper01_auxiliary_function_green_representation} and $\phi^{(x)}_{\mathcal I}$ in~\eqref{Paper01_number_particles_position_x_number_mutations_in_I}.
\end{proof}
Having proved our two preliminary lemmas, we now outline the overall proof strategy for the remainder of this section.
We will prove Theorem~\ref{Paper01_thm_existence_uniqueness_IPS_spatial_muller} by establishing that there exists $\mathfrak{A} \subseteq \mathscr{C}_{b}(\mathcal{S}, \mathbb{R})$ such that $\mathfrak{A}$ and the sequence of processes $\Big((\eta^{n}(t))_{t \geq 0}\Big)_{n \in \mathbb{N}}$ satisfy Assumptions~\ref{Paper01_assumption_candidate_set_functions_domain_generator} and~\ref{Paper01_assumption_tightness_sequence_processes}, so that Theorem~\ref{Paper01_thm_existence_uniqueness_IPS_spatial_muller} will follow from Theorem~\ref{Paper01_general_thm_convergence_feller_non_locally_compact}, and Theorem~\ref{Paper01_thm_moments_estimates_spatial_muller_ratchet} will then follow from Proposition~\ref{Paper01_bound_total_mass} and Theorem~\ref{Paper01_thm_existence_uniqueness_IPS_spatial_muller}. Recall the definition of $\mathscr{C}_{b,*}^{\textrm{cyl}}(\mathcal{S}, \mathbb{R})$ from~\eqref{Paper01_general_definition_cylindrical_function}, and observe that by Lemma~\ref{Paper01_characterisation_lipschitz_functions_semi_metric}(v), $\mathscr{C}_{b,*}^{\textrm{cyl}}(\mathcal{S}, \mathbb{R}) \subseteq \mathscr{C}_*(\mathcal{S}, \mathbb{R})$. 
Recall that $\mathcal{L}$ is the linear operator defined in~\eqref{Paper01_infinitesimal_generator} and~\eqref{Paper01_generator_foutel_etheridge_model}.
To carry out this programme, we will establish the following statements:
\begin{enumerate}[(i)]
    \item $\mathscr{C}^{\textrm{cyl}}_{b,*}(\mathcal{S}, \mathbb{R})$ satisfies Assumption~\ref{Paper01_assumption_candidate_set_functions_domain_generator}, and for any $\phi \in  \mathscr{C}^{\textrm{cyl}}_{b,*}(\mathcal{S}, \mathbb{R})$ and any $\boldsymbol{\xi} \in \mathcal{S}$,
    \begin{equation*}
        \lim_{t \rightarrow 0^{+}} \frac{(P^{n}_{t}\phi)(\boldsymbol{\xi}) - \phi(\boldsymbol{\xi})}{t} = (\mathcal{L}^{n}\phi)(\boldsymbol{\xi}) \; \forall n \in \mathbb{N} \quad \textrm{and} \quad \lim_{n \rightarrow \infty} (\mathcal{L}^{n}\phi)(\boldsymbol{\xi}) = (\mathcal{L}\phi)(\boldsymbol{\xi});
    \end{equation*}
    \item For any $\boldsymbol{\eta} \in \mathcal{S}_{0}$, conditioning on $\eta^{n}(0) = \boldsymbol{\eta}$ $\forall \, n \in \mathbb{N}$, the sequence of processes $(\eta^{n})_{n \in \mathbb{N}}$ is tight in $\mathscr{D}([0, \infty), \mathcal{S})$;
    \item For any $T \geq 0$, any compact subset $\mathcal{K} \subset \mathcal{S}$ and any $\varepsilon > 0$, there exists a compact subset $\mathscr{K}' \defeq \mathscr{K}'(T, \mathcal{K}, \varepsilon) \subset \mathcal{S}$ such that
        \begin{equation*}
            \inf_{\boldsymbol{\xi} \in \mathcal{K}} \; \inf_{n \in \mathbb{N}} \; \inf_{t \in [0,T]} \; \mathbb{P}_{\boldsymbol{\xi}}\left(\eta^{n}(t) \in \mathscr{K}'\right) \geq 1 - \varepsilon;
        \end{equation*}
    \item For any $\phi \in \mathscr{C}^{\textrm{cyl}}_{b,*}(\mathcal{S}, \mathbb{R})$, any $t \geq 0$ and any $\boldsymbol{\xi} \in \mathcal{S}$, the sequence  $\left((P^{n}_{t}\phi)(\boldsymbol{\xi})\right)_{n \in \mathbb{N}}$ converges as $n \rightarrow \infty$. Moreover, this convergence holds uniformly on compact time intervals and compact subsets of $\mathcal{S}$.
\end{enumerate}

Our first step will be to prove statement~(i), i.e.~to show that $\mathscr{C}^{\textrm{cyl}}_{b,*}(\mathcal{S}, \mathbb{R})$ satisfies
Assumption~\ref{Paper01_assumption_candidate_set_functions_domain_generator} and Assumption~\ref{Paper01_assumption_tightness_sequence_processes}(i).  We will need the following result.

\begin{lemma} \label{Paper01_convergence_infinitesimal_generator_lipischitz_semi_metric}
    For any $\phi \in \mathscr{C}_{*}(\mathcal{S}, \mathbb{R})$ and any compact set $\mathcal{K} \subset \mathcal{S}$, 
    \begin{equation} \label{P01:eq_simple_unif_bound_L_n_compacts}
        \sup_{n \in \mathbb{N}} \; \sup_{\boldsymbol{\xi} \in \mathcal{K}} \vert \mathcal{L}^n\phi(\boldsymbol{\xi}) \vert < \infty.
    \end{equation}
    Moreover, the series defining $\mathcal{L}\phi$ in~\eqref{Paper01_infinitesimal_generator} converge absolutely and uniformly on compact subsets of~$\mathcal{S}$.
    Furthermore, for any $\phi \in \mathscr{C}_{*}(\mathcal{S}, \mathbb{R})$ and any sequence $(\boldsymbol{\xi}^{n})_{n \in \mathbb{N}} \subset \mathcal{S}$ such that $\lim_{n \rightarrow \infty} \boldsymbol{\xi}^{n} = \boldsymbol{\xi}$, the following limit holds:
    \begin{equation} \label{Paper01_convergence_sequence_generators_action_bizarre_set_functions}
        \lim_{n \rightarrow \infty} \,  \Big\vert (\mathcal{L}^{n}\phi)(\boldsymbol{\xi}^{(n)}) - (\mathcal{L}\phi)(\boldsymbol{\xi}) \Big\vert = 0.
    \end{equation}
\end{lemma}

\begin{proof}
Throughout the proof, we write $\boldsymbol{\xi}\in \mathcal S$ as a shorthand for $\boldsymbol{\xi}=(\xi_i(x))_{x\in L^{-1}\mathbb Z ,\, i\in \mathbb N_0}\in \mathcal S$.
   Let $\mathcal{K} \subset \mathcal{S}$ be a compact set. Then by Proposition~\ref{Paper01_topological_properties_state_space}(i), $\mathcal K$ is $\vert \vert \vert \cdot \vert \vert \vert_{\mathcal S}$-bounded, and so, by~\eqref{Paper01_definition_state_space_formal}, in particular there exists $C_{\mathcal K}<\infty$ such that 
    \begin{equation} \label{eq:KisSbounded}
        \|\xi(y)\|_{\ell_1}\le C_{\mathcal K}(1+|y|)^2 \quad \forall \boldsymbol{\xi}\in \mathcal K \text{ and }y\in L^{-1}\mathbb Z.
    \end{equation}
   Recall the definition of $(F^{b}_{k})_{k \in \mathbb{N}_{0}}$ and $(F^{d}_{k})_{k \in \mathbb{N}_{0}}$ from~\eqref{Paper01_discrete_birth_rate} and~\eqref{Paper01_discrete_death_rate} respectively. For $A \subseteq L^{-1}\mathbb Z$ and $B \in \mathbb{N}_0$, we observe that
   by the definition of $\mathscr{C}_*(\mathcal{S}, \mathbb{R})$ in~\eqref{Paper01_definition_Lipschitz_function}, by~\eqref{Paper01_discrete_birth_rate},~\eqref{Paper01_discrete_death_rate} and since $s_k\le 1$ $\forall k\in \mathbb N_0$ by Assumption~\ref{Paper01_assumption_fitness_sequence}, for $\phi\in \mathscr{C}_*(\mathcal{S}, \mathbb{R})$, $n\in \mathbb N$ and $\boldsymbol{\xi}\in \mathcal{K}$,
    \begin{equation} \label{Paper01_simple_bound_gen_compact}
    \begin{aligned}
        & \sum_{x \in A} \; \sum_{j = B}^{\infty} \left(\Big\vert\phi\Big(\boldsymbol{\xi} + \boldsymbol{e}_{j}^{(x)}\Big) - \phi(\boldsymbol{\xi})\Big\vert F^{b}_{j}(\xi(x)) + \Big\vert\phi\Big(\boldsymbol{\xi} - \boldsymbol{e}_{j}^{(x)}\Big) - \phi(\boldsymbol{\xi})\Big\vert F^{d}_{j}(\xi(x))\right)
        \\ & \quad \lesssim_{\phi} \sum_{x \in L^{-1}\mathbb{Z}} \; \sum_{j = (B-1) \, \vee \, 0}^{\infty} \; \frac{\xi_{j}(x) \Big(q_{+}^{N}(\|\xi(x)\|_{\ell_1}) + q_{-}^{N}(\|\xi(x)\|_{\ell_1})\Big)}{(1 + \vert x \vert)^{2(1 + \deg q_-)}} \\
          & \quad \lesssim_{q_+, q_-, N}  \sum_{x \in L^{-1}\mathbb{Z}} \; \frac{\|\xi(x)\|_{\ell_1}^{\deg q_-}+ \|\xi(x)\|_{\ell_1}^{\deg q_+}+1}{(1 + \vert x \vert)^{2\deg q_-}}\sum_{j = (B-1)\vee 0}^{\infty} \; \frac{\xi_{j}(x) }{(1 + \vert x \vert)^{2}} \\
           & \quad \lesssim_{\mathcal K}  \sum_{x \in L^{-1}\mathbb{Z}} \; \sum_{j = (B-1) \,  \vee \,  0}^{\infty} \; \frac{\xi_{j}(x) }{(1 + \vert x \vert)^{2}},
    \end{aligned}
    \end{equation}
    where the second inequality follows from~\eqref{Paper01_scaled_polynomials_carrying_capacity}, and the last inequality follows by~\eqref{eq:KisSbounded} and since $\deg q_{+} < \deg q_{-}$ by Assumption~\ref{Paper01_assumption_polynomials}. Moreover, also by the definition of $\mathscr{C}_*(\mathcal{S}, \mathbb{R})$ in~\eqref{Paper01_definition_Lipschitz_function} and then by the fact that by Assumption~\ref{Paper01_assumption_polynomials}, $\deg q_- > 0$, for $\phi \in \mathscr{C}_*(\mathcal{S}, \mathbb{R})$, $n \in \mathbb{N}$ and $\boldsymbol{\xi}  \in \mathcal{K}$,
    \begin{equation} \label{Paper01_bound_mig_high_mutat_compact}
    \begin{aligned}
        & \sum_{x \in A} \; \sum_{j = B}^\infty \xi_{j}(x) \left\vert\phi\Big(\boldsymbol{\xi} + \boldsymbol{e}_{j}^{(x+L^{-1})} - \boldsymbol{e}_{j}^{(x)}\Big) + \phi\Big(\boldsymbol{\xi} + \boldsymbol{e}_{j}^{(x - L^{-1})} - \boldsymbol{e}_{j}^{(x)}\Big) - 2\phi(\boldsymbol{\xi}) \right\vert \\ & \quad \lesssim_{\phi} \sum_{x \in L^{-1}\mathbb{Z}} \; \sum_{j = B}^{\infty} \; \frac{\xi_{j}(x) }{(1 + \vert x \vert)^{2(1 + \deg q_-)}} \\ & \quad \leq \sum_{x \in L^{-1}\mathbb{Z}} \; \sum_{j = B}^{\infty} \; \frac{\xi_{j}(x) }{(1 + \vert x \vert)^{2}}.
    \end{aligned}
    \end{equation}
   
   Then,~\eqref{P01:eq_simple_unif_bound_L_n_compacts} follows from applying the definition of $\mathcal{L}^n$ in~\eqref{Paper01_generator_foutel_etheridge_model_restriction_n} and~\eqref{Paper01_infinitesimal_generator_restriction_n}, and the identities~\eqref{Paper01_simple_bound_gen_compact},~\eqref{Paper01_bound_mig_high_mutat_compact},~\eqref{Paper01_definition_state_space_formal}, and the fact that since $\mathcal{K} \subset \mathcal{S}$ is compact, by Proposition~\ref{Paper01_topological_properties_state_space}(i), $\mathcal K$ is $\vert \vert \vert \cdot \vert \vert \vert_{\mathcal S}$-bounded.
   
   Now, take $(a_n)_{n\in \mathbb N}\subset [0,\infty)$ and $(B_n)_{n\in \mathbb N}\subseteq \mathbb N$ increasing sequences such that $a_n\to \infty$ as $n \rightarrow \infty$ and $B_{n} \rightarrow \infty$ as $n \rightarrow \infty$. Let $A_n=L^{-1}\mathbb Z \cap [-a_n,a_n]$ for $n\in \mathbb N$. Then, by using~\eqref{Paper01_easy_limit_bound_high_mutation_load_space} in Proposition~\ref{Paper01_topological_properties_state_space} and~\eqref{Paper01_simple_bound_gen_compact}, it follows that for any $\phi\in \mathscr{C}_*(\mathcal{S}, \mathbb{R})$, the following limits hold:
    \begin{equation} \label{Paper01_intermediate_step_showing_control_number_mutations_strange_compacts}
    \begin{aligned}
        \lim_{n \rightarrow \infty} \; \sup_{\boldsymbol{\xi} \in \mathcal{K}} \; \sum_{x \in A_n} \; \sum_{j = B_{n}}^{\infty} \Big\vert\phi\Big(\boldsymbol{\xi} + \boldsymbol{e}_{j}^{(x)}\Big) - \phi(\boldsymbol{\xi})\Big\vert F^{b}_{j}(\xi(x)) & =0, \\ \lim_{n \rightarrow \infty} \; \sup_{\boldsymbol{\xi} \in \mathcal{K}} \; \sum_{x \in A_n} \; \sum_{j = B_{n}}^{\infty} \Big\vert\phi\Big(\boldsymbol{\xi} - \boldsymbol{e}_{j}^{(x)}\Big) - \phi(\boldsymbol{\xi})\Big\vert F^{d}_{j}(\xi(x)) & = 0.
    \end{aligned}
    \end{equation}
    and by~\eqref{Paper01_bound_mig_high_mutat_compact},
    \begin{equation} \label{Paper01_converging_migration_generator_large_mutations}
    \lim_{n \rightarrow \infty} \; \sup_{\boldsymbol{\xi} \in \mathcal{K}} 
        \sum_{x \in A_n} \; \sum_{j = B_n}^{\infty} \xi_{j}(x) \left\vert\phi\Big(\boldsymbol{\xi} + \boldsymbol{e}_{j}^{(x+L^{-1})} - \boldsymbol{e}_{j}^{(x)}\Big) + \phi\Big(\boldsymbol{\xi} + \boldsymbol{e}_{j}^{(x - L^{-1})} - \boldsymbol{e}_{j}^{(x)}\Big) - 2\phi(\boldsymbol{\xi}) \right\vert \\ 
         = 0.
    \end{equation}
     Moreover, for $\phi \in \mathscr{C}_*(\mathcal{S}, \mathbb{R})$, by the definition of $\mathscr{C}_{*}(\mathcal{S}, \mathbb{R})$ in~\eqref{Paper01_definition_Lipschitz_function}, we have
    \begin{equation} \label{Paper01_converging_migration_generator}
    \begin{aligned}
        & \limsup_{n \rightarrow \infty} \; \sup_{\boldsymbol{\xi} \in \mathcal{K}} 
        \sum_{x \in L^{-1}\mathbb{Z}\setminus A_n} \; \sum_{k = 0}^{\infty} \xi_{k}(x) \left\vert\phi\Big(\boldsymbol{\xi} + \boldsymbol{e}_{k}^{(x+L^{-1})} - \boldsymbol{e}_{k}^{(x)}\Big) + \phi\Big(\boldsymbol{\xi} + \boldsymbol{e}_{k}^{(x - L^{-1})} - \boldsymbol{e}_{k}^{(x)}\Big) - 2\phi(\boldsymbol{\xi}) \right\vert \\ 
        & \quad \lesssim_{\phi} \limsup_{n \rightarrow \infty} \; \sup_{\boldsymbol{\xi} \in \mathcal{K}} \; \sum_{x \in L^{-1}\mathbb{Z}\setminus A_n} \; \frac{\vert \vert \xi(x) \vert \vert_{\ell_1}}{(1 + \vert x \vert)^{2(1 + \deg q_{-})}} \\ & \quad = 0,
    \end{aligned}
    \end{equation}
    where the last identity follows from Proposition~\ref{Paper01_topological_properties_state_space}(ii), since $a_n\to \infty$ as $n\to \infty$. Recall the definitions of $(F^b_k)_{k \in \mathbb{N}_{0}}$ and $(F^d_k)_{k \in \mathbb{N}_{0}}$ from~\eqref{Paper01_discrete_birth_rate} and~\eqref{Paper01_discrete_death_rate}, and the definition of $\bar{q}^{N}$ in~\eqref{Paper01_bar_polynomial}. By the definition of~$\mathscr{C}_*(\mathcal{S}, \mathbb{R})$, we also have that for $\phi\in \mathscr{C}_*(\mathcal{S}, \mathbb{R})$, $n\in \mathbb N$ and $\boldsymbol{\xi} \in \mathcal{K}$,
    \begin{equation} \label{Paper01_converging_reaction_generator}
    \begin{aligned}
        & \sum_{x \in L^{-1}\mathbb{Z}\setminus A_n} \; \sum_{k = 0}^{\infty} \left(\Big\vert\phi\Big(\boldsymbol{\xi} + \boldsymbol{e}_{k}^{(x)}\Big) - \phi(\boldsymbol{\xi})\Big\vert F^{b}_{k}(\xi(x)) + \Big\vert\phi\Big(\boldsymbol{\xi} - \boldsymbol{e}_{k}^{(x)}\Big) - \phi(\boldsymbol{\xi})\Big\vert F^{d}_{k}(\xi(x))\right) \\ 
        & \quad \lesssim_{\phi} \sum_{x \in L^{-1}\mathbb{Z}\setminus A_n} \frac{1}{(1 + \vert x \vert)^{2(1 + \deg q_-)}} \sum_{k = 0}^{\infty} \Big(F^{b}_{k}(\xi(x)) + F^{d}_{k}(\xi(x))\Big) \\
        & \quad \leq \sum_{x \in L^{-1}\mathbb{Z}\setminus A_n} \frac{\bar{q}^{N}(\vert \vert \xi(x) \vert \vert_{\ell_1})}{(1 + \vert x \vert)^{2(1 + \deg q_-)}} \\ 
        & \quad \lesssim_{q_{+},q_{-},N} \sum_{x \in L^{-1}\mathbb{Z}\setminus A_n} \; \frac{\vert \vert \xi(x) \vert \vert_{\ell_1}^{1+\deg q_-}+1}{(1 + \vert x \vert)^{2(1 + \deg q_-)}},
    \end{aligned}
    \end{equation}
    where we used~\eqref{Paper01_simple_bound_total_birth_death_rates_deterministic} in the second inequality, and in the last inequality we used the definition of~$\bar{q}^N$ in~\eqref{Paper01_bar_polynomial} and the fact that, by Assumption~\ref{Paper01_assumption_polynomials}, $0 \leq \deg q_+ < \deg q_-$. 
    Then by applying~\eqref{eq:KisSbounded} to bound $\|\xi(x)\|_{\ell_1}^{\deg q_-}$ and then applying  Proposition~\ref{Paper01_topological_properties_state_space}(ii) to~\eqref{Paper01_converging_reaction_generator}, we conclude that for $\phi\in \mathscr{C}_*(\mathcal{S}, \mathbb{R})$,
    \begin{equation} \label{Paper01_converging_reaction_generator_ii}
        \lim_{n \rightarrow \infty} \; \sup_{\boldsymbol{\xi} \in \mathcal{K}} \; \sum_{x \in L^{-1}\mathbb{Z}\setminus A_n} \; \sum_{k = 0}^{\infty} \left(\Big\vert\phi\Big(\boldsymbol{\xi} + \boldsymbol{e}_{k}^{(x)}\Big) - \phi(\boldsymbol{\xi})\Big\vert F^{b}_{k}(\xi(x)) + \Big\vert\phi\Big(\boldsymbol{\xi} - \boldsymbol{e}_{k}^{(x)}\Big) - \phi(\boldsymbol{\xi})\Big\vert F^{d}_{k}(\xi(x))\right) = 0.
    \end{equation}
    It follows from the definition of $\mathcal{L}$ in~\eqref{Paper01_generator_foutel_etheridge_model} and~\eqref{Paper01_infinitesimal_generator}, and from~\eqref{Paper01_intermediate_step_showing_control_number_mutations_strange_compacts},~\eqref{Paper01_converging_migration_generator_large_mutations},~\eqref{Paper01_converging_migration_generator} and~\eqref{Paper01_converging_reaction_generator_ii}, that the series defining $\mathcal{L}\phi$ converge absolutely and uniformly on compact subsets of $\mathcal{S}$ for any $\phi \in \mathscr{C}_{*}(\mathcal{S}, \mathbb{R})$.

    It remains to establish~\eqref{Paper01_convergence_sequence_generators_action_bizarre_set_functions}. Observe that if $(\boldsymbol{\xi}^n)_{n \in \mathbb{N}} \subset \mathcal{S}$ is such that $\lim_{n \rightarrow \infty} \boldsymbol{\xi}^n = \boldsymbol{\xi} \in \mathcal{S}$, then the set $\mathcal{K} \defeq \{\boldsymbol{\xi}^{n}: \, n \in \mathbb{N}\} \cup \{\boldsymbol{\xi}\}$ is compact. Therefore,~\eqref{Paper01_convergence_sequence_generators_action_bizarre_set_functions} follows from the definitions of $\mathcal{L}$ in~\eqref{Paper01_generator_foutel_etheridge_model} and~\eqref{Paper01_infinitesimal_generator}, and of $\mathcal{L}^{n}$ in~\eqref{Paper01_generator_foutel_etheridge_model_restriction_n} and~\eqref{Paper01_infinitesimal_generator_restriction_n}, and then by~\eqref{Paper01_scaling_parameters_for_discrete_approximation} and using~\eqref{Paper01_intermediate_step_showing_control_number_mutations_strange_compacts},~\eqref{Paper01_converging_migration_generator} and~\eqref{Paper01_converging_reaction_generator_ii}; this completes the proof.
\end{proof}
We can now establish statement~(i) of our programme.
\begin{corollary} \label{cor:cylinderfnsassumptions}
The set $ \mathscr{C}^{\textrm{cyl}}_{b,*}(\mathcal{S}, \mathbb{R})$, the generator $\mathcal L$ and the sequence of processes $(\eta^n)_{n\in \mathbb N}$ satisfy Assumption~\ref{Paper01_assumption_candidate_set_functions_domain_generator} and Assumption~\ref{Paper01_assumption_tightness_sequence_processes}(i).
\end{corollary}
\begin{proof}
    By Lemma~\ref{Paper01_characterisation_lipschitz_functions_semi_metric} in the appendix and Lemma~\ref{Paper01_convergence_infinitesimal_generator_lipischitz_semi_metric}, Assumption~\ref{Paper01_assumption_candidate_set_functions_domain_generator} is satisfied.
    By combining Lemma~\ref{Paper01_characterisation_lipschitz_functions_semi_metric}(v), Lemma~\ref{Paper01_auxilliary_lemma_easy_martingale_representation} (and taking expectations on both sides of~\eqref{Paper01_local_martingale_problem_discrete_approximation}) and Lemma~\ref{Paper01_convergence_infinitesimal_generator_lipischitz_semi_metric} with standard results from the theory of local martingale problems (see e.g.~\cite[Theorem~6.4.1]{ethier2009markov}), Assumption~\ref{Paper01_assumption_tightness_sequence_processes}(i) is satisfied.
\end{proof}
To prove statement~(ii), we will need the following result, which can be thought of as a stronger version of estimate~\eqref{Paper01_moments_estimates_uniform_space} in Proposition~\ref{Paper01_bound_total_mass}. For $x \in L^{-1}\mathbb{Z}$ and $\mathcal{I} \subseteq \mathbb{N}_{0}$, recall the definition of $\phi^{(x)}_{\mathcal{I}}$ in~\eqref{Paper01_number_particles_position_x_number_mutations_in_I}.

\begin{lemma} \label{Paper01_supremum_moments_discrete_particle_system}
    For any $p \in \mathbb{N}$ and $\boldsymbol{\eta} \in \mathcal{S}_{0}$, there exists a non-decreasing function $C_{N,p,\boldsymbol{\eta}}: [0,\infty) \rightarrow [0,\infty)$, such that for all $T \geq 0$,
    \begin{equation} \label{Paper01_uniform_estimate_discrete_process_any_moments_finite_time}
        \sup_{n \in \mathbb{N}} \; \sup_{x \in L^{-1}\mathbb{Z}} \mathbb{E}_{\boldsymbol{\eta}}\left[\sup_{t \leq T} \vert \vert \eta^{n}(t,x) \vert \vert^{p}_{\ell_{1}} \right] \leq C_{N,p,\boldsymbol{\eta}}(T).
    \end{equation}
Moreover, for any $\boldsymbol{\eta} \in \mathcal{S}_{0}$ and $T \geq 0$,
    \begin{equation} \label{Paper01_supremum_control_spatial_dispersion_mass}
        \lim_{R \rightarrow \infty} \; \sup_{n \in \mathbb{N}} \; \mathbb{E}_{\boldsymbol{\eta}}\Bigg[\sup_{t \leq T} \, \sum_{\{x \in L^{-1}\mathbb{Z}: \, \vert x \vert \geq R\}} \frac{\vert \vert \eta^{n}(t,x) \vert \vert_{\ell_{1}}}{(1 + \vert x \vert)^{2}}\Bigg] = 0.
    \end{equation}
    Also, for any $\boldsymbol{\eta} \in \mathcal{S}_{0}$ and $T \geq 0$,
    \begin{equation} \label{Paper01_local_time_bound_large_mutation_load}
        \lim_{k \rightarrow \infty} \; \sup_{n \in \mathbb{N}} \; \sup_{x \in L^{-1}\mathbb{Z}} \, \sup_{t \leq T} \, \mathbb{E}_{\boldsymbol{\eta}}\Bigg[ \sum_{j = k}^{\infty} \eta^{n}_{j}(t,x)\Bigg] = 0.
    \end{equation}
    Furthermore, for any $\boldsymbol{\eta} \in \mathcal{S}_{0}$ and $T \geq 0$,
    \begin{equation} \label{Paper01_supremum_time_bound_large_mutation_load}
        \lim_{k \rightarrow \infty} \; \sup_{n \in \mathbb{N}} \; \sup_{x \in L^{-1}\mathbb{Z}} \, \mathbb{E}_{\boldsymbol{\eta}}\Bigg[\sup_{t \leq T} \sum_{j = k}^{\infty} \eta^{n}_{j}(t,x)\Bigg] = 0.
    \end{equation}
\end{lemma}

\begin{proof}

We will divide the proof into steps corresponding to each of the estimates~\eqref{Paper01_uniform_estimate_discrete_process_any_moments_finite_time}-\eqref{Paper01_supremum_time_bound_large_mutation_load}.

\medskip

\noindent \underline{Step~(1): Proof of~\eqref{Paper01_uniform_estimate_discrete_process_any_moments_finite_time}:}

\medskip

As observed after~\eqref{Paper01_number_particles_position_x_number_mutations_in_I}, for any $x\in L^{-1}\mathbb Z$ we have $\phi^{(x)}_{\mathbb{N}_{0}} \in \mathscr{C}_*(\mathcal{S}, \mathbb{R})$. Hence, by Lemma~\ref{Paper01_auxilliary_lemma_easy_martingale_representation}, for any $\boldsymbol{\eta}\in \mathcal S$, conditioning on $\eta^n(0)=\boldsymbol{\eta}$,
the process $(M^{n,(x),\boldsymbol{\eta}}(t))_{t \geq 0}$ given by, for any $T \geq 0$,
      \begin{equation} \label{Paper01_martingale_formulation_local_moments}
    \begin{aligned}
        \vert \vert \eta^{n}(T,x) \vert \vert_{\ell_{1}} =  \vert \vert {\eta}(x) \vert \vert_{\ell_{1}} + M^{n,(x),\boldsymbol{\eta}}(T) + \int_{0}^{T} (\mathcal{L}^{n} \phi_{\mathbb N_0}^{(x)})(\eta^{n}(t-)) \, dt
    \end{aligned}
    \end{equation}
    is a square-integrable càdlàg martingale. For $T \geq 0$, $p \in \mathbb{N}$ and $x \in L^{-1}\mathbb{Z}$, by~\eqref{Paper01_martingale_formulation_local_moments}, we have
    \begin{equation} \label{Paper01_equivalence_bound_supremum_moments_number_particles}
    \begin{aligned}
    \sup_{t \leq T} \vert \vert \eta^{n}(t,x) \vert \vert_{\ell_{1}}^{p} & \lesssim_{p} \vert \vert {\eta}(x) \vert \vert_{\ell_{1}}^{p} + \sup_{t \leq T} \left\vert M^{n,(x),\boldsymbol{\eta}}(t) \right\vert^{p} + \sup_{t \leq T} \left(\int_{0}^{t}
    \mathcal{L}^{n}\phi_{\mathbb N_0}^{(x)}(\eta^{n}(\tau -))
    \, d\tau\right)^{p} \\ & \leq \; \vert \vert {\eta}(x) \vert \vert_{\ell_{1}}^{p} + \sup_{t \leq T} \left\vert M^{n,(x),\boldsymbol{\eta}}(t) \right\vert^{p} + \left(\int_{0}^{T} \left\vert \mathcal{L}^{n}\phi_{\mathbb N_0}^{(x)}(\eta^{n}(t -)) \right\vert \, dt \right)^{p}.
    \end{aligned}
    \end{equation}
    We now bound the expectation of each of the terms on the right-hand side of~\eqref{Paper01_equivalence_bound_supremum_moments_number_particles} separately. Tackling the integral term first, observe that by Jensen's inequality, and then by the definition of~$\mathcal{L}^{n}$ in~\eqref{Paper01_infinitesimal_generator_restriction_n} and~\eqref{Paper01_generator_foutel_etheridge_model_restriction_n},
    \begin{equation} \label{Paper01_bound_moment_FV_process}
    \begin{aligned}
     & \left(\int_{0}^{T} \left\vert \mathcal{L}^{n}\phi_{\mathbb N_0}^{(x)}(\eta^{n}(t -)) \right\vert \, dt \right)^{p} \\ 
     &\lesssim_{T,p} \int_{0}^{T} \left\vert \mathcal{L}^{n}\phi_{\mathbb N_0}^{(x)}(\eta^{n}(t -)) \right\vert^{p} \, dt \\ 
     &\lesssim_{p} m^{p}\int_{0}^{T} \left(\vert \vert \eta^{n}(t-, x+L^{-1}) \vert \vert_{\ell_{1}}^{p} + \vert \vert \eta^{n}(t-, x - L^{-1}) \vert \vert_{\ell_{1}}^{p} + \vert \vert \eta^{n}(t-, x) \vert \vert_{\ell_{1}}^{p}\right) \, dt \\ & \quad \quad \quad \quad + \int_{0}^{T} \Bigg(\sum_{k = 0}^{\infty} \left(F^{b}_{k}\left(\eta^{n}(t-,x)\right) + F^{d}_{k}\left(\eta^{n}(t-,x)\right)\right) \Bigg)^{p}\, dt \\ 
     &\lesssim_m \int_{0}^{T} \left( \vert \vert \eta^{n}(t-, x+L^{-1}) \vert \vert_{\ell_{1}}^{p} + \vert \vert \eta^{n}(t-, x - L^{-1}) \vert \vert_{\ell_{1}}^{p} + \vert \vert \eta^{n}(t-, x) \vert \vert_{\ell_{1}}^{p} + \bar{q}^{N}(\vert \vert \eta^{n}(t-, x) \vert \vert_{\ell_{1}})^{p}\right) dt,
    \end{aligned}
    \end{equation}
    where~$(F^{b}_{k})_{k \in \mathbb{N}_{0}}$ and~$(F^{d}_{k})_{k \in \mathbb{N}_{0}}$ are given by~\eqref{Paper01_discrete_birth_rate} and~\eqref{Paper01_discrete_death_rate} respectively, and we used~\eqref{Paper01_simple_bound_total_birth_death_rates_deterministic} in the last inequality. By taking expectations on both sides of~\eqref{Paper01_bound_moment_FV_process}, and applying estimate~\eqref{Paper01_moments_estimates_uniform_space} in Proposition~\ref{Paper01_bound_total_mass}, we conclude that for $p \in \mathbb{N}$ and $\boldsymbol{\eta}\in \mathcal S_0$, there exists a non-decreasing function $C^{(1)}_{N,p,\boldsymbol{\eta}}: [0, \infty) \rightarrow [0, \infty)$ such that for $T \geq 0$,
    \begin{equation} \label{Paper01_bound_estimate_FV_moments_intermediate}
        \sup_{n \in \mathbb{N}} \, \sup_{x \in L^{-1}\mathbb{Z}} \, \mathbb{E}_{\boldsymbol{\eta}}\Bigg[\left(\int_{0}^{T} \left\vert \mathcal{L}^{n}\phi_{\mathbb N_0}^{(x)}(\eta^{n}(t -)) \right\vert \, dt \right)^{p}\Bigg] \le C^{(1)}_{N,p,\boldsymbol{\eta}}(T).
    \end{equation}
    
    It remains to tackle the càdlàg martingale term on the right-hand side of~\eqref{Paper01_equivalence_bound_supremum_moments_number_particles}. 
    Observe that as a corollary of estimate~\eqref{Paper01_discrete_approximation_does_not_explode} in Proposition~\ref{Paper01_bound_total_mass} and~\eqref{Paper01_martingale_formulation_local_moments}, we have for $T\ge 0$ and $\boldsymbol{\eta}\in \mathcal S$, conditioning on $\eta^n(0)=\boldsymbol{\eta}$, almost surely,
    \begin{equation} \label{Paper01_bound_size_jumps}
        \sup_{t \leq T} \left\vert M^{n,(x),\boldsymbol{\eta}}(t) - M^{n,(x),\boldsymbol{\eta}}(t-) \right\vert =  \sup_{t \leq T} \left\vert \vert \vert \eta^{n}(t,x) \vert \vert_{\ell_1} - \vert \vert \eta^{n}(t-,x) \vert \vert_{\ell_1}  \right\vert \leq 1.
    \end{equation}
    Moreover, by using~\eqref{Paper01_predictable_bracket_process_martingale_lipschitz_functions_discrete_IPS} in Lemma~\ref{Paper01_auxilliary_lemma_easy_martingale_representation}, we have
     \begin{equation} \label{Paper01_predictable_bracket_process_local_number_particles}
    \begin{aligned}
     & \left\langle M^{n,(x),\boldsymbol{\eta}} \right\rangle (T) \\ & \leq \frac{m}{2} \int_{0}^{T}  \left(\vert \vert \eta^{n}(t-, x+L^{-1}) \vert \vert_{\ell_{1}} + \vert \vert \eta^{n}(t-, x - L^{-1}) \vert \vert_{\ell_{1}} + 2\vert \vert \eta^{n}(t-, x) \vert \vert_{\ell_{1}}\right) \, dt \\ & \quad + \int_{0}^{T} \sum_{k = 0}^{\infty} \left(F^{b}_{k}\left(\eta^{n}(t-,x)\right) + F^{d}_{k}\left(\eta^{n}(t-,x)\right)\right) \, dt \\ 
     & \lesssim_{m} \int_{0}^{T} \left(\vert \vert \eta^{n}(t-, x+L^{-1}) \vert \vert_{\ell_{1}} + \vert \vert \eta^{n}(t-, x - L^{-1}) \vert \vert_{\ell_{1}} + 2\vert \vert \eta^{n}(t-, x) \vert \vert_{\ell_{1}} + \bar{q}^{N}(\vert \vert \eta^{n}(t-, x) \vert \vert_{\ell_{1}})\right) \, dt,
     \end{aligned}
     \end{equation}
     where we used~\eqref{Paper01_simple_bound_total_birth_death_rates_deterministic} in the last inequality. For $p\ge 2$, by taking both sides of~\eqref{Paper01_predictable_bracket_process_local_number_particles} to the power of $p/2$, 
     and then using Jensen's inequality and estimate~\eqref{Paper01_moments_estimates_uniform_space} again,
     we conclude that for $p\ge 2$ and $\boldsymbol{\eta}\in \mathcal S_0$ there exists a non-decreasing function $C^{(2)}_{N,p,\boldsymbol{\eta}}: [0, \infty) \rightarrow [0, \infty)$ such that for $T \geq 0$,
    \begin{equation} \label{Paper01_bound_predictable_moments_martingale}
        \sup_{n \in \mathbb{N}} \, \sup_{x \in L^{-1}\mathbb{Z}} \, \mathbb{E}_{\boldsymbol{\eta}}\left[\left(\left\langle M^{n,(x),\boldsymbol{\eta}}\right\rangle(T)\right)^{p/2}\right] \leq C^{(2)}_{N,p,\boldsymbol{\eta}}(T). 
    \end{equation}
    Therefore, by the form of the BDG inequality stated in~\eqref{Paper01_BDG_inequality_cadlag_version}, and since, conditioning on $\eta^n(0)=\boldsymbol{\eta}$, $M^{n,(x),\boldsymbol{\eta}}(0) = 0$ and $M^{n,(x),\boldsymbol{\eta}}$ is a square-integrable càdlàg martingale, combining~\eqref{Paper01_bound_size_jumps} and~\eqref{Paper01_bound_predictable_moments_martingale} we have
    \begin{equation} \label{Paper01_bound_sup_moments_finite_times_martingale_final_term}
        \sup_{n \in \mathbb{N}} \; \sup_{x \in L^{-1}\mathbb{Z}} \mathbb{E}_{\boldsymbol{\eta}}\Big[\sup_{t \leq T} \left\vert M^{n,(x),\boldsymbol{\eta}}(t) \right\vert^{p} \Big] \lesssim_{p} 1 + C^{(2)}_{N,p,\boldsymbol{\eta}}(T).
    \end{equation}
    The estimate~\eqref{Paper01_uniform_estimate_discrete_process_any_moments_finite_time} follows by taking expectations on both sides of~\eqref{Paper01_equivalence_bound_supremum_moments_number_particles}, and then by applying~\eqref{Paper01_bound_estimate_FV_moments_intermediate} and~\eqref{Paper01_bound_sup_moments_finite_times_martingale_final_term}.

    \medskip

    \noindent \underline{Step (2): Proof of~\eqref{Paper01_supremum_control_spatial_dispersion_mass}:}

    \medskip

    Take $\boldsymbol{\eta} \in \mathcal{S}_{0}$, $T \geq 0$ and $R > 0$. Then
    \begin{equation} \label{Paper01_auxliary_estimate_spatial_direction}
    \begin{aligned}
    &  \sup_{n\in \mathbb N} \; \mathbb{E}_{\boldsymbol{\eta}}\Bigg[\sup_{t \leq T} \, \sum_{\substack{\{x \in L^{-1}\mathbb{Z}: \, \vert x \vert \geq R\}}} \frac{\vert \vert \eta^{n}(t,x) \vert \vert_{\ell_{1}}}{(1 + \vert x \vert)^{2}} \Bigg] \\ 
    & \quad \leq  \sum_{\substack{\{x \in L^{-1}\mathbb{Z}: \, \vert x \vert \geq R\}}} \frac{1}{(1 + \vert x \vert)^{2}} \sup_{n\in \mathbb N} \mathbb{E}_{\boldsymbol{\eta}}\Bigg[\sup_{t \leq T} \, \vert \vert \eta^{n}(t,x) \vert \vert_{\ell_{1}} \Bigg] \\ 
    & \quad \le C_{N,1,\boldsymbol{\eta}}(T) \sum_{\{\substack{x \in L^{-1}\mathbb{Z}: \, \vert x \vert \geq R\}}} \frac{1}{(1 + \vert x \vert)^{2}},
\end{aligned}
\end{equation}
where we used~\eqref{Paper01_uniform_estimate_discrete_process_any_moments_finite_time} in the last inequality. The statement~\eqref{Paper01_supremum_control_spatial_dispersion_mass} follows by taking the limit as $R \rightarrow \infty$ on both sides of~\eqref{Paper01_auxliary_estimate_spatial_direction}.

\medskip

\noindent \underline{Step (3): Proof of~\eqref{Paper01_local_time_bound_large_mutation_load}:}

\medskip

For $k \in \mathbb{N}_{0}$, let $\mathcal{I}_{k} \defeq \{j \in \mathbb{N}_{0}: \, j \geq k\}$, and let $\phi^{(x)}_{\mathcal{I}_{k}}: \mathcal{S} \rightarrow \mathbb{N}_{0}$ be given by~\eqref{Paper01_number_particles_position_x_number_mutations_in_I}, and recall the definition of $(X^{n}(t))_{t \geq 0}$ from after~\eqref{Paper01_number_particles_position_x_number_mutations_in_I}. By Lemma~\ref{Paper01_greens_function_representation}, for $\boldsymbol{\eta} = (\eta_{k}(x))_{k \in \mathbb{N}_{0}, \, x \in L^{-1}\mathbb{Z}}\in \mathcal{S}_{0}$, $n \in \mathbb{N}$, $k \in \mathbb{N}$, $x \in \Lambda_n$, $T \geq 0$ and $t \in [0,T]$,
\begin{equation} \label{Paper01_martingale_formulation_particles_with_mutations_discrete_IPS}
\begin{aligned}
    \mathbb{E}_{\boldsymbol{\eta}}\Bigg[\sum_{j = k}^{\infty} \eta^{n}_{j}(t,x)\Bigg] & = \sum_{y \in \Lambda_n} \mathbb{P}_{x}(X^{n}(t) = y) \sum_{j = k}^{\infty} \eta_{j}(y) \\ & \quad \quad + \sum_{y \in \Lambda_n} \, \int_{0}^{t} \mathbb{P}_{x}(X^{n}(t-t') = y) \mathbb{E}_{\boldsymbol{\eta}}\Big[(\mathcal{L}^{n}_{r}\phi^{(y)}_{\mathcal{I}_k})(\eta^{n}(t'-))\Big] \, dt'.
    \end{aligned}
    \end{equation}
    We now bound the two terms on the right-hand side of~\eqref{Paper01_martingale_formulation_particles_with_mutations_discrete_IPS} separately. For the first term, we have
    \begin{equation*}
    \begin{aligned}
       \sum_{y \in \Lambda_n} \mathbb{P}_{x}(X^{n}(t) = y) \sum_{j = k}^{\infty} \eta_{j}(y) \leq   \sum_{y \in \Lambda_n} \mathbb{P}_{x}(X^{n}(t) = y) \sup_{z \in L^{-1}\mathbb{Z}} \; \sum_{j = k}^{\infty} \eta_{j}(z) = \sup_{z \in L^{-1}\mathbb{Z}} \; \sum_{j = k}^{\infty} \eta_{j}(z).
    \end{aligned}
    \end{equation*}
    Since $\boldsymbol{\eta} \in \mathcal{S}_{0}$, by the definition of $\mathcal{S}_{0}$ in~\eqref{Paper01_set_initial_configurations}, it follows that there exists $k^{*} = k^{*}(\boldsymbol{\eta}) \in \mathbb{N}_{0}$ such that for $k \geq k^{*}$,
    \begin{equation} \label{Paper01_bound_type_space_first_sum_initial_particles}
        \sup_{n \in \mathbb{N}} \; \sup_{x \in \Lambda_n} \; \sup_{t \leq T} \; \sum_{y \in \Lambda_n} \mathbb{P}_{x}(X^{n}(t) = y) \sum_{j = k}^{\infty} \eta_{j}(y) = 0 \quad \forall k \geq k^{*}.
    \end{equation}
   We now bound the integral term on the right-hand side of~\eqref{Paper01_martingale_formulation_particles_with_mutations_discrete_IPS}. Recalling the definition of $\mathcal{L}^{n}_{r}$ from~\eqref{Paper01_infinitesimal_generator_restriction_n}, observe that
   \begin{equation} \label{Paper01_basic_application_generator_reaction}
   \begin{aligned}
       & \sum_{y \in \Lambda_n} \, \int_{0}^{t} \mathbb{P}_{x}(X^{n}(t-t') = y) \mathbb{E}_{\boldsymbol{\eta}}\Big[(\mathcal{L}^{n}_{r}\phi^{(y)}_{\mathcal{I}_k})(\eta^{n}(t'-))\Big] \, dt' \\ & \quad  = \int_{0}^{t} \sum_{y \in \Lambda_n} \mathbb{P}_{x}(X^{n}(t-t') = y) \mathbb{E}_{\boldsymbol{\eta}}\Bigg[\mathds 1_{\{k\le K_n\}}\sum_{j = k}^{K_n} F^{b}_{j}(\eta^{n}(t'-,y)) - \sum_{j = k}^{\infty} F^{d}_{j}(\eta^{n}(t'-,y))\Bigg] \, dt',
    \end{aligned}
   \end{equation}
   where $(F^{b}_{j})_{j \in \mathbb{N}_{0}}$ and $(F^{d}_{j})_{j \in \mathbb{N}_{0}}$ are given by~\eqref{Paper01_discrete_birth_rate} and~\eqref{Paper01_discrete_death_rate} respectively. By~\eqref{Paper01_discrete_birth_rate} and since $(s_{j})_{j \in \mathbb{N}_{0}}$ is monotonically non-increasing by Assumption~\ref{Paper01_assumption_fitness_sequence}, for ${u} = (u_{j})_{j \in \mathbb{N}_{0}} \in \ell_1^+$,
   \begin{equation} \label{Paper01_bound_birth_rate_high_mutation_load}
       \mathds 1_{\{k\le K_n\}}\sum_{j = k}^{K_n} F^{b}_j(u) \leq     \sum_{j = k}^{\infty} F^{b}_j(u) \leq s_{k-1} \vert \vert u \vert \vert_{\ell_{1}}q^{N}_{+}(\vert \vert u \vert \vert_{\ell_{1}}), 
        \end{equation}
   where $q^{N}_{+}$ is the polynomial given by~\eqref{Paper01_scaled_polynomials_carrying_capacity}. Moreover, for ${u} \in \ell_1^+$ and $j \in \mathbb{N}_{0}$, $F^{d}_{j}(u) \geq 0$ by~\eqref{Paper01_discrete_death_rate} and Assumption~\ref{Paper01_assumption_polynomials}. Therefore, by applying~\eqref{Paper01_bound_birth_rate_high_mutation_load} to~\eqref{Paper01_basic_application_generator_reaction}, we can write
    \begin{equation} \label{Paper01_almost_final_step_control_high_mutation_load}
    \begin{aligned}
         & \sup_{n \in \mathbb{N}}\; \sup_{x \in \Lambda_n}\; \sup_{t \leq T} \sum_{y \in \Lambda_n} \, \int_{0}^{t} \mathbb{P}_{x}(X^{n}(t-t') = y) \mathbb{E}_{\boldsymbol{\eta}}\Big[(\mathcal{L}^{n}_{r}\phi^{(y)}_{\mathcal{I}_k})(\eta^{n}(t'-))\Big] \, dt' \\ & \quad \leq \sup_{t\le T}\int_{0}^{t} \sup_{n \in \mathbb{N}}\; \sup_{x \in \Lambda_n}\sum_{y \in \Lambda_n} \mathbb{P}_{x}(X^n(t-t') = y)  \mathbb{E}_{\boldsymbol{\eta}}\Big[s_{k-1}\vert \vert \eta^{n}(t'-,y) \vert \vert_{\ell_{1}}q^{N}_{+}(\vert \vert \eta^{n}(t'-,y) \vert \vert_{\ell_{1}})\Big] \, dt' \\ & \quad \lesssim_{T, \boldsymbol{\eta},N,q_{+},q_-} s_{k-1},
    \end{aligned}
    \end{equation}
    where for the second estimate we used the fact by Assumption~\ref{Paper01_assumption_polynomials} and~\eqref{Paper01_scaled_polynomials_carrying_capacity}, $q^{N}_{+}$ is a non-negative polynomial, and applied estimate~\eqref{Paper01_moments_estimates_uniform_space} in Proposition~\ref{Paper01_bound_total_mass}. Since by Assumption~\ref{Paper01_assumption_fitness_sequence}, $s_{k} \rightarrow 0$ as $k \rightarrow \infty$, by taking the limit as $k \rightarrow \infty$ on both sides of~\eqref{Paper01_almost_final_step_control_high_mutation_load}, we get
    \begin{equation} \label{Paper01_bound_integral_term_discrete_IPS_type_space}
       \lim_{k \rightarrow \infty} \, \sup_{n \in \mathbb{N}}\; \sup_{x \in \Lambda_n}\; \sup_{t \leq T} \sum_{y \in \Lambda_n} \, \int_{0}^{t} \mathbb{P}_{x}(X^{n}(t-t') = y) \mathbb{E}_{\boldsymbol{\eta}}\Big[(\mathcal{L}^{n}_{r}\phi^{(y)}_{\mathcal{I}_k})(\eta^{n}(t'-))\Big] \, dt' = 0.
    \end{equation}
    Estimate~\eqref{Paper01_local_time_bound_large_mutation_load} follows by applying~\eqref{Paper01_bound_type_space_first_sum_initial_particles} and~\eqref{Paper01_bound_integral_term_discrete_IPS_type_space} to~\eqref{Paper01_martingale_formulation_particles_with_mutations_discrete_IPS}, and since particles outside $\Lambda_n$ are frozen in the process $\eta^n$.

    \medskip

    \noindent \underline{Step (4): Proof of~\eqref{Paper01_supremum_time_bound_large_mutation_load}:}

    \medskip

    For $k \in \mathbb{N}_{0}$, recall the definition of $\mathcal{I}_{k} \subseteq \mathbb{N}_{0}$ from the beginning of Step~(3), and for $x \in L^{-1}\mathbb{Z}$, let $\phi^{(x)}_{\mathcal{I}_k}$ be given by~\eqref{Paper01_number_particles_position_x_number_mutations_in_I}. As observed after~\eqref{Paper01_number_particles_position_x_number_mutations_in_I}, we have $\phi^{(x)}_{\mathcal{I}_k} \in \mathscr{C}_*(\mathcal{S}, \mathbb{R})$. Hence, by Lemma~\ref{Paper01_auxilliary_lemma_easy_martingale_representation}, for $\boldsymbol \eta = (\eta_{k}(x))_{k \in \mathbb{N}_{0}, \, x \in L^{-1}\mathbb Z} \in \mathcal{S}_0$, conditioning on $\eta^n(0)=\boldsymbol{\eta}$, there exists a càdlàg martingale $(M^{n,x,\boldsymbol{\eta}}_{\mathcal{I}_k}(t))_{t \geq 0}$ such that for $t \geq 0$,
\begin{equation*}
    \sum_{j = k}^{\infty} \eta^{n}_j(t,x) = \sum_{j = k}^{\infty} \eta_j(x) + M^{n,x,\boldsymbol{\eta}}_{\mathcal{I}_k}(t) + \int_{0}^{t} \left(\mathcal{L}^{n}\phi^{(x)}_{\mathcal I_k}\right)(\eta^{n}(t'-)) \, dt'.
\end{equation*}
Hence, for $T \geq 0$,
\begin{equation} \label{Paper01_martingale_problem_formulation_number_particles_high_mutations}
    \sup_{t \leq T} \sum_{j = k}^{\infty} \eta^{n}_j(t,x) \leq  \sum_{j = k}^{\infty} \eta_j(x) + \sup_{t \leq T} \Big\vert M^{n,x,\boldsymbol{\eta}}_{\mathcal{I}_k}(t) \Big\vert + \int_{0}^{T} \Big\vert\left(\mathcal{L}^{n}\phi^{(x)}_{\mathcal I_k}\right)(\eta^{n}(t'-))\Big\vert \, dt'.
\end{equation}
We now bound the terms on the right-hand side of~\eqref{Paper01_martingale_problem_formulation_number_particles_high_mutations} separately. For the martingale term, recall the definitions of $(F^b_j)_{j \in \mathbb{N}_{0}}$ and $(F^d_j)_{j \in \mathbb{N}_{0}}$ in~\eqref{Paper01_discrete_birth_rate} and~\eqref{Paper01_discrete_death_rate} respectively. By~\eqref{Paper01_predictable_bracket_process_martingale_lipschitz_functions_discrete_IPS} in Lemma~\ref{Paper01_auxilliary_lemma_easy_martingale_representation}, we have that for $T \geq 0$,
\begin{equation} \label{Paper01_martingale_bound_particles_carrying_high_number_mutations}
\begin{aligned}
\left\langle M^{n,x,\boldsymbol{\eta}}_{\mathcal{I}_k} \right\rangle(T) & \leq \frac{m}{2} \int_{0}^{T}  \sum_{j = k}^{\infty} \left(\eta^{n}_j(t'-, x+L^{-1})  + \eta^{n}_j(t'-, x - L^{-1}) + 2\eta^{n}_j(t'-, x) \right) \, dt' \\ & \quad \quad + \int_{0}^{T} \sum_{j = k}^{\infty} \left(F^{b}_{j}\left(\eta^{n}(t'-,x)\right) + F^{d}_{j}\left(\eta^{n}(t'-,x)\right)\right) \, dt'.
\end{aligned}
\end{equation}
We now bound the expectation of the terms on the right-hand side of~\eqref{Paper01_martingale_bound_particles_carrying_high_number_mutations} separately. For the first term, by~\eqref{Paper01_local_time_bound_large_mutation_load}, for $k \in \mathbb{N}_{0}$, there exists $c_{k,\boldsymbol{\eta}}^{(m)}: [0, \infty) \rightarrow [0, \infty)$ such that for any $T \geq 0$,
$$\lim_{k \rightarrow \infty} c_{k,\boldsymbol{\eta}}^{(m)}(T) = 0,$$
and
\begin{equation} \label{Paper01_intermediate_step_migration_term_control_mutation_load_uniform_time}
    \sup_{n \in \mathbb{N}} \, \sup_{x \in L^{-1}\mathbb{Z}} \, \frac{m}{2} \mathbb{E}_{\boldsymbol{\eta}}\Bigg[\int_{0}^{T}  \sum_{j = k}^{\infty} \left(\eta^{n}_j(t'-, x+L^{-1})  + \eta^{n}_j(t'-, x - L^{-1}) + 2\eta^{n}_j(t'-, x) \right) \, dt'\Bigg] \leq c_{k,\boldsymbol{\eta}}^{(m)}(T).
\end{equation}
For the second term on the right-hand side of~\eqref{Paper01_martingale_bound_particles_carrying_high_number_mutations}, we observe that since $s_{j} \leq 1$ for every $j \in \mathbb{N}_{0}$ by Assumption~\ref{Paper01_assumption_fitness_sequence}, and then by~\eqref{Paper01_discrete_birth_rate} and~\eqref{Paper01_discrete_death_rate}, we have that for any $u = (u_{j})_{j \in \mathbb{N}_{0}} \in \ell_{1}^{+}$ and $k \in \mathbb{N}$,
     \begin{equation} \label{Paper01_simple_bound_total_birth_death_rates_deterministic_modified_starting_point}
         \sum_{j = k}^{\infty} (F^{b}_{k}(u) + F^{d}_{k}(u)) \leq \Bigg(\sum_{j = k-1}^{\infty} u_j\Bigg)\Big({q}^{N}_+(\vert \vert u \vert \vert_{\ell_{1}}) + {q}^{N}_-(\vert \vert u \vert \vert_{\ell_{1}})\Big),
     \end{equation}
where $q^N_{+}, q^N_{-}$ are the non-negative polynomials given by~\eqref{Paper01_scaled_polynomials_carrying_capacity}. Moreover, observe that for any $a \geq 0$, $u=(u_j)_{j\in \mathbb N_0} \in \ell_1^+$ and $k\in \mathbb N$,
\begin{equation} \label{Paper01_especial_observation_cauchy_schwarz}
\begin{aligned}
   \Bigg(\sum_{j = k-1}^{\infty} u_j\Bigg) \vert \vert u \vert \vert_{\ell_{1}}^{a} =   \Bigg(\sum_{j = k-1}^{\infty} u_j\Bigg)^{1/2} \Bigg(\sum_{j = k-1}^{\infty} u_j\Bigg)^{1/2} \vert \vert u \vert \vert_{\ell_{1}}^{a}  \leq \Bigg(\sum_{j = k-1}^{\infty} u_j\Bigg)^{1/2} \vert \vert u \vert \vert_{\ell_{1}}^{\frac 12 +a}.
\end{aligned}
\end{equation}
Hence, by using~\eqref{Paper01_simple_bound_total_birth_death_rates_deterministic_modified_starting_point}, the fact that by Assumption~\ref{Paper01_assumption_polynomials}, $q_{+}$ and $q_{-}$ are non-negative polynomials, and then~\eqref{Paper01_especial_observation_cauchy_schwarz} and the Cauchy-Schwarz inequality, we conclude that there exists $c^{(r)}_{k,\boldsymbol{\eta}}: [0, \infty) \rightarrow [0, \infty)$ with $\lim_{k \rightarrow \infty} c^{(r)}_{k,\boldsymbol{\eta}}(T) = 0$ for any $T \geq 0$, such that
\begin{equation} \label{Paper01_intermediate_step_reaction_term_control_mutation_load_uniform_time}
\begin{aligned}
    & \sup_{n \in \mathbb{N}} \, \sup_{x \in L^{-1}\mathbb{Z}} \mathbb{E}_{\boldsymbol{\eta}}\Bigg[\int_{0}^{T} \sum_{j = k}^{\infty} \left(F^{b}_{j}\left(\eta^{n}(t'-,x)\right) + F^{d}_{j}\left(\eta^{n}(t'-,x)\right)\right) \, dt'\Bigg] \\ & \quad \lesssim_{q_{+},q_{-},N} \int_{0}^{T} \sup_{n \in \mathbb{N}} \, \sup_{x \in L^{-1}\mathbb{Z}} \sum_{a = 0}^{\deg q_-}\mathbb{E}_{\boldsymbol{\eta}}\Bigg[\sum_{j = k-1}^{\infty} \eta^{n}_j(t'-,x)\Bigg]^{1/2} \mathbb{E}_{\boldsymbol{\eta}}\Big[\vert \vert \eta^n(t'-,x) \vert \vert_{\ell_1}^{1 + 2a}\Big]^{1/2} \, dt' \\ & \quad \leq c^{(r)}_{k,\boldsymbol{\eta}}(T),
\end{aligned}
\end{equation}
where in the last inequality we used~\eqref{Paper01_local_time_bound_large_mutation_load} and estimate~\eqref{Paper01_moments_estimates_uniform_space} in Proposition~\ref{Paper01_bound_total_mass}. By taking expectations on both sides of~\eqref{Paper01_martingale_bound_particles_carrying_high_number_mutations}, and then by applying~\eqref{Paper01_intermediate_step_migration_term_control_mutation_load_uniform_time} and~\eqref{Paper01_intermediate_step_reaction_term_control_mutation_load_uniform_time},
\begin{equation*}
\begin{aligned}
    \sup_{n \in \mathbb{N}} \; \sup_{x \in L^{-1}\mathbb{Z}} \mathbb{E}\left[\left\langle M^{n,x,\boldsymbol{\eta}}_{\mathcal{I}_k} \right\rangle(T)\right] \lesssim_{q_+,q_-,N} c^{(m)}_{k,\boldsymbol{\eta}}(T) + c^{(r)}_{k,\boldsymbol{\eta}}(T),
\end{aligned}
\end{equation*}
where $c^{(m)}_{k,\boldsymbol{\eta}}(T)$ and $c^{(r)}_{k,\boldsymbol{\eta}}(T)$ are such that $\lim_{k \rightarrow \infty} (c^{(m)}_{k,\boldsymbol{\eta}}(T) + c^{(r)}_{k,\boldsymbol{\eta}}(T)) = 0$. Hence, by Doob's and Jensen's inequalities, we have
\begin{equation} \label{Paper01_bound_martingale_type_space_supremum}
    \lim_{k \rightarrow \infty} \; \sup_{n \in \mathbb{N}} \; \sup_{x \in L^{-1}\mathbb{Z}} \mathbb{E}\left[\sup_{t \leq T} \left\vert  M^{n,x,\boldsymbol{\eta}}_{\mathcal{I}_k}(t) \right\vert \right] = 0. 
\end{equation}

We now bound the integral term on the right-hand side of~\eqref{Paper01_martingale_problem_formulation_number_particles_high_mutations}. By the definition of $\mathcal{L}^{n}$ in~\eqref{Paper01_infinitesimal_generator_restriction_n} and~\eqref{Paper01_generator_foutel_etheridge_model_restriction_n}, and by the triangle inequality, 
\begin{equation} \label{Paper01_bound_FV_type_space_supremum_almost_final_step}
\begin{aligned}
    & \int_{0}^{T} \Big\vert\left(\mathcal{L}^{n}\phi^{(x)}_{\mathcal I_k}\right)(\eta^{n}(t'-))\Big\vert \, dt' \\ & \quad \leq \frac{m}{2} \int_{0}^{T}  \sum_{j = k}^{\infty} \left(\eta^{n}_j(t'-, x+L^{-1})  + \eta^{n}_j(t'-, x - L^{-1}) + 2\eta^{n}_j(t'-, x) \right) \, dt' \\ & \quad \quad \quad + \int_{0}^{T} \sum_{j = k}^{\infty} \left(F^{b}_{j}\left(\eta^{n}(t'-,x)\right) + F^{d}_{j}\left(\eta^{n}(t'-,x)\right)\right) \, dt';
\end{aligned}
\end{equation}
note that the right-hand side of~\eqref{Paper01_bound_FV_type_space_supremum_almost_final_step} is the same as the right-hand side of~\eqref{Paper01_martingale_bound_particles_carrying_high_number_mutations}. Therefore, by taking expectations on both sides of~\eqref{Paper01_bound_FV_type_space_supremum_almost_final_step}, and then by applying~\eqref{Paper01_intermediate_step_migration_term_control_mutation_load_uniform_time} and~\eqref{Paper01_intermediate_step_reaction_term_control_mutation_load_uniform_time}, we get
\begin{equation} \label{Paper01_bound_FV_type_space_supremum}
    \lim_{k \rightarrow \infty} \; \sup_{n \in \mathbb{N}} \; \sup_{x \in L^{-1}\mathbb{Z}} \mathbb{E}\Bigg[\int_{0}^{T} \Big\vert\left(\mathcal{L}^{n}\phi^{(x)}_{\mathcal I_k}\right)(\eta^{n}(t'-))\Big\vert \, dt' \Bigg] = 0. 
\end{equation}
The result~\eqref{Paper01_supremum_time_bound_large_mutation_load} follows by taking expectations on both sides of~\eqref{Paper01_martingale_problem_formulation_number_particles_high_mutations}, and then using~\eqref{Paper01_bound_martingale_type_space_supremum},~\eqref{Paper01_bound_FV_type_space_supremum}, and the definition of $\mathcal{S}_{0}$ in~\eqref{Paper01_set_initial_configurations}.
\end{proof}

We will now proceed to prove tightness of the sequence of processes $(\eta^{n})_{n \in \mathbb{N}}$. We will first establish a strong compact containment condition for the sequence $(\eta^{n})_{n \in \mathbb{N}}$.

\begin{lemma} \label{Paper01_further_characterisation_relatively_compactness_discrete_IPS}
For any $\boldsymbol{\eta} \in \mathcal{S}_{0}$, conditioning on $\eta^{n}(0) = \boldsymbol{\eta}$ for every $n \in \mathbb{N}$, for any $\varepsilon > 0$ and $T \geq 0$, there exists a compact subset $\mathscr{K} = \mathscr{K}(\boldsymbol{\eta}, T, \varepsilon) \subset \mathcal{S}$ such that
    \begin{equation*}
        \inf_{n \in \mathbb{N}} \; \mathbb{P}_{\boldsymbol{\eta}}\Big(\eta^{n}(t) \in \mathscr{K} \; \forall t \in [0,T]\Big) \geq 1 - \varepsilon.
    \end{equation*}
\end{lemma}

\begin{proof}
The result will follow from Lemma~\ref{Paper01_supremum_moments_discrete_particle_system}, the characterisation of compact subsets of $\mathcal{S}$ in Proposition~\ref{Paper01_topological_properties_state_space}, and a standard application of Markov's inequality (see e.g.~\cite[Lemma~4.1]{arrigoni2003deterministic}). For $n \in \mathbb{N}$, $T \geq 0$, $\alpha > 0$, $(\alpha_{i})_{i \in \mathbb{N}} \subset (0, \infty)$ such that $\lim_{i \rightarrow \infty} \, \alpha_{i} = 0$, and sequences $(R_{i})_{i \in \mathbb{N}}  \subset (0, \infty)$ and $(k_{i})_{i \in \mathbb{N}} \subset \mathbb{N}$ such that $\lim_{i \rightarrow \infty} \, R_{i} = \lim_{i \rightarrow \infty} \, k_{i} = \infty$, we define the events
\begin{equation*}
\begin{aligned}
    \mathscr{A}_{\textrm{tot}}^{n}(T, \alpha) & \defeq \Bigg\{\sum_{x \in L^{-1}\mathbb{Z}} \, \frac{\vert \vert \eta^{n}(t,x) \vert \vert_{\ell_{1}}}{(1 + \vert x \vert)^{2}} \leq \alpha \; \forall \, t \in [0,T] \Bigg\}, \\ \mathscr{A}_{\textrm{spat}}^{n}(T,  (\alpha_{i})_{i \in \mathbb{N}}, (R_{i})_{i \in \mathbb{N}}) & \defeq \Bigg\{\sum_{\substack{\{x \in L^{-1}\mathbb{Z}: \, \vert x \vert \geq R_{i}\}}} \, \frac{\vert \vert \eta^{n}(t,x) \vert \vert_{\ell_{1}}}{(1 + \vert x \vert)^{2}} \leq \alpha_{i} \; \forall \, i \in \mathbb{N} \textrm{ and } \forall \, t \in [0,T] \Bigg\}, \\ \mathscr{A}_{\textrm{mut}}^{n}(T,  (\alpha_{i})_{i \in \mathbb{N}}, (k_{i})_{i \in \mathbb{N}}) & \defeq \Bigg\{\sum_{\substack{x \in L^{-1}\mathbb{Z}}} \, \sum_{j = k_{i}}^{\infty} \, \frac{\eta^{n}_{j}(t,x)}{(1 + \vert x \vert)^{2}} \leq \alpha_{i} \; \forall \, i \in \mathbb{N} \textrm{ and } \forall \, t \in [0,T] \Bigg\}.
\end{aligned}
\end{equation*}
By the characterisation of compact subsets of $\mathcal{S}$ given by Proposition~\ref{Paper01_topological_properties_state_space}, there exists a compact subset $\mathscr{K}  = \mathscr{K}(\alpha, (\alpha_{i})_{i \in \mathbb{N}}, (R_{i})_{i \in \mathbb{N}}, (k_{i})_{i \in \mathbb{N}}) \subset \mathcal{S}$, which does not depend on $n$, such that
for each $n\in \mathbb N$, on the event
$$\mathscr{A}_{\textrm{tot}}^{n}(T, \alpha) \cap  \mathscr{A}_{\textrm{spat}}^{n}(T,  (\alpha_{i})_{i \in \mathbb{N}}, (R_{i})_{i \in \mathbb{N}}) \cap \mathscr{A}_{\textrm{mut}}^{n}(T, (\alpha_{i})_{i \in \mathbb{N}}, (k_{i})_{i \in \mathbb{N}}),$$
we have $\eta^{n}(t) \in \mathscr{K}$ for all $t \in [0,T]$. By estimates~\eqref{Paper01_uniform_estimate_discrete_process_any_moments_finite_time},~\eqref{Paper01_supremum_control_spatial_dispersion_mass} and~\eqref{Paper01_supremum_time_bound_large_mutation_load} in Lemma~\ref{Paper01_supremum_moments_discrete_particle_system}, for any $\boldsymbol{\eta} \in \mathcal{S}_{0}$, $T \geq 0$ and $\varepsilon > 0$, it is possible to take $\alpha$, $(\alpha_{i})_{i \in \mathbb{N}}$, $(R_{i})_{i \in \mathbb{N}}$ and $(k_{i})_{i \in \mathbb{N}}$ such that
$\lim_{i \rightarrow \infty} \, \alpha_{i} = 0$, $\lim_{i \rightarrow \infty} \, R_{i} = \lim_{i \rightarrow \infty} \, k_{i} = \infty$ and
\begin{equation*}
    \inf_{n \in \mathbb{N}} \; \mathbb{P}_{\boldsymbol{\eta}}\Big(\mathscr{A}_{\textrm{tot}}^{n}(T, \alpha) \cap  \mathscr{A}_{\textrm{spat}}^{n}(T, (\alpha_{i})_{i \in \mathbb{N}}, (R_{i})_{i \in \mathbb{N}}) \cap \mathscr{A}_{\textrm{mut}}^{n}(T, (\alpha_{i})_{i \in \mathbb{N}}, (k_{i})_{i \in \mathbb{N}})\Big) \geq 1 - \varepsilon,
\end{equation*}
which completes the proof.
\end{proof}
In order to be able to characterise the limiting process in terms of the solution of a martingale problem later on in the proof, in our proof of tightness we will establish both tightness of evaluations and uniform integrability of the quadratic variation of a class of martingales. Recall the definitions of $\mathscr{C}_{b,*}^{\textrm{cyl}}(\mathcal{S}, \mathbb{R})$ and $\mathscr{C}_*(\mathcal{S}, \mathbb{R})$ in~\eqref{Paper01_general_definition_cylindrical_function} and~\eqref{Paper01_definition_Lipschitz_function}, respectively, and that by Lemma~\ref{Paper01_characterisation_lipschitz_functions_semi_metric}(v), we have $\mathscr{C}_{b,*}^{\textrm{cyl}}(\mathcal{S}, \mathbb{R}) \subseteq \mathscr{C}_*(\mathcal{S}, \mathbb{R})$. Also, for $\phi \in \mathscr{C}_*(\mathcal{S}, \mathbb{R})$, $n \in \mathbb{N}$ and $\boldsymbol{\eta} \in \mathcal{S}_{0}$, recall the definition of $(M^{n,\phi,\boldsymbol{\eta}}(t))_{t \geq 0}$ from~\eqref{Paper01_local_martingale_problem_discrete_approximation}, and recall that by Lemma~\ref{Paper01_auxilliary_lemma_easy_martingale_representation}, conditioning on $\eta^n(0)=\boldsymbol{\eta}$,
$(M^{n,\phi,\boldsymbol{\eta}}(t))_{t \geq 0}$ is a square-integrable càdlàg martingale.

\begin{proposition} \label{Paper01_tightness_evaluations_discrete_approximating_IPS}
    For $T > 0$, $\phi \in \mathscr{C}_*(\mathcal{S}, \mathbb{R})$, $\boldsymbol{\eta} \in \mathcal{S}_{0}$, conditioning on $\eta^{n}(0) = \boldsymbol{\eta}$ for every $n \in \mathbb{N}$, the sequences $\left(\left(\phi(\eta^{n}(t)\right)_{t \in [0,T]}\right)_{n \in \mathbb{N}}$ and $\left((M^{n,\phi,\boldsymbol{\eta}}(t))_{t \in [0,T]}\right)_{n \in \mathbb{N}}$ are both tight in $\mathscr{D}([0,T], \mathbb{R})$, and the sequence of real-valued random variables $\left([ M^{n,\phi,\boldsymbol{\eta}}](T)\right)_{n \in \mathbb{N}}$ is uniformly integrable.
\end{proposition}

\begin{proof}
We will organise the proof into steps to improve the readability.

\medskip

\noindent \underline{Step $(1)$: Tightness of the sequence of martingales $(M^{n,\phi,\boldsymbol{\eta}})_{n \in \mathbb{N}}$}

\medskip

Take $\boldsymbol{\eta}=(\eta_k(x))_{x\in L^{-1}\mathbb Z, \, k\in \mathbb N_0} \in \mathcal{S}_{0}$, $\phi \in  \mathscr{C}_*(\mathcal{S}, \mathbb{R})$ and $T > 0$. To prove tightness of the sequence $\left((M^{n,\phi,\boldsymbol{\eta}}(t))_{t \in [0,T]}\right)_{n \in \mathbb{N}}$, we must establish the following conditions (see for instance~\cite[Theorem~3.8.6(c)]{ethier2009markov}):
\begin{enumerate}[(i)]
    \item \underline{Compact containment condition:} For $\varepsilon > 0$, there exists $\mathscr{K} = \mathscr{K}(\varepsilon, \boldsymbol{\eta}, \phi, T) \subset \mathbb{R}$, $\mathscr{K}$~compact, such that
    \begin{equation*}
        \inf_{n \in \mathbb{N}} \mathbb{P}_{\boldsymbol{\eta}}\Big(M^{n,\phi,\boldsymbol{\eta}}(t) \in \mathscr{K} \; \forall \, t \in [0,T]\Big) \geq 1 - \varepsilon.
    \end{equation*}
    \item \underline{Aldous' criterion~\cite{aldous1978stopping}:} Let $\mathscr{T}(T)$ denote the family of $\{\mathcal{F}^{\eta^n}_{t}\}_{t \geq 0}$-stopping times bounded by $T$. Then, for any $\varepsilon > 0$, 
    \begin{equation} \label{Paper01_aldous_criterion_martingale_term_discrete_process}
    \lim_{\gamma \rightarrow 0^{+}} \; \limsup_{n \rightarrow \infty} \; \sup_{\substack{\tau \in \mathscr{T}(T), \\ 0 < t \leq \gamma}} \mathbb{P}_{\boldsymbol{\eta}}\left(\left\vert M^{n,\phi,\boldsymbol{\eta}}((\tau + t) \wedge T) - M^{n,\phi,\boldsymbol{\eta}}(\tau) \right\vert > \varepsilon\right) = 0.
    \end{equation}
\end{enumerate}

We will first establish condition~(i), i.e.~the compact containment condition. By Markov's inequality, it will suffice to establish that
    \begin{equation} \label{Paper01_trivial_martingale_lipschitz_bound}
        \sup_{n \in \mathbb{N}} \, \mathbb{E}_{\boldsymbol{\eta}}\left[\sup_{t \leq T} \left(M^{n,\phi,\boldsymbol{\eta}}(t)\right)^{2}\right] < \infty.
    \end{equation}
    By Doob's maximal inequality, for every $n \in \mathbb{N}$, we have
    \begin{equation*}
    \mathbb{E}_{\boldsymbol{\eta}}\left[\sup_{t \leq T} \left(M^{n,\phi,\boldsymbol{\eta}}(t)\right)^{2}\right] \leq 4 \mathbb{E}_{\boldsymbol{\eta}}\left[\left\langle M^{n,\phi,\boldsymbol{\eta}}\right\rangle(T)\right],
    \end{equation*}
    so that~\eqref{Paper01_trivial_martingale_lipschitz_bound} will be established after proving that
    \begin{equation} \label{Paper01_trivial_martingale_lipschitz_bound_ii}
        \sup_{n \in \mathbb{N}} \, \mathbb{E}_{\boldsymbol{\eta}}\left[\left\langle M^{n,\phi,\boldsymbol{\eta}}\right\rangle(T)\right] < \infty.
    \end{equation}
Since $\phi \in \mathscr{C}_*(\mathcal{S}, \mathbb{R})$, by the definition of $\mathscr{C}_*(\mathcal{S}, \mathbb{R})$ in~\eqref{Paper01_definition_Lipschitz_function} and by~\eqref{Paper01_predictable_bracket_process_martingale_lipschitz_functions_discrete_IPS} in Lemma~\ref{Paper01_auxilliary_lemma_easy_martingale_representation}, we conclude that there exists $C_{\phi} > 0$ such that
    \begin{equation} \label{Paper01_bound_first_step_uniform_integrability_martingale_lipschitz}
    \begin{aligned}
       \left\langle M^{n,\phi,\boldsymbol{\eta}}\right\rangle(T) &\leq C_\phi\sum_{\substack{x \in \Lambda_n}} \, \sum_{k = 0}^{\infty} \int_{0}^{T} \,  \frac{(m \vee 1) \left(\eta^{n}_{k}(t-,x)+ F^{b}_{k}(\eta^{n}(t-,x)) + F^{d}_{k}(\eta^{n}(t-,x))\right)}{(1 + \vert x \vert)^{4(1 + \deg q_-)}} \, dt \\ &\leq C_\phi \sum_{\substack{x \in \Lambda_n}} \, \int_{0}^{T}  \frac{(m \vee 1) \left(\vert \vert \eta^{n}(t-,x) \vert \vert_{\ell_1} + \bar{q}^{N}(\vert \vert \eta^{n}(t-,x) \vert \vert_{\ell_1})\right)}{(1 + \vert x \vert)^{4(1 + \deg q_-)}} \, dt,
    \end{aligned}
    \end{equation}
    where $\bar{q}^{N}$ is defined in~\eqref{Paper01_bar_polynomial}, and we used~\eqref{Paper01_simple_bound_total_birth_death_rates_deterministic} in the last inequality. Estimate~\eqref{Paper01_trivial_martingale_lipschitz_bound_ii} then follows by taking expectations on both sides of~\eqref{Paper01_bound_first_step_uniform_integrability_martingale_lipschitz}, and by applying estimate~\eqref{Paper01_moments_estimates_uniform_space} in Proposition~\ref{Paper01_bound_total_mass} and the fact that $\bar{q}^{N}$ is a polynomial.
    
    It remains to establish condition~(ii), i.e.~Aldous' criterion. By the optional sampling theorem (see for instance~\cite[Theorem~2.2.13]{ethier2009markov}), for any $\tau \in \mathscr{T}(T)$, $(M^{n,\phi,\boldsymbol{\eta}}(\tau+t))_{t \geq 0}$ is a càdlàg martingale with respect to the filtration $\{\mathcal{F}^{\eta^n}_{\tau + t}\}_{t \geq 0}$. We will prove~\eqref{Paper01_aldous_criterion_martingale_term_discrete_process} by establishing bounds on the expectation of the process $$\left(\left\langle M^{n,\phi,\boldsymbol{\eta}} \right\rangle(\tau + t) - \left\langle M^{n,\phi,\boldsymbol{\eta}} \right\rangle(\tau) \right)_{t \geq 0}.$$
    Since $\phi \in \mathscr{C}_*\left(\mathcal{S}, \mathbb{R}\right)$, by the definition of $\mathscr{C}_*\left(\mathcal{S}, \mathbb{R}\right)$ in~\eqref{Paper01_definition_Lipschitz_function} and by~\eqref{Paper01_predictable_bracket_process_martingale_lipschitz_functions_discrete_IPS}, and then by using~\eqref{Paper01_simple_bound_total_birth_death_rates_deterministic} as in the derivation of~\eqref{Paper01_bound_first_step_uniform_integrability_martingale_lipschitz}, there exists $C_\phi > 0$ such that for $\tau \in \mathscr{T}(T)$, $\gamma > 0$ and $t \in (0, \gamma]$,
    \begin{equation} \label{Paper01_intermediate_step_bound_martingale_part}
    \begin{aligned}
        & \mathbb{E}_{\boldsymbol{\eta}}\left[\left\langle M^{n,\phi,\boldsymbol{\eta}} \right\rangle((\tau + t) \wedge T) - \left\langle M^{n,\phi,\boldsymbol{\eta}} \right\rangle(\tau) \right] \\ 
        & \quad  \leq \mathbb{E}_{\boldsymbol{\eta}}\Bigg[ \sum_{\substack{x \in \Lambda_n}} \, \int_{\tau}^{(\tau + t) \wedge T}  \frac{ C_\phi(m \vee 1) \left(\vert \vert \eta^{n}(t'-,x) \vert \vert_{\ell_1} + \bar{q}^{N}(\vert \vert \eta^{n}(t'-,x) \vert \vert_{\ell_1})\right)}{(1 + \vert x \vert)^{4(1 + \deg q_-)}} \, dt'\Bigg] \\ 
        & \quad \leq \sum_{\substack{x \in L^{-1}\mathbb{Z}}} \frac{C_\phi(m \vee 1) \mathbb{E}_{\boldsymbol{\eta}}\left[\displaystyle\int_{0}^{T} \left(\vert \vert \eta^{n}(t'-,x) \vert \vert_{\ell_1} + \bar{q}^{N}(\vert \vert \eta^{n}(t'-,x) \vert \vert_{\ell_1})\right) \cdot \mathds{1}_{\{\tau \leq t' \leq \tau + t\}} dt'\right]}{(1+\vert x \vert)^{4(1 + \deg q_-)}}
        \\ & \quad \leq \sum_{\substack{x \in L^{-1}\mathbb{Z}}} \frac{C_\phi(m \vee 1)}{(1 + \vert x \vert)^{4}} \\ & \quad \quad \quad \quad \; \cdot \mathbb{E}_{\boldsymbol{\eta}}\Bigg[\displaystyle \Bigg(\int_{0}^{T} \left(\vert \vert \eta^{n}(t'-,x) \vert \vert_{\ell_1} + \bar{q}^{N}(\vert \vert \eta^{n}(t'-,x) \vert \vert_{\ell_1})\right)^{2} dt' \Bigg)^{1/2} \Bigg(\int_{0}^{T} \mathds{1}_{\{\tau \leq t' \leq \tau + t\}} \, dt'\Bigg)^{1/2}\Bigg] \\ & \quad \leq \sum_{\substack{x \in L^{-1}\mathbb{Z}}} \frac{t^{1/2}C_\phi(m \vee 1)}{(1 + \vert x \vert)^{4}} \Bigg(\int_{0}^{T} \mathbb{E}_{\boldsymbol{\eta}}\Bigg[ \left(\vert \vert \eta^{n}(t'-,x) \vert \vert_{\ell_1} + \bar{q}^{N}(\vert \vert \eta^{n}(t'-,x) \vert \vert_{\ell_1})\right)^{2} \Big] dt' \Bigg)^{1/2}
         \\ & \quad \lesssim_{\boldsymbol{\eta},m,L,T,\phi,N,q_-,q_+} \gamma^{1/2},
    \end{aligned}
    \end{equation}
    where in the third inequality we applied the Cauchy-Schwarz inequality, in the fourth inequality we used Jensen's inequality, and in the last inequality we used estimate~\eqref{Paper01_moments_estimates_uniform_space} in Proposition~\ref{Paper01_bound_total_mass} and the fact that $\bar{q}^{N}$ defined in~\eqref{Paper01_bar_polynomial} is a polynomial, together with the fact that $t \in (0, \gamma]$. We can then obtain~\eqref{Paper01_aldous_criterion_martingale_term_discrete_process} by combining~\eqref{Paper01_intermediate_step_bound_martingale_part} with Markov’s inequality, applied to the squared increment of  $M^{n,\phi,\boldsymbol{\eta}}$ on the interval from $\tau$ to $(\tau + t) \wedge T$, completing Step~(1) of the proof.

    \medskip

    \noindent \underline{Step $(2)$: Uniform integrability of  $\left([M^{n,\phi,\boldsymbol{\eta}}](T)\right)_{n \in \mathbb{N}}$}

    \medskip

    Take $\boldsymbol{\eta}\in \mathcal S_0$, $\phi \in \mathscr{C}_*(\mathcal{S}, \mathbb{R})$ and $T>0$.
    To prove that the sequence $\left([M^{n,\phi,\boldsymbol{\eta}}](T)\right)_{n \in \mathbb{N}}$ is uniformly integrable, 
    it will suffice to establish that
    \begin{equation} \label{Paper01_uniform_integrability_quadratic_variation_claim}
        \sup_{n\in \mathbb N} \, \mathbb{E}_{\boldsymbol{\eta}}\left[[ M^{n,\phi,\boldsymbol{\eta}} ]^{2}(T)\right] < \infty.
    \end{equation}
    Recall from before~\eqref{Paper01_Poisson_process_representation_discrete_approximation} that the collections $\{\mathcal{Q}^{x,x + z}_{k}: \, x \in L^{-1}\mathbb{Z}, \, z \in \{ -L^{-1}, L^{-1}\}, \, k \in \mathbb{N}_{0}\}$, $\left\{\mathcal{R}^{x}_{k}: \, x \in L^{-1}\mathbb{Z}, \, k \in \mathbb{N}_{0}\right\}$ and $\left\{\mathcal{D}^{x}_{k}: \, x \in L^{-1}\mathbb{Z}, \, k \in \mathbb{N}_{0}\right\}$ are families of i.i.d.~Poisson random measures on $[0, \infty) \times [0, \infty)$, each with intensity measure given by the Lebesgue measure $\lambda$. 
    Since $\phi \in \mathscr{C}_*(\mathcal{S}, \mathbb{R})$, by~\eqref{Paper01_definition_quadratic_variation_martingale_problem} in Lemma~\ref{Paper01_auxilliary_lemma_easy_martingale_representation} and by~\eqref{Paper01_definition_Lipschitz_function}, there exists $C_\phi > 0$ such that for each $n\in \mathbb N$,
    \begin{equation} \label{Paper01_first_immediate_bound_quadratic_variation}
    \begin{aligned}
     & \left[M^{n,\phi,\boldsymbol{\eta}}\right](T) \\ 
     & \quad \leq C_\phi\sum_{x \in \Lambda_n} \, \sum_{k = 0}^{\infty} \, \sum_{z \in \{-L^{-1}, L^{-1}\}} \int_{[0,T]\times [0,\infty)} \frac{1}{(1 + \vert x \vert)^{4(1 + \deg q_-)}}\cdot \mathds{1}_{\left\{y \leq \frac{m}{2}\eta^{n}_{k}(t-,x)\right\}} \mathcal{Q}^{x,x+z}_{k}(dt \times dy) \\ 
     & \; \quad \quad   + C_\phi \sum_{x \in \Lambda_n} \, \sum_{k = 0}^{K_n} \, \int_{[0,T]\times [0,\infty)} \frac{1}{(1 + \vert x \vert)^{4(1 + \deg q_-)}} \cdot \mathds{1}_{\left\{y \leq F^{b}_{k}(\eta^{n}(t-,x))\right\}} \mathcal{R}^{x}_{k}(dt \times dy) \\ & \; \quad \quad + C_\phi\sum_{x \in \Lambda_n} \, \sum_{k = 0}^{\infty} \, \int_{[0,T]\times [0,\infty)} \frac{1}{(1 + \vert x \vert)^{4(1 + \deg q_-)}} \cdot \mathds{1}_{\left\{y \leq F^{d}_{k}(\eta^{n}(t-,x))\right\}} \mathcal{D}^{x}_{k}(dt \times dy).
    \end{aligned}
    \end{equation}
    By the independence of the families $\{\mathcal{Q}^{x,x + z}_{k}\}_{x,z,k}$, $\left\{\mathcal{R}^{x}_{k}\right\}_{x,k}$ and $\left\{\mathcal{D}^{x}_{k}\right\}_{x,k}$, and the thinning property of Poisson random measures, together with~\eqref{Paper01_simple_bound_total_birth_death_rates_deterministic} and~\eqref{Paper01_first_immediate_bound_quadratic_variation}, we can construct a countable family $\{\mathcal{P}^{x}: \, x \in L^{-1}\mathbb{Z}\}$ of i.i.d.~Poisson random measures on $[0, \infty) \times [0, \infty)$ with intensity measure given by the Lebesgue measure $\lambda$, on the same probability space as $(\eta^{n}(t))_{t \geq 0}$ such that the following inequality holds almost surely:
    \begin{align} \label{Paper01_useful_bound_quadratic_variation_uniform_integrability}
        & \left[M^{n,\phi,\boldsymbol{\eta}}\right](T) \notag \\ 
        & \quad \leq C_\phi\sum_{x \in \Lambda_n} \frac{1}{(1 + \vert x \vert)^{4(1 + \deg q_-)}} \int_{[0,T]\times [0,\infty)} \mathds{1}_{\left\{y \leq m\vert \vert \eta^{n}(t-,x) \vert \vert_{\ell_{1}} + \bar{q}^{N}(\vert \vert \eta^{n}(t-,x) \vert \vert_{\ell_{1}})\right\}} \mathcal{P}^{x}(dt \times dy),
    \end{align}
    where $\bar{q}^{N}$ is the polynomial defined in~\eqref{Paper01_bar_polynomial}.
    By taking the square of both sides of~\eqref{Paper01_useful_bound_quadratic_variation_uniform_integrability}, we obtain
    \begin{align} \label{Paper01_useful_bound_quadratic_variation_uniform_integrability_ii}
    \left[M^{n,\phi,\boldsymbol{\eta}}\right]^{2}(T) & \leq \sum_{\substack{x,z \in \Lambda_n}} \, \frac{C_\phi^2}{(1 + \vert x \vert)^{4(1 + \deg q_-)}(1 + \vert z \vert)^{4(1 + \deg q_-)}} \notag \\ 
    & \quad \quad \quad \quad \quad \quad \quad \quad \cdot \Bigg( \int_{[0,T]\times [0,\infty)} \mathds{1}_{\left\{y \leq m\vert \vert \eta^{n}(t-,x) \vert \vert_{\ell_{1}} + \bar{q}^{N}(\vert \vert \eta^{n}(t-,x) \vert \vert_{\ell_{1}})\right\}} \mathcal{P}^{x}(dt \times dy)\Bigg) \notag  \\ 
    & \quad \quad \quad \quad \quad \quad \quad \quad \cdot \Bigg( \int_{[0,T]\times [0,\infty)} \mathds{1}_{\left\{y \leq m\vert \vert \eta^{n}(t-,z) \vert \vert_{\ell_{1}} + \bar{q}^{N}(\vert \vert \eta^{n}(t-,z) \vert \vert_{\ell_{1}})\right\}} \mathcal{P}^{z}(dt \times dy)\Bigg). 
    \end{align}
    Since
    \begin{equation*}
        \sum_{x,z \in L^{-1}\mathbb{Z}} \, \frac{1}{(1 + \vert x \vert)^{4(1 + \deg q_-)}(1 + \vert z \vert)^{4(1 + \deg q_-)}} < \infty,
    \end{equation*}
    estimate~\eqref{Paper01_uniform_integrability_quadratic_variation_claim} will follow from~\eqref{Paper01_useful_bound_quadratic_variation_uniform_integrability_ii} after establishing that
    \begin{equation} \label{Paper01_target_estimate_quadratic_variation}
    \begin{aligned}
        \sup_{n \in \mathbb{N}} \, \sup_{x,z \in L^{-1}\mathbb{Z}} \mathbb{E}_{\boldsymbol{\eta}}\Bigg[&\Bigg( \int_{[0,T]\times [0,\infty)} \mathds{1}_{\left\{y \leq m\vert \vert \eta^{n}(t-,x) \vert \vert_{\ell_{1}} + \bar{q}^{N}(\vert \vert \eta^{n}(t-,x) \vert \vert_{\ell_{1}})\right\}} \mathcal{P}^{x}(dt \times dy)\Bigg) \\ & \quad \cdot \Bigg( \int_{[0,T]\times [0,\infty)} \mathds{1}_{\left\{y \leq m\vert \vert \eta^{n}(t-,z) \vert \vert_{\ell_{1}} + \bar{q}^{N}(\vert \vert \eta^{n}(t-,z) \vert \vert_{\ell_{1}})\right\}} \mathcal{P}^{z}(dt \times dy)\Bigg)\Bigg] < \infty.
    \end{aligned}
    \end{equation}
    To prove~\eqref{Paper01_target_estimate_quadratic_variation}, we first observe that by the Cauchy-Schwarz inequality,
    \begin{equation} \label{Paper01_target_estimate_quadratic_variation_ii}
    \begin{aligned}
       & \sup_{n \in \mathbb{N}} \, \sup_{x,z \in L^{-1}\mathbb{Z}} \mathbb{E}_{\boldsymbol{\eta}}\Bigg[\Bigg( \int_{[0,T]\times [0,\infty)} \mathds{1}_{\left\{y \leq m\vert \vert \eta^{n}(t-,x) \vert \vert_{\ell_{1}} + \bar{q}^{N}(\vert \vert \eta^{n}(t-,x) \vert \vert_{\ell_{1}})\right\}} \mathcal{P}^{x}(dt \times dy)\Bigg) \\ 
       & \quad \quad \quad \quad \quad \quad \quad \quad \cdot \Bigg( \int_{[0,T]\times [0,\infty)} \mathds{1}_{\left\{y \leq m\vert \vert \eta^{n}(t-,z) \vert \vert_{\ell_{1}} + \bar{q}^{N}(\vert \vert \eta^{n}(t-,z) \vert \vert_{\ell_{1}})\right\}} \mathcal{P}^{z}(dt \times dy)\Bigg)\Bigg] \\ & \quad \quad  \leq \sup_{n \in \mathbb{N}} \, \sup_{x \in L^{-1}\mathbb{Z}} \mathbb{E}_{\boldsymbol{\eta}}\Bigg[\Bigg( \int_{[0,T]\times [0,\infty)} \mathds{1}_{\left\{y \leq m\vert \vert \eta^{n}(t-,x) \vert \vert_{\ell_{1}} + \bar{q}^{N}(\vert \vert \eta^{n}(t-,x) \vert \vert_{\ell_{1}})\right\}} \mathcal{P}^{x}(dt \times dy)\Bigg)^{2}\Bigg].
    \end{aligned}
    \end{equation}
    To bound the right-hand side of~\eqref{Paper01_target_estimate_quadratic_variation_ii}, we apply Itô's isometry (see Lemma~\ref{Paper01_easy_inequality_square_integral_poisson_point_process}) to obtain
    \begin{align*}
        & \sup_{n \in \mathbb{N}} \, \sup_{x \in L^{-1}\mathbb{Z}} \mathbb{E}_{\boldsymbol{\eta}}\Bigg[\Bigg( \int_{[0,T]\times [0,\infty)} \mathds{1}_{\left\{y \leq m\vert \vert \eta^{n}(t-,x) \vert \vert_{\ell_{1}} + \bar{q}^{N}(\vert \vert \eta^{n}(t-,x) \vert \vert_{\ell_{1}})\right\}} \mathcal{P}^{x}(dt \times dy)\Bigg)^{2}\Bigg] \\ 
        & \quad \leq  2 \sup_{n \in \mathbb{N}} \, \sup_{x \in L^{-1}\mathbb{Z}} \mathbb{E}_{\boldsymbol{\eta}}\Bigg[\int_{0}^{T} \left( m\vert \vert \eta^{n}(t-,x) \vert \vert_{\ell_{1}} + \bar{q}^{N}(\vert \vert \eta^{n}(t-,x) \vert \vert_{\ell_{1}})\right) dt\Bigg] \\ 
        & \quad \quad + 2 \sup_{n \in \mathbb{N}} \, \sup_{x \in L^{-1}\mathbb{Z}} \mathbb{E}_{\boldsymbol{\eta}}\Bigg[\left(\int_{0}^{T} \left( m\vert \vert \eta^{n}(t-,x) \vert \vert_{\ell_{1}} + \bar{q}^{N}(\vert\vert \eta^{n}(t-,x) \vert \vert_{\ell_{1}})\right) dt\right)^{2}\Bigg] \\ 
        & \quad \leq 2T \sup_{n \in \mathbb{N}} \, \sup_{x \in L^{-1}\mathbb{Z}} \sup_{t \leq T} \mathbb{E}_{\boldsymbol{\eta}}\Big[m\vert \vert \eta^{n}(t-,x) \vert \vert_{\ell_{1}} + \bar{q}^{N}(\vert\vert \eta^{n}(t-,x) \vert \vert_{\ell_{1}})\Big] \\ 
        & \quad \quad + 2T \sup_{n \in \mathbb{N}} \, \sup_{x \in L^{-1}\mathbb{Z}} \sup_{t \leq T} \mathbb{E}_{\boldsymbol{\eta}}\Big[\Big(m\vert \vert \eta^{n}(t-,x) \vert \vert_{\ell_{1}} + \bar{q}^{N}(\vert\vert \eta^{n}(t-,x) \vert \vert_{\ell_{1}})\Big)^{2}\Big]\\
        & \quad <\infty,
    \end{align*}
    where in the second inequality we used the Cauchy-Schwarz inequality,
    and for the last line we used estimate~\eqref{Paper01_moments_estimates_uniform_space} in Proposition~\ref{Paper01_bound_total_mass} and the fact that $\bar{q}^{N}$ is a polynomial by~\eqref{Paper01_bar_polynomial}.
By~\eqref{Paper01_target_estimate_quadratic_variation_ii}, this establishes~\eqref{Paper01_target_estimate_quadratic_variation} and completes the proof of Step~(2). 
    
    \medskip

    \noindent \underline{Step~$(3)$: Tightness of the sequence of càdlàg processes $\left(\phi(\eta^{n})\right)_{n \in \mathbb{N}}$}

    \medskip

   Let $\boldsymbol{\eta} \in \mathcal{S}_{0}$, $\phi \in  \mathscr{C}_*(\mathcal{S}, \mathbb{R})$ and $T > 0$. As in Step~(1) of the proof, to prove tightness of the sequence $\Big((\phi(\eta^{n}(t)))_{t \in [0,T]}\Big)_{n \in \mathbb{N}}$ in $\mathscr{D}([0,T], \mathbb{R})$, conditioned on $\eta^{n}(0) = \boldsymbol{\eta}$ for every $n \in \mathbb{N}$, we must establish the following conditions (see for instance~\cite[Theorem~3.8.6(c)]{ethier2009markov}):
\begin{enumerate}[(i)]
    \item \underline{Compact containment condition:} For $\varepsilon > 0$, there exists $\mathscr{K} = \mathscr{K}(\varepsilon, \boldsymbol{\eta}, \phi, T) \subset \mathbb{R}$, $\mathscr{K}$~compact, such that
    \begin{equation*}
        \inf_{n \in \mathbb{N}} \mathbb{P}_{\boldsymbol{\eta}}\Big(\phi(\eta^{n}(t)) \in \mathscr{K} \; \forall \, t \in [0,T]\Big) \geq 1 - \varepsilon.
    \end{equation*}
    \item \underline{Aldous' criterion:} Let $\mathscr{T}(T)$ denote the family of $\{\mathcal{F}^{\eta^n}_{t}\}_{t \geq 0}$-stopping times bounded by $T$. Then, for any $\varepsilon > 0$, 
    \begin{equation} \label{Paper01_aldous_criterion_lipschitz_functions_discrete_process}
    \lim_{\gamma \rightarrow 0^{+}} \; \limsup_{n \rightarrow \infty} \; \sup_{\substack{\tau \in \mathscr{T}(T), \\ 0 < t \leq \gamma}} \mathbb{P}_{\boldsymbol{\eta}}\Big(\left\vert \phi(\eta^{n}((\tau + t) \wedge T)) - \phi(\eta^{n}(\tau)) \right\vert > \varepsilon\Big) = 0.
    \end{equation}
\end{enumerate}
   
   Since the continuous image of a compact set is compact, condition~(i) follows from Lemma~\ref{Paper01_further_characterisation_relatively_compactness_discrete_IPS} and the fact that $\phi$ is continuous. It remains to establish condition~(ii), i.e.~\eqref{Paper01_aldous_criterion_lipschitz_functions_discrete_process}. By the definition of $M^{n,\phi,\boldsymbol{\eta}}$ in~\eqref{Paper01_local_martingale_problem_discrete_approximation}, for $n \in \mathbb{N}$, $\tau \in \mathscr{T}(T)$, $\gamma > 0$ and $t \in (0, \gamma]$, we have
    \begin{equation} \label{Paper01_intermediate_step_tightness_evaluations}
    \begin{aligned}
        & \left\vert \phi(\eta^{n}((\tau + t) \wedge T)) - \phi(\eta^{n}(\tau)) \right\vert \\ & \quad \leq \left\vert M^{n,\phi,\boldsymbol{\eta}}((\tau + t) \wedge T) - M^{n,\phi,\boldsymbol{\eta}}(\tau) \right\vert + \int_{\tau}^{(\tau + t) \wedge T}\left\vert (\mathcal{L}^{n}  \phi)(\eta^{n}(t'-))\right\vert \, dt'.
    \end{aligned}
    \end{equation}
   Hence, by~\eqref{Paper01_intermediate_step_tightness_evaluations},~\eqref{Paper01_aldous_criterion_martingale_term_discrete_process} and Markov's inequality,~\eqref{Paper01_aldous_criterion_lipschitz_functions_discrete_process} will be proved after establishing that
    \begin{equation} \label{Paper01_target_estimate_aldous_test_functions}
    \begin{aligned}
        \lim_{\gamma \rightarrow 0^{+}} \, \limsup_{n \rightarrow \infty} \, \sup_{\substack{\tau \in \mathscr{T}(T), \\ 0 < t \leq \gamma}} \mathbb{E}_{\boldsymbol{\eta}}\Bigg[\int_{\tau}^{(\tau+t) \wedge T} \left\vert (\mathcal{L}^{n}  \phi)(\eta^{n}(t'-))\right\vert \, dt' \Bigg] = 0.
    \end{aligned}
    \end{equation}
    Recall the definition of $\mathcal{L}^{n}$ from~\eqref{Paper01_generator_foutel_etheridge_model_restriction_n} and~\eqref{Paper01_infinitesimal_generator_restriction_n}, and the definition of $(F^{b}_{k})_{k \in \mathbb{N}_{0}}$ and~$(F^{d}_{k})_{k \in \mathbb{N}_{0}}$ from~\eqref{Paper01_discrete_birth_rate} and~\eqref{Paper01_discrete_death_rate} respectively. Then, by the triangle inequality and the definition of $\mathscr{C}_*(\mathcal{S}, \mathbb{R})$ in~\eqref{Paper01_definition_Lipschitz_function}, we conclude that for $n \in \mathbb{N}$, $\tau \in \mathscr{T}(T)$, $\gamma > 0$ and $t \in (0, \gamma]$,
    \begin{equation*}
    \begin{aligned}
        & \mathbb{E}_{\boldsymbol{\eta}}\Bigg[\int_{\tau}^{(\tau+ t) \wedge T} \left\vert (\mathcal{L}^{n}  \phi)(\eta^{n}(t'-))\right\vert \, dt' \Bigg] \\ 
        & \quad \lesssim_{\phi} \sum_{x \in L^{-1}\mathbb{Z}} \frac{1}{(1 + \vert x \vert)^{2(1 + \deg q_-)}}\mathbb{E}_{\boldsymbol{\eta}}\Bigg[\int_{0}^{T} \sum_{k = 0}^{\infty} \left( m \eta^{n}_{k}(t'-,x) + F^{b}_{k}(\eta^{n}(t'-,x)) + F^{d}_{k}(\eta^{n}(t'-,x)) \right)\\[-3mm]
        & \quad \quad \quad \quad \quad \quad \quad \quad \quad \quad \quad \quad \quad \quad \quad \quad \quad \quad \quad \quad \quad \quad \quad \quad \quad \quad \quad \quad \quad \quad \quad \cdot \mathds{1}_{\{\tau \leq t' \leq (\tau + t) \wedge T\}}\, dt'\Bigg] \\ 
        & \quad \leq \sum_{x \in L^{-1}\mathbb{Z}}\frac{1}{(1 + \vert x \vert)^{2(1 + \deg q_-)}}\, \mathbb{E}_{\boldsymbol{\eta}}\Bigg[\int_{0}^{T} \Big(m\vert \vert \eta^{n}(t'-,x) \vert \vert_{\ell_{1}} + \bar{q}^{N}(\vert \vert \eta^{n}(t'-,x) \vert \vert_{\ell_{1}})\Big) \\[-3mm] 
        & \quad \quad \quad \quad \quad \quad \quad \quad \quad \quad \quad \quad \quad \quad \quad \quad \quad \quad \quad \quad \quad \quad \quad \quad \quad \quad \quad \quad \quad \quad \quad \quad \quad \cdot \mathds{1}_{\{\tau \leq t' \leq \tau + t\}} \, dt' \Bigg] \\ 
        & \quad \leq \sum_{x \in L^{-1}\mathbb{Z}}\frac{1}{(1 + \vert x \vert)^{2(1 + \deg q_-)}}\, \mathbb{E}_{\boldsymbol{\eta}}\Bigg[\Bigg(\int_{0}^{T} \Big(m\vert \vert \eta^{n}(t'-,x) \vert \vert_{\ell_{1}} + \bar{q}^{N}(\vert \vert \eta^{n}(t'-,x) \vert \vert_{\ell_{1}})\Big)^{2} \, dt'\Bigg)^{1/2} \\ & \quad \quad \quad \quad \quad \quad \quad \quad \quad \quad \quad \quad \quad \quad \quad \quad \quad \quad \quad \quad \quad \quad \quad \quad \quad \quad \quad \quad \cdot \Bigg(\int_{0}^{T} \mathds{1}_{\{\tau \leq t' \leq \tau + t\}} \, dt'\Bigg)^{1/2}\Bigg] \\ & \quad \lesssim_{L,q_+,q_-,N,m,\boldsymbol{\eta}, T} \gamma^{1/2},
    \end{aligned}
    \end{equation*}
    where in the second inequality we used~\eqref{Paper01_simple_bound_total_birth_death_rates_deterministic}, in the third inequality we used the Cauchy-Schwarz inequality, and in the last inequality we used Jensen's inequality, the fact that $\bar{q}^{N}$ is a polynomial (defined in~\eqref{Paper01_bar_polynomial}) and estimate~\eqref{Paper01_moments_estimates_uniform_space} in Proposition~\ref{Paper01_bound_total_mass} (in the same way as in the last two lines of~\eqref{Paper01_intermediate_step_bound_martingale_part}). 
    This establishes~\eqref{Paper01_target_estimate_aldous_test_functions}, which completes the proof.
\end{proof}

We can now prove statement~(ii) of our programme, i.e.~we prove tightness of the sequence of processes $(\eta^{n})_{n \in \mathbb{N}}$, conditioned on $\eta^{n}(0) = \boldsymbol{\eta} \in \mathcal{S}_{0}$ for every $n \in \mathbb{N}$, in~$\mathscr{D}([0,\infty),\mathcal{S})$; we also characterise any subsequential limit via a martingale problem.

\begin{proposition} \label{Paper01_existence_solution_martingale_problem}
    For any $\boldsymbol{\eta} \in \mathcal{S}_{0}$, conditioning on $\eta^{n}(0) = \boldsymbol{\eta}$ for every $n \in \mathbb{N}$, the sequence $(\eta^{n})_{n \in \mathbb{N}}$ is tight in $\mathscr{D}([0,\infty), \mathcal{S})$ with respect to the $J_1$-topology. Moreover, for any subsequential limit $(\eta(t))_{t \geq 0}$ and any test function $\phi \in \mathscr{C}_*(\mathcal{S}, \mathbb{R})$, the process $(M^{\phi,\boldsymbol{\eta}}(t))_{t \geq 0}$ given by
    \begin{equation*}
        M^{\phi,\boldsymbol{\eta}}(T) \defeq \phi(\eta(T)) -\phi(\boldsymbol{\eta}) - \int_{0}^{T} (\mathcal{L}\phi)(\eta(t-)) \, dt \quad \forall T \geq 0
    \end{equation*}
    is a square-integrable càdlàg martingale.
\end{proposition}

\begin{proof}
    For $\boldsymbol{\eta} \in \mathcal{S}_{0}$, we will first establish tightness of the sequence $(\eta^{n})_{n \in \mathbb{N}}$ conditioned on $\eta^{n}(0) = \boldsymbol{\eta}$ for every $n \in \mathbb{N}$, by applying the Jakubowski criterion~\cite[Theorem~3.1]{jakubowski1986skorokhod}. According to this criterion, it suffices to verify two conditions:
    \begin{enumerate}[(i)]
        \item \underline{Strong compact containment condition:} For any $T \geq 0$ and $\varepsilon > 0$, there exists a compact subset $\mathscr{K} = \mathscr{K}({\varepsilon,\boldsymbol{\eta}, T}) \subset  \mathcal{S}$ such that
    \begin{equation*}
        \inf_{n \in \mathbb{N}} \; \mathbb{P}_{\boldsymbol{\eta}}\left(\eta^{n}(t) \in \mathscr{K} \; \forall t \in [0,T]\right) \geq 1 - \varepsilon;
    \end{equation*}
        \item \underline{Tightness of evaluations:} There is a subset $\mathfrak{A} \subset \mathscr{C}(\mathcal{S}, \mathbb{R})$ which is an algebra over $\mathbb{R}$ and separates points of $\mathcal{S}$, such that for any $\phi \in \mathfrak{A}$ and any $T \geq 0$, the sequence of $\mathbb{R}$-valued càdlàg processes $\left((\phi(\eta^{n}(t)))_{t \geq 0}\right)_{n \in \mathbb{N}}$ is tight in $\mathscr{D}([0,T], \mathbb{R})$.
    \end{enumerate}
   Condition~(i) follows from Lemma~\ref{Paper01_further_characterisation_relatively_compactness_discrete_IPS}, while condition~(ii) follows from Proposition~\ref{Paper01_tightness_evaluations_discrete_approximating_IPS} and from the fact that by Lemma~\ref{Paper01_characterisation_lipschitz_functions_semi_metric}, $\mathscr{C}^{\textrm{cyl}}_{b,*}(\mathcal{S}, \mathbb{R}) \subseteq \mathscr{C}(\mathcal{S}, \mathbb{R})$ is an algebra over $\mathbb{R}$ which separates the points of $\mathcal{S}$. Therefore, by Jakubowski's criterion, under $\mathbb{P}_{\boldsymbol{\eta}}$, the sequence of stochastic processes $(\eta^{n})_{n \in \mathbb{N}}$ is tight in~$\mathscr{D}([0, \infty), \mathcal{S})$.
    
    Let $(\eta^{n_i})_{i \in \mathbb{N}}$ be a subsequence that weakly converges to an $\mathcal{S}$-valued càdlàg process $(\eta(t))_{t \geq 0}$. It is immediate that $\eta(0) = \boldsymbol{\eta}$ almost surely. For $\phi \in \mathscr{C}_*(\mathcal{S}, \mathbb{R})$ and $n \in \mathbb{N}$, recall the definition of the square-integrable càdlàg martingale $(M^{n,\phi,\boldsymbol{\eta}}(t))_{t \geq 0}$ from~\eqref{Paper01_local_martingale_problem_discrete_approximation}. 
    Fix $T>0$.
    By the tightness of $(\eta^{n})_{n \in \mathbb{N}}$, and then by the uniform integrability of the sequence $([M^{n,\phi,\boldsymbol{\eta}}](T))_{n \in \mathbb{N}}$ and the tightness of $((M^{n,\phi,\boldsymbol{\eta}}(t))_{t \in [0,T]})_{n \in \mathbb{N}}$ from Proposition~\ref{Paper01_tightness_evaluations_discrete_approximating_IPS}, by the continuous mapping theorem and by standard results in stochastic analysis (see for instance~\cite[Proposition~4.3]{katzenberger1991solutions}), the sequence $((\phi(\eta^{n_i}(t)), M^{n_i,\phi,\boldsymbol{\eta}}(t))_{t \in [0,T]})_{i \in \mathbb{N}}$ weakly converges as $i \rightarrow \infty$ to a pair of real valued càdlàg processes $(\phi(\eta(t)), M^{\phi,\boldsymbol{\eta}}(t))_{t \in [0,T]}$, where $(M^{\phi,\boldsymbol{\eta}}(t))_{t \in [0,T]}$ is a square-integrable càdlàg martingale with respect to the filtration $\{\mathcal{F}^{\eta}_{t}\}_{t \in [0,T]}$. By Skorokhod's representation theorem (see e.g.~\cite[Theorem~3.1.7]{ethier2009markov}), we can construct the convergent subsequences and the limiting processes $\eta$ and $M^{\phi,\boldsymbol{\eta}}$ on the same probability space $(\Omega, \mathbb{P}_{\boldsymbol{\eta}}, \mathcal{F})$ in such a way that convergence occurs almost surely. Hence, by taking the limit as $i \rightarrow \infty$ on both sides of~\eqref{Paper01_local_martingale_problem_discrete_approximation} with $n=n_i$, we have
    \begin{equation} \label{Paper01_small_step_formulation_limiting_process_martingale_problem}
        \phi(\eta(t)) = \phi(\boldsymbol{\eta}) + M^{\phi,\boldsymbol{\eta}}(t) + \lim_{i \rightarrow \infty} \int_{0}^{t} (\mathcal{L}^{n_i}\phi)(\eta^{n_i}(t'-)) \, dt' \; \forall \, t \in [0,T].
    \end{equation}
    For $\omega \in \Omega$, let $(\eta^{n_i}(t))(\omega)$ (respectively $(\eta(t))(\omega)$) denote the configuration of the realisation of $\eta^{n_i}$ (respectively $\eta$) at time $t \geq 0$ associated to $\omega$. By the characterisation of the $J_1$-topology in $\mathscr{D}([0, T], \mathcal{S})$ (see e.g.~\cite[Theorem~3.7.2]{ethier2009markov}), under this construction, for almost every $\omega \in \Omega$, there exists a compact set $\mathcal{K} = \mathcal{K}(\omega) \subset \mathcal{S}$ such that
    \begin{equation*}
        (\eta^{n_i}(t))(\omega) \in \mathscr{K} \quad \forall \, i \in \mathbb{N} \textrm{ and } \forall \, t \in [0,T].
    \end{equation*}
    Then, by estimate~\eqref{P01:eq_simple_unif_bound_L_n_compacts} in Lemma~\ref{Paper01_convergence_infinitesimal_generator_lipischitz_semi_metric}, we can use dominated convergence, and therefore~\eqref{Paper01_convergence_sequence_generators_action_bizarre_set_functions} in Lemma~\ref{Paper01_convergence_infinitesimal_generator_lipischitz_semi_metric} implies that under this construction, the following limit holds almost surely for $t\in [0,T]$:
    \begin{equation*}
    \begin{aligned}
        \phi(\eta(t)) & = \phi(\boldsymbol{\eta}) + M^{\phi,\boldsymbol{\eta}}(t) + \int_{0}^{t} \lim_{i \rightarrow \infty} \left((\mathcal{L}^{n_i}\phi)(\eta^{n_i}(t'-))\right) \, dt' \\ & = \phi(\boldsymbol{\eta}) + M^{\phi,\boldsymbol{\eta}}(t) + \int_{0}^{t} (\mathcal{L}\phi)(\eta(t'-)) \, dt',
    \end{aligned}
    \end{equation*}
    which completes the proof.
\end{proof}

We now proceed to statement~(iii) of our programme, namely to establish tightness of the one-dimensional distributions for general initial conditions. A difficulty that must be addressed is that, when conditioning on $\eta^{n}(0) = \boldsymbol{\eta} \in \mathcal{S} \setminus \mathcal{S}_{0}$ for every $n \in \mathbb{N}$, estimate~\eqref{Paper01_moments_estimates_uniform_space} from Proposition~\ref{Paper01_bound_total_mass} is no longer available. Instead, we must rely on the cruder estimate~\eqref{Paper01_crude_estimate_moments_general_configuration}, which does not hold uniformly in $n \in \mathbb{N}$. To prove statement~(iii), we require estimates for the local number of particles that hold uniformly on compact subsets of $\mathcal{S}$. With this goal in mind, we will now show that, when the process is started from some configuration in a compact subset of $\mathcal{S}$,  the moments of the local particle numbers grow at most polynomially with distance from the origin. We will use the characterisation of compact subsets of $\mathcal{S}$ given in Proposition~\ref{Paper01_topological_properties_state_space}.

\begin{proposition} \label{Paper01_estimates_holding_uniformly_compact_subsets_S_N}
    For $p \in \mathbb{N}$, let $\psi^{(p)}: \mathcal{S} \rightarrow [0, \infty)$ be given by
\begin{equation} \label{Paper01_ell_p_norm_S_N}
    \psi^{(p)}(\boldsymbol{\xi}) \defeq \sum_{x \in L^{-1}\mathbb{Z}} \; \frac{\vert \vert \xi(x) \vert \vert_{\ell_{1}}^{p}}{(1 + \vert x \vert)^{2p}} \quad \forall  \boldsymbol{\xi} = (\xi_k(x))_{k \in \mathbb{N}_{0}, \, x \in L^{-1}\mathbb{Z}} \in \mathcal{S}.
\end{equation}
For a compact subset $\mathcal{K} \subset \mathcal{S}$, $T \geq 0$ and $p \in \mathbb{N}$, 
    \begin{equation} \label{Paper01_uniform_control_compacts_total_mass}
        \sup_{n \in \mathbb{N}} \; \sup_{\boldsymbol{\eta} \in \mathcal{K}} \; \sup_{t \in [0,T]} \; \mathbb{E}_{\boldsymbol{\eta}} \left[\psi^{(p)}(\eta^{n}(t))\right] < \infty.
    \end{equation}
    \ignore{Also, there exists $C_{T, \mathcal{K}} > 0$ and $p \geq 1$ such that, for any $x \in L^{-1}\mathbb{Z}$, 
    \begin{equation} \label{Paper01_uniform_control_sup_over_finite_time_number_particles_polynomial_function_deme}
        \sup_{n \in \mathbb{N}} \; \sup_{\boldsymbol{\eta} \in \mathcal{K}} \; \mathbb{E}_{\boldsymbol{\xi}}\left[\sup_{t \leq T} \; \vert \vert \eta^{n}(t,x) \vert \vert_{\ell_{1}}\right] \leq C_{T, \mathcal{K}} (1 + \vert x \vert)^{p}.
    \end{equation}}
    Moreover, for $x \in L^{-1}\mathbb{Z}$, 
     \begin{equation} \label{Paper01_uniform_control_compacts_type_space}
        \lim_{k \rightarrow \infty} \;  \sup_{n \in \mathbb{N}} \; \sup_{\boldsymbol{\eta} \in \mathcal{K}} \; \sup_{t \in [0,T]} \; \mathbb{E}_{\boldsymbol{\eta}} \Bigg[\sum_{j = k}^{\infty}\eta^{n}_{j}(t, x)\Bigg] = 0.
    \end{equation}
    Furthermore,
    \begin{equation} \label{Paper01_uniform_control_compacts_spatial_tail}
        \lim_{R \rightarrow \infty} \; \sup_{n \in \mathbb{N}} \, \sup_{\boldsymbol{\eta} \in \mathcal{K}} \; \sup_{t \in [0,T]} \; \mathbb{E}_{\boldsymbol{\eta}} \Bigg[\sum_{\substack{\{x \in L^{-1}\mathbb{Z}: \, \vert x \vert \geq R\}}} \frac{\vert \vert \eta^{n}(t,x) \vert \vert_{\ell_{1}}}{(1 + \vert x \vert)^{2}}\Bigg] = 0.
    \end{equation}
    In particular, for any compact subset $\mathcal{K} \subset \mathcal{S}$, $T \geq 0$ and $\varepsilon > 0$, there exists a compact subset $\mathscr{K}' = \mathscr{K}'(\mathcal{K}, T, \varepsilon) \subset \mathcal{S}$ such that
    \begin{equation} \label{Paper01_tightness_one_dimensional_distributions}
       \inf_{n \in \mathbb{N}} \; \inf_{\boldsymbol{\eta} \in \mathcal{K}} \; \inf_{t \in [0,T]} \; \mathbb{P}_{\boldsymbol{\eta}}\Big(\eta^{n}(t) \in \mathscr{K}'\Big) \geq 1 - \varepsilon.
    \end{equation}
\end{proposition}

\begin{proof}
     We will divide the proof into steps corresponding to each of the estimates~\eqref{Paper01_uniform_control_compacts_total_mass}-\eqref{Paper01_tightness_one_dimensional_distributions}.

    \medskip

    \noindent \underline{Step~$(1)$: uniform bounds on $\psi^{(p)}$}

    \medskip

   We start by proving~\eqref{Paper01_uniform_control_compacts_total_mass}; take $\mathcal K\subset \mathcal S$ compact and $T\ge 0$. Observe that by~\eqref{Paper01_definition_state_space_formal}, for any $\boldsymbol{\xi} \in \mathcal{S}$, the sequence $(\vert \vert \xi(x) \vert \vert_{\ell_{1}}  (1 + \vert x \vert)^{-2})_{x \in L^{-1}\mathbb{Z}}$ can be seen as an element of $\ell_{1}(L^{-1}\mathbb{Z})$, and since $\ell_{1}$ is continuously embedded in $\ell_{p}$ for $p \in (1, \infty)$, we conclude that the map $\psi^{(p)} : \mathcal{S} \rightarrow [0, \infty)$ given by~\eqref{Paper01_ell_p_norm_S_N} is continuous for $p \in \mathbb{N}$. Since $\mathcal{K}$ is compact, the continuity of $\psi^{(p)}$ implies that
    \begin{equation} \label{Paper01_immediate_estimate_continuity_quasi_moments_total_mass}
        \sup_{\boldsymbol{\xi} \in \mathcal{K}} \; \psi^{(p)}(\boldsymbol{\xi}) < \infty.
    \end{equation}
    Moreover, for $n \in \mathbb{N}$ and $\boldsymbol{\eta} \in \mathcal{K}$ fixed, let $(\eta^{n}(t))_{n \in \mathbb{N}}$ be the Markov process with infinitesimal generator $\mathcal{L}^{n}$ conditioned on $\eta^{n}(0) = \boldsymbol{\eta}$. Then, by Proposition~\ref{Paper01_bound_total_mass} and standard results in stochastic analysis (see~e.g.~\cite[Theorem~6.4.1]{ethier2009markov}), for $p \in \mathbb{N}$, there exists a càdlàg local martingale $M^{n,(p)}$ such that for any $t \in [0,T]$,
    \begin{equation} \label{Paper01_semi_martingale_formulation_ell_p_norm}
        \psi^{(p)}(\eta^{n}(t)) = \psi^{(p)}(\boldsymbol{\eta}) + M^{n,(p)}(t) + \int_{0}^{t} (\mathcal{L}^{n} \psi^{(p)})(\eta^{n}(t'-)) \, dt'.
    \end{equation}
    We will first establish that $(M^{n,(p)}(t))_{t \in [0,T]}$ is a càdlàg martingale. Observe that for $p \geq 1$, there exists $c_p > 0$ such that
    \begin{equation} \label{Paper01_elementary_inequality_difference_moments}
        1\vee (pa^{p-1})\le (a + 1)^{p} - a^{p} \leq c_p(a^{p-1} + 1) \quad \forall \, a \in [0, \infty).
    \end{equation}
    Recall the definition of $\boldsymbol{e}^{(x)}_{k}$ from before~\eqref{Paper01_scaled_polynomials_carrying_capacity}. 
We will show that~\eqref{Paper01_easy_criterion_martingale_feller_process_N_d} holds with $\phi=\psi^{(p)}$.
By~\eqref{Paper01_ell_p_norm_S_N} and then by~\eqref{Paper01_elementary_inequality_difference_moments}, we have that for $\boldsymbol{\zeta} =(\zeta_k(x))_{k\in \mathbb N_0,\, x\in L^{-1}\mathbb Z}\in \mathcal{S}$, $x \in L^{-1}\mathbb{Z}$ and $k \in \mathbb{N}_0$,
    \begin{equation} \label{Paper01_bound_norm_reaction_action_total_moment}
    \begin{aligned}
    & \Big\vert \psi^{(p)}\Big(\boldsymbol{\zeta} + \boldsymbol{e}^{(x)}_{k}\Big) - \psi^{(p)}(\boldsymbol{\zeta}) \Big\vert ^2 +  \mathds{1}_{\{\zeta_{k}(x)>0\}}\Big\vert \psi^{(p)}\Big(\boldsymbol{\zeta} - \boldsymbol{e}^{(x)}_{k}\Big) - \psi^{(p)}(\boldsymbol{\zeta}) \Big\vert ^2 \\ 
    & \quad = \frac{\Big((\vert \vert \zeta(x) \vert \vert_{\ell_1} + 1)^{p} - \vert \vert \zeta(x) \vert \vert_{\ell_1}^{p}\Big)^2 + \mathds{1}_{\{\zeta_k(x) > 0\}} \Big(\vert \vert \zeta(x) \vert \vert_{\ell_1}^{p} - (\vert \vert \zeta(x) \vert \vert_{\ell_1} - 1)^{p}\Big)^2}{(1 + \vert x \vert)^{4p}} \\ 
    & \quad \leq c_p^2 \frac{(\vert \vert \zeta(x) \vert \vert_{\ell_1}^{p-1}+1)^2 + \mathds{1}_{\{\zeta_k(x) > 0\}} ((\vert \vert \zeta(x) \vert \vert_{\ell_1} - 1)^{p-1} + 1)^2}{(1 + \vert x \vert)^{4p}} \\ 
    & \quad \leq \frac{2c^2_p(\vert \vert \zeta(x) \vert \vert_{\ell_1}^{p-1} + 1)^2}{(1 + \vert x \vert)^{4p}},
    \end{aligned}
    \end{equation}
    where for the last inequality we used that $\vert \vert \zeta(x) \vert \vert_{\ell_1} \in \mathbb{N}_{0}$. Moreover, for $z \in \{-L^{-1},L^{-1}\}$,
    by~\eqref{Paper01_ell_p_norm_S_N}, and then by~\eqref{Paper01_elementary_inequality_difference_moments} and since $\|\zeta(x)\|_{\ell_1}\in \mathbb N_0$,
    \begin{equation} \label{Paper01_bound_norm_migration_action_total_moment}
    \begin{aligned}
        & \mathds{1}_{\{\zeta_{k}(x)>0\}} \Big\vert \psi^{(p)}\Big(\boldsymbol{\zeta} + \boldsymbol{e}^{(x+ z)}_{k} - \boldsymbol{e}^{(x)}_{k}\Big) - \psi^{(p)}(\boldsymbol{\xi}) \Big\vert \\ 
        & \quad \leq \frac{(\vert \vert \zeta(x + z) \vert \vert_{\ell_1} + 1)^{p} - \vert \vert \zeta(x + z) \vert \vert_{\ell_1}^{p}}{(1 + \vert x + z \vert)^{2p}} + \mathds{1}_{\{\zeta_{k}(x)>0\}} \frac{\vert \vert \zeta(x) \vert \vert_{\ell_1}^{p} - (\vert \vert \zeta(x) \vert \vert_{\ell_1} - 1)^{p}}{(1 + \vert x\vert)^{2p}} \\ & \quad \leq \frac{c_p(\vert \vert \zeta(x + z) \vert \vert_{\ell_1}^{p-1} + 1)}{(1 + \vert x + z \vert)^{2p}} + \frac{c_p(\vert \vert \zeta(x) \vert \vert_{\ell_1}^{p-1} + 1)}{(1 + \vert x \vert)^{2p}}.
    \end{aligned}
    \end{equation}
    Recall the definitions of $(F^{b}_{k})_{k \in \mathbb{N}_{0}}$ and $(F^{d}_{k})_{k \in \mathbb{N}_{0}}$ from~\eqref{Paper01_discrete_birth_rate} and~\eqref{Paper01_discrete_death_rate}, respectively, and the definition of $\mathcal{L}^{n}$ from~\eqref{Paper01_generator_foutel_etheridge_model_restriction_n} and~\eqref{Paper01_infinitesimal_generator_restriction_n}. To bound the left-hand side of~\eqref{Paper01_easy_criterion_martingale_feller_process_N_d} with $\phi=\psi^{(p)}$, by~\eqref{Paper01_bound_norm_reaction_action_total_moment} and~\eqref{Paper01_bound_norm_migration_action_total_moment} we conclude that for $\boldsymbol{\eta} \in \mathcal{S}$,
    \begin{equation} \label{Paper01_intermediate_step_proof_martingale_property_associated_moments_total_mass}
    \begin{aligned}
    &E^n(\psi^{(p)},\boldsymbol{\eta},T)\\
    & \quad \leq \sum_{x \in \Lambda_n} \, \sum_{z \in \{-L^{-1}, L^{-1}\}} \int_{0}^{T} \Bigg(\frac{mc_p^2}{2} \\ 
    & \quad \quad \quad \quad \quad \quad \quad  \cdot \mathbb{E}_{\boldsymbol{\eta}}\Bigg[\vert \vert \eta^{n}(t-,x) \vert \vert_{\ell_1} \Bigg(\frac{(\vert \vert \eta^{n}(t-,x + z) \vert \vert_{\ell_1}^{p-1} + 1)}{(1 + \vert x + z \vert)^{2p}} + \frac{(\vert \vert \eta^{n}(t-,x) \vert \vert_{\ell_1}^{p-1} + 1)}{(1 + \vert x \vert)^{2p}} \Bigg)^{2}\Bigg]\Bigg) dt \\ 
    & \quad \quad + \sum_{x \in \Lambda_n} \sum_{k = 0}^{\infty} 4c_p^2\int_{0}^{T} \mathbb{E}_{\boldsymbol{\eta}}\Bigg[\frac{(\vert \vert \eta^{n}(t-,x) \vert \vert_{\ell_1}^{p-1} + 1)^{2}}{(1 + \vert x \vert)^{4p}}\Big(F^{b}_{k}(\eta^{n}(t-,x)) + F^{d}_{k}(\eta^{n}(t-,x))\Big)\Bigg] dt \\ 
    & \quad \leq \sum_{x \in \Lambda_n} \, \sum_{z \in \{-L^{-1}, L^{-1}\}} mc_p^2 \int_{0}^{T} \mathbb{E}_{\boldsymbol{\eta}}\Bigg[\vert \vert \eta^{n}(t-,x) \vert \vert_{\ell_1} \frac{(\vert \vert \eta^{n}(t-,x + z) \vert \vert_{\ell_1}^{p-1} + 1)^{2}}{(1 + \vert x + z \vert)^{4p}} \Bigg] dt \\ & \quad \quad + \sum_{x \in \Lambda_n} \, \sum_{z \in \{-L^{-1}, L^{-1}\}} mc_p^2 \int_{0}^{T} \mathbb{E}_{\boldsymbol{\eta}}\Bigg[\vert \vert \eta^{n}(t-,x) \vert \vert_{\ell_1} \frac{(\vert \vert \eta^{n}(t-,x) \vert \vert_{\ell_1}^{p-1} + 1)^{2}}{(1 + \vert x \vert)^{4p}} \Bigg] dt \\ & \quad \quad + \sum_{x \in \Lambda_n} 4c_p^2\int_{0}^{T} \mathbb{E}_{\boldsymbol{\eta}}\Bigg[\frac{(\vert \vert \eta^{n}(t-,x) \vert \vert_{\ell_1}^{p-1} + 1)^{2}}{(1 + \vert x \vert)^{4p}}\bar{q}^{N}(\vert \vert \eta^{n}(t-,x) \vert \vert_{\ell_1})\Bigg] dt,
\end{aligned}
\end{equation}
where for the last inequality we used~\eqref{Paper01_simple_bound_total_birth_death_rates_deterministic} and the elementary inequality
$(a+b)^{2} \leq 2(a^2 + b^2)$ $\forall a,b \in \mathbb{R}.$
We now bound each term on the right-hand side of~\eqref{Paper01_intermediate_step_proof_martingale_property_associated_moments_total_mass}. For the last two terms, since $\Lambda_n$ (defined after~\eqref{Paper01_scaling_parameters_for_discrete_approximation}) is a finite set, and since~$\bar{q}^{N}$ (defined in~\eqref{Paper01_bar_polynomial}) is a polynomial, by estimate~\eqref{Paper01_crude_estimate_moments_general_configuration} in Proposition~\ref{Paper01_bound_total_mass}, there exists $c^{(1)}_{n,p,T,\boldsymbol{\eta}} > 0$ such that
\begin{equation} \label{Paper01_intermediate_step_proof_martingale_property_associated_moments_total_mass_i}
\begin{aligned}
    & \sum_{x \in \Lambda_n} \, \sum_{z \in \{-L^{-1}, L^{-1}\}} mc_p^2 \int_{0}^{T} \mathbb{E}_{\boldsymbol{\eta}}\Bigg[\vert \vert \eta^{n}(t-,x) \vert \vert_{\ell_1} \frac{(\vert \vert \eta^{n}(t-,x) \vert \vert_{\ell_1}^{p-1} + 1)^{2}}{(1 + \vert x \vert)^{4p}} \Bigg] dt \\ & \quad + \sum_{x \in \Lambda_n} 4c_p^2\int_{0}^{T} \mathbb{E}_{\boldsymbol{\eta}}\Bigg[\frac{(\vert \vert \eta^{n}(t-,x) \vert \vert_{\ell_1}^{p-1} + 1)^{2}}{(1 + \vert x \vert)^{4p}}\bar{q}^{N}(\vert \vert \eta^{n}(t-,x) \vert \vert_{\ell_1})\Bigg] dt \\ & \quad \leq c^{(1)}_{n,p,T,\boldsymbol{\eta}}.
\end{aligned}
\end{equation}
To bound the first term on the right-hand side of~\eqref{Paper01_intermediate_step_proof_martingale_property_associated_moments_total_mass}, we use the Cauchy-Schwarz inequality, and then the fact $\Lambda_n$ is a finite set and~\eqref{Paper01_crude_estimate_moments_general_configuration} to conclude that there exists $c^{(2)}_{n,p,T,\boldsymbol{\eta}} > 0$ such that
\begin{equation} \label{Paper01_intermediate_step_proof_martingale_property_associated_moments_total_mass_ii}
\begin{aligned}
    & \sum_{x \in \Lambda_n} \, \sum_{z \in \{-L^{-1}, L^{-1}\}} mc_p^2 \int_{0}^{T} \mathbb{E}_{\boldsymbol{\eta}}\Bigg[\vert \vert \eta^{n}(t-,x) \vert \vert_{\ell_1} \frac{(\vert \vert \eta^{n}(t-,x + z) \vert \vert_{\ell_1}^{p-1} + 1)^{2}}{(1 + \vert x + z \vert)^{4p}} \Bigg] dt \\ & \quad \leq \sum_{x \in \Lambda_n} \, \sum_{z \in \{-L^{-1}, L^{-1}\}} mc_p^2 \int_{0}^{T} \mathbb{E}_{\boldsymbol{\eta}}\Big[\vert \vert \eta^{n}(t-,x) \vert \vert_{\ell_1}^{2}\Big]^{1/2} \mathbb{E}_{\boldsymbol{\eta}}\Bigg[\frac{(\vert \vert \eta^{n}(t-,x + z) \vert \vert_{\ell_1}^{p-1} + 1)^{4}}{(1 + \vert x + z \vert)^{8p}} \Bigg]^{1/2} dt \\ & \quad \leq c^{(2)}_{n,p,T,\boldsymbol{\eta}}.
\end{aligned}
\end{equation}
    Hence, as explained after~\eqref{Paper01_easy_criterion_martingale_feller_process_N_d}, by applying~\eqref{Paper01_intermediate_step_proof_martingale_property_associated_moments_total_mass_i} and~\eqref{Paper01_intermediate_step_proof_martingale_property_associated_moments_total_mass_ii} to~\eqref{Paper01_intermediate_step_proof_martingale_property_associated_moments_total_mass}, we conclude that the process $(M^{n,(p)}(t))_{t \in [0,T]}$ defined in~\eqref{Paper01_semi_martingale_formulation_ell_p_norm} is a square-integrable càdlàg martingale. Thus, by taking expectations on both sides of~\eqref{Paper01_semi_martingale_formulation_ell_p_norm}, we conclude that for $t \in [0,T]$,
    \begin{equation} \label{Paper01_first_step_control_total_mass_on_compacts}
    \begin{aligned}
         \mathbb{E}_{\boldsymbol{\eta}}\left[\psi^{(p)}(\eta^{n}(t))\right] & = \psi^{(p)}(\boldsymbol{\eta}) + \mathbb{E}_{\boldsymbol{\eta}}\Bigg[\int_{0}^{t} (\mathcal{L}^{n} \psi^{(p)})(\eta^{n}(t'-)) \, dt' \Bigg] \\ & = \psi^{(p)}(\boldsymbol{\eta}) + \frac{m}{2} \int_{0}^{t}  \mathbb{E}_{\boldsymbol{\eta}}\left[ (\mathcal{L}^{n}_{m} \psi^{(p)})(\eta^{n}(t'-)) \right] \, dt' + \int_{0}^{t} \mathbb{E}_{\boldsymbol{\eta}}\left[ (\mathcal{L}^{n}_{r} \psi^{(p)})(\eta^{n}(t'-)) \right] \, dt',
    \end{aligned}
    \end{equation}
    where $\mathcal{L}^{n}_{m}$ and $\mathcal{L}^{n}_{r}$ are given by~\eqref{Paper01_infinitesimal_generator_restriction_n}. Observe that by the definition of $\psi^{(p)}$ in~\eqref{Paper01_ell_p_norm_S_N}, and by applying~\eqref{Paper01_elementary_inequality_difference_moments},
    we have that for $\boldsymbol{\zeta} \in \mathcal{S}$, $x \in L^{-1}\mathbb{Z}$, $k \in \mathbb{N}_{0}$, $p \in \mathbb{N}$ and $z \in \{-L^{-1}, L^{-1}\}$,
    \begin{align}
        \psi^{(p)}\Big(\boldsymbol{\zeta} + \boldsymbol{e}^{(x)}_{k}\Big) - \psi^{(p)}(\boldsymbol{\zeta}) & \leq \frac{c_p(\vert \vert \zeta(x) \vert \vert_{\ell_1}^{p-1} + 1)}{(1 + \vert x \vert)^{2p}}, \label{Paper01_action_birth_ell_p}\\ 
        \mathds{1}_{\{\zeta_k(x) > 0\}}\Big(\psi^{(p)}\Big(\boldsymbol{\zeta} - \boldsymbol{e}^{(x)}_{k}\Big) - \psi^{(p)}(\boldsymbol{\zeta})\Big) & \leq - \mathds{1}_{\{\zeta_{k}(x)>0\}}\frac{(p(\vert \vert \zeta(x) \vert \vert_{\ell_1} - 1)^{p-1}) \vee 1}{(1 + \vert x \vert)^{2p}}, \label{Paper01_action_death_ell_p} \\ 
        \mathds{1}_{\{\zeta_{k}(x)>0\}} \Big(\psi^{(p)}\Big(\boldsymbol{\zeta} + \boldsymbol{e}^{(x+ z)}_{k} - \boldsymbol{e}^{(x)}_{k}\Big) - \psi^{(p)}(\boldsymbol{\xi}) \Big) &\leq \frac{c_p(\vert \vert \zeta(x + z) \vert \vert_{\ell_1}^{p-1} + 1)}{(1 + \vert x + z \vert)^{2p}} \label{Paper01_action_migration_ell_p_norm_only_positive_stuff}.
    \end{align}
    Recall the definitions of~$(F^{b}_{k})_{k \in \mathbb{N}_0}$ and~$(F^{d}_{k})_{k \in \mathbb{N}_0}$ in~\eqref{Paper01_discrete_birth_rate} and~\eqref{Paper01_discrete_death_rate}, and that the sequence $(K_{n})_{n \in \mathbb{N}} \subset \mathbb{N}$ satisfies~\eqref{Paper01_scaling_parameters_for_discrete_approximation}. Notice that since $s_{k} \leq 1$ for every $k \in \mathbb{N}_{0}$ by Assumption~\ref{Paper01_assumption_fitness_sequence}, and then since~$0 \leq \deg q_{+} < \deg q_-$ by Assumption~\ref{Paper01_assumption_polynomials}, we can conclude that there exists $c^{(3)}_{p,q_{+},q_{-},N} \ge 0$ such that for any $u = (u_{k})_{k \in \mathbb{N}_{0}} \in \ell_1^{+}$ and any~$n \in \mathbb{N}$,
    \begin{equation} \label{Paper01_pre_control_reaction_term_total_mass_on_compacts}
    \begin{aligned}
        & \sum_{k = 0}^{K_n} c_p(\vert \vert u \vert \vert_{\ell_1}^{p-1} + 1)F^b_k(u)  - \sum_{k = 0}^{\infty} ((p(\vert \vert u \vert \vert_{\ell_1} - 1)^{p-1})\vee 1)F^d_k(u) \\ 
        & \quad \leq \sum_{k = 0}^{\infty} \Big(c_p(\vert \vert u \vert \vert_{\ell_1}^{p-1} + 1)F^{b}_k(u) - ((p(\vert \vert u \vert \vert_{\ell_1} - 1)^{p-1})\vee 1)F^{d}_k(u) \Big) \\ 
        & \quad \leq c_p(\vert \vert u \vert \vert_{\ell_1}^{p-1} + 1)\vert \vert u \vert \vert_{\ell_1} q_{+}^{N}(\vert \vert u \vert \vert_{\ell_1}) - ((p(\vert \vert u \vert \vert_{\ell_1} - 1)^{p-1} )\vee 1)\vert \vert u \vert \vert_{\ell_1} q^{N}_{-}(\vert \vert u \vert \vert_{\ell_1})  \\ 
        & \quad \leq c^{(3)}_{p,q_{+},q_{-},N}.
    \end{aligned}
    \end{equation}
    We can now bound the integral terms on the right-hand side of~\eqref{Paper01_first_step_control_total_mass_on_compacts}. Starting with the second integral term, by the definition of $\mathcal{L}^{n}_{r}$ in~\eqref{Paper01_infinitesimal_generator_restriction_n} and by~\eqref{Paper01_action_birth_ell_p} and~\eqref{Paper01_action_death_ell_p}, we conclude that there exists $c^{(4)}_{p,q_+,q_-,N,L} \ge 0$ such that for any~$n \in \mathbb{N}$, $\boldsymbol{\eta} \in \mathcal{S}$, $t\in [0,T]$ and $p\ge 1$,
    \begin{equation} \label{Paper01_control_reaction_term_total_mass_on_compacts}
    \begin{aligned}
         & \int_{0}^{t} \mathbb{E}_{\boldsymbol{\eta}}\left[ (\mathcal{L}^{n}_{r} \psi^{(p)})(\eta^{n}(t'-)) \right] dt' \\ & \quad \leq \int_{0}^{t} \sum_{x \in \Lambda_n} \frac{1}{(1 + \vert x \vert)^{2p}} \mathbb{E_{\boldsymbol{\eta}}}\Bigg[\sum_{k = 0}^{K_n} c_p(\vert \vert \eta^{n}(t'-,x) \vert \vert_{\ell_1}^{p-1} + 1)F^b_k(\eta^{n}(t'-,x)) \\ & \quad \quad \quad \quad \quad \quad \quad \quad \quad \quad \quad \quad \quad - \sum_{k = 0}^{\infty} ((p(\vert \vert \eta^{n}(t'-,x) \vert \vert_{\ell_1} - 1)^{p-1})\vee 1)F^d_k(\eta^{n}(t'-,x))\Bigg] \, dt' \\ & \quad \leq t \sum_{x \in \Lambda_n} \frac{c^{(3)}_{p,q_+,q_-,N}}{(1 + \vert x \vert)^{2p}} \\ & \quad \leq tc^{(4)}_{p,q_+,q_-,N,L},
    \end{aligned}
    \end{equation}
    where in the second inequality we used~\eqref{Paper01_pre_control_reaction_term_total_mass_on_compacts}.

    We now bound the first integral term on the right-hand side of~\eqref{Paper01_first_step_control_total_mass_on_compacts}. Observe that for any $p > 1$, any non-negative sequence $(a_{x})_{x \in L^{-1}\mathbb{Z}}$, any $n \in \mathbb{N}$ and any $z \in \{-L^{-1}, L^{-1}\}$, an application of H\"older's inequality with conjugate exponents $p$ and $p/(p-1)$ yields
    \begin{equation} \label{Paper01_holder_control_migration_moments-total_norm}
    \begin{aligned}
        \sum_{x \in \Lambda_n} \frac{a_{x}a_{x + z}^{p-1}}{(1 + \vert x+z \vert)^{2p}} & \leq \sum_{x \in L^{-1}\mathbb{Z}} \frac{a_{x}a_{x + z}^{p-1}}{(1 + \vert x+z \vert)^{2}(1 + \vert x+z \vert)^{2(p-1)}} \\ 
        & \leq \Bigg( \sum_{x \in L^{-1}\mathbb{Z}} \frac{a_{x}^{p}}{(1 + \vert x+z \vert)^{2p}}\Bigg)^{1/p}\Bigg(\sum_{x \in L^{-1}\mathbb{Z}} \frac{a_{x+z}^{p}}{(1 + \vert x+z \vert)^{2p}}\Bigg)^{(p-1)/p} \\ 
        & \leq \sup_{y \in L^{-1}\mathbb{Z}} \Bigg(\frac{1 + \vert y \vert }{1 + \vert y+z \vert}\Bigg)^{2}\sum_{x \in L^{-1}\mathbb{Z}} \frac{a_{x}^{p}}{(1 + \vert x \vert)^{2p}}.
    \end{aligned}
    \end{equation}
    We also have
    \begin{equation} \label{eq:changeybyL-1}
        \sup_{y \in L^{-1}\mathbb{Z}} \left(\Bigg(\frac{1 + \vert y \vert}{1 + \vert y+L^{-1} \vert}\Bigg)^{2} + \Bigg(\frac{1 + \vert y \vert}{1 + \vert y-L^{-1} \vert}\Bigg)^{2} \right)< \infty.
    \end{equation}
    Therefore, by the definition of $\mathcal{L}^{n}_{m}$ in~\eqref{Paper01_infinitesimal_generator_restriction_n}, and by~\eqref{Paper01_action_migration_ell_p_norm_only_positive_stuff}, and then in the last inequality, 
    using the fact that for $u \in \ell_{1} \cap (\mathbb{N}_{0})^{\mathbb{N}_{0}}$ we have $\vert \vert u \vert \vert_{\ell_1} \leq \vert \vert u \vert \vert_{\ell_1}^{p}$ for the first term, 
    applying~\eqref{Paper01_holder_control_migration_moments-total_norm} to the second term in the case $p>1$, and using~\eqref{eq:changeybyL-1} and the definition of $\psi^{(p)}$ in~\eqref{Paper01_ell_p_norm_S_N}, we conclude that for $p\ge 1$, there exists $c^{(5)}_{p,m,L} \ge 0$ such that for any $n \in \mathbb{N}$, $\boldsymbol{\eta} \in \mathcal{S}$ and $t\in [0,T]$,
\begin{equation} \label{Paper01_control_migration_term_total_mass_on_compacts}
    \begin{aligned}
        & \frac{m}{2} \int_{0}^{t}  \mathbb{E}_{\boldsymbol{\eta}}\left[ (\mathcal{L}^{n}_{m} \psi^{(p)})(\eta^{n}(t'-)) \right] \, dt' \\ & \quad  \leq \frac{mc_p}{2} \sum_{z \in \{-L^{-1}, L^{-1}\}} \int_{0}^{t}\mathbb{E}_{\boldsymbol{\eta}}\Bigg[\sum_{x \in \Lambda_n} \frac{\vert \vert \eta^{n}(t'-,x) \vert \vert_{\ell_1}(\vert \vert \eta^{n}(t'-,x + z) \vert \vert_{\ell_1}^{p-1} + 1)}{(1 + \vert x + z \vert)^{2p}} \Bigg] \, dt' \\ 
        & \quad = \frac{mc_p}2  \sum_{z \in \{-L^{-1}, L^{-1}\}}\int_0^t \mathbb{E}_{\boldsymbol{\eta}}\Bigg[\sum_{x \in \Lambda_n} \frac{\vert \vert \eta^{n}(t'-,x) \vert \vert_{\ell_1}}{(1 + \vert x + z \vert)^{2p}} \Bigg] \, dt' \\ & \quad \quad + \frac{mc_p}{2} \sum_{z \in \{-L^{-1}, L^{-1}\}} \int_{0}^{t}\mathbb{E}_{\boldsymbol{\eta}}\Bigg[\sum_{x \in \Lambda_n} \frac{\vert \vert \eta^{n}(t'-,x) \vert \vert_{\ell_1}\vert \vert \eta^{n}(t'-,x + z) \vert \vert_{\ell_1}^{p-1}}{(1 + \vert x + z \vert)^{2p}} \Bigg] \, dt' \\ & \quad \leq c^{(5)}_{p,m,L} \int_{0}^{t}\mathbb{E}_{\boldsymbol{\eta}}\Big[ \psi^{(p)}(\eta^{n}(t'-)) \Big] \, dt'.
    \end{aligned}
    \end{equation}
    Applying~\eqref{Paper01_control_reaction_term_total_mass_on_compacts} and~\eqref{Paper01_control_migration_term_total_mass_on_compacts} to~\eqref{Paper01_first_step_control_total_mass_on_compacts}, we conclude that for any $n\in \mathbb N$ and $\boldsymbol{\eta} \in \mathcal{K}$,
    \begin{equation} \label{Paper01_almost_final_step_polynomial_growth}
        \mathbb{E}_{\boldsymbol{\eta}}\left[\psi^{(p)}(\eta^{n}(t))\right] \leq \sup_{\boldsymbol{\xi} \in \mathcal{K}} \; \psi^{(p)}(\boldsymbol{\xi})+ c^{(4)}_{p,q_+,q_-,N,L}t + c^{(5)}_{p,m,L}\int_{0}^{t}    \mathbb{E}_{\boldsymbol{\eta}}\left[\psi^{(p)}(\eta^{n}(t'-))\right] \, dt' \; \forall \, t \in [0,T].
    \end{equation}
    By the definition of $\psi^{(p)}$ in~\eqref{Paper01_ell_p_norm_S_N} and estimate~\eqref{Paper01_crude_estimate_moments_general_configuration} in Proposition~\ref{Paper01_bound_total_mass}, the map $t \mapsto   \mathbb{E}_{\boldsymbol{\eta}}\left[\psi^{(p)}(\eta^{n}(t)\right]$ is well defined and integrable over $[0,T]$ for any $n,p \in \mathbb{N}$ and any $\boldsymbol{\eta} \in \mathcal{S}$. Hence, using~\eqref{Paper01_immediate_estimate_continuity_quasi_moments_total_mass}, we can apply Gr\"{o}nwall's lemma to~\eqref{Paper01_almost_final_step_polynomial_growth}, which completes the proof of estimate~\eqref{Paper01_uniform_control_compacts_total_mass}.
    \ignore{\medskip

    \noindent \underline{Step $(2)$: local bounds on the number of particles over finite time intervals}

    \medskip

    To prove estimate~\eqref{Paper01_uniform_control_sup_over_finite_time_number_particles_polynomial_function_deme}, it is enough to repeat the same argument we used in the proof of Lemma~\ref{Paper01_supremum_moments_discrete_particle_system}, with the difference that here we apply estimate~\eqref{Paper01_uniform_control_compacts_total_mass} instead of Proposition~\ref{Paper01_bound_total_mass}, which explains why estimate~\eqref{Paper01_uniform_control_sup_over_finite_time_number_particles_polynomial_function_deme} is not uniform on $x \in L^{-1}\mathbb{Z}$. Since the proof follows along the same lines, we omit the details.}
    
    \medskip

    \noindent \underline{Step $(2)$: local bounds on the number of particles carrying a high number of mutations}

    \medskip

    We now proceed to the proof of estimate~\eqref{Paper01_uniform_control_compacts_type_space}; take $\mathcal K \subset \mathcal S$ compact and $T>0$. For $k \in \mathbb{N}_{0}$, let $\mathcal{I}_{k} \defeq \{j \in \mathbb{N}_{0}: \, j \geq k\}$. For $x \in L^{-1}\mathbb{Z}$, recall the definition of $\phi^{(x)}_{\mathcal{I}_{k}}$ in~\eqref{Paper01_number_particles_position_x_number_mutations_in_I}, and, for $n \in \mathbb{N}$, let $(X^{n}(t))_{t \geq 0}$ be a simple symmetric random walk with reflecting boundaries on $\Lambda_n$, as defined before~\eqref{eq:RWinitcond}. By Lemma~\ref{Paper01_greens_function_representation}, for $t \in [0,T]$, $n \in \mathbb{N}$, $k\in \mathbb N_0$, $x\in \Lambda_n$ and $\boldsymbol{\eta} \in \mathcal{K}$, we have
    \begin{equation} \label{Paper01_green_function_representation_local_density_higher_number_mutations}
    \begin{aligned}
        \mathbb{E}_{\boldsymbol{\eta}}\Bigg[\sum_{j = k}^{\infty} \eta^{n}_j(t,x)\Bigg] & = \sum_{\substack{y \in \Lambda_{n}}} \; \mathbb{P}_{x}(X^{n}(t) = y) \phi_{\mathcal I_k}^{(y)}(\boldsymbol{\eta}) \\ & \quad \quad + \int_{0}^{t} \sum_{\substack{y \in  \Lambda_{n}}} \; \mathbb{P}_{x}(X^{n}(t - t') = y) \mathbb{E}_{\boldsymbol{\eta}}\left[(\mathcal{L}^{n}_{r}\phi_{\mathcal I_k}^{(y)})(\eta^{n}(t'-))\right] \, dt'.
    \end{aligned}
    \end{equation}
    We now bound the terms on the right-hand side of~\eqref{Paper01_green_function_representation_local_density_higher_number_mutations} separately. For the first term, we will show that for any $x\in L^{-1}\mathbb Z$,
    \begin{equation} \label{Paper01_limit_migration_number_mutations_uniform_compact_sets}
        \lim_{k \rightarrow \infty} \; \sup_{n \in \mathbb{N}} \; \sup_{\boldsymbol{\eta} \in \mathcal{K}} \; \sup_{t \in [0,T]} \; \sum_{\substack{y \in \Lambda_n}} \; \mathbb{P}_{x}(X^{n}(t) = y) \phi_{\mathcal I_k}^{(y)}(\boldsymbol{\eta}) = 0.
    \end{equation}
    For $\varepsilon > 0$ and  $x\in L^{-1}\mathbb Z$, let $R_{\varepsilon,x} = R_{\varepsilon,x}(T, \mathcal{K},m,L) > e^2(mT + 1)L^{-1}$ be some large positive number, to be specified later. Observe that for $n\in \mathbb N$, $k\in \mathbb N_0$, $\boldsymbol{\eta}\in \mathcal K$ and $t\in [0,T]$,
    \begin{equation} \label{Paper01_first_steps_local_control_mutations_on_compacts}
    \begin{aligned}
       & \sum_{\substack{y \in \Lambda_n}} \; \mathbb{P}_{x}\Big(X^n(t) = y\Big) \phi_{\mathcal I_k}^{(y)}(\boldsymbol{\eta}) \\ & \quad = \sum_{\substack{\{y \in \Lambda_n: \, \vert y - x \vert \leq R_{\varepsilon,x}\}}} \; \mathbb{P}_{x}(X^n(t) = y) \phi_{\mathcal I_k}^{(y)}(\boldsymbol{\eta}) + \sum_{\substack{\{y \in \Lambda_n: \, \vert y - x \vert > R_{\varepsilon,x}\}}} \; \mathbb{P}_{x}(X^n(t) = y) \phi_{\mathcal I_k}^{(y)}(\boldsymbol{\eta}) \\ 
       & \quad \leq \sup_{\boldsymbol{\xi} \in \mathcal{K}} \; \sup_{\substack{\{y \in L^{-1}\mathbb{Z}: \, \vert y - x \vert \leq R_{\varepsilon,x}\}}} \phi_{\mathcal I_k}^{(y)}(\boldsymbol{\xi})  \\ & \quad \quad \quad + \sum_{\substack{\{y \in \Lambda_n: \, \vert y - x \vert > R_{\varepsilon,x}\}}} \; \mathbb{P}_{x}(X^n(t) = y)(1 + \vert y \vert)^{2} \Bigg(\sup_{\boldsymbol{\xi}=(\xi(z))_{z\in L^{-1}\mathbb Z} \in \mathcal{K}} \, \sup_{z \in L^{-1}\mathbb{Z}} \, \frac{\vert \vert \xi(z) \vert \vert_{\ell_1}}{(1 + \vert z \vert)^{2}}\Bigg) \\ & \quad \leq \sup_{\boldsymbol{\xi} \in \mathcal{K}} \; \sup_{\substack{\{y \in L^{-1}\mathbb{Z}: \, \vert y - x \vert \leq R_{\varepsilon,x}\}}} \phi_{\mathcal I_k}^{(y)}(\boldsymbol{\xi})  \\ & \quad \quad \quad + \sup_{\boldsymbol{\xi} \in \mathcal{K}} \vert \vert \vert \boldsymbol{\xi} \vert \vert \vert_{\mathcal{S}}\sum_{\substack{\{y \in \Lambda_n: \, \vert y - x \vert > R_{\varepsilon,x}\}}} \mathbb{P}_{x}(X^n(t) = y)(1 + \vert y \vert)^{2},
    \end{aligned}
    \end{equation}
    where the last inequality follows from~\eqref{Paper01_definition_state_space_formal}.
    We now bound the terms on the right-hand side of~\eqref{Paper01_first_steps_local_control_mutations_on_compacts} separately. For the second term, since we chose $R_{\varepsilon,x} > e^2(mT + 1)L^{-1}$, a standard large deviations estimate for continuous-time random walks (see Corollary~\ref{Paper01_corollary_lemma_tail_continuous_time_rw}) yields
    \begin{equation} \label{Paper01_intermediate_step_local_bound_high_number_mutations_i}
     \sum_{\substack{\{y \in \Lambda_n: \, \vert y - x \vert > R_{\varepsilon,x}\}}} \mathbb{P}_{x}(X^n(t) = y)(1 + \vert y \vert)^{2}  \leq   \sum_{\substack{\{y \in L^{-1}\mathbb{Z}: \, \vert y - x \vert > R_{\varepsilon,x}\}}} \; e^{-L \vert y - x\vert}(1 + \vert y \vert)^{2}.
    \end{equation}
    Since $\mathcal{K} \subset \mathcal{S}$ is compact, by Proposition~\ref{Paper01_topological_properties_state_space}(i) we have $\sup_{\boldsymbol{\xi} \in \mathcal{K}} \vert \vert \vert \boldsymbol{\xi} \vert \vert \vert_{\mathcal{S}} < \infty$, and so we can choose $R_{\varepsilon,x} > 0$ sufficiently large that~\eqref{Paper01_intermediate_step_local_bound_high_number_mutations_i} yields
    \begin{equation} \label{Paper01_intermediate_step_local_bound_high_number_mutations_ii}
        \sup_{\boldsymbol{\xi} \in \mathcal{K}} \vert \vert \vert \boldsymbol{\xi} \vert \vert \vert_{\mathcal{S}} \sum_{\substack{\{y \in \Lambda_n: \, \vert y - x \vert > R_{\varepsilon,x}\}}} \mathbb{P}_{x}(X^n(t) = y)(1 + \vert y \vert)^{2} < \varepsilon.
    \end{equation}
    Fixing such an $R_{\varepsilon,x}$, observe that since $\mathcal{K}$ is compact, by Proposition~\ref{Paper01_topological_properties_state_space}(iii), there exists $k_{\varepsilon, x} \in \mathbb{N}$ such that for $k \geq k_{\varepsilon, x}$,
    \begin{equation} \label{Paper01_intermediate_step_local_bound_high_number_mutations_iii}
       \sup_{\boldsymbol{\xi} \in \mathcal{K}} \; \sup_{\substack{\{y \in L^{-1}\mathbb{Z}: \, \vert y - x \vert \leq R_{\varepsilon,x}\}}} \phi_{\mathcal I_k}^{(y)}(\boldsymbol{\xi})  = 0.
    \end{equation}
    By applying~\eqref{Paper01_intermediate_step_local_bound_high_number_mutations_ii} and~\eqref{Paper01_intermediate_step_local_bound_high_number_mutations_iii} to~\eqref{Paper01_first_steps_local_control_mutations_on_compacts}, we conclude that for $k \geq k_{\varepsilon, x}$,
    \begin{equation*}
        \sup_{n \in \mathbb{N}} \; \sup_{\boldsymbol{\eta} \in \mathcal{K}} \; \sup_{t \in [0,T]} \; \sum_{\substack{y \in \Lambda_n}} \; \mathbb{P}_{x}(X^{n}(t) = y) \phi_{\mathcal I_k}^{(y)}(\boldsymbol{\eta}) < \varepsilon.
    \end{equation*}
    Since $\varepsilon > 0$ was arbitrary, the limit in~\eqref{Paper01_limit_migration_number_mutations_uniform_compact_sets} holds.
    
    For the second term on the right-hand side of~\eqref{Paper01_green_function_representation_local_density_higher_number_mutations}, we first notice that by~\eqref{Paper01_ell_p_norm_S_N} and~\eqref{Paper01_uniform_control_compacts_total_mass}, for every $p \in \mathbb{N}_{0}$, there exists $c^{(6)}_{p,\mathcal{K},T} \ge 0$ such that
    \begin{equation} \label{Paper01_simple_consequence_estimate_total_mass}
        \sup_{t \in [0,T]} \, \sup_{n \in \mathbb{N}} \, \sup_{\boldsymbol{\xi} \in \mathcal{K}} \, \mathbb{E}_{\boldsymbol{\xi}}\Big[\vert \vert \eta^n(t,y) \vert \vert_{\ell_1}^{p}\Big] \leq c^{(6)}_{p,\mathcal{K},T} (1 + \vert y \vert)^{2p} \quad \forall \, y \in L^{-1}\mathbb{Z}. 
    \end{equation}
    By combining~\eqref{Paper01_basic_application_generator_reaction} with~\eqref{Paper01_bound_birth_rate_high_mutation_load} from the proof of Lemma~\ref{Paper01_supremum_moments_discrete_particle_system}, we conclude that there exists $c^{(7)}_{q_+,N} \ge 0$ such that for any $n,k \in \mathbb{N}$, any $\boldsymbol{\eta} \in \mathcal{K}$ and any $t \in [0,T]$,
    \begin{equation} \label{Paper01_intermediate_step_local_bound_high_number_mutations_iv}
    \begin{aligned}
        & \int_{0}^{t} \sum_{\substack{y \in  \Lambda_{n}}} \; \mathbb{P}_{x}(X^{n}(t - t') = y) \mathbb{E}_{\boldsymbol{\eta}}\left[(\mathcal{L}^{n}_{r}\phi_{\mathcal I_k}^{(y)})(\eta^{n}(t'-))\right] \, dt' \\ & \quad \leq \int_{0}^{t} \sum_{\substack{y \in  \Lambda_{n}}} \; \mathbb{P}_{x}(X^{n}(t - t') = y) \mathbb{E}_{\boldsymbol{\eta}}\left[s_{k-1} \vert \vert \eta^{n}(t'-,y) \vert \vert_{\ell_{1}}q^{N}_{+}(\vert \vert \eta^{n}(t'-,y) \vert \vert_{\ell_{1}})\right] \, dt' \\ & \quad \le s_{k-1} \int_{0}^{t} \sum_{\substack{y \in  \Lambda_{n}}} \; \mathbb{P}_{x}(X^{n}(t - t') = y)c^{(7)}_{q_+,N} \sum_{p = 1}^{1 + \deg q_+}\mathbb{E}_{\boldsymbol{\eta}}\Big[ \vert \vert \eta^{n}(t'-,y) \vert \vert_{\ell_{1}}^{p}\Big] \, dt' \\ & \quad \leq s_{k-1}c^{(7)}_{q_+,N} c^{(6)}_{p,\mathcal{K},T} \int_{0}^{t} \sum_{\substack{y \in  \Lambda_{n}}} \; \mathbb{P}_{x}(X^{n}(t - t') = y) \sum_{p = 1}^{1 + \deg q_+} (1 + \vert y \vert)^{2p} \, dt' \\ & \quad \leq s_{k-1}c^{(7)}_{q_+,N}c^{(6)}_{p,\mathcal{K},T}(1 + \deg q_+) \int_0^t \mathbb{E}_{x}[(1 + \vert X^n(t - t') \vert)^{2(1+\deg q_+)}] \, dt',
    \end{aligned}
    \end{equation}
    where in the third inequality we used~\eqref{Paper01_simple_consequence_estimate_total_mass}, and in the last inequality we used the fact that $1 + \vert y \vert \geq 1$ $\forall y\in L^{-1}\mathbb Z$.
    Since $\lim_{k \rightarrow \infty} s_k = 0$ by Assumption~\ref{Paper01_assumption_fitness_sequence}, and since the moments of $((X^{n}(t))_{t \in [0,T]})_{n\in \mathbb N}$ under $\mathbb P_x$ are uniformly bounded (e.g.~as an easy consequence of Corollary~\ref{Paper01_corollary_lemma_tail_continuous_time_rw}(ii)), we conclude from~\eqref{Paper01_intermediate_step_local_bound_high_number_mutations_iv} that
    \begin{equation} \label{Paper01_second_important_limit_local_control_particles_mutations_compacts}
        \lim_{k \rightarrow \infty} \; \sup_{n \in \mathbb{N}}\; \sup_{\boldsymbol{\eta \in \mathcal{K}}} \; \sup_{t \in [0,T]} \; \int_{0}^{t} \sum_{{y \in \Lambda_n}} \; \mathbb{P}_{x}(X^{n}(t - t') = y) \mathbb{E}_{\boldsymbol{\eta}}\left[(\mathcal{L}^{n}_{r}\phi_{\mathcal I_k}^{(y)})(\eta^{n}(t'-))\right] \, dt' = 0.
    \end{equation}
    Thus, applying~\eqref{Paper01_limit_migration_number_mutations_uniform_compact_sets} and~\eqref{Paper01_second_important_limit_local_control_particles_mutations_compacts} to~\eqref{Paper01_green_function_representation_local_density_higher_number_mutations}, and using Proposition~\ref{Paper01_topological_properties_state_space}(iii) for the case of $n\in \mathbb N$ such that $x\notin \Lambda_n$, we conclude that~\eqref{Paper01_uniform_control_compacts_type_space} holds, which completes Step~(2).

    \medskip

    \noindent \underline{Step~$(3)$: Uniform control on the spatial spread of mass}

    \medskip

    We now establish~\eqref{Paper01_uniform_control_compacts_spatial_tail}; take $\mathcal K\subset \mathcal S$ compact and $T>0$. For $x \in L^{-1}\mathbb{Z}$, recall the definition of $\phi^{(x)}_{\mathbb{N}_{0}}$ in~\eqref{Paper01_number_particles_position_x_number_mutations_in_I}, and, for $n \in \mathbb{N}$, let $(X^{n}(t))_{t \geq 0}$ be a simple symmetric random walk with reflecting boundaries on $\Lambda_n$ (as defined before~\eqref{eq:RWinitcond}).
    Recall the definitions of $(F^b_k)_{k \in \mathbb{N}_{0}}$ and $(F^d_k)_{k \in \mathbb{N}_{0}}$ from~\eqref{Paper01_discrete_birth_rate} and~\eqref{Paper01_discrete_death_rate}, respectively. By the definition of $\mathcal{L}^{n}_{r}$ in~\eqref{Paper01_infinitesimal_generator_restriction_n}, and then since $s_k\le 1$ $\forall k\in \mathbb N_0$ by Assumption~\ref{Paper01_assumption_fitness_sequence}, and finally since $0\le \deg q_+<\deg q_-$ by Assumption~\ref{Paper01_assumption_polynomials},
    we conclude that there exists $c^{(8)}_{q_+,q_-,N} \ge 0$ such that for any $\boldsymbol{\xi} =(\xi(x))_{x\in L^{-1}\mathbb Z}\in \mathcal{S}$ and~$y \in \Lambda_n$,
    \begin{equation} \label{Paper01_rewriting_trivial_consequence_deaths_prevail}
    \begin{aligned}
        (\mathcal{L}^{n}_{r}\phi^{(y)}_{\mathbb{N}_{0}})(\boldsymbol{\xi}) = \sum_{j = 0}^{K_n} F^b_j(\xi(y)) - \sum_{j = 0}^{\infty} F^d_j(\xi(y)) 
        \le \|\xi(y)\|_{\ell_1}(q_+^N(\|\xi(y)\|_{\ell_1})-q_-^N(\|\xi(y)\|_{\ell_1}))
        \leq c^{(8)}_{q_+,q_-,N}.
    \end{aligned}
    \end{equation}
    By Lemma~\ref{Paper01_greens_function_representation} and~\eqref{Paper01_rewriting_trivial_consequence_deaths_prevail} (and recalling that $\eta^n(t,x)=\eta^n(0,x)$ for $x\notin \Lambda_n$ and $t\ge 0$), for $t \in [0,T]$, $n \in \mathbb{N}$, $x\in L^{-1}\mathbb Z$ and $\boldsymbol{\eta}=(\eta(y))_{y\in L^{-1}\mathbb Z} \in \mathcal{K}$, we have
    \begin{equation*}
    \begin{aligned}
        \mathbb{E}_{\boldsymbol{\eta}}\Big[\vert \vert \eta^{n}(t,x) \vert \vert_{\ell_1}\Big] \leq \sum_{y \in L^{-1}\mathbb Z} \; \mathbb{P}_{x}(X^{n}(t) = y) \vert \vert \eta(y) \vert \vert_{\ell_1} + c^{(8)}_{q_+,q_-,N}t.
    \end{aligned}
    \end{equation*}
    Hence, for $R > 0$, $n\in \mathbb N$, $\boldsymbol{\eta}=(\eta(y))_{y\in L^{-1}\mathbb Z} \in \mathcal{K}$ and $t \in [0,T]$,
    \begin{equation} \label{Paper01_rewriting_number_of_particles_for_bound_spatial_dispertion}
    \begin{aligned}
        & \mathbb{E}_{\boldsymbol{\eta}} \Bigg[\sum_{\{x \in L^{-1}\mathbb{Z}: \, \vert x \vert \geq R\}} \frac{\vert \vert \eta^{n}(t,x) \vert \vert_{\ell_{1}}}{(1 + \vert x \vert)^{2}}\Bigg] \\ & \quad \leq \sum_{\substack{\{x \in L^{-1}\mathbb{Z}: \, \vert x \vert \geq R\}}} \frac{c^{(8)}_{q_{+},q_{-},N}t}{(1 + \vert x \vert)^{2}} + \sum_{\substack{\{x \in L^{-1}\mathbb{Z}: \, \vert x \vert \geq R\}}} \; \sum_{y \in L^{-1}\mathbb Z} \mathbb{P}_{x}(X^n(t) = y)\frac{\vert \vert \eta(y) \vert \vert_{\ell_{1}}}{(1 + \vert x \vert)^{2}}.
    \end{aligned}
    \end{equation}
    Observe that
    \begin{equation} \label{Paper01_first_limit_spatial_bound_compacts}
    \begin{aligned}
        \lim_{R \rightarrow \infty} \sum_{\substack{\{x \in L^{-1}\mathbb{Z}: \, \vert x \vert \geq R\}}} \frac{c^{(8)}_{q_{+},q_{-},N}T}{(1 + \vert x \vert)^{2}} = 0.
    \end{aligned}
    \end{equation}
    By~\eqref{Paper01_rewriting_number_of_particles_for_bound_spatial_dispertion} and~\eqref{Paper01_first_limit_spatial_bound_compacts}, the proof of~\eqref{Paper01_uniform_control_compacts_spatial_tail} will follow once we establish that 
        \begin{equation}
        \label{Paper01_equivalent_formulation_limit_spatial_tails}
        \lim_{R \rightarrow \infty} \; \sup_{n \in \mathbb{N}} \; \sup_{\boldsymbol{\eta}=(\eta(y))_{y\in L^{-1}\mathbb Z} \in \mathcal{K}} \; \sup_{t \in [0,T]} \; \sum_{\substack{\{x \in L^{-1}\mathbb{Z}: \, \vert x \vert \geq R\}}} \; \sum_{y \in L^{-1}\mathbb{Z}} \mathbb{P}_{x}(X^n(t) = y)\frac{\vert \vert \eta(y) \vert \vert_{\ell_{1}}}{(1 + \vert x \vert)^{2}} = 0.
    \end{equation}
  Note that by the triangle inequality,
  \begin{equation} \label{Paper01_application_triangle_inequality_spatial_dispersion}
      \frac{1}{1 + \vert x \vert} \leq \frac{1 + \vert y - x \vert}{1 + \vert y \vert} \quad \forall \, x,y \in L^{-1}\mathbb{Z}.
  \end{equation}
 For $\varepsilon > 0$, let $r_{\varepsilon}$ and $R_{\varepsilon}$ be positive numbers, to be specified later, such that $$R_{\varepsilon} > r_{\varepsilon} > e^2(mT+1)L^{-1},$$ 
 and take $R \ge R_\varepsilon$. By~\eqref{Paper01_application_triangle_inequality_spatial_dispersion}, for $n\in \mathbb N$, $\boldsymbol{\eta}=(\eta(y))_{y\in L^{-1}\mathbb Z} \in \mathcal{K}$ and $t\in [0,T]$,
    \begin{equation} \label{Paper01_intermediate_step_control_spatial_tail}
    \begin{aligned}
        & \sum_{\substack{\{x \in L^{-1}\mathbb{Z}: \, \vert x \vert \geq R\}}} \; \sum_{y \in L^{-1}\mathbb{Z}} \mathbb{P}_{x}(X^n(t) = y)\frac{\vert \vert \eta(y) \vert \vert_{\ell_{1}}}{(1 + \vert x \vert)^{2}} \\ & \quad \leq \sum_{\{\substack{x \in L^{-1}\mathbb{Z}: \, \vert x \vert \geq R\}}} \; \sum_{\substack{\{y \in L^{-1}\mathbb{Z}: \, \vert y - x \vert > r_{\varepsilon}\}}} \; \mathbb{P}_{x}(X^n(t) = y) \cdot {(1 + \vert y - x \vert)^{2}} \cdot\frac{\vert \vert \eta(y) \vert \vert_{\ell_{1}}}{(1 + \vert y \vert)^{2}} \\ & \quad \quad \quad + \sum_{\substack{\{x \in L^{-1}\mathbb{Z}: \, \vert x \vert \geq R\}}} \; \sum_{\substack{\{y \in L^{-1}\mathbb{Z}: \, \vert y - x \vert \leq r_{\varepsilon}\}}} \; \mathbb{P}_{x}(X^n(t) = y ) \cdot {(1 + \vert y - x \vert)^{2}} \cdot\frac{\vert \vert \eta(y) \vert \vert_{\ell_{1}}}{(1 + \vert y \vert)^{2}}.
    \end{aligned}
    \end{equation}
    We now bound the two terms on the right-hand side of~\eqref{Paper01_intermediate_step_control_spatial_tail} separately. For the first term, recall from before~\eqref{Paper01_intermediate_step_control_spatial_tail} that $r_\varepsilon > e^2(mT + 1)L^{-1}$, and so by Corollary~\ref{Paper01_corollary_lemma_tail_continuous_time_rw},
    \begin{equation*}
    \begin{aligned}
        & \sum_{\substack{\{x \in L^{-1}\mathbb{Z}: \, \vert x \vert \geq R\}}} \; \sum_{\substack{\{y \in L^{-1}\mathbb{Z}: \, \vert y - x \vert > r_{\varepsilon}\}}} \; \mathbb{P}_{x}\Big(X^n(t) = y\Big) \cdot {(1 + \vert y - x \vert)^{2}} \cdot\frac{\vert \vert \eta(y) \vert \vert_{\ell_{1}}}{(1 + \vert y \vert)^{2}} \\ & \quad \leq \sum_{\substack{x \in L^{-1}\mathbb{Z}}} \; \sum_{\substack{\{y \in L^{-1}\mathbb{Z}: \, \vert y - x \vert > r_{\varepsilon}\}}} e^{-L\vert y - x \vert} (1 + \vert y - x \vert)^2 \frac{\vert \vert \eta(y) \vert \vert_{\ell_1}}{(1 + \vert y \vert)^2} \\ & \quad = \sum_{y \in L^{-1}\mathbb{Z}} \; \sum_{\{z \in L^{-1}\mathbb{Z}: \, \vert z \vert > r_\varepsilon\}} \; e^{-L\vert z \vert} (1 + \vert z \vert)^2 \frac{\vert \vert \eta(y) \vert \vert_{\ell_1}}{(1 + \vert y \vert)^2} \\ & \quad \leq \sup_{\boldsymbol{\xi} \in \mathcal{K}} \vert \vert \vert \boldsymbol{\xi} \vert \vert \vert_{\mathcal{S}} \sum_{\{z \in L^{-1}\mathbb{Z}: \, \vert z \vert > r_\varepsilon\}} \; e^{-L\vert z \vert} (1 + \vert z \vert)^2,
    \end{aligned}
    \end{equation*}
    where for the last inequality we used the definition of $\vert \vert \vert \cdot \vert \vert \vert_\mathcal{S}$ in~\eqref{Paper01_definition_state_space_formal}. Since $\mathcal{K} \subset \mathcal{S}$ is compact, by Proposition~\ref{Paper01_topological_properties_state_space}(i) we have $\sup_{\boldsymbol{\xi} \in \mathcal{K}} \vert \vert \vert \boldsymbol{\xi} \vert \vert \vert_{\mathcal{S}} < \infty$, and therefore we can pick $r_{\varepsilon} > 0$ sufficiently large that for any $R\ge R_\varepsilon >r_\varepsilon$,
    \begin{equation} \label{Paper01_first_bound_term_limit_spatial_tail_compacts}
        \sup_{n\in \mathbb N}\sup_{\boldsymbol{\eta} \in \mathcal{K}} \; \sup_{t \in [0,T]} \; \sum_{\substack{\{x \in L^{-1}\mathbb{Z}: \, \vert x \vert \geq R\}}} \; \sum_{\substack{\{y \in L^{-1}\mathbb{Z}: \, \vert y - x \vert > r_{\varepsilon}\}}} \; \mathbb{P}_{x}\Big(X^n(t) = y\Big) \cdot {(1 + \vert y - x \vert)^{2}} \cdot\frac{\vert \vert \eta(y) \vert \vert_{\ell_{1}}}{(1 + \vert y \vert)^{2}} \leq \frac{\varepsilon}{2}.
    \end{equation}
    For the second term on the right-hand side of~\eqref{Paper01_intermediate_step_control_spatial_tail}, we observe that for any $x,y \in L^{-1}\mathbb{Z}$ with $\vert x\vert \geq R > R_{\varepsilon}$ and $\vert y - x \vert \leq r_{\varepsilon}$, we have $\vert y \vert \geq R-r_{\varepsilon}$. Moreover, for any $y \in L^{-1}\mathbb{Z}$, we have  $\#(L^{-1}\mathbb{Z} \, \cap \, [y -r_{\varepsilon}, y + r_{\varepsilon}])\le 2\lceil r_{\varepsilon}L \rceil + 1
    $. Hence,
    \begin{equation}
    \label{Paper01_second_bound_term_limit_spatial_tail_compacts_pre_step}
    \begin{aligned}
        & \sum_{\substack{\{x \in L^{-1}\mathbb{Z}: \, \vert x \vert \geq R\}}} \; \sum_{\substack{\{y \in L^{-1}\mathbb{Z}: \, \vert y - x \vert \leq r_{\varepsilon}\}}} \; \mathbb{P}_{x}(X^n(t) = y) \cdot {(1 + \vert y - x \vert)^{2}} \cdot\frac{\vert \vert \eta(y) \vert \vert_{\ell_{1}}}{(1 + \vert y \vert)^{2}} \\ & \quad \leq \sum_{\substack{\{y \in L^{-1}\mathbb{Z}: \, \vert y \vert \geq R - r_{\varepsilon}\}}} \; \sum_{\substack{\{x \in L^{-1}\mathbb{Z}: \, \vert y - x \vert \leq r_{\varepsilon}\}}} \; \mathbb{P}_{x}(X^n(t) = y) \cdot {(1 + \vert y - x \vert)^{2}} \cdot\frac{\vert \vert \eta(y) \vert \vert_{\ell_{1}}}{(1 + \vert y \vert)^{2}} \\ & \quad \leq (2\lceil r_{\varepsilon}L \rceil + 1)(1+r_{\varepsilon})^{2} \sum_{\substack{\{y \in L^{-1}\mathbb{Z}: \, \vert y \vert \geq R - r_{\varepsilon}\}}} \frac{\vert \vert \eta(y) \vert \vert_{\ell_{1}}}{(1 + \vert y \vert)^{2}} .
    \end{aligned}
    \end{equation}
    By Proposition~\ref{Paper01_topological_properties_state_space}(ii), we can pick $R_{\varepsilon} > 0$ sufficiently large that for any $\boldsymbol{\xi}=(\xi(y))_{y\in L^{-1}\mathbb Z}\in \mathcal K$ and any $R \geq R_{\varepsilon}$,
    \begin{equation} \label{Paper01_second_bound_term_limit_spatial_tail_compacts}
        (2\lceil r_{\varepsilon}L \rceil + 1)(1+r_{\varepsilon})^{2}  \sum_{\substack{\{y \in L^{-1}\mathbb{Z}: \, \vert y \vert \geq R - r_{\varepsilon}\}}} \frac{\vert \vert \xi(y) \vert \vert_{\ell_{1}}}{(1 + \vert y \vert)^{2}} \leq \frac{\varepsilon}{2}.
    \end{equation}
    Therefore, by applying~\eqref{Paper01_first_bound_term_limit_spatial_tail_compacts},~\eqref{Paper01_second_bound_term_limit_spatial_tail_compacts_pre_step} and~\eqref{Paper01_second_bound_term_limit_spatial_tail_compacts} to~\eqref{Paper01_intermediate_step_control_spatial_tail}, we conclude that for $\varepsilon > 0$, there exists $R_{\varepsilon} > 0$ such that for any $R > R_{\varepsilon}$,
    \begin{equation*}
        \sup_{n \in \mathbb{N}} \; \sup_{\boldsymbol{\eta} \in \mathcal{K}} \; \sup_{t \in [0,T]} \; \sum_{\substack{\{x \in L^{-1}\mathbb{Z}: \, \vert x \vert \geq R\}}} \; \sum_{y \in L^{-1}\mathbb{Z}} \mathbb{P}_{x}(X^n(t) = y)\frac{\vert \vert \eta(y) \vert \vert_{\ell_{1}}}{(1 + \vert x \vert)^{2}} \leq \varepsilon.
    \end{equation*}
    Since $\varepsilon > 0$ was arbitrary, we conclude that the limit in~\eqref{Paper01_equivalent_formulation_limit_spatial_tails} holds, which completes the proof of~\eqref{Paper01_uniform_control_compacts_spatial_tail}.

    \medskip

    \noindent \underline{Step $(4)$: Tightness of one-dimensional distributions}

    \medskip

    In order to prove that for any $\varepsilon > 0$, $\mathcal K\subset \mathcal S$ compact and $T\ge 0$ there exists a compact subset $\mathscr{K}' \subset \mathcal{S}$ such that~\eqref{Paper01_tightness_one_dimensional_distributions} holds, it is enough to repeat the same argument we used in the proof of the strong compact containment condition in Lemma~\ref{Paper01_further_characterisation_relatively_compactness_discrete_IPS}, with the difference that here we apply estimates~\eqref{Paper01_uniform_control_compacts_total_mass},~\eqref{Paper01_uniform_control_compacts_type_space} and~\eqref{Paper01_uniform_control_compacts_spatial_tail} of this result. Since the proof follows along the same lines, we omit the details.
\end{proof}

We are now ready to move forward to statement~(iv) of our programme, i.e.~to the proof that for any $\phi \in \mathscr{C}^{\textrm{cyl}}_{b,*}(\mathcal{S}, \mathbb{R})$, any $T \geq 0$ and any $\boldsymbol{\eta} \in \mathcal{S}$, $((P^{n}_{T}\phi)(\boldsymbol{\eta}))_{n \in \mathbb{N}}$ is a Cauchy sequence in $\mathbb{R}$. It will be convenient to define, for any $\boldsymbol{\eta}^{(1)}=(\eta^{(1)}_k(x))_{k\in \mathbb N_0,\, x\in L^{-1}\mathbb Z}$, $\boldsymbol{\eta}^{(2)} =(\eta^{(2)}_k(x))_{k\in \mathbb N_0,\, x\in L^{-1}\mathbb Z}\in \mathcal{S}$, the \emph{difference set} between $\boldsymbol{\eta}^{(1)}$ and $\boldsymbol{\eta}^{(2)}$ by letting
\begin{equation} \label{Paper01_difference_set_reference_modified}
    \Delta(\boldsymbol{\eta}^{(1)}, \boldsymbol{\eta}^{(2)}) \defeq \left\{(x,k) \in L^{-1}\mathbb{Z} \times \mathbb{N}_{0}: \; \eta^{(1)}_{k}(x) \neq \eta_{k}^{(2)}(x)\right\}.
\end{equation}
It will also be convenient to define the \emph{spatial projection of the difference set} by
\begin{equation} \label{Paper01_spatial_projection_difference_set_reference_modified}
    \Delta_{\textrm{sp}}(\boldsymbol{\eta}^{(1)}, \boldsymbol{\eta}^{(2)}) \defeq \left\{x \in L^{-1}\mathbb{Z}: \;\exists k \in \mathbb{N}_{0} \textrm{ such that } (x,k) \in  \Delta(\boldsymbol{\eta}^{(1)}, \boldsymbol{\eta}^{(2)})\right\}.
\end{equation}
Recall that as defined at the end of Section~\ref{Paper01_introduction}, $\# \Delta(\boldsymbol{\eta}^{(1)}, \boldsymbol{\eta}^{(2)})$ denotes the number of elements of $ \Delta(\boldsymbol{\eta}^{(1)}, \boldsymbol{\eta}^{(2)})$. In order to prove statement~(iv) of our programme, we will need to construct, for each $n \in \mathbb{N}$, realisations of the Markov process associated to $\mathcal{L}^{n}$ starting from different configurations on the same probability space in such a way that after some finite period of time, the realisations agree with each other with high probability on a large spatial box around the origin. More concretely, we will need to prove the following result.

\begin{proposition}[Bounds on the spread of infection] \label{Paper01_proposition_bound_spread_infection}
    For any $r > 0$ and $T > 0$, there exist $C_{1} = C_{1}(T,m,L,N,r) > 0$, $C_{2} = C_{2}(T,m,L,N)>0$ and $R_{c}  = R_{c}(T,r,m,L) > 0$ such that the following holds. For any $n\in \mathbb N$, for any $R > R_{c}$, and for any $\boldsymbol{\eta}^{(1)}, \boldsymbol{\eta}^{(2)} \in \mathcal{S}$ with $\# \Delta(\boldsymbol{\eta}^{(1)}, \boldsymbol{\eta}^{(2)}) < \infty$ and $$ \Delta_{\textrm{sp}}(\boldsymbol{\eta}^{(1)}, \boldsymbol{\eta}^{(2)}) \subset (- \infty, - R] \, \cup \, [R, \infty),$$
    for any $n \in \mathbb{N}$, realisations $(\eta^{n,(1)}(t))_{t \geq 0}$ and $(\eta^{n,(2)}(t))_{t \geq 0}$ of the Markov process associated to the infinitesimal generator $\mathcal{L}^{n}$ can be constructed on the same probability space in such a way that $\eta^{n,(1)}(0) = \boldsymbol{\eta}^{(1)}$ and $\eta^{n,(2)}(0) = \boldsymbol{\eta}^{(2)}$ almost surely, and
    \begin{equation} \label{eq:boundspreadinfect}
        \sup_{t \in [0,T]} \; \mathbb{P}\Big(\Delta_{\textrm{sp}}(\eta^{n,(1)}(t), \eta^{n,(2)}(t))\,  \cap  \, [-r,r] \neq \emptyset \Big) \leq C_{1}e^{-C_2 R}\sum_{(x,k) \in  \Delta(\boldsymbol{\eta}^{(1)}, \boldsymbol{\eta}^{(2)})} \Big\vert \eta^{(1)}_{k}(x) - \eta^{(2)}_{k}(x) \Big\vert .
    \end{equation}
\end{proposition}

Proposition~\ref{Paper01_proposition_bound_spread_infection} will be proved in Section~\ref{Paper01_section_spread_infection} by constructing a particle system that encodes the difference between the realisations in terms of infected and partially recovered particles. Our next result is a direct consequence of Proposition~\ref{Paper01_proposition_bound_spread_infection}. For any $r > 0$ and $K \in \mathbb{N}$, let 
\begin{equation} \label{Paper01_definition_cylindrical_functions}
\begin{aligned}
    & \mathscr{C}^{\textrm{cyl}}_{b,r,K}(\mathcal{S}, \mathbb{R}) \\ & \quad \defeq \Big\{\phi \in \mathscr{C}^{\textrm{cyl}}_{b,*}(\mathcal{S}, \mathbb{R}): \; \exists J \in \mathbb{N}, \, (x_{1}, \ldots, x_{J}) \in (L^{-1}\mathbb{Z} \cap [-r,r])^{J}, \, (k_{1}, \ldots, k_{J}) \in \left([K]_{0}\right)^{J}  \\ & \quad \quad \quad \quad \quad \quad \quad \quad \quad\quad \quad \textrm{and } \varphi \in \mathscr{C}_b(\mathbb{N}_0^{J}, \mathbb{R}) \textrm{ s.t. } \phi(\boldsymbol{\eta}) = \varphi(\eta_{k_{1}}(x_{1}), \ldots, \eta_{k_{J}}(x_{J}))\; \forall \boldsymbol{\eta} \in \mathcal{S}\Big\}.
\end{aligned}
\end{equation}
Observe that by~\eqref{Paper01_general_definition_cylindrical_function},
\begin{equation*}
\mathscr{C}^{\textrm{cyl}}_{b,*}(\mathcal{S}, \mathbb{R}) = \bigcup_{r > 0} \, \bigcup_{K \in \mathbb{N}} \, \mathscr{C}^{\textrm{cyl}}_{b,r,K}(\mathcal{S}, \mathbb{R}).
\end{equation*}

\begin{corollary} \label{Paper01_control_distance_expectation_cylindrical_functions}
    For any $r > 0$ and $T > 0$, there exist $C_{1} = C_{1}(T,m,L,N,r) > 0$, $C_{2} = C_{2}(T,m,L,N)$ and $R_{c}  = R_{c}(T,r,m,L) > 0$ such that the following holds. For any $R > R_{c}$, any $\boldsymbol{\eta}^{(1)}, \boldsymbol{\eta}^{(2)} \in \mathcal{S}$ with $\# \Delta(\boldsymbol{\eta}^{(1)}, \boldsymbol{\eta}^{(2)}) < \infty$ and $ \Delta_{\textrm{sp}}(\boldsymbol{\eta}^{(1)}, \boldsymbol{\eta}^{(2)}) \subset (- \infty, - R] \, \cup \, [R, \infty)$, any $K \in \mathbb{N}$ and any $\phi \in \mathscr{C}_{b,r,K}^{\textrm{cyl}}(\mathcal{S}, \mathbb{R})$,
    \begin{equation*}
        \sup_{n \in \mathbb{N}} \; \sup_{t \in [0,T]} \; \vert (P^{n}_{t}\phi)(\boldsymbol{\eta}^{(1)}) - (P^{n}_{t}\phi)(\boldsymbol{\eta}^{(2)}) \vert \leq 2 C_{1}e^{-C_2 R}\vert \vert \phi \vert \vert_{L_{\infty}(\mathcal{S}; \mathbb{R})} \sum_{(x,k) \in  \Delta(\boldsymbol{\eta}^{(1)}, \boldsymbol{\eta}^{(2)})} \Big\vert \eta^{(1)}_{k}(x) - \eta^{(2)}_{k}(x) \Big\vert .
    \end{equation*}
\end{corollary}

\begin{proof}
Take $n\in \mathbb N$.
    Using the construction in the statement of Proposition~\ref{Paper01_proposition_bound_spread_infection}, let $(\eta^{n,(1)}(t))_{t \geq 0}$ and $(\eta^{n,(2)}(t))_{t \geq 0}$ be realisations of the Markov process associated to the infinitesimal generator $\mathcal{L}^{n}$ on the same probability space with $\eta^{n,(1)}(0) = \boldsymbol{\eta}^{(1)}$ and $\eta^{n,(2)}(0) = \boldsymbol{\eta}^{(2)}$ and such that~\eqref{eq:boundspreadinfect} holds. Observe that, for any $K>0$, any $\phi \in \mathscr{C}_{b,r,K}^{\textrm{cyl}}(\mathcal{S}, \mathbb{R})$ and any $t\in [0,T]$, on the event 
    \begin{equation*}
        \mathscr{A}^{n} \defeq \left\{\Delta_{\textrm{sp}}(\eta^{n,(1)}(t), \eta^{n,(2)}(t))\,  \cap  \, [-r,r] = \emptyset\right\}
    \end{equation*}
    we have $\phi(\eta^{n,(1)}(t)) = \phi(\eta^{n,(2)}(t))$. Therefore,
    \begin{equation*}
    \begin{aligned}
        \vert (P^{n}_{t}\phi)(\boldsymbol{\eta}^{(1)}) - (P^{n}_{t}\phi)(\boldsymbol{\eta}^{(2)}) \vert \leq 2 \vert \vert \phi \vert \vert_{L_{\infty}(\mathcal{S}; \mathbb{R})} (1 -  \mathbb{P}(\mathscr{A}^{n})).
    \end{aligned}
    \end{equation*}
    The result then follows from Proposition~\ref{Paper01_proposition_bound_spread_infection}.
\end{proof}

We are now ready to finish the proof of statement~(iv) of our programme.

\begin{proposition} \label{Paper01_convergence_one_dimensional_distributions_test_function-muller_ratchet}
    For any $T \geq 0$, any $\phi \in \mathscr{C}_{b,*}^{\textrm{cyl}}(\mathcal{S}, \mathbb{R})$, any compact set $\mathcal{K} \subset \mathcal{S}$ and any $\varepsilon > 0$, there exists $n(T, \phi, \mathcal{K}, \varepsilon) \in \mathbb{N}$ such that
        \begin{equation*}
            \sup_{n_{1}, n_{2} \geq n(T, \phi, \mathcal{K}, \varepsilon)} \; \sup_{\boldsymbol{\eta} \in \mathcal{K}} \; \sup_{t \in [0,T]} \; \left\vert (P^{n_{1}}_{t}\phi)(\boldsymbol{\eta}) - (P^{n_{2}}_{t}\phi)(\boldsymbol{\eta}) \right\vert \leq \varepsilon.
        \end{equation*}
\end{proposition}

\begin{proof}
    Fix $\mathcal{K} \subset \mathcal{S}$ compact, $T \geq 0$, $r>0$, $K'\in \mathbb N$ and $\phi \in \mathscr{C}^{\textrm{cyl}}_{b,r,K'}(\mathcal{S}, \mathbb{R})$. By Proposition~\ref{Paper01_bound_total_mass}, the integration by parts formula stated in~\cite[Lemma~1.6.2]{ethier2009markov} and standard results about (strongly) Feller process taking values in $(\mathbb{N}_0)^d$ (see~\cite[Theorem~6.4.1]{ethier2009markov}), we have that for any $n_{1}, n_{2} \in \mathbb{N}$, $\boldsymbol{\eta} \in \mathcal{K}$ and $t \in [0,T]$,
    \begin{equation} \label{Paper01_integration_parts_formula_semigroup}
    \begin{aligned}
        (P^{n_{1}}_{t}\phi)(\boldsymbol{\eta}) -  (P^{n_{2}}_{t}\phi)(\boldsymbol{\eta}) & = \int_{0}^{t} \Big(P^{n_{1}}_{t-t'}(\mathcal{L}^{n_{1}} - \mathcal{L}^{n_{2}})P^{n_{2}}_{t'} \phi\Big)(\boldsymbol{\eta}) \, dt' \\ & = \int_{0}^{t} \mathbb{E}_{\boldsymbol{\eta}}\Big[\left(\mathcal{L}^{n_{1}} - \mathcal{L}^{n_{2}}\right)(P^{n_{2}}_{t'} \phi)(\eta^{n_{1}}(t - t'))\Big] \, dt'.
    \end{aligned}
    \end{equation}
   Recall the sequences of birth and death rates $(F^{b}_{k})_{k \in \mathbb{N}_{0}}$ and $(F^{d}_{k})_{k \in \mathbb{N}_{0}}$ given by~\eqref{Paper01_discrete_birth_rate} and~\eqref{Paper01_discrete_death_rate} respectively, and recall the definition of $\boldsymbol{e}^{(x)}_{k}$ from before~\eqref{Paper01_scaled_polynomials_carrying_capacity}. 
   Recall from after~\eqref{Paper01_scaling_parameters_for_discrete_approximation} that $\Lambda_n = L^{-1}\mathbb{Z} \cap [- \lambda_n, \lambda_n]$, where $(\lambda_n)_{n \in \mathbb{N}}$ is increasing with $\lim_{n \rightarrow \infty} \lambda_n = \infty$, and that $(K_n)_{n\in \mathbb N}$ is increasing with $\lim_{n\to \infty}K_n=\infty$. 
   By the definition of the sequence of generators $(\mathcal{L}^{n})_{n \in \mathbb{N}}$ in~\eqref{Paper01_generator_foutel_etheridge_model_restriction_n}, note that for any $n_1,n_2\in \mathbb N$ with $n_1>n_2$, $\boldsymbol{\xi}=(\xi_j(x))_{j\in \mathbb N_0,\, x\in L^{-1}\mathbb Z} \in \mathcal{S}$ and $0\le t' \le t$,
  \begin{align} \label{Paper01_first_step_proof_convergence_one_dimensional_distributions}
         & \left(\mathcal{L}^{n_{1}} - \mathcal{L}^{n_{2}}\right)(P^{n_{2}}_{t'} \phi)(\boldsymbol{\xi}) \notag \\ 
         & \quad = \frac{m}{2} \sum_{x \in \Lambda_{n_1} \setminus \Lambda_{n_2}} \; \sum_{j = 0}^{\infty}\sum_{z\in \{-L^{-1},L^{-1}\}}\mathds 1_{\{x+z\in \Lambda_{n_1}\}} \xi_{j}(x) \Big(P^{n_{2}}_{t'} \phi\Big(\boldsymbol{\xi} + \boldsymbol{e}^{(x + z)}_{j} - \boldsymbol{e}^{(x)}_{j}\Big)  - P^{n_{2}}_{t'} \phi(\boldsymbol{\xi} )\Big) \notag \\ 
         & \quad \quad + \frac{m}{2} \sum_{x \in \Lambda_{n_2}} \; \sum_{j = 0}^{\infty}\sum_{z\in \{-L^{-1},L^{-1}\}}\mathds 1_{\{x+z\notin \Lambda_{n_2}\}} \xi_{j}(x) \Big(P^{n_{2}}_{t'} \phi\Big(\boldsymbol{\xi} + \boldsymbol{e}^{(x + z)}_{j} - \boldsymbol{e}^{(x)}_{j}\Big)  - P^{n_{2}}_{t'} \phi(\boldsymbol{\xi} )\Big) \notag \\ 
         &\quad \quad + \sum_{x \in \Lambda_{n_{2}}} \; \sum_{j = K_{n_{2}} + 1}^{K_{n_{1}}} F^{b}_{j}(\xi(x))\Big(P^{n_{2}}_{t'} \phi\Big(\boldsymbol{\xi} + \boldsymbol{e}^{(x)}_{j}\Big) - P^{n_{2}}_{t'} \phi(\boldsymbol{\xi})\Big) \notag \\ 
         & \quad \quad + \sum_{\substack{x \in \Lambda_{n_{1}} \setminus \Lambda_{n_{2}}}} \; \sum_{j = 0}^{K_{n_{1}}} F^{b}_{j}(\xi( x))\Big(P^{n_{2}}_{t'} \phi\Big(\boldsymbol{\xi} + \boldsymbol{e}^{(x)}_{j}\Big) - P^{n_{2}}_{t'} \phi(\boldsymbol{\xi})\Big) \notag \\ 
         & \quad \quad + \sum_{\substack{x \in \Lambda_{n_{1}} \setminus \Lambda_{n_{2}}}} \; \sum_{j = 0}^{\infty} F^{d}_{j}(\xi( x))\Big(P^{n_{2}}_{t'} \phi\Big(\boldsymbol{\xi} - \boldsymbol{e}^{(x)}_{j}\Big) - P^{n_{2}}_{t'} \phi(\boldsymbol{\xi})\Big).
   \end{align}
   Then, by applying Corollary~\ref{Paper01_control_distance_expectation_cylindrical_functions} to~\eqref{Paper01_first_step_proof_convergence_one_dimensional_distributions}, we conclude that there exist
$R_c=R_c(T,r,m,L,N)>0$, $C_1=C_1(T,m,L,N,r)>0$, $C_2=C_2(T,m,L,N)>0$
such that for $R>R_c$, for $n_{2}$ sufficiently large that $\lambda_{n_2}-L^{-1} > R>R_c$ and $n_1>n_2$, for $0\le t'\le t\le T$ and $\boldsymbol{\xi}=(\xi_j(x))_{j\in \mathbb N_0,\, x\in L^{-1}\mathbb Z}\in \mathcal S$,
   \begin{equation} \label{Paper01_second_step_proof_convergence_one_dimensional_distributions}
   \begin{aligned}
           & \Big\vert \left(\mathcal{L}^{n_{1}} - \mathcal{L}^{n_{2}}\right)(P^{n_{2}}_{t'} \phi)(\boldsymbol{\xi})\Big\vert \\ 
            & \quad \le 4 m C_1 \vert \vert \phi \vert \vert_{L_{\infty}(\mathcal{S}; \mathbb{R})} \sum_{\substack{\{x \in L^{-1}\mathbb{Z}: \, \lambda_{n_{2}}-L^{-1} < \vert x \vert \leq \lambda_{n_{1}}}\}} \; \sum_{j = 0}^{\infty} \xi_{j}(x) e^{-C_2 (\vert x \vert-L^{-1})} \\
           & \quad \quad + 2 \vert \vert \phi \vert \vert_{L_{\infty}(\mathcal{S}; \mathbb{R})} \sum_{\substack{\{x \in L^{-1}\mathbb{Z}: \, \vert x \vert \leq R}\}} \; \sum_{j = K_{n_{2}} + 1}^{K_{n_{1}}}F^{b}_{j}(\xi(x)) \\ 
           & \quad \quad \; + 2 C_1 \vert \vert \phi \vert \vert_{L_{\infty}(\mathcal{S}; \mathbb{R})} \sum_{\substack{\{x \in L^{-1}\mathbb{Z}: \, R < \vert x \vert \leq \lambda_{n_{2}}}\}} \; \sum_{j = K_{n_{2}} + 1}^{K_{n_{1}}}F^{b}_{j}(\xi(x)) e^{-C_2 \vert x \vert} \\ 
           & \quad \quad \; + 2 C_1 \vert \vert \phi \vert \vert_{L_{\infty}(\mathcal{S}; \mathbb{R})} \sum_{\substack{\{x \in L^{-1}\mathbb{Z}: \, \lambda_{n_{2}} < \vert x \vert \leq \lambda_{n_{1}}}\}} \; e^{-C_2 \vert x \vert}\Bigg(\sum_{j = 0}^{K_{n_{1}}} F^{b}_{j}(\xi(x)) + \sum_{j = 0}^{\infty} F^{d}_{j}(\xi(x)) \Bigg) \\ 
           & \quad \leq 2 s_{K_{n_2}} \vert \vert \phi \vert \vert_{L_{\infty}(\mathcal{S}; \mathbb{R})}  \sum_{\substack{\{x \in L^{-1}\mathbb{Z}: \, \vert x \vert \leq R\}}}  \vert \vert \xi(x) \vert \vert_{\ell_1}q^{N}_{+}(\vert \vert \xi(x) \vert \vert_{\ell_1})
           \\ 
           & \quad \quad \; + 2 C_1 \vert \vert \phi \vert \vert_{L_{\infty}(\mathcal{S}; \mathbb{R})} \sum_{\substack{\{x \in L^{-1}\mathbb{Z}: \,  \vert x \vert > R\}}} \Big(2e^{C_2L^{-1}} m\vert \vert \xi(x) \vert \vert_{\ell_{1}} + \bar{q}^{N}(\vert \vert \xi(x) \vert \vert_{\ell_{1}})\Big) e^{-C_2 \vert x \vert},
   \end{aligned}
   \end{equation}
   where in the last inequality we used~\eqref{Paper01_bound_birth_rate_high_mutation_load} and~\eqref{Paper01_simple_bound_total_birth_death_rates_deterministic}; recall that $\bar{q}^{N}$ is the polynomial defined in~\eqref{Paper01_bar_polynomial}. By~\eqref{Paper01_definition_state_space_formal}, for any $\boldsymbol{\xi} \in \mathcal{S}$, the map $x \mapsto \vert \vert \xi(x) \vert \vert_{\ell_1}$ grows at most polynomially in $|x|$. Hence, the right-hand side of~\eqref{Paper01_second_step_proof_convergence_one_dimensional_distributions} is finite for any $\boldsymbol{\xi} \in \mathcal{S}$.
   
   Now fix $\varepsilon > 0$, and take $R_{\varepsilon} = R_\varepsilon (r,\phi, \mathcal{K}, T) > R_c$ to be specified later. For $n_{2}$ sufficiently large that $\lambda_{n_{2}} > R_{\varepsilon}+L^{-1}$, by applying~\eqref{Paper01_second_step_proof_convergence_one_dimensional_distributions}, and then Fubini's theorem, to~\eqref{Paper01_integration_parts_formula_semigroup}, we obtain that for $n_1>n_2$ and $t\in [0,T]$,
   \begin{equation} \label{Paper01_third_step_proof_convergence_one_dimensional_distributions}
    \begin{aligned}
        & \Big\vert (P^{n_{1}}_{t}\phi)(\boldsymbol{\eta}) -  (P^{n_{2}}_{t}\phi)(\boldsymbol{\eta}) \Big\vert \\
        & \quad \leq 2 s_{K_{n_2}} \vert \vert \phi \vert \vert_{L_{\infty}(\mathcal{S}; \mathbb{R})}  \int_{0}^{t} \sum_{\substack{\{x \in L^{-1}\mathbb{Z}: \, \vert x \vert \leq R_{\varepsilon}\}}} \;  \mathbb{E}_{\boldsymbol{\eta}}\Big[\vert \vert \eta^{n_1}(t - t', x) \vert \vert_{\ell_1}q^{N}_{+}(\vert \vert \eta^{n_1}(t-t',x) \vert \vert_{\ell_1})\Big] \, dt'\\ & \quad \quad \;  + 2 C_1 \vert \vert \phi \vert \vert_{L_{\infty}(\mathcal{S}; \mathbb{R})} \int_{0}^{t} \sum_{\substack{\{x \in L^{-1}\mathbb{Z}: \,  \vert x \vert > R_\varepsilon\}}} e^{-C_2 \vert x \vert} \\ & \quad \quad \quad \quad \quad \quad \quad \quad \quad \quad \quad \quad \quad \quad \cdot \mathbb{E}_{\boldsymbol{\eta}}\Big[2e^{C_2 L^{-1}} m\vert \vert \eta^{n_1}(t-t', x) \vert \vert_{\ell_{1}} + \bar{q}^{N}(\vert \vert \eta^{n_1}(t-t',x) \vert \vert_{\ell_{1}})\Big]  dt'.
    \end{aligned}
   \end{equation}
   We now bound the terms on the right-hand side of~\eqref{Paper01_third_step_proof_convergence_one_dimensional_distributions} separately. For the second term, since $0 \leq \deg q_{+} < \deg q_-$ by Assumption~\ref{Paper01_assumption_polynomials}, 
   and using the definition of $\bar{q}^{N}$ in~\eqref{Paper01_bar_polynomial},
   by~\eqref{Paper01_uniform_control_compacts_total_mass} in Proposition~\ref{Paper01_estimates_holding_uniformly_compact_subsets_S_N}
   combined with the definition of $\psi^{(p)}$ in~\eqref{Paper01_ell_p_norm_S_N},
   there exists $c^{(1)}_{q_+,q_-,N,\mathcal{K},T,m,L} > 0$ such that for $n_1>n_2$,
    \begin{align} \label{eq:supetasuptbound}
        & \sup_{\boldsymbol{\eta} \in \mathcal{K}} \; \sup_{t \in [0,T]} \int_{0}^{t} \sum_{\substack{\{x \in L^{-1}\mathbb{Z}: \,  \vert x \vert > R_\varepsilon\}}} e^{-C_2 \vert x \vert} \mathbb{E}_{\boldsymbol{\eta}}\Big[2e^{C_2 L^{-1}} m\vert \vert \eta^{n_1}(t-t', x) \vert \vert_{\ell_{1}} + \bar{q}^{N}(\vert \vert \eta^{n_1}(t-t',x) \vert \vert_{\ell_{1}})\Big] dt' \notag \\ & \quad \leq c^{(1)}_{q_+,q_-,N,\mathcal{K},T,m,L} \int_{0}^{T} \sum_{\substack{\{x \in L^{-1}\mathbb{Z}: \,  \vert x \vert > R_\varepsilon\}}} e^{-C_2 \vert x \vert} \sum_{p = 0}^{1 + \deg q_-} (1 + \vert x \vert)^{2p} dt'.
    \end{align}
   We can then take $R_\varepsilon > 0$ sufficiently large that
   \begin{equation} \label{Paper01_third_step_proof_convergence_one_dimensional_distributions_i}
        c^{(1)}_{q_+,q_-,N,\mathcal{K},T,m,L} \int_{0}^{T} \sum_{\substack{\{x \in L^{-1}\mathbb{Z}: \,  \vert x \vert > R_\varepsilon\}}} e^{-C_2 \vert x \vert} \sum_{p = 0}^{1 + \deg q_-} (1 + \vert x \vert)^{2p} dt' \leq (2C_1 \vert \vert \phi \vert \vert_{L_{\infty}(\mathcal{S}; \mathbb{R})})^{-1} \frac{\varepsilon}{2}.
   \end{equation}
   Fix $R_\varepsilon > R_c$ such that~\eqref{Paper01_third_step_proof_convergence_one_dimensional_distributions_i} holds. Then, for the first term on the right-hand side of~\eqref{Paper01_third_step_proof_convergence_one_dimensional_distributions}, we observe that $$\# (L^{-1}\mathbb{Z} \cap [-R_\varepsilon, R_\varepsilon]) \leq 2\lfloor L R_\varepsilon \rfloor + 1.$$
   Hence, since $q^{N}_{+}$ is a polynomial by~\eqref{Paper01_scaled_polynomials_carrying_capacity}, 
   and again using~\eqref{Paper01_uniform_control_compacts_total_mass} in Proposition~\ref{Paper01_estimates_holding_uniformly_compact_subsets_S_N}
   combined with~\eqref{Paper01_ell_p_norm_S_N},
   there exists $c^{(2)}_{q_+,N,\mathcal{K},T} > 0$ such that
   \begin{equation*}
   \begin{aligned}
       &  \sup_{\boldsymbol{\eta} \in \mathcal{K}} \; \sup_{t \in [0,T]} \int_{0}^{t} \sum_{\substack{\{x \in L^{-1}\mathbb{Z}: \, \vert x \vert \leq R_{\varepsilon}\}}} \;  \mathbb{E}_{\boldsymbol{\eta}}\Big[\vert \vert \eta^{n_1}(t - t', x) \vert \vert_{\ell_1}q^{N}_{+}(\vert \vert \eta^{n_1}(t-t',x) \vert \vert_{\ell_1})\Big] \, dt' \\ & \quad \leq c^{(2)}_{q_+,N,\mathcal{K},T} T (2\lfloor L R_\varepsilon \rfloor + 1) \sum_{p = 0}^{1 + \deg q_+} (1 + R_\varepsilon)^{2p}.
   \end{aligned}
   \end{equation*}
   Since $s_{k} \rightarrow 0$ as $k \rightarrow \infty$ by Assumption~\ref{Paper01_assumption_fitness_sequence}, and since $K_{n} \rightarrow \infty$ and $\lambda_n\to \infty$ as $n \rightarrow \infty$ by~\eqref{Paper01_scaling_parameters_for_discrete_approximation}, there exists $n(T,r,\phi, \mathcal{K}, \varepsilon) \in \mathbb{N}$ such that for $n_{1} > n_{2} \geq n(T,r,\phi, \mathcal{K}, \varepsilon)$ we have $\lambda_{n_{2}} > R_{\varepsilon}+L^{-1}$ and
   \begin{equation} \label{Paper01_third_step_proof_convergence_one_dimensional_distributions_ii}
   \begin{aligned}
        & \sup_{\boldsymbol{\eta} \in \mathcal{K}} \, \sup_{t \in [0,T]} 2 s_{K_{n_2}} \vert \vert \phi \vert \vert_{L_{\infty}(\mathcal{S}; \mathbb{R})}  \int_{0}^{t} \sum_{\substack{\{x \in L^{-1}\mathbb{Z}: \, \vert x \vert \leq R_{\varepsilon}\}}} \;  \mathbb{E}_{\boldsymbol{\eta}}\Big[\vert \vert \eta^{n_1}(t - t', x) \vert \vert_{\ell_1}q^{N}_{+}(\vert \vert \eta^{n_1}(t-t',x) \vert \vert_{\ell_1})\Big] \, dt' \\ & \quad \leq \frac{\varepsilon}{2}.
   \end{aligned}
   \end{equation}
   Thus, the proof is completed by applying~\eqref{eq:supetasuptbound},~\eqref{Paper01_third_step_proof_convergence_one_dimensional_distributions_i} and~\eqref{Paper01_third_step_proof_convergence_one_dimensional_distributions_ii} to~\eqref{Paper01_third_step_proof_convergence_one_dimensional_distributions}. 
\end{proof}

We are finally ready to prove Theorem~\ref{Paper01_thm_existence_uniqueness_IPS_spatial_muller}.

\begin{proof}[Proof of Theorem~\ref{Paper01_thm_existence_uniqueness_IPS_spatial_muller}]
    By Proposition~\ref{Paper01_bound_total_mass} and Corollary~\ref{cor:cylinderfnsassumptions},
    the set $\mathscr{C}_{b,*}^{\textrm{cyl}}(\mathcal{S}, \mathbb{R})$, the generator $\mathcal L$ and the sequence of Markov processes $(\eta^n)_{n\in \mathbb N}$ satisfy
    Assumptions~\ref{Paper01_assumption_candidate_set_functions_domain_generator} and~\ref{Paper01_assumption_tightness_sequence_processes}(i). 
Moreover, the set $\mathscr{C}_{b,*}^{\textrm{cyl}}(\mathcal{S}, \mathbb{R})$ and the sequence
$(\eta^{n})_{n \in \mathbb{N}}$ satisfy Assumption~\ref{Paper01_assumption_tightness_sequence_processes}(ii) by Proposition~\ref{Paper01_existence_solution_martingale_problem}, Assumption~\ref{Paper01_assumption_tightness_sequence_processes}(iii) by estimate~\eqref{Paper01_tightness_one_dimensional_distributions} of Proposition~\ref{Paper01_estimates_holding_uniformly_compact_subsets_S_N}, and Assumption~\ref{Paper01_assumption_tightness_sequence_processes}(iv) by Proposition~\ref{Paper01_convergence_one_dimensional_distributions_test_function-muller_ratchet}.
    
    Therefore, Theorem~\ref{Paper01_general_thm_convergence_feller_non_locally_compact} implies
    the claimed properties of the semigroup $\{P_t\}_{t\ge 0}$ given by
    $P_T \phi(\boldsymbol{\eta}):=\lim_{n\to \infty}(P^n_T \phi)(\boldsymbol{\eta})$, and that 
    for $\boldsymbol{\eta}\in \mathcal S_0$, conditioning on $\eta^n(0)=\boldsymbol{\eta}$ for every $n\in \mathbb N$,
    the sequence of Markov processes $(\eta^{n})_{n \in \mathbb{N}}$ converges weakly with respect to the $J_1$-topology on $\mathscr{D}([0, \infty),\mathcal{S})$ as $n \rightarrow \infty$ to a càdlàg Markov process $(\eta(t))_{t \geq 0}$ with $\eta(0)=\boldsymbol{\eta}$ almost surely and with finite dimensional distributions given by $\{P_t\}_{t\ge 0}$, and $(\eta(t))_{t \geq 0}$ is a strong Markov process with respect to its right-continuous natural filtration. Observe that the limiting semigroup $\{P_t\}_{t \geq 0}$ does not depend on the particular choice of sequences $(\lambda_n)_{n \in \mathbb{N}}$ and $(K_n)_{n \in \mathbb{N}}$ satisfying~\eqref{Paper01_scaling_parameters_for_discrete_approximation}. Indeed, our convergence result implies that for any two such choices, one may interleave the corresponding sequences and still obtain the same limit. To complete our argument, we observe that Proposition~\ref{Paper01_existence_solution_martingale_problem} implies that $(\eta(t))_{t \geq 0}$ solves the desired martingale problem.
\end{proof}

As explained in Section~\ref{Paper01_model_description}, we call the process $(\eta(t))_{t \geq 0}$ the spatial Muller's ratchet. We now prove Theorem~\ref{Paper01_thm_moments_estimates_spatial_muller_ratchet} regarding moment bounds for the spatial Muller's ratchet.

\begin{proof}[Proof of Theorem~\ref{Paper01_thm_moments_estimates_spatial_muller_ratchet}]
    Take $\boldsymbol{\eta} \in \mathcal{S}_{0}$, and let the sequence of Markov processes $\left((\eta^{n}(t))_{t \geq 0}\right)_{n \in \mathbb{N}}$ be such that for every $n \in \mathbb{N}$, $(\eta^{n}(t))_{t \geq 0}$ is associated to the infinitesimal generator $\mathcal{L}^{n}$ and $\eta^{n}(0) = \boldsymbol{\eta}$ almost surely. By Theorem~\ref{Paper01_thm_existence_uniqueness_IPS_spatial_muller}, $\left((\eta^{n}(t))_{t \geq 0}\right)_{n \in \mathbb{N}}$ converges weakly with respect to the $J_1$-topology on $\mathscr{D}([0, \infty),\mathcal{S})$ as $n \rightarrow \infty$ to the spatial Muller's ratchet $(\eta(t))_{t \geq 0}$. Since $(\mathcal{S}, d_{\mathcal{S}})$ is a complete and separable metric space, by Skorokhod's representation theorem, it is possible to construct the sequence of processes $\left((\eta^{n}(t))_{t \geq 0}\right)_{n \in \mathbb{N}}$ and $(\eta(t))_{t \geq 0}$ on the same probability space in such a way that the convergence holds almost surely. For any $p \in \mathbb{N}$, any $T \geq 0$ and any $x\in L^{-1}\mathbb Z$, we then have
    \begin{equation} \label{Paper01_intermediate_step_control_moments_real_mullers_ratchet}
         \mathbb{E}\left[\vert \vert \eta(T,x) \vert \vert_{\ell_{1}}^{p}\right]  =  \mathbb{E}\left[\lim_{n \rightarrow \infty} \; \vert \vert \eta^{n}(T,x) \vert \vert_{\ell_{1}}^{p}\right]  \leq \liminf_{n \rightarrow \infty} \; \mathbb{E}\left[ \vert \vert \eta^{n}(T,x) \vert \vert_{\ell_{1}}^{p}\right],
    \end{equation}
    where we used Fatou's lemma for the last inequality. Recall that for double sequences of positive real numbers $(a^{n}(x))_{x \in L^{-1}\mathbb{Z}, n \in \mathbb{N}}$, the following elementary inequality holds:
    \begin{equation} \label{Paper01_interchange_sup_liminf}
        \sup_{x \in L^{-1}\mathbb{Z}} \; \liminf_{n \rightarrow \infty} \; a^{n}(x) \leq \liminf_{n \rightarrow \infty} \;   \sup_{x \in L^{-1}\mathbb{Z}} \; a^{n}(x).
    \end{equation}
    Applying~\eqref{Paper01_interchange_sup_liminf} to~\eqref{Paper01_intermediate_step_control_moments_real_mullers_ratchet}, we get that for $p\in \mathbb N$ and $T\ge 0$,
    \begin{equation*}
         \sup_{x \in L^{-1}\mathbb{Z}} \; \mathbb{E}\left[\vert \vert \eta(T,x) \vert \vert_{\ell_{1}}^{p}\right] \leq \liminf_{n \rightarrow \infty} \; \sup_{x \in L^{-1}\mathbb{Z}} \; \mathbb{E}\left[ \vert \vert \eta^{n}(T,x) \vert \vert_{\ell_{1}}^{p}\right].
    \end{equation*}
    Then, by estimate~\eqref{Paper01_moments_estimates_uniform_space} in Proposition~\ref{Paper01_bound_total_mass}, the result follows.
\end{proof}

\section{Spread of infection in an expanding population} \label{Paper01_section_spread_infection}

In this section, we will prove Proposition~\ref{Paper01_proposition_bound_spread_infection}, i.e.~we will show that for every $n \in \mathbb{N}$ and $\boldsymbol{\eta}^{(1)}, \boldsymbol{\eta}^{(2)} \in \mathcal{S}$ sufficiently close to each other in a suitable sense, it is possible to construct realisations (on the same probability space) of the Markov process associated to $\mathcal{L}^{n}$ starting from the configurations $\boldsymbol{\eta}^{(1)}$ and $\boldsymbol{\eta}^{(2)}$, in such a way that after a fixed time, with high probability, the two processes agree with each other at all demes in a box around the origin. To prove this result, we will encode the evolution of the difference between the two processes in terms of the spread of an infection, in which `susceptible' particles are contained in both processes, and `infected' and `partially recovered' particles are only contained in one process or the other. 

Recall from~\eqref{Paper01_difference_set_reference_modified} that the difference set between $\boldsymbol{\eta}^{(1)} = (\eta^{(1)}_{k}(x))_{x \in L^{-1}\mathbb{Z}, \, k \in \mathbb{N}_{0}} \in \mathcal{S}$ and $\boldsymbol{\eta}^{(2)} = (\eta^{(2)}_{k}(x))_{x \in L^{-1}\mathbb{Z}, \, k \in \mathbb{N}_{0}} \in \mathcal{S}$ is given by
\begin{equation*}
    \Delta(\boldsymbol{\eta}^{(1)}, \boldsymbol{\eta}^{(2)}) \defeq \left\{(x,k) \in L^{-1}\mathbb{Z} \times \mathbb{N}_{0}: \; \eta^{(1)}_{k}(x) \neq \eta^{(2)}_{k}(x)\right\}.
\end{equation*}
Take $n \in \mathbb{N}$ and 
two configurations $\boldsymbol{\eta}^{(1)}, \boldsymbol{\eta}^{(2)} \in \mathcal{S}$ such that $\#  \Delta(\boldsymbol{\eta}^{(1)}, \boldsymbol{\eta}^{(2)})  < \infty$.
We will define a coupling
$\left(\eta^{n,(1)}(t), \eta^{n,(2)}(t) \right)_{t \geq 0}$
in such a way that for each $i\in \{1,2\}$, the process
$(\eta^{n,(i)}(t))_{t\ge 0}$ is
a realisation of the Markov process associated with the infinitesimal generator $\mathcal{L}^{n}$ defined in~\eqref{Paper01_generator_foutel_etheridge_model_restriction_n} and~\eqref{Paper01_infinitesimal_generator_restriction_n}, 
with $\eta^{n,(i)}(0) = \boldsymbol{\eta}^{(i)}$ almost surely.

The coupling will be encoded by giving each particle a `class' in addition to its `type' (the number of mutations that it carries) and spatial location. There will be five different classes of particles:
\begin{itemize}
    \item Class 0 (\emph{susceptible}) particles;
    \item Class 1 and class 2 (\emph{infected}) particles;
    \item Class $1*$ and class $2*$ (\emph{partially recovered}) particles. 
\end{itemize}
The class 0 particles will be shared by both processes $(\eta^{n,(1)}(t))_{t\ge 0}$ and $(\eta^{n,(2)}(t))_{t\ge 0}$, while the class 1 and class $1*$ particles will be exclusive to the process $(\eta^{n,(1)}(t))_{t\ge 0}$, and the class 2 and class $2*$ particles will be exclusive to the process $(\eta^{n,(2)}(t))_{t\ge 0}$. Additionally, each partially recovered class $1*$ (respectively, class $2*$) particle will be uniquely associated with a class $2*$ (respectively, class $1*$) particle, referred to as its \emph{dual particle}, with the same spatial location (but not necessarily with the same number of mutations). A class~$1*$ particle and its class~$2*$ dual particle will be referred to as a \emph{dual pair} of particles.

We will construct the coupling in such a way that class $1*$ and class $2*$ particles generate new infected particles at a lower rate than class $1$ and class $2$ particles. Moreover, the construction will be such that for $i\in \{1,2\}$, whenever a class $i$ particle ‘transmits' the infection to a class 0 particle at a deme with a large number of particles, it partially recovers and becomes a class~$i*$ particle. 

For each $j\in \{0,1,1*,2,2*\}$, $k\in \mathbb N_0$, $t\ge 0$ and $x\in L^{-1}\mathbb Z$, we write
$\eta^{n,j}_k(t,x)$ for the number of particles in the coupling at time $t$ at deme~$x$ with class $j$, and carrying exactly~$k$ mutations.
We write $\eta^{n,j}(t,x)=(\eta^{n,j}_k(t,x))_{k\in \mathbb N_0}$ and $\eta^{n,j}(t)=(\eta^{n,j}(t,x))_{x\in L^{-1}\mathbb Z}$.
For $t \geq 0$, $i\in \{1,2\}$, $k\in \mathbb N_0$ and $x\in L^{-1}\mathbb Z$, we set 
\begin{equation} \label{Paper01_definition_number_particles_composition_reference_modified_processes}
\eta^{n,(i)}_k(t,x)=\eta^{n,0}_k(t,x)+\eta^{n,i}_k(t,x)+\eta^{n,i*}_k(t,x)
\end{equation}
and $\eta^{n,(i)}(t,x)=(\eta^{n,(i)}_k(t,x))_{k\in \mathbb N_0}$. For~$i \in \{1,2\}$, $t \geq 0$ and $x \in L^{-1}\mathbb{Z}$, we define the event that there are more class $i$ particles than class $3-i$ particles at deme $x$ at time $t$ as
\begin{equation} \label{Paper01_event_more_particles_process_i_than_3_minus_i}
    A^{i}(t,x) \defeq \Big\{\vert \vert \eta^{n,i}(t,x) \vert \vert_{\ell_{1}} > \vert \vert \eta^{n,3 - i}(t,x) \vert \vert_{\ell_{1}}\Big\}.
\end{equation}
We will now define the process 
\begin{equation*}
    {\left(\eta^{n,\textrm{(coupling)}}(t)\right)_{t \geq 0} \defeq \left(\eta^{n,0}(t), \eta^{n,1}(t), \eta^{n,2}(t),\eta^{n,1*}(t),\eta^{n,2*}(t)\right)_{t \geq 0}}
\end{equation*}
as an $\mathcal{S}^{5}$-valued càdlàg Markov process in such a way that almost surely, for all $x \in L^{-1}\mathbb{Z}$ and $t \geq 0$,
\begin{equation} \label{Paper01_equations_colouring_pink_purple_processes}
    \vert \vert \eta^{{n},1*}(t,x) \vert \vert_{\ell_{1}} = \vert \vert \eta^{{n},2*}(t,x) \vert \vert_{\ell_{1}},
\end{equation}
because each class~$1*$ particle has a unique dual pair class~$2*$ particle at the same deme, and vice versa.
We now begin to define the coupling. First we assign classes to particles in the initial configurations $\boldsymbol{\eta}^{(1)}$ and $\boldsymbol{\eta}^{(2)}$ as follows:
\begin{itemize}
    \item For each $(x,k) \in \Delta(\boldsymbol{\eta}^{(1)}, \boldsymbol{\eta}^{(2)})$ with $\eta_{k}^{(1)}(x) > \eta_{k}^{(2)}(x)$, we assign class $1$ to exactly $\eta_{k}^{(1)}(x) - \eta_{k}^{(2)}(x)$ of the particles in $\boldsymbol{\eta}^{(1)}$ at deme $x$ carrying $k$ mutations. All other particles in $\boldsymbol{\eta}^{(1)}$ are assigned class 0.
    \item Similarly, for each $(x,k) \in \Delta(\boldsymbol{\eta}^{(1)}, \boldsymbol{\eta}^{(2)})$ with $\eta_{k}^{(2)}(x) > \eta_{k}^{(1)}(x)$, we assign class $2$ to exactly $\eta_{k}^{(2)}(x) - \eta_{k}^{(1)}(x)$ of the particles in $\boldsymbol{\eta}^{(2)}$ at deme $x$ carrying $k$ mutations. All other particles in $\boldsymbol{\eta}^{(2)}$ are assigned class 0.
\end{itemize}
In other words, for all $x \in L^{-1}\mathbb{Z}$ and $k \in \mathbb{N}_{0}$,
we set
\begin{equation} \label{Paper01_initial_condition_coloured_process_infection}
\left\{\begin{array}{ccl}
   \eta^{n,0}_{k}(0,x) & = & \eta_{k}^{(1)}(x) \wedge \eta_{k}^{(2)}(x), \\
    \eta^{n,1}_{k}(0,x) & = & \left(\eta_{k}^{(1)}(x) - \eta_{k}^{(2)}(x)\right)^{+}, \\  
    \eta^{n,2}_{k}(0,x) & = & \left(\eta_{k}^{(2)}(x) - \eta_{k}^{(1)}(x)\right)^{+}, \\ 
    \eta^{n,1*}_{k}(0,x) & = & 0, \\ 
    \eta^{n,2*}_{k}(0,x) & = & 0.
\end{array}\right.
\end{equation}
The coupling process $\left(\eta^{n,\textrm{(coupling)}}(t)\right)_{t \geq 0}$ then evolves as follows.
\begin{itemize}
\item \textit{Migration of particles.}
\begin{itemize}
    \item \textit{Susceptible and infected particles}: For each $t \geq 0$ and $x \in \Lambda_n$, each class $0$, $1$ or $2$ particle at deme $x$ at time $t-$ independently attempts to make a jump at rate $m$ to a uniformly chosen deme in $\{x - L^{-1}, \, x + L^{-1}\}$; if this new deme is outside $\Lambda_n$ then the particle does not jump.
    \item \textit{Partially recovered particles}: For each $t \geq 0$ and $x \in \Lambda_n$, each dual pair of class $1*$ and class $2*$ particles at deme $x$ at time $t-$ independently attempts to make a jump at rate $m$ to a uniformly chosen deme in $\{x - L^{-1}, \, x + L^{-1}\}$; if this new deme is outside $\Lambda_n$ then the dual pair does not jump.
    \end{itemize}
 \item \textit{Reproduction of particles.} In each of the cases below, whenever a new particle is produced, say at time $t$ and carrying $k$ mutations, then with probability $\mu$, a mutation occurs at time~$t$, and the new particle now carries $k+1$ mutations. In the case $k+1>K_n$, this newly mutated particle is immediately deleted from the process.
 \begin{itemize}
    \item \textit{Infected and partially recovered particles}: For each $i\in \{1,2\}$, for $t \geq 0$, $x \in \Lambda_n$ and $k\in \mathbb N_0$, each class $i$ or class $i*$ particle at deme $x$ at time $t-$ carrying exactly $k$ mutations independently reproduces at rate 
    \begin{equation} \label{Paper01_eq:betainfdefn}
        \beta^i_{\textrm{infected}}(t,x,k):=\mathds 1_{\{k\le K_n\}}s_{k}\displaystyle q_{+}^{N}\left({\| \eta^{n,(i)}(t-,x) \|_{\ell_{1}}}\right).
    \end{equation}
After such an event, a new class $i$ particle carrying $k$ mutations is produced at deme $x$.
    \item \textit{Uninfected reproduction of susceptible particles}: For $t \geq 0$, $x \in \Lambda_n$ and $k\in \mathbb N_0$, each class $0$ particle at deme $x$ at time $t-$ carrying exactly $k$ mutations independently reproduces at rate
    \begin{equation} \label{Paper01_eq:betamindefn}
        \beta_{\min}(t,x,k):=\mathds 1_{\{k\le K_n\}} s_{k}\displaystyle \bigg(q_{+}^{N}\left({\| \eta^{n,(1)}(t-,x) \|_{\ell_{1}}}\right) \wedge q_{+}^{N}\left({\| \eta^{n,(2)}(t-,x) \|_{\ell_{1}}}\right)\bigg).
    \end{equation}
    After such an event, a new class $0$ particle carrying $k$ mutations is produced at deme $x$. 
     \item \textit{Infected reproduction of susceptible particles}:  
     Recalling~\eqref{Paper01_event_more_particles_process_i_than_3_minus_i}, for each choice of ordered pair $(i_1,i_2) \in \{1,2\}^{2}$, for $t \geq 0$ and $x \in \Lambda_n$, each class $i_1$ particle at deme $x$ at time~$t-$ independently \emph{induces} an infected reproduction event at rate
    \begin{align} \label{Paper01_eq:betainducedefn}
    \beta_{\textrm{induce}}^{i_1,i_2}(t,x)&:= \mathds 1_{A^{i_1}(t-,x)} \| \eta^{n,i_1}(t-,x) \|_{\ell_{1}}^{-1} \left(\sum_{k = 0}^{K_n} s_{k} \eta^{n,0}_{k}(t-,x) \right)\notag \\
    & \qquad \cdot \left(q_{+}^{N}\left({\| \eta^{n,(i_2)}(t-,x) \|_{\ell_{1}}}\right)-q_{+}^{N}\left({\| \eta^{n,(3-i_2)}(t-,x) \|_{\ell_{1}}}\right)\right)^{+}.
    \end{align}
    After such an event, we choose $\kappa \in [K_{n}]_{0}$ in such a way that for $k\in [K_{n}]_{0}$, the event $\{\kappa=k\}$ occurs with probability
    $$ \frac{s_{k}\eta^{n,0}_k(t-,x)}{\sum_{k' = 0}^{K_n} s_{k'} \eta^{n,0}_{k'}(t-,x) }.$$
    Then a new class $i_2$ particle carrying $\kappa$ mutations is produced at deme $x$. 
\end{itemize}
\item \textit{Death of particles.}
\begin{itemize}
    \item \textit{Death of infected particles}: For each $i\in \{1,2\}$, for $t \geq 0$ and $x \in \Lambda_n$, each class $i$ particle at deme $x$ at time $t-$ independently dies at rate
    \begin{equation} \label{Paper01_eq:deltainfdefn}
        \delta_{\textrm{infected}}^i(t,x):=q^{N}_{-}\left({\vert \vert \eta^{n,(i)}(t-,x) \vert \vert_{\ell_{1}}}\right).
    \end{equation}
    After such an event, the particle is removed from the process.
    \item \textit{Death of susceptible particles}: For~$t \geq 0$ and $x \in \Lambda_n$, each class $0$ particle at deme $x$ at time~$t-$ independently dies at rate
    \begin{equation} \label{Paper01_eq:deltamindefn}
       \delta_{\textrm{min}}(t,x):=q_{-}^{N}\left({\vert \vert \eta^{n,(1)}(t-,x) \vert \vert_{\ell_{1}}}\right) \wedge q_{-}^{N}\left({\vert \vert \eta^{n,(2)}(t-,x) \vert \vert_{\ell_{1}}}\right).
    \end{equation}
    After such an event, the particle is removed from the process.
     \item \textit{Transmission of infection to susceptible particles}: Recalling~\eqref{Paper01_event_more_particles_process_i_than_3_minus_i},
     for each choice of ordered pair $(i_1,i_2) \in \{1,2\}^{2}$,
     for $t \geq 0$ and $x \in \Lambda_n$, each class $i_1$ particle at deme $x$ at time~$t-$ independently has a \emph{transmission} event at rate
     \begin{align} \label{Paper01_eq:deltatransmitdefn}
         \delta_{\textrm{transmit}}^{i_1,i_2}(t,x)&:=\mathds 1_{A^{i_1}(t-,x)}
         \frac{\| \eta^{n,0}(t-,x) \|_{\ell_{1}}}{\| \eta^{n,i_1}(t-,x) \|_{\ell_{1}}} \notag \\
         &\qquad \cdot \left(q_{-}^{N}\left({\| \eta^{n,(i_2)}(t-,x) \|_{\ell_{1}}}\right)-q_{-}^{N}\left({\| \eta^{n,(3-i_2)}(t-,x) \|_{\ell_{1}}}\right)\right)^{+}.
     \end{align}
    After such an event, an \emph{infection} occurs: a particle is chosen uniformly at random from the class $0$ particles at deme $x$ at time $t-$; if $i_1 \neq i_2$, this particle is given class $3 - i_2$.
    Otherwise, if $i_1 = i_2$, a \emph{partial recovery} occurs: this particle is given class $(3 - i_2)*$ and becomes a dual pair with the class $i_1$ particle that transmitted the infection (which now becomes class $i_1*$). 
    \item \textit{Death of dual pairs of particles}: For $t \geq 0$ and $x \in \Lambda_n$, each dual pair of class $1*$ and class $2*$ particles at deme $x$ at time $t-$ dies independently at rate $\delta_{\min}(t,x)$ defined in~\eqref{Paper01_eq:deltamindefn}. 
    After such an event, the dual pair of particles are removed from the process.

    \item \textit{Death of partially recovered particles with reinfection or replacement}: 
    Recalling~\eqref{Paper01_event_more_particles_process_i_than_3_minus_i}, for each choice of ordered pair $(i_1,i_2) \in \{1,2\}^{2}$, for $t \geq 0$ and $x \in \Lambda_n$, each class $i_2*$ particle at time~$t-$ at deme $x$ independently dies at rate 
     \begin{equation} \label{Paper01_eq:deltapartialdefn}
         \delta_{\textrm{partial}}^{i_1,i_2}(t,x):=\mathds 1 _{A^{i_1}(t-,x)} \left(q_{-}^{N}\left({\| \eta^{n,(i_2)}(t-,x) \|_{\ell_{1}}}\right)-q_{-}^{N}\left({\| \eta^{n,(3-i_2)}(t-,x) \|_{\ell_{1}}}\right)\right)^{+}.
     \end{equation}
     After such an event, the class $i_2*$ particle (call this particle $p$) is removed from the process. Let $p'$ denote the class $(3 - i_2)*$ dual particle of particle $p$ at time $t-$.
     If $i_1\neq i_2$,
     then a \emph{reinfection} occurs: particle $p'$ becomes class $3 - i_2$.
     If instead $i_1=i_2$,
     then a \emph{replacement} occurs: we choose a class $i_2$ particle at deme $x$ at time $t-$ uniformly at random (note that there must be at least one), and this particle becomes the new dual particle of particle $p'$, and becomes class $i_2*$.
    \end{itemize}
\end{itemize}

For the sake of completeness, we state the following result. Recall the definition of the infinitesimal generator~$\mathcal{L}^{n}$ in~\eqref{Paper01_generator_foutel_etheridge_model_restriction_n} and~\eqref{Paper01_infinitesimal_generator_restriction_n}.

\begin{lemma} \label{Paper01_coupling_lemma_between_coloured_and_uncoloured_processes}
For $\boldsymbol{\eta}^{(1)}, \boldsymbol{\eta}^{(2)} \in \mathcal{S}$ such that $\# \Delta(\boldsymbol{\eta}^{(1)}, \boldsymbol{\eta}^{(2)}) < \infty$,
   under the construction above, the Markov processes ${\left(\eta^{n,(1)}(t), \eta^{n,(2)}(t)\right)_{t \geq 0}}$ and
   \begin{equation*}
    {\left(\eta^{n,\textrm{(coupling)}}(t)\right)_{t \geq 0} \defeq \left(\eta^{n,0}(t), \eta^{n,1}(t), \eta^{n,2}(t),\eta^{n,1*}(t),\eta^{n,2*}(t)\right)_{t \geq 0}}
\end{equation*} can be defined on the same probability space in such a way that
     for each $i\in \{1,2\}$, the process
$(\eta^{n,(i)}(t))_{t\ge 0}$ is
a realisation of the Markov process associated with the infinitesimal generator $\mathcal{L}^{n}$
with $\eta^{n,(i)}(0) = \boldsymbol{\eta}^{(i)}$ almost surely, and almost surely the following relation holds for all $t \ge 0$ and~$i \in \{1,2\}$:
    \begin{equation*}
    \eta^{n,(i)}(t) = \eta^{n,0}(t) + \eta^{n,i}(t) + \eta^{n,i*}(t).
\end{equation*}
\end{lemma}

The proof is omitted since it follows directly from computing and comparing the action of the corresponding infinitesimal generators. 
For $\boldsymbol{\eta}^{(1)}, \boldsymbol{\eta}^{(2)} \in \mathcal{S}$ with $\# \Delta(\boldsymbol{\eta}^{(1)}, \boldsymbol{\eta}^{(2)}) < \infty$,
we write $\mathbb{P}_{\boldsymbol{\eta}^{(1)}, \boldsymbol{\eta}^{(2)}}$ to denote the probability measure under which $\eta^{n,(1)}(0) = \boldsymbol{\eta}^{(1)}$ and $\eta^{n,(2)}(0) = \boldsymbol{\eta}^{(2)}$, and~$\mathbb{E}_{\boldsymbol{\eta}^{(1)}, \boldsymbol{\eta}^{(2)}}$ to denote the corresponding expectation.

\begin{remark}  We collect here some observations regarding the definition of the process $\eta^{n,\textrm{(coupling)}}$.

\begin{enumerate}[(a)]
    \item   Note that in the process~$\eta^{n,\textrm{(coupling)}}$, a dual pair of class~$1*$ and class~$2*$ cannot be merged into a single class~$0$ particle because the two particles may carry different numbers of mutations.

    \item If the process $\eta^{n,\textrm{(coupling)}}$ were defined with only susceptible and infected particles, it would not be possible to have enough control on the rate at which new infected particles are produced. The construction involving partially recovered particles will be crucial in allowing us to prove Proposition~\ref{Paper01_proposition_bound_spread_infection}.
\end{enumerate}
\end{remark}

Recall that $q_{+}, q_{-}: [0,\infty) \rightarrow [0,\infty)$ satisfy Assumption~\ref{Paper01_assumption_polynomials}, i.e.~that $q_+$ and $q_-$ are non-negative polynomials with $0 \leq \deg q_{+} < \deg q_{-}$, and therefore the leading coefficient of~$q_{+}$ must be non-negative, and the leading coefficient of~$q_{-}$ must be strictly positive.
It follows that for any $\epsilon>0$, there exists $U_\epsilon$ such that the following holds.
\begin{definition} \label{Paper01_definition_cutoff_density_parameter}
    For each $\epsilon>0$, define $U_{\epsilon} > 1$ sufficiently large that the following conditions hold:
    \begin{enumerate}[(i)]
    \item The polynomials $q_{+}$, $q_{-}$ and $q'_+$ are non-negative and monotonically non-decreasing on $[U_{\epsilon}, \infty)$.
    \item For every $0 \le U \leq U_\epsilon$,
    \begin{equation*}
        q_{+}(U) \le q_{+}(U_{\epsilon}), \quad q'_{+}(U) \le q'_{+}(U_{\epsilon}), \quad  \quad \textrm{and} \quad  q_{-}(U) \le q_{-}(U_{\epsilon}).
    \end{equation*}
    \item For every $U\ge U_\epsilon$,
    \begin{equation*}
        \frac{q_{+}(U)+ Uq'_+(U)}{q_{-}(U)} \leq \epsilon.
    \end{equation*}
\end{enumerate}
\end{definition}

We can now make the following observations for demes with high particle density.

\begin{lemma}\label{Paper01_lem:cutoffconseq}
    For $\epsilon>0$, $t\ge 0$, and $x\in \Lambda_n$ with $\|\eta^{n,(i)}(t-,x)\|_{\ell_1}\ge NU_\epsilon$ for some $i\in \{1,2\}$, the following holds:
    \begin{enumerate}[(i)]
        \item If an infected (class~$1$ or~$2$) particle at deme $x$ has a transmission event at time $t$, then a partial recovery occurs. If a partially recovered (class~$1*$ or~$2*$) particle at deme $x$ dies at time $t$ and its dual particle does not, then a replacement occurs (i.e.~a reinfection does not occur).
        Moreover, for each $j\in \{1,2\}$, if a class $j$ particle at deme $x$ induces an infected reproduction event at time $t$, then the new particle added at deme~$x$ is a class $j$ particle. 
        \item For any $k\in \mathbb N_0$,
        \begin{align}
            \frac{\beta^i_{\textrm{infected}}(t,x,k)+\beta^{i,i}_{\textrm{induce}}(t,x)}{\beta^i_{\textrm{infected}}(t,x,k)+\beta^{i,i}_{\textrm{induce}}(t,x)+ \delta^i_{\textrm{infected}}(t,x) + \delta^{i,i}_{\textrm{transmit}}(t,x)}&\leq \epsilon \label{Paper01_estimate_rate_giving_birth_or_inducing_reprod_event_high_density} \\
        \text{and}\qquad
        \frac{\beta^i_{\textrm{infected}}(t,x,k)}{\beta^i_{\textrm{infected}}(t,x,k)+\delta_{\textrm{min}}(t,x)+\delta^{i,i}_{\textrm{partial}}(t,x)}&\leq \epsilon \label{Paper01_estimate_rate_only_giving_birth_high_density}. 
        \end{align}
    \end{enumerate}
\end{lemma}

\begin{proof}
    We begin by proving~(i).
    By the definition of the process $\eta^{n,(\textrm{coupling)}}$, it suffices to show that for each $j\in \{1,2\}$ we have
    \begin{equation} \label{eq:cutofflemclaim}
        \delta^{j,3-j}_{\textrm{transmit}}(t,x)=0, \quad \delta^{j,3-j}_{\textrm{partial}}(t,x)=0, \quad \text{and}\quad \beta^{j,3-j}_{\textrm{induce}}(t,x)=0.
    \end{equation}
    Take $j\in \{1,2\}$ and suppose first that 
    \begin{equation} \label{eq:spsmorej}
        \vert \vert \eta^{n,j}(t-,x) \vert \vert_{\ell_{1}} \ge \vert \vert \eta^{n,3-j}(t-,x) \vert \vert_{\ell_{1}}.
    \end{equation}
    By~\eqref{Paper01_definition_number_particles_composition_reference_modified_processes} and~\eqref{Paper01_equations_colouring_pink_purple_processes}, it follows that
     \[
     \| \eta^{n,(j)}(t-,x)\|_{\ell_{1}} \ge \vert \vert \eta^{n,(3-j)}(t-,x) \vert \vert_{\ell_{1}}.
     \]
     In particular, by the assumption at the start of the statement, we have $\vert \vert \eta^{n,(j)}(t-,x) \vert \vert_{\ell_1} \geq NU_\epsilon$. 
     Therefore, by the definition of $q^N_{+}$ and $q^N_{-}$ in~\eqref{Paper01_scaled_polynomials_carrying_capacity}  and by Definition~\ref{Paper01_definition_cutoff_density_parameter}(i) and~(ii), we have
     \begin{equation} \label{Paper01_immediate_conclusion_high_density_particles}
     \begin{aligned}
         & q_{-}^{N}\left({\| \eta^{n,(j)}(t-,x) \|_{\ell_{1}}}\right) \ge q_{-}^{N}\left({\| \eta^{n,(3-j)}(t-,x) \|_{\ell_{1}}}\right)\\
         \text{and} \quad &q_{+}^{N}\left({\| \eta^{n,(j)}(t-,x) \|_{\ell_{1}}}\right)\ge q_{+}^{N}\left({\| \eta^{n,(3-j)}(t-,x) \|_{\ell_{1}}}\right),
    \end{aligned}
     \end{equation}
     and the claim~\eqref{eq:cutofflemclaim} follows from\eqref{Paper01_eq:deltatransmitdefn},~\eqref{Paper01_eq:deltapartialdefn} and~\eqref{Paper01_eq:betainducedefn}.
     If instead~\eqref{eq:spsmorej} does not hold, then by~\eqref{Paper01_event_more_particles_process_i_than_3_minus_i} we have 
     $\mathds 1_{A^j(t-,x)}=0$, and~\eqref{eq:cutofflemclaim} follows.

     It remains to prove (ii); we begin by establishing~\eqref{Paper01_estimate_rate_giving_birth_or_inducing_reprod_event_high_density}. Recall from Assumption~\ref{Paper01_assumption_fitness_sequence} that $s_k\le 1$ $\forall k\in \mathbb N_0$.
     First suppose that $\vert \vert \eta^{n,i}(t-,x) \vert \vert_{\ell_1} \le \vert \vert \eta^{n,3-i}(t-,x) \vert \vert_{\ell_1}$.
     Then for $k\in \mathbb N_0$, we have $\beta^i_{\textrm{infected}}(t,x,k)\le q_{+}^{N}({\| \eta^{n,(i)}(t-,x) \|_{\ell_{1}}})$ by~\eqref{Paper01_eq:betainfdefn}, $\beta^{i,i}_{\textrm{induce}}(t,x)=0$ by~\eqref{Paper01_eq:betainducedefn} and~\eqref{Paper01_event_more_particles_process_i_than_3_minus_i}, and
     $\delta^i_{\textrm{infected}}(t,x)= q_{-}^{N}({\| \eta^{n,(i)}(t-,x) \|_{\ell_{1}}})$ by~\eqref{Paper01_eq:deltainfdefn}.
     Since $\| \eta^{n,(i)}(t-,x) \|_{\ell_{1}}\ge N U_\epsilon$ by the assumption at the start of the statement, it follows from Definition~\ref{Paper01_definition_cutoff_density_parameter}(i) and~(iii) that~\eqref{Paper01_estimate_rate_giving_birth_or_inducing_reprod_event_high_density} holds. 
     
     Now suppose instead that $\vert \vert \eta^{n,i}(t-,x) \vert \vert_{\ell_1} > \vert \vert \eta^{n,3-i}(t-,x) \vert \vert_{\ell_1}$. Note that by~\eqref{Paper01_definition_number_particles_composition_reference_modified_processes} and~\eqref{Paper01_equations_colouring_pink_purple_processes}, we have
     \begin{equation} \label{Paper01_eq:diffofetas}
      0 < \|\eta^{n,(i)}(t-,x) \|_{\ell_{1}}-\| \eta^{n,(3-i)}(t-,x) \|_{\ell_{1}}
     \le \| \eta^{n,i}(t-,x) \|_{\ell_{1}}.
     \end{equation}
     Also, since $\vert \vert \eta^{n,(i)}(t-,x) \vert \vert_{\ell_{1}} \geq NU_{\epsilon}$, we have by Definition~\ref{Paper01_definition_cutoff_density_parameter}(i) and~(ii) that
     \begin{equation} \label{Paper01_trivial_inequality_derivative_birth_polynomial}
         q'_{+}\left(\frac{u}{N}\right) \leq q'_{+}\left(\frac{\vert \vert \eta^{n,(i)}(t-,x) \vert \vert_{\ell_{1}}}{N}\right) \quad \forall \, u \in \Big[0, \vert \vert \eta^{n,(i)}(t-,x) \vert \vert_{\ell_{1}}\Big].
     \end{equation}
     Therefore, by~\eqref{Paper01_eq:betainducedefn}, and by~\eqref{Paper01_eq:diffofetas},~\eqref{Paper01_scaled_polynomials_carrying_capacity} and since $s_{k} \leq 1$ for all~$k \in \mathbb{N}_{0}$ (by Assumption~\ref{Paper01_assumption_fitness_sequence}) in the first inequality, and then in the second inequality by~\eqref{Paper01_eq:diffofetas} and~\eqref{Paper01_trivial_inequality_derivative_birth_polynomial},
     we have
     \begin{align*}
    \beta_{\textrm{induce}}^{i,i}(t,x)
    &\le \mathds{1}_{A^{i}(t-,x)}\| \eta^{n,i}(t-,x) \|_{\ell_{1}}^{-1} \|\eta^{n,0}(t-,x) \|_{\ell_1} 
    \int_{\| \eta^{n,(3-i)}(t-,x) \|_{\ell_{1}}}^{\| \eta^{n,(i)}(t-,x) \|_{\ell_{1}}} N^{-1} \left(q'_+\left(\frac{u}{N}\right)\right)^+ du\\
    &\le N^{-1}\|\eta^{n,0}(t-,x) \|_{\ell_1} 
    q'_+\left(\frac 1 N\| \eta^{n,(i)}(t-,x) \|_{\ell_{1}}\right).
    \end{align*}
    By~\eqref{Paper01_eq:betainfdefn},~\eqref{Paper01_eq:deltainfdefn} and~\eqref{Paper01_eq:deltatransmitdefn}, and since $s_{k} \leq 1$ for all~$k \in \mathbb{N}_{0}$, it follows that for $k \in \mathbb N_0$,
    \begin{align*}
        &\frac{\beta^i_{\textrm{infected}}(t,x,k)+\beta^{i,i}_{\textrm{induce}}(t,x)}{\beta^i_{\textrm{infected}}(t,x,k)+\beta^{i,i}_{\textrm{induce}}(t,x)+\delta^i_{\textrm{infected}}(t,x)+\delta^{i,i}_{\textrm{transmit}}(t,x)}\\
        &\le \frac{q^N_+(\| \eta^{n,(i)}(t-,x) \|_{\ell_{1}})+N^{-1}\|\eta^{n,0}(t-,x) \|_{\ell_1} 
    q'_+(\tfrac 1 N\| \eta^{n,(i)}(t-,x) \|_{\ell_{1}})}{q^N_-(\| \eta^{n,(i)}(t-,x) \|_{\ell_{1}})}\\
    & \leq \epsilon,
    \end{align*}
    where the last inequality follows by Definition~\ref{Paper01_definition_cutoff_density_parameter}(iii)
    and since $\vert \vert \eta^{n,(i)}(t-,x) \vert \vert_{\ell_{1}} \geq NU_{\epsilon}$. 
    This completes the proof of~\eqref{Paper01_estimate_rate_giving_birth_or_inducing_reprod_event_high_density}.

    To prove~\eqref{Paper01_estimate_rate_only_giving_birth_high_density}, note that by~\eqref{Paper01_eq:deltamindefn} and~\eqref{Paper01_eq:deltapartialdefn}, using~\eqref{Paper01_immediate_conclusion_high_density_particles} again in the case 
    $\|\eta^{n,3-i}(t-,x) \|_{\ell_{1}}\ge \| \eta^{n,i}(t-,x) \|_{\ell_{1}}$ and using that $\mathds 1_{A^i(t-,x)}=1$ otherwise, we have
    \[
    \delta_{\textrm{min}}(t,x)+\delta^{i,i}_{\textrm{partial}}(t,x)=
    q^N_-(\| \eta^{n,(i)}(t-,x) \|_{\ell_{1}}).
    \]
    Therefore, for $k\in \mathbb N_0$, by~\eqref{Paper01_eq:betainfdefn}, and then by Definition~\ref{Paper01_definition_cutoff_density_parameter}(iii),
    \[
        \frac{\beta^i_{\textrm{infected}}(t,x,k)}{\beta^i_{\textrm{infected}}(t,x,k)+\delta_{\textrm{min}}(t,x)+\delta^{i,i}_{\textrm{partial}}(t,x)}
        \le \frac{q^N_+(\| \eta^{n,(i)}(t-,x) \|_{\ell_{1}})}{q^N_-(\| \eta^{n,(i)}(t-,x) \|_{\ell_{1}})} \leq \epsilon,
    \]
    which completes the proof.
\end{proof}

It will be convenient to assign labels to infected and partially recovered particles in the process $\eta^{n,\textrm{(coupling)}}$.
For an infected (class $1$ or class $2$) particle, we say that the particle \emph{generates} another infected or partially recovered particle when the particle reproduces, or induces an infected reproduction event, or has a transmission event (when a class $0$ particle becomes either infected or partially recovered).
For a partially recovered (class $1*$ or class $2*$) particle, we say that the particle \emph{generates} another infected or partially recovered particle when the particle reproduces.
We define a label set
\begin{equation} \label{Paper01_ulam-harris-scheme}
    \mathcal{U} \defeq  \mathbb{N} \times \bigcup_{i = 0}^{\infty} \mathbb{N}^{i},
\end{equation}
and we label infected and partially recovered particles according to an Ulam-Harris labelling scheme in the following sense. 
The infected and partially recovered particles in $\eta^{n,(\textrm{coupling})}$ at time $0$ are given some arbitrary ordering; let $\mathcal N_0$ denote the total number of these particles.
Then for each $i \in [\mathcal N_0]$, the $i$th particle in this ordering is given label $i$.
When a particle is given a label, the particle
keeps the same label until it dies.
For each $u\in \mathcal U$ and $j\in \mathbb N$, the $j$th particle generated by the particle with label $u$ is given label $uj$.

We note that a particle’s label remains unchanged if its status shifts between infected and partially recovered. Moreover, labels are not assigned to susceptible (class $0$) particles that remain susceptible until their death.

We now make some more definitions related to the labelling of the particles.

\begin{itemize}
\item For $\ell\in \mathbb N_0$, $u=u_0u_1\ldots u_\ell\in \mathcal U$ and $i\in [\ell]_{0}$, we let
$u|_i=u_0u_1\ldots u_i$. We say that $u|_i$ is an {\textit{ancestor}} of $u$ in $\mathcal U$.
\item For $\ell\in \mathbb N_0$ and $u=u_0u_1\ldots u_\ell\in \mathcal U$, we let
\begin{equation} \label{Paper01_length_branch_ulam_harris_tree}
    |u|:=\mathds 1 _{\{\ell\ge 1\}}\sum_{i=1}^\ell u_i.
\end{equation}
\item For $\ell\in \mathbb N_0$ and $u=u_0 u_1 \ldots u_\ell\in \mathcal U$, we write
\begin{equation} \label{Paper01_set_labels_particles_generated_same_ancestral_path}
    \mathcal D_u :=
    \begin{cases}
        \emptyset \quad &\text{if }\ell=0,\\
    \{u\vert_{i}j:i\in [\ell -1]_0, j\in [u_{i+1}]\} \quad &\text{if }\ell\ge 1,
    \end{cases}
\end{equation}
for the set of labels of certain particles generated by ancestors of $u$. Note that $\# \mathcal D_u=|u|$.
\item For $t\ge 0$, we let $\mathcal N^n(t)\subset \mathcal U$ denote the set of labels of particles in $\eta^{n,(\textrm{coupling})}$ that are alive and either infected or partially recovered (class $1,1*,2$ or $2*$) at time $t$.
\item For $t\ge 0$, we let $\mathcal I^n(t)\subset \mathcal U$ denote the set of labels of particles in $\eta^{n,(\textrm{coupling})}$ that are alive and infected (class $1$ or $2$) at time $t$.
\item For $t\ge 0$, we let $\mathcal P^n(t)\subset \mathcal U$ denote the set of labels of particles in $\eta^{n,(\textrm{coupling})}$ that are alive and partially recovered (class $1*$ or $2*$) at time $t$.
\item For $u\in \mathcal U$, we let $\tau^{n,\textrm{birth}}_u$ denote the time at which the particle with label $u$ first becomes an infected or partially recovered particle, and we let $\tau^{n,\textrm{death}}_u$ denote the time at which it dies, i.e.
\begin{equation*}
\begin{aligned}
    \tau^{n,\textrm{birth}}_u &\defeq \inf(\{t\ge 0:u\in \mathcal N^n(t)\} \cup \{\infty\}),  \\
    \text{and}\quad 
    \tau^{n,\textrm{death}}_u &\defeq \inf(\{t\ge  \tau^{n,\textrm{birth}}_u :u \notin \mathcal N^n(t)\} \cup \{\infty\}). 
\end{aligned}
\end{equation*}
    
    \item For $u\in \mathcal U$ and $t\ge 0$ with $u\in \mathcal N^n(t)$,
we let $Z_u^n(t)\in L^{-1}\mathbb Z$ denote the deme of the particle with label $u$ at time $t$.
For $\ell \in \mathbb N_0$ and $u = u_{0}u_{1}\ldots u_{\ell} \in \mathcal{U}$ with $\tau^{n,\textrm{birth}}_{u} <\infty$, for $t \in [0,\tau^{n,\textrm{birth}}_u)$, we let $Z_u^n(t)\in L^{-1}\mathbb Z$ denote the deme of the most recent ancestor (in $\mathcal U$) of particle $u$ that is alive at time $t$, i.e.
\begin{equation} \label{eq:Zbulletlist}
Z_u^n(t):= Z_{u |_{i_{\max}(t)}}^n(t), \quad \text{where }
    i_{\max}(t) \defeq \max\{i\in [\ell]_0:u|_i \in \mathcal N^n(t)\}.
\end{equation}
\item For $u\in \mathcal U$ with $\tau^{n,\textrm{birth}}_u<\infty$,
we let $\gamma^{n,\textrm{mut}}_u$ denote the number of mutations carried by the particle with label $u$.
\item For $u\in \mathcal U$ with $\tau^{n,\textrm{birth}}_u<\infty$, and $t\ge 0$, we let $J^n_u(t)$ denote the number of jumps (migration events) made by the particle with label $u$ and its ancestors (in $\mathcal{U}$) up to time~$t$, i.e.
\begin{equation} \label{eq:Jbulletlist}
J^{n}_u(t) \defeq \#\Big\{s\in [0,t] \cap [0,\tau^{n,\textrm{death}}_u): Z_{u}^n(s)\neq Z_{u}^n(s-)\Big\}.
\end{equation}
\item For $\ell\in \mathbb N_0$ and $u=u_0 u_1 \ldots u_\ell\in \mathcal U$ with $\tau^{n,\textrm{birth}}_{u} < \infty$, we let $R^n_u$ denote the number of times that the ancestors (in $\mathcal U$) of the particle with label $u$ are reinfected during the time intervals $(\tau^{n,\textrm{birth}}_{u|_i},\tau^{n,\textrm{birth}}_{u|_{i+1}})$ for $i\in [\ell-1]_0$, i.e.
\begin{equation} \label{eq:Rbulletlist}
    R^n_u \defeq \mathds 1_{\{\ell\ge 1\}}\sum_{i=0}^{\ell-1}\#\Big\{s\in (\tau^{n,\textrm{birth}}_{u|_i},\tau^{n,\textrm{death}}_{u|_i}\wedge \tau^{n,\textrm{birth}}_{u|_{i+1}}): u|_i\in \mathcal I^n(s) \textrm{ and } u|_i\in \mathcal P^n(s-)\Big\}.
\end{equation}
\item For $u\in \mathcal U$ with $\tau^{n,\textrm{birth}}_{u} < \infty$, let $\gamma^{n,\textrm{class}}_u$ indicate whether the particle with label $u$ is in classes $1$ and $1*$ at its birth time, or in classes $2$ and $2*$, i.e.~for $i \in \{1,2\}$,
\[
\gamma^{n,\textrm{class}}_u \defeq i \text{ if the particle with label }u\text{ is class }i \text{ or } i*\text{ at time }\tau^{n,\textrm{birth}}_u.
\]
\item For $u\in \mathcal U$ and $j\in \mathbb N$, we let $\gamma^{n,\textrm{infected},-}_{uj}$ and $\gamma^{n,\textrm{infected},+}_{uj}$ indicate whether or not the particles $u$ and $uj$ are infected at times $\tau^{n,\textrm{birth}}_{uj}-$ and $\tau^{n,\textrm{birth}}_{uj}$ respectively, i.e.
\begin{align} \label{eq:gammaIbulletlist}
\gamma^{n,\textrm{infected},-}_{uj} &\defeq
\mathds 1 _{\{\tau^{n,\textrm{birth}}_{uj}<\infty\}\cap \{u\in \mathcal I^n(\tau^{n,\textrm{birth}}_{uj}-)\}} \notag\\
\text{and}\quad
\gamma^{n,\textrm{infected},+}_{uj} &\defeq
\mathds 1 _{\{\tau^{n,\textrm{birth}}_{uj}<\infty\}\cap\{uj\in \mathcal I^n(\tau^{n,\textrm{birth}}_{uj})\}}.
\end{align}
\item For $u\in \mathcal U$, $j\in \mathbb N$ and $\epsilon>0$, we let $\gamma^{n,\textrm{density},\epsilon}_{uj}$ indicate whether or not particle $u$ is in a high density deme (as determined by the quantity $U_\epsilon$ introduced in Definition~\ref{Paper01_definition_cutoff_density_parameter}) at time $\tau^{n,\textrm{birth}}_{uj}-$, i.e.
\begin{equation} \label{eq:gammaDbulletlist}
\gamma^{n,\textrm{density},\epsilon}_{uj} \defeq
\mathds{1}_{\{\tau^{n,\textrm{birth}}_{uj}<\infty\} \cap A^{\epsilon}_{u,j}},  
\end{equation}
where
\[
A^{\epsilon}_{u,j} \defeq \{\|\eta^{n,(\gamma^{n,\textrm{class}}_{u})}(\tau^{n,\textrm{birth}}_{uj}-, Z^n_{u}(\tau^{n,\textrm{birth}}_{uj}-))\|_{\ell_1}\ge N U_\epsilon\}.
\]
\end{itemize}

Now take $\epsilon>0$, $\ell\in \mathbb N_0$, $u=u_0u_1\ldots u_\ell \in \mathcal U$, $J\in \mathbb N_0$,
$R\in \mathbb N_0$,
$\gamma^-=(\gamma^-_{v})_{v\in \mathcal D_u} \in \{0,1\}^{\mathcal D_u}$,
$\gamma^+=(\gamma^+_{v})_{v\in \mathcal D_u} \in \{0,1\}^{\mathcal D_u}$,
$d=(d_{v})_{v\in \mathcal D_u} \in \{0,1\}^{\mathcal D_u}$
and $T\ge 0$.
Then define the event
\begin{equation} \label{Paper01_event_genealogical_path_considering_the_number_reinfection_events}
\begin{aligned}
    E^{n,\epsilon}_{u,J,R,\gamma^-,\gamma^+,d}(T)
    &=\Big\{
    u\in \mathcal N^n(T),\; J^n_{u}(T)\ge J,\; R^n_{u}= R,
    \gamma^{n,\textrm{infected},-}_{v}=\gamma^-_{v} \; \forall v \in \mathcal D_u,\\
    &\hspace{1cm}
    \gamma^{n,\textrm{infected},+}_{v}=\gamma^+_{v} \; \forall v \in \mathcal D_u,\; \gamma^{n,\textrm{density},\epsilon}_{v}=d_{v} \; \forall v\in \mathcal D_u
    \Big\}.
\end{aligned}
\end{equation}
We now use a construction of the process $\eta^{n,(\textrm{coupling})}$ that depends on~$u$
to bound the probability of the event $E^{n,\epsilon}_{u,J,R,\gamma^-,\gamma^+,d}(T)$. Recall the definition of $U_{\epsilon}$ from Definition~\ref{Paper01_definition_cutoff_density_parameter}.

\begin{lemma} \label{Paper01_lem:boundPE}
For $\alpha\ge 0$, let $Z_\alpha$ denote a Poisson random variable with mean $\alpha$. 
For $\epsilon>0$, let
\begin{equation} \label{Paper01_eq:Vepsdefn}
    V^N_\epsilon := (2+NU_\epsilon)(q_+(U_\epsilon)+q_-(U_\epsilon)).
\end{equation}
Then for $\boldsymbol{\eta}^{(1)},\boldsymbol{\eta}^{(2)} \in \mathcal{S}$ with $\#\Delta (\boldsymbol{\eta}^{(1)},\boldsymbol{\eta}^{(2)})<\infty$, for
$\epsilon>0$, $\ell\in \mathbb N_0$, $u=u_0u_1\ldots u_\ell \in \mathcal U$, $J\in \mathbb N_0$,
$R\in \mathbb N_0$,
$\gamma^-=(\gamma^-_{v})_{v\in \mathcal D_u} \in \{0,1\}^{\mathcal D_u}$,
$\gamma^+=(\gamma^+_{v})_{v\in \mathcal D_u} \in \{0,1\}^{\mathcal D_u}$,
$d=(d_{v})_{v\in \mathcal D_u} \in \{0,1\}^{\mathcal D_u}$
and $T\ge 0$,
\begin{equation} \label{Paper01_upper_bound_probability_one_infected_particle_living_too_far_from_original_point}
\begin{aligned}
    \mathbb P_{\boldsymbol{\eta}^{(1)},\boldsymbol{\eta}^{(2)} }\left(E^{n,\epsilon}_{u,J,R,\gamma^-,\gamma^+,d}(T) \right)
    &\le
    \mathbb P\left(Z_{mT}\ge J\right)
    \mathbb P\left(Z_{V_\epsilon^N T}\ge R\right)\\
    &\qquad \cdot \mathbb P\left(Z_{V^N_\epsilon T}\ge \#\{v\in \mathcal D_u: d_{v}=0\}\right)
    \epsilon^{\#\{v\in \mathcal D_u:d_{v}=1,(\gamma^-_{v},\gamma^+_{v})\neq (1,0)\}}.
\end{aligned}
\end{equation}
Moreover, if
\begin{equation} \label{Paper01_deterministic_condition_path_genealogical_infection}
    2\Big(\#\{v\in \mathcal D_u: d_{v}=0\}+\#\{v\in \mathcal D_u:d_{v}=1,(\gamma^-_{v},\gamma^+_{v})\neq (1,0)\}\Big)+R<|u|-1,
\end{equation}
then 
\begin{equation} \label{Paper01_event_cannot_happen_if_number_reinfections_too_low}
    \mathbb P_{\boldsymbol{\eta}^{(1)},\boldsymbol{\eta}^{(2)} }\left(E^{n,\epsilon}_{u,J,R,\gamma^-,\gamma^+,d}(T) \right)=0.
\end{equation}
\end{lemma}

\begin{proof}
    Let $(X^{\textrm{jump}}(t),t\ge 0)$, $(X^{\textrm{low}}(t),t\ge 0)$, $(X^{\textrm{high}}(t),t\ge 0)$, and $(X^{\textrm{reinf}}(t),t\ge 0)$ denote i.i.d.~Poisson processes with rate $1$.
    Let $(B^\epsilon_i)_{i=1}^\infty$ denote an independent sequence of i.i.d.~Bernoulli($\epsilon$) random variables. 
    We will now construct a realisation of $\eta^{n,(\textrm{coupling})}$ up to a suitable stopping time, using these Poisson processes and random variables to determine (some of) the behaviour of the particles with labels $u|_0$, $u|_1,\ldots , u|_\ell=u$. 
    
    It will be useful for the construction to define, for any $t\ge 0$, the most recent ancestor of particle $u$ that was born before time $t$, i.e.~for $t\ge 0$ we let
\[
p(t):=u|_{\ell(t)}, \quad \text{where }\ell(t)=\max\Big(\Big\{\ell' \in [\ell]_{0} : \tau^{n,\textrm{birth}}_{u|_{\ell'}}\le t\Big\}\cup \{0\}\Big).
\]
We let $\tau^*$ denote the last time at which the particle with label $u$ or an ancestor (in $\mathcal{U}$) of particle $u$ is alive, i.e.~we let
\[
\tau^*:=\sup \Big\{t\ge 0:p(t)\in \mathcal N^n(t)\Big\}. 
\]
Note that if $u_0>\#\mathcal N^n(0)$, then by our choice of labelling after~\eqref{Paper01_ulam-harris-scheme} we have $\tau^*=0$. 
We also let $i(t)$ indicate whether the particle with label $p(t)$ is in classes $1$ and $1*$, or in classes $2$ and $2*$, and we let $k(t)$ denote the number of mutations that it carries at time $t$, i.e.~for $t \in [0, \tau^{*})$, we let
\begin{equation} \label{eq:starinEeventpf}
   i(t)=\gamma^{n,\textrm{class}}_{p(t)}\quad \text{and} \quad
k(t)=\gamma^{n,\textrm{mut}}_{p(t)}. 
\end{equation}

We now define time changes of the Poisson processes $X^{\textrm{jump}}$, $X^{\textrm{low}}$, $X^{\textrm{high}}$ and $X^{\textrm{reinf}}$. Firstly, for $t\ge 0$, let
\begin{equation} \label{eq:1Eeventpf}
    Y^{\textrm{jump}}(t)=X^{\textrm{jump}}(m t).
\end{equation}
This process will determine the jump times of the particle with label $p(t-)$.
Then, for $t\in [0,\tau^*)$, define the event that the particle with label $p(t)$ is at a ‘low density' deme at time $t$ by writing 
\begin{equation} \label{eq:**Eeventpf}
    E^{\textrm{low}}(t)=\Big\{\|\eta^{n,(i(t))}(t, Z^n_{p(t)}(t))\|_{\ell_1}< N U_\epsilon\Big\}.
\end{equation}
Recalling~\eqref{Paper01_eq:Vepsdefn},
for $t\in [0,\tau^*]$, let
\begin{equation} \label{eq:2Eeventpf}
    Y^{\textrm{low}}(t)=X^{\textrm{low}}\left(V^N_\epsilon\int_0^t \mathds{1}_{E^{\textrm{low}}(s-)}ds
\right).
\end{equation}
The jump times of $Y^{\textrm{low}}$ will be times at which the particle with label $p(t-)$ is at a ‘low density' deme and a reproduction or death event generated by particle $p(t-)$ may occur.
For $i\in \{1,2\}$, $t\ge 0$, $x\in \Lambda_n$ and $k\in \mathbb N_0$, let
\begin{align} \label{eq:daggerEeventpf}
    \beta^{\mathcal I, i}(t,x,k)
    &=\beta^i_{\textrm{infected}}(t,x,k)+\beta^{i,i}_{\textrm{induce}}(t,x) +\delta^i_{\textrm{infected}}(t,x)+\delta^{i,i}_{\textrm{transmit}}(t,x), \notag \\
    \text{and }\quad
    \beta^{\mathcal P, i}(t,x,k)
    &=\beta^i_{\textrm{infected}}(t,x,k)+\delta_{\textrm{min}}(t,x)+\delta^{i,i}_{\textrm{partial}}(t,x).
\end{align}
Using Lemma~\ref{Paper01_lem:cutoffconseq}(i), if $\vert \vert \eta^{n,(i)}(t-,x) \vert \vert_{\ell_{1}} \geq NU_{\epsilon}$, then the quantity $\beta^{\mathcal{I},i}(t,x,k)$ is the rate at which a class~$i$ particle at deme~$x$ at time~$t$ carrying~$k$ mutations reproduces, induces an infected reproduction event, dies, or has a transmission event, and the quantity $\beta^{\mathcal{P},i}(t,x,k)$ is the rate at which a class~$i*$ particle at deme $x$ at time $t$ reproduces or dies.
For $t\in [0,\tau^*]$, let
\begin{align} \label{eq:3Eeventpf}
Y^{\textrm{high}}(t)&=X^{\textrm{high}}\Bigg(\int_0^t 
\Big(\beta^{\mathcal I, i(s-)}(s,Z^n_{p(s-)}(s-),k(s-))\mathds 1_{\{p(s-)\in \mathcal I^n(s-)\}} \notag \\
&\hspace{2.5cm} +\beta^{\mathcal P, i(s-)}(s,Z^n_{p(s-)}(s-),k(s-))\mathds 1_{\{p(s-)\in \mathcal P^n(s-)\}}\Big) \notag \\
&\hspace{5.5cm}
\cdot \mathds{1}_{\left(E^{\textrm{low}}(s-)\right)^c}\mathds{1}_{\{Z^n_{p(s-)}(s-)\in \Lambda_n\}}ds
\Bigg).
\end{align}
The jump times of $Y^{\textrm{high}}$ will be the times at which the particle with label $p(t-)$ is at a high density deme and generates a reproduction or death event.
Finally, for $t\ge 0$, let
\begin{equation} \label{eq:4Eeventpf}
  Y^{\textrm{reinf}}(t)=X^{\textrm{reinf}}\left(q_-(U_\epsilon)t
\right).  
\end{equation}
The jump times of $Y^{\textrm{reinf}}$ will be times at which the particle with label $p(t-)$ may be reinfected (from partially recovered to infected). 

We will use the following bounds on the total rate of events at `low density' demes.
Suppose $i\in \{1,2\}$, $t\ge 0$ and $x\in \Lambda_n$ with $\|\eta^{n,(i)}(t-,x)\|_{\ell_1}\le N U_\epsilon.$
By~\eqref{Paper01_event_more_particles_process_i_than_3_minus_i},~\eqref{Paper01_definition_number_particles_composition_reference_modified_processes} and~\eqref{Paper01_equations_colouring_pink_purple_processes}, on the event $A^i(t-,x)$ we also have $\|\eta^{n,(3-i)}(t-,x)\|_{\ell_1}\le N U_\epsilon.$
Therefore, for $k\in \mathbb N_0$ and $i^+,i^-\in \{1,2\}$, by Definition~\ref{Paper01_definition_cutoff_density_parameter}(ii) combined with~\eqref{Paper01_eq:betainfdefn},~\eqref{Paper01_eq:betainducedefn},~\eqref{Paper01_eq:deltainfdefn},~\eqref{Paper01_eq:deltatransmitdefn} and since $s_k\le 1$ by Assumption~\ref{Paper01_assumption_fitness_sequence},
\begin{align} \label{eq:1lowdensitybound}
    &\beta^i_{\textrm{infected}}(t,x,k)+\beta^{i,i^+}_{\textrm{induce}}(t,x)+\delta^i_{\textrm{infected}}(t,x)+\delta^{i,i^-}_{\textrm{transmit}}(t,x) \notag \\
    &\quad \le q_+(U_\epsilon)+NU_\epsilon q_+ (U_\epsilon)+q_-(U_\epsilon)+NU_\epsilon q_-(U_\epsilon) \notag \\
    &\quad \le V^N_\epsilon,
\end{align}
where the last line follows from~\eqref{Paper01_eq:Vepsdefn}.
Moreover, for $k\in \mathbb N_0$ and $i^\circ \in \{1,2\}$, by Definition~\ref{Paper01_definition_cutoff_density_parameter}(ii) combined with~\eqref{Paper01_eq:betainfdefn},~\eqref{Paper01_eq:deltamindefn},~\eqref{Paper01_eq:deltapartialdefn} and since $s_k\le 1$,
\begin{align} \label{eq:2lowdensitybound}
    \beta^i_{\textrm{infected}}(t,x,k)+\delta_{\textrm{min}}(t,x)+\delta^{i^\circ,i}_{\textrm{partial}}(t,x)
    &\le q_+(U_\epsilon)+q_-(U_\epsilon)+q_-(U_\epsilon) \notag \\
    &\le V^N_\epsilon,
\end{align}
where again the last line follows from~\eqref{Paper01_eq:Vepsdefn}.

We define $\eta^{n,(\textrm{coupling})}(0)$ as in~\eqref{Paper01_initial_condition_coloured_process_infection}.
Then we construct the process $\eta^{n,(\textrm{coupling})}$ on the time interval $[0,T\wedge \tau^*]$ as follows. 

As in the original definition of $\eta^{n,\textrm{(coupling)}}$, whenever a new particle is added to the process, say carrying $k$ mutations, then with probability $\mu$, a mutation occurs and the new particle carries $k+1$ mutations; if $k+1>K_n$ then the particle is immediately deleted. 
\begin{itemize}

\item \textit{Migration of particles.}
At a jump time $t$ in the process $Y^{\textrm{jump}}$ (defined in~\eqref{eq:1Eeventpf}), if $Z^n_{p(t-)}(t-)\in \Lambda_n$ then the particle with label $p(t-)$
makes a jump to a uniformly chosen neighbouring deme in $L^{-1}\mathbb Z$ (with reflection at the boundary of $\Lambda_n$). 
If particle $p(t-)$ has a dual particle at time $t-$, then this dual particle makes the same jump at time $t$.

All other susceptible and infected particles at demes in $\Lambda_n$ make jumps independently at rate $m$ to uniformly chosen neighbouring demes (with reflection at the boundary of $\Lambda_n$). All other dual pairs of partially recovered particles at demes in $ \Lambda_n$ independently make jumps at rate $m$ to uniformly chosen neighbouring demes (with reflection at the boundary of $\Lambda_n$).
\item \textit{Reproduction and death events generated by particle $p(t-)$ at low density demes.}
Suppose $t$ is a jump time in the process $Y^{\textrm{low}}$ (defined in~\eqref{eq:2Eeventpf}). Recalling~\eqref{eq:starinEeventpf}, let $x=Z^n_{p(t-)}(t-)$, $i=i(t-)$ and $k=k(t-)$,
and take $i^\circ,i^-,i^+\in \{1,2\}$ such that
\begin{align*}
\|\eta^{n,i^\circ}(t-,x)\|_{\ell_1} &\ge \|\eta^{n,3-i^\circ}(t-,x)\|_{\ell_1},\\
   q^N_-(\|\eta^{n,(i^-)}(t-,x)\|_{\ell_1}) &\ge q^N_-(\|\eta^{n,(3-i^-)}(t-,x)\|_{\ell_1}),\\
   \text{and}\quad q^N_+(\|\eta^{n,(i^+)}(t-,x)\|_{\ell_1})& \ge q^N_+(\|\eta^{n,(3-i^+)}(t-,x)\|_{\ell_1}).
\end{align*}
Suppose $x\in \Lambda_n$ (otherwise, nothing happens).
If $p(t-)\in \mathcal I^n(t-)$, then (using~\eqref{eq:1lowdensitybound} and~\eqref{eq:**Eeventpf})
\begin{itemize}
    \item with probability $\beta^i_{\textrm{infected}}(t,x,k)/V_\epsilon^N$, the particle $p(t-)$ reproduces, adding a new class $i$ particle carrying $k$ mutations at deme $x$.
    \item with probability $\beta^{i,i^+}_{\textrm{induce}}(t,x)/V_\epsilon^N$, the particle $p(t-)$ induces an infected reproduction event, adding a new class $i^+$ particle at deme $x$ carrying $\kappa$ mutations, where, for each $k' \in [K_{n}]_{0}$, the event $\{\kappa=k'\}$ occurs with probability
    $$ \frac{s_{k'}\eta^{n,0}_{k'}(t-,x)}{\sum_{k'' = 0}^{K_n} s_{k''} \eta^{n,0}_{k''}(t-,x) }.$$
    \item with probability $\delta^i_{\textrm{infected}}(t,x)/V^N_\epsilon$, the particle $p(t-)$ dies.
    \item with probability $\delta^{i,i^-}_{\textrm{transmit}}(t,x)/V^N_\epsilon$, the particle $p(t-)$ has a transmission event: If $i\neq i^-$, then a particle chosen uniformly at random from the class $0$ particles at deme $x$ at time $t-$ is given class $3-i^-$.
    If instead $i=i^-$, then a particle chosen uniformly at random from the class $0$ particles at deme $x$ at time $t-$ is given class $(3-i^-)*$ and becomes a dual pair with particle $p(t-)$, which becomes class $i*$.
    \item otherwise, nothing happens.   
\end{itemize}
If instead $p(t-)\in \mathcal P^n(t-)$, then (using~\eqref{eq:2lowdensitybound} and~\eqref{eq:**Eeventpf})
\begin{itemize}
    \item with probability $\beta^i_{\textrm{infected}}(t,x,k)/V_\epsilon^N$, the particle $p(t-)$ reproduces, adding a new class $i$ particle carrying $k$ mutations at deme $x$.
    \item with probability $\delta_{\textrm{min}}(t,x)/V_\epsilon^N$, the particle $p(t-)$ and its dual particle both die.
    \item with probability $\delta^{i^\circ,i}_{\textrm{partial}}(t,x)/V_\epsilon^N$, particle $p(t-)$ dies and is removed from the process.
    Let $p'$ denote the class $(3-i)*$ dual particle of particle $p(t-)$ at time $t-$.
     If $i^\circ\neq i$,
     then a reinfection occurs: particle $p'$ becomes class $3-i$.
     If instead $i^\circ=i$,
     then a replacement occurs: choose a class $i$ particle at deme $x$ uniformly at random, and this particle becomes the new dual particle of particle $p'$, and becomes class $i*$.
    \item otherwise, nothing happens.
\end{itemize}

\item \textit{Reproduction and death events generated by particle $p(t-)$ at high density demes.}
Suppose $t$ is a jump time in the process $Y^{\textrm{high}}$ (defined in~\eqref{eq:3Eeventpf}). Recalling~\eqref{eq:starinEeventpf}, let $x=Z^n_{p(t-)}(t-)$, $i=i(t-)$ and $k=k(t-)$.

If $p(t-)\in \mathcal I^n(t-)$ and $B^\epsilon_{Y^{\textrm{high}}(t)}=1$, then (using~\eqref{Paper01_estimate_rate_giving_birth_or_inducing_reprod_event_high_density} in Lemma~\ref{Paper01_lem:cutoffconseq},~\eqref{eq:**Eeventpf} and~\eqref{eq:daggerEeventpf})
\begin{itemize}
    \item with probability $\beta^i_{\textrm{infected}}(t,x,k)/(\epsilon \beta^{\mathcal I,i}(t,x,k))$, the particle $p(t-)$ reproduces, adding a new class $i$ particle carrying $k$ mutations at deme $x$.
    \item with probability $\beta^{i,i}_{\textrm{induce}}(t,x)/(\epsilon \beta^{\mathcal I,i}(t,x,k))$, the particle $p(t-)$ induces an infected reproduction event, adding a new class $i$ particle at deme $x$ carrying $\kappa$ mutations, where for each $k'\in [K_n]_0$, the event $\{\kappa=k'\}$ occurs with probability
    $$ \frac{s_{k'}\eta^{n,0}_{k'}(t-,x)}{\sum_{k'' = 0}^{K_n} s_{k''} \eta^{n,0}_{k''}(t-,x) }.$$
    \item with probability 
    \[\left(1-\frac{\beta^i_{\textrm{infected}}(t,x,k)+\beta^{i,i}_{\textrm{induce}}(t,x)}{\epsilon \beta^{\mathcal I,i}(t,x,k)}\right)\frac{\delta^i_{\textrm{infected}}(t,x)}{\delta^i_{\textrm{infected}}(t,x)+\delta^{i,i}_{\textrm{transmit}}(t,x)},\]
    the particle $p(t-)$ dies.
    \item with probability
    \[\left(1-\frac{\beta^i_{\textrm{infected}}(t,x,k)+\beta^{i,i}_{\textrm{induce}}(t,x)}{\epsilon \beta^{\mathcal I,i}(t,x,k)}\right)\frac{\delta^{i,i}_{\textrm{transmit}}(t,x)}{\delta^i_{\textrm{infected}}(t,x)+\delta^{i,i}_{\textrm{transmit}}(t,x)},\]
    the particle $p(t-)$ has a transmission event with partial recovery: a particle chosen uniformly at random from the class $0$ particles at deme $x$ at time $t-$ is given class $(3-i)*$ and becomes a dual pair with particle $p(t-)$, which becomes class $i*$.
\end{itemize}
If $p(t-)\in \mathcal I^n(t-)$ and $B^\epsilon_{Y^{\textrm{high}}(t)}=0$, then
\begin{itemize}
        \item with probability $\delta^i_{\textrm{infected}}(t,x)/(\delta^i_{\textrm{infected}}(t,x)+\delta^{i,i}_{\textrm{transmit}}(t,x))$, the particle $p(t-)$ dies.
    \item with probability $\delta^{i,i}_{\textrm{transmit}}(t,x)/(\delta^i_{\textrm{infected}}(t,x)+\delta^{i,i}_{\textrm{transmit}}(t,x))$, the particle $p(t-)$ has a transmission event with partial recovery: a particle chosen uniformly at random from the class $0$ particles at deme $x$ at time $t-$ is given class $(3-i)*$ and becomes a dual pair with particle $p(t-)$, which becomes class $i*$.
    \end{itemize}   
If $p(t-)\in \mathcal P^n(t-)$ and $B^\epsilon_{Y^{\textrm{high}}(t)}=1$, then (using~\eqref{Paper01_estimate_rate_only_giving_birth_high_density} in Lemma~\ref{Paper01_lem:cutoffconseq},~\eqref{eq:**Eeventpf} and~\eqref{eq:daggerEeventpf})
\begin{itemize}
    \item with probability $\beta^i_{\textrm{infected}}(t,x,k)/(\epsilon \beta^{\mathcal P,i}(t,x,k))$, the particle $p(t-)$ reproduces, adding a new class $i$ particle carrying $k$ mutations at deme $x$.
    \item with probability 
    \[\left(1-\frac{\beta^i_{\textrm{infected}}(t,x,k)}{\epsilon \beta^{\mathcal P,i}(t,x,k)}\right)\frac{\delta_{\textrm{min}}(t,x)}{\delta_{\textrm{min}}(t,x)+\delta^{i,i}_{\textrm{partial}}(t,x)},\]
    the particle $p(t-)$ and its dual particle both die.
    \item with probability
    \[\left(1-\frac{\beta^i_{\textrm{infected}}(t,x,k)}{\epsilon \beta^{\mathcal P,i}(t,x,k)}\right)\frac{\delta^{i,i}_{\textrm{partial}}(t,x)}{\delta_{\textrm{min}}(t,x)+\delta^{i,i}_{\textrm{partial}}(t,x)},\]
    particle $p(t-)$ dies and is removed from the process.
    Let $p'$ denote the class $(3-i)*$ dual particle of particle $p(t-)$ at time $t-$.
     A class $i$ particle at deme $x$ chosen uniformly at random becomes the new dual particle of particle $p'$, and becomes class $i*$.
\end{itemize}
If $p(t-)\in \mathcal P^n(t-)$ and $B^\epsilon_{Y^{\textrm{high}}(t)}=0$, then
\begin{itemize}
        \item with probability $\delta_{\textrm{min}}(t,x)/(\delta_{\textrm{min}}(t,x)+\delta^{i,i}_{\textrm{partial}}(t,x))$, the particle $p(t-)$ and its dual particle both die.
    \item with probability $\delta^{i,i}_{\textrm{partial}}(t,x)/(\delta_{\textrm{min}}(t,x)+\delta^{i,i}_{\textrm{partial}}(t,x))$, the particle $p(t-)$ dies and is removed from the process.
    Let $p'$ denote the class $(3-i)*$ dual particle of particle $p(t-)$ at time $t-$.
     A class $i$ particle at deme $x$ chosen uniformly at random becomes the new dual particle of particle $p'$, and becomes class $i*$.
    \end{itemize}   

\item \textit{Reinfection events of particle $p(t-)$.}
Suppose $t$ is a jump time in the process $Y^{\textrm{reinf}}$ (defined in~\eqref{eq:4Eeventpf}). Recalling~\eqref{eq:starinEeventpf},
let $i=i(t-)$ and $x=Z^n_{p(t-)}(t-)$.
If $p(t-)\in \mathcal P^n(t-)$ and $x \in  \Lambda_n$, then let $p'$ denote the dual particle of $p(t-)$ at time $t-$;
with probability
$
\delta^{i,3-i}_{\textrm{partial}}(t,x)/q_-(U_\epsilon),
$
particle $p'$ dies and particle $p(t-)$ is reinfected, i.e.~becomes class $i$.
Otherwise, nothing happens.
(Here, we used Lemma~\ref{Paper01_lem:cutoffconseq}(i), which implies that $\delta^{i,3-i}_{\textrm{partial}}(t,x)=0$ if $\|\eta^{n,(3-i)}(t-,x)\|_{\ell_1}\ge NU_\epsilon$, and~\eqref{Paper01_eq:deltapartialdefn} combined with Definition~\ref{Paper01_definition_cutoff_density_parameter}(ii).)

\item \textit{Other events for the dual particle of particle $p(t-)$.}
\begin{itemize}
    \item At rate
    $\delta^{3-i(t-),3-i(t-)}_{\textrm{partial}}(t,Z^n_{p(t-)}(t-))\mathds{1}_{\{p(t-)\in \mathcal P^n(t-)\}}\mathds{1}_{\{Z^n_{p(t-)}(t-)\in \Lambda_n\}}$,
    the dual particle of particle $p(t-)$ dies and a class
$3-i(t-)$
    particle, chosen uniformly at random from the class
$3-i(t-)$
    particles at deme $x$ at time $t-$, becomes the new dual particle of particle $p(t-)$, and becomes class
    $(3-i(t-))*$.
    \item For $t\ge 0$ such that $p(t)\in \mathcal P^n(t)$, let $k'(t)$ denote the number of mutations carried by the dual particle of particle $p(t)$ at time $t$.
    
    At rate 
    $\beta^{i(t-)}_{\textrm{infected}}(t,Z^n_{p(t-)}(t-),k'(t-))\mathds 1_{\{p(t-)\in \mathcal P^n(t-)\}}\mathds{1}_{\{Z^n_{p(t-)}(t-)\in \Lambda_n\}},$
    the dual particle of particle $p(t-)$ reproduces, adding a new class $i(t-)$ particle carrying $k'(t-)$ mutations at deme $Z^n_{p(t-)}(t-)$.
\end{itemize}

\item \textit{Reproduction and death events for other particles.}
    Other particles  and dual pairs (particles that are not $p(t-)$ or its dual particle) reproduce, induce reproduction events, die, and have transmission events independently, as in the original definition of the process.
\end{itemize}

Using Lemma~\ref{Paper01_lem:cutoffconseq}(i) to see that all possible reproduction and death events generated by $p(t-)$ at high density demes are included in the construction, we see that this process is a realisation of $\eta^{n,(\textrm{coupling})}$ on the time interval $[0,T\wedge \tau^*]$.

Now note that under this construction, on the event $\{u\in \mathcal N^n(T)\}$ we have $\tau^*\ge T$ and so, by~\eqref{eq:Jbulletlist} and~\eqref{eq:Zbulletlist}, and then by our construction,
\begin{equation} \label{Paper01_eq:Jnbound}
    J^n_u(T)=  \#\{t\le T:Z^n_{p(t)}(t-)\neq Z^n_{p(t)}(t)\}\le Y^{\textrm{jump}}(T)=X^{\textrm{jump}}(mT),
\end{equation}
by~\eqref{eq:1Eeventpf}.
Moreover, by our construction of the reinfection times of particle $p(t-)$, we have
\[
\{t\in [0,\tau^*\wedge T]:p(t)=p(t-), p(t)\in \mathcal I^n(t), p(t)\in \mathcal P^n(t-)\}
\subseteq \{t\in [0,T]:Y^{\textrm{reinf}}(t)\neq Y^{\textrm{reinf}}(t-)\},
\]
and therefore on the event $\{u\in \mathcal N^n(T)\}$, by~\eqref{eq:Rbulletlist} and then~\eqref{eq:4Eeventpf} we have
\begin{equation} \label{Paper01_eq:Rnbound}
    R^n_u\le Y^{\textrm{reinf}}(T)=X^{\textrm{reinf}}(q_-(U_\epsilon)T).
\end{equation}
Recall the definition of $\mathcal D_u$ in~\eqref{Paper01_set_labels_particles_generated_same_ancestral_path};
on the event $\{u\in \mathcal N^n(T)\}$, by~\eqref{eq:gammaDbulletlist}, and then by our construction of reproduction events generated by particle $p(t-)$ at low density demes, and finally by~\eqref{eq:2Eeventpf},
we have
\begin{align} \label{Paper01_eq:Xlowbound}
    &\#\Big\{v\in \mathcal D_u: \gamma^{n,\textrm{density},\epsilon}_{v}=0\Big\} \notag \\
    &\le \# \Big\{t\le T:\|\eta^{n,(\gamma^{n,\textrm{class}}_{p(t-)})}(t-, Z^n_{p(t-)}(t-))\|_{\ell_1}< N U_\epsilon
    \text{ and }t=\tau^{n,\textrm{birth}}_{p(t-)j}\text{ for some }j\in \mathbb N \Big\}\notag \\
    &\le Y^{\textrm{low}}(T)\notag \\
    &\le X^{\textrm{low}}(V^N_\epsilon T).
\end{align}
By our construction of reproduction and death events generated by particle $p(t-)$ at high density demes,
if $t\le \tau^*$ is a jump time in $Y^{\textrm{high}}$, then $p(t-)$ either dies or generates a new particle at time $t$.
Therefore,
on the event $\{u\in \mathcal N^n(T)\}$, using~\eqref{eq:3Eeventpf},~\eqref{eq:**Eeventpf} and~\eqref{eq:gammaDbulletlist}
we can write
\begin{equation} \label{eq:*Yhigh}
    \{t\in [0,\tau^{n,\textrm{birth}}_u]:Y^{\textrm{high}}(t)\neq Y^{\textrm{high}}(t-)\}
    =\{ \tau^{n,\textrm{birth}}_{v}: v \in \mathcal D_u, \gamma^{n,\textrm{density},\epsilon}_{v}=1\}.
\end{equation}
For $v = u\vert_{i}j \in \mathcal D_u$ with $d_{v}=1$, we denote
\begin{equation} \label{eq:stardvtotal}
    d^{\textrm{total}}_{v} \defeq \sum_{w \in \{u\vert_{i'}j'\in \mathcal D_u:i'<i\}\cup \{u\vert_{i}j':j'\le j\}}d_{w}.
\end{equation}
Then, recalling~\eqref{Paper01_event_genealogical_path_considering_the_number_reinfection_events} and using~\eqref{eq:*Yhigh},
on the event $E^{n,\epsilon}_{u,J,R,\gamma^-,\gamma^+,d}(T)$,
for $v \in \mathcal D_u$ with $d_{v}=1$, we have
\[
Y^{\textrm{high}}(\tau^{n,\textrm{birth}}_{v})=d^{\textrm{total}}_{v}.
\]
Moreover, on the event $E^{n,\epsilon}_{u,J,R,\gamma^-,\gamma^+,d}(T)$, for $v=u|_i j \in \mathcal D_u$ with $d_{v}=1$ and $(\gamma^-_{v},\gamma^+_{v})\neq (1,0)$, we must have
$B^\epsilon_{Y^{\textrm{high}}(\tau^{n,\textrm{birth}}_{v})}=1$
(because particle $u\vert_i$ cannot die before generating $v$, and if $\tau$ is a jump time for $Y^{\textrm{high}}$ with $B^\epsilon_{Y^{\textrm{high}}(\tau)}=0$ then by our construction, either particle $p(\tau-)$ dies at time $\tau$ or $p(\tau-)\in \mathcal I^n(\tau-)$ and $p(\tau)\in \mathcal P^n(\tau)$, which is not possible by~\eqref{eq:gammaIbulletlist}). 
Therefore, on the event $E^{n,\epsilon}_{u,J,R,\gamma^-,\gamma^+,d}(T)$, we have
\begin{equation} \label{Paper01_eq:B1highdensity}
    B^\epsilon_{d^{\textrm{total}}_{v}}=1\; \forall v\in \mathcal D_u\text{ with }d_{v}=1\text{ and }(\gamma^-_{v},\gamma^+_{v})\neq (1,0).
\end{equation}
Combining~\eqref{Paper01_eq:Jnbound},~\eqref{Paper01_eq:Rnbound},~\eqref{Paper01_eq:Xlowbound} and~\eqref{Paper01_eq:B1highdensity},
we now have that
\begin{equation} \label{Paper01_alternative_formulation_event}
\begin{aligned}
&E^{n,\epsilon}_{u,J,R,\gamma^-,\gamma^+,d}(T)\\
&\quad \subseteq 
\{X^{\textrm{jump}}(mT)\ge J\}
\cap \{X^{\textrm{reinf}}(q_-(U_\epsilon)T)\ge R\}\\
&\qquad \cap \{ X^{\textrm{low}}(V^N_\epsilon T)\ge \#\{v \in \mathcal D_u: d_{v}=0\}\}
\cap \bigcap_{\{v\in \mathcal D_u : d_{v}=1\text{ and }(\gamma^-_{v},\gamma^+_{v})\neq (1,0)\}}\{B^\epsilon_{d^{\textrm{total}}_{v}}=1\}.
\end{aligned}
\end{equation}
Note that by~\eqref{eq:stardvtotal}, for $v \neq v' \in \mathcal D_u$ with $d_{v}=1=d_{v'}$, we have $d^{\textrm{total}}_{v}\neq d^{\textrm{total}}_{v'}$.
Since, by construction at the start of the proof, $(X^{\textrm{jump}}(t),t\ge 0)$, $(X^{\textrm{low}}(t),t\ge 0)$, $(X^{\textrm{high}}(t),t\ge 0)$, and $(X^{\textrm{reinf}}(t),t\ge 0)$ are i.i.d.~rate $1$ Poisson processes, and $(B^\epsilon_i)_{i=1}^\infty$ is an independent sequence of i.i.d.~Bernoulli($\epsilon$) random variables,~\eqref{Paper01_upper_bound_probability_one_infected_particle_living_too_far_from_original_point} follows directly from~\eqref{Paper01_alternative_formulation_event}.

It remains to establish that~\eqref{Paper01_deterministic_condition_path_genealogical_infection} implies~\eqref{Paper01_event_cannot_happen_if_number_reinfections_too_low}. On the event $E^{n,\epsilon}_{u,J,R,\gamma^-,\gamma^+,d}(T)$,
by~\eqref{eq:gammaIbulletlist}, for
each consecutive pair of times $s_1<s_2$ in the set
\[
\{\tau^{n,\textrm{birth}}_{v}:v \in \mathcal D_u, (\gamma^-_{v},\gamma^+_{v})=(1,0)\},
\]
we have $p(s_1)\in \mathcal P^n(s_1)$ and $p(s_2-)\in \mathcal I^n(s_2-)$, and so
since $p(t)\in \mathcal N^n(t)$ $\forall t\in (s_1,s_2)$,
there exists a time $t\in (s_1,s_2)$ such that $p(t-)\in \mathcal P^n(t-)$ and $p(t)\in \mathcal I^n(t)$.
Hence, since $R^n_u=R$ and using~\eqref{eq:Rbulletlist} and~\eqref{eq:gammaIbulletlist}, we must have
\begin{align*}
\#\{v\in \mathcal D_u:(\gamma^-_{v},\gamma^+_{v})=(1,0)\}-1
&\le R+\#\{v\in \mathcal D_u:(\gamma^-_{v},\gamma^+_{v})=(0,1)\}\\
&\le R+\#\mathcal D_u-\#\{v\in \mathcal D_u:(\gamma^-_{v},\gamma^+_{v})=(1,0)\}.
\end{align*}
Therefore,
\begin{equation} \label{Paper01_lower_bound_number_reinfection_events_that_must_happen}
\begin{aligned}
    R&\ge 2\#\{v\in \mathcal D_u:(\gamma^-_{v},\gamma^+_{v})=(1,0)\}-\#\mathcal D_u-1\\
    &\ge 2\Big(\# \mathcal D_u-\#\{v\in \mathcal D_u:d_{v}=0\}-\#\{v\in \mathcal D_u:d_{v}=1,(\gamma^-_{v},\gamma^+_{v})\neq (1,0)\}\Big)-\# \mathcal D_u-1\\
    &=\# \mathcal D_u -2\Big(\#\{v\in \mathcal D_u:d_{v}=0\}+\#\{v\in \mathcal D_u:d_{v}=1,(\gamma^-_{v},\gamma^+_{v})\neq (1,0)\}\Big)-1.
\end{aligned}
\end{equation}
As observed after~\eqref{Paper01_set_labels_particles_generated_same_ancestral_path}, we have $\# \mathcal{D}_{u} = \vert u \vert$, and so by~\eqref{Paper01_lower_bound_number_reinfection_events_that_must_happen}, on the event $E^{n,\epsilon}_{u,J,R,\gamma^-,\gamma^+,d}(T)$, we cannot have~\eqref{Paper01_deterministic_condition_path_genealogical_infection}, which completes the proof.
\end{proof}

For $\ell\in \mathbb N_0$, $u=u_0u_1\ldots u_\ell \in \mathcal U$, $J\in \mathbb N_0$,
and $T\ge 0$,
define the event
\begin{align} \label{Paper01_event_infected_partially_recovered_particle_travelled_long_distance}
    E^{n}_{u,J}(T)
    &=\Big\{
    u\in \mathcal N^n(T),\; J^n_{u}(T)\ge J
    \Big\}.
\end{align}
Using Lemma~\ref{Paper01_lem:boundPE} and a union bound, we can bound the probability of the event $E_{u,J}^n$.
Recall the definition of $|u|$ in~\eqref{Paper01_length_branch_ulam_harris_tree}.
\begin{corollary} \label{Paper01_corollary_almost_final_bound_spread_infection}
For $\boldsymbol{\eta}^{(1)},\boldsymbol{\eta}^{(2)} \in \mathcal{S}$ with $\#\Delta(\boldsymbol{\eta}^{(1)},\boldsymbol{\eta}^{(2)})<\infty$, $\delta \in (0,1)$ and $T^*\ge 0$, there exists $C_{\delta,N,q_-,q_+,T^*}>0$ such that for $n\in \mathbb N$, $u\in \mathcal U$, $J\in \mathbb N$ and $T\in [0,T^*]$,
    \[
\mathbb P_{\boldsymbol{\eta}^{(1)},\boldsymbol{\eta}^{(2)}}\left(E^{n}_{u,J}(T) \right)\le  C_{\delta,N,q_-,q_+,T^*}\delta^{|u|}\left( e^{J\left(1-\log\left(\frac J {mT}\right)\right)}\mathds 1_{\{J\ge mT\}}+\mathds 1_{\{J< mT\}}\right).   
    \]
\end{corollary}
\begin{proof}
    For $\epsilon\in (0,1]$, 
$\gamma^-=(\gamma^-_{v})_{v\in \mathcal D_u} \in \{0,1\}^{\mathcal D_u}$,
$\gamma^+=(\gamma^+_{v})_{v\in \mathcal D_u} \in \{0,1\}^{\mathcal D_u}$, and
$d=(d_{v})_{v \in \mathcal D_u} \in \{0,1\}^{\mathcal D_u}$,
define the event
\begin{equation} \label{Paper01_event_genealogical_path_without_considering_the_number_reinfection_events}
    E^{n,\epsilon}_{u,J,\gamma^-,\gamma^+,d}(T)
    =\bigcup_{R\in \mathbb N_0}E^{n,\epsilon}_{u,J,R,\gamma^-,\gamma^+,d}(T).
\end{equation}
For $\alpha \ge 0$, let $Z_\alpha$ denote a Poisson random variable with mean $\alpha$. Let
\begin{equation} \label{eq:*D0D1defn}
   D_0:=\#\{v \in \mathcal D_u: d_{v}=0\}
\quad \text{and}\quad 
D_1:=\#\{v \in \mathcal D_u:d_{v}=1,(\gamma^-_{v},\gamma^+_{v})\neq (1,0)\}, 
\end{equation}
let $R_0=(|u|-1-2(D_0+D_1))^+$,
and recall the definition of $V^N_\epsilon$ in~\eqref{Paper01_eq:Vepsdefn}. 
By a union bound, and then by~\eqref{Paper01_upper_bound_probability_one_infected_particle_living_too_far_from_original_point} and~\eqref{Paper01_event_cannot_happen_if_number_reinfections_too_low} in Lemma~\ref{Paper01_lem:boundPE}, we have
\begin{align} \label{eq:1probEunionbd}
    \mathbb P_{\boldsymbol{\eta}^{(1)},\boldsymbol{\eta}^{(2)}}\left(E^{n,\epsilon}_{u,J,\gamma^-,\gamma^+,d}(T)
    \right)
    &\le \sum_{R=0}^\infty\mathbb P_{\boldsymbol{\eta}^{(1)},\boldsymbol{\eta}^{(2)}}\left(E^{n,\epsilon}_{u,J,R,\gamma^-,\gamma^+,d}(T) \right) \notag \\
    &\le \mathbb P\left(Z_{mT}\ge J\right)\mathbb P\left(Z_{V^N_\epsilon T}\ge D_0\right)
    \epsilon^{D_1}
    \sum_{R=R_0}^\infty\mathbb P\left(Z_{V_\epsilon^N T}\ge R\right).
\end{align}
Then to bound the sum on the right-hand side,
\begin{align} \label{eq:2sumprobZ}
    \sum_{R=R_0}^\infty \mathbb P\left(Z_{V_\epsilon^NT}\ge R\right)
    &=\mathbb E\left[ (Z_{V_\epsilon^NT} -R_0+1) \mathds 1_{\{ Z_{V_\epsilon^NT} \ge R_0 \}}
    \right] \notag \\
    &\le \mathbb E\left[ (Z_{V_\epsilon^NT}+1 )^2\right]^{1/2} \mathbb P\left(Z_{V_\epsilon^N T} \ge R_0\right)^{1/2},
\end{align}
where the second inequality follows by the Cauchy-Schwarz inequality.
Note that by~\eqref{eq:*D0D1defn} we have
\[
\max\left( D_0, D_1, R_0  \right)\ge \tfrac 16 |u|-1.
\]
Therefore, by~\eqref{eq:1probEunionbd} and~\eqref{eq:2sumprobZ} we can write
\begin{equation} \label{eq:3probEintermsZ}
\begin{aligned}
    & \mathbb P_{\boldsymbol{\eta}^{(1)},\boldsymbol{\eta}^{(2)}}\left(E^{n,\epsilon}_{u,J,\gamma^-,\gamma^+,d}(T)
    \right) \\ 
    & \quad \le \mathbb P\left(Z_{mT}\ge J\right)\mathbb E\left[ (Z_{V_\epsilon^N T}+1 )^2\right]^{1/2} \left(\mathbb P\left(Z_{V_\epsilon^N T}\ge \tfrac 16 |u|-1\right)^{1/2} \vee 
    \left(\epsilon^{\frac 16 |u|-1}\right)\right).
\end{aligned}
\end{equation}
By the Chernoff bound in Lemma~\ref{Paper01_lemma_upper_tail_poisson}, 
\begin{align*}
 \mathbb P\left(Z_{V_\epsilon^N T}\ge \tfrac 16 |u|-1\right)^{1/2}
 &\le \mathds 1_{\{\frac 16 |u|-1\ge V^N_\epsilon T\}}(e^{1-\log((\frac 16 |u|-1)/(V^N_\epsilon T))})^{\frac 12 (\frac 16 |u|-1)}
 +\mathds 1_{\{\frac 16 |u|-1< V^N_\epsilon T\}}\\
 &\le \mathds 1_{\{|u|\ge C^{(1)}_{\epsilon,N,q_-,q_+,T^*}\}}\epsilon^{\frac 16 |u|-1}
 +\mathds 1_{\{|u|< C^{(1)}_{\epsilon,N,q_-,q_+,T^*}\}}
\end{align*}
for some $C^{(1)}_{\epsilon,N,q_-,q_+,T^*}>0$ sufficiently large.
Therefore, by~\eqref{eq:3probEintermsZ}, there exists $C^{(2)}_{\epsilon,N,q_-,q_+,T^*}>0$ such that
\[
\mathbb P_{\boldsymbol{\eta}^{(1)},\boldsymbol{\eta}^{(2)}}\left(E^{n,\epsilon}_{u,J,\gamma^-,\gamma^+,d}(T)
    \right)
    \le \mathbb P\left(Z_{mT}\ge J\right)C^{(2)}_{\epsilon,N,q_-,q_+,T^*}
    \epsilon^{\frac 16 |u|}.
\]
By~\eqref{Paper01_event_infected_partially_recovered_particle_travelled_long_distance},~\eqref{Paper01_event_genealogical_path_considering_the_number_reinfection_events},~\eqref{Paper01_event_genealogical_path_without_considering_the_number_reinfection_events} and a union bound
over all possible choices of $\gamma^-$, $\gamma^+$ and $d$, 
and using that $\#\mathcal D_u=|u|$ as observed after~\eqref{Paper01_set_labels_particles_generated_same_ancestral_path},
it follows that 
    \begin{align*}
    \mathbb P_{\boldsymbol{\eta}^{(1)},\boldsymbol{\eta}^{(2)}}\left(E^{n}_{u,J}(T)
    \right)
    &\le 2^{3|u|}\mathbb P\left(Z_{mT}\ge J\right)C^{(2)}_{\epsilon,N,q_-,q_+,T^*}
    \epsilon^{\frac 16 |u|}.
\end{align*}
By setting $\epsilon=(\delta/8)^6$ so that $2^{3|u|}\epsilon^{\frac 16 |u|}=\delta^{|u|}$,
and applying Lemma~\ref{Paper01_lemma_upper_tail_poisson}
again to bound $\mathbb P(Z_{mT}\ge J)$, the result follows.
\end{proof}

We are now finally ready to prove Proposition~\ref{Paper01_proposition_bound_spread_infection}.

\begin{proof}[Proof of Proposition~\ref{Paper01_proposition_bound_spread_infection}]
    Fix $r > 0$ and $T > 0$, and take $R > r$  and $\boldsymbol{\eta}^{(1)}=(\eta^{(1)}_k(x))_{k\in \mathbb N_0,\, x\in L^{-1}\mathbb Z},$ $ \boldsymbol{\eta}^{(2)}=(\eta^{(2)}_k(x))_{k\in \mathbb N_0,\, x\in L^{-1}\mathbb Z} \in \mathcal{S}$ with $\# \Delta(\boldsymbol{\eta}^{(1)}, \boldsymbol{\eta}^{(2)}) < \infty$ and
    \begin{equation} \label{eq:starDeltaspbd}
        \Delta_{\textrm{sp}}(\boldsymbol{\eta}^{(1)}, \boldsymbol{\eta}^{(2)}) \subset (- \infty, - R] \, \cup \, [R, \infty),
    \end{equation}
    where $\Delta(\boldsymbol{\eta}^{(1)}, \boldsymbol{\eta}^{(2)})$ is defined in~\eqref{Paper01_difference_set_reference_modified}, and $\Delta_{\textrm{sp}}(\boldsymbol{\eta}^{(1)}, \boldsymbol{\eta}^{(2)})$ is defined in~\eqref{Paper01_spatial_projection_difference_set_reference_modified}. 
    By Lemma~\ref{Paper01_coupling_lemma_between_coloured_and_uncoloured_processes}, under $\mathbb P_{\boldsymbol{\eta}^{(1)},\boldsymbol{\eta}^{(2)}}$, the processes $(\eta^{n,(1)}(t))_{t\ge 0}$ and $(\eta^{n,(2)}(t))_{t\ge 0}$ are realisations of the Markov process associated to $\mathcal L^n$ with $\eta^{n,(1)}(0)=\boldsymbol{\eta}^{(1)}$ and $\eta^{n,(2)}(0)=\boldsymbol{\eta}^{(2)}$.
    Let
    \begin{equation} \label{eq:*Delta0defn}
        \Delta_0=\sum_{(x,k)\in \Delta(\boldsymbol{\eta}^{(1)},\boldsymbol{\eta}^{(2)})}|\eta^{(1)}_k(x)-\eta^{(2)}_k(x)|,
    \end{equation}
    and recall from~\eqref{Paper01_initial_condition_coloured_process_infection} that under $\mathbb P_{\boldsymbol{\eta}^{(1)},\boldsymbol{\eta}^{(2)}}$, almost surely the total number of infected and partially recovered particles at time $0$ is given by $\Delta_0$.
    Therefore, under the labelling defined after~\eqref{Paper01_ulam-harris-scheme}, under $\mathbb P_{\boldsymbol{\eta}^{(1)},\boldsymbol{\eta}^{(2)}}$ almost surely
    \begin{equation} \label{eq:1Nntisin}
        \mathcal N^n(t)\subseteq \{u=u_0 u_1 \ldots u_{\ell}\in \mathcal U: \ell \in \mathbb N_0, u_0\le \Delta_0 \} \quad \forall t\ge 0.
    \end{equation}
    Moreover, for $t\ge 0$, by~\eqref{Paper01_spatial_projection_difference_set_reference_modified} and~\eqref{Paper01_definition_number_particles_composition_reference_modified_processes}, and since $\mathcal N^n(t)$ contains all the class $1$, $1*$, $2$ and $2*$ particles alive at time $t$, and then by~\eqref{eq:Zbulletlist},~\eqref{eq:starDeltaspbd} and~\eqref{Paper01_initial_condition_coloured_process_infection}, and finally by~\eqref{eq:Jbulletlist},~\eqref{Paper01_event_infected_partially_recovered_particle_travelled_long_distance} and since migration jumps have size $L^{-1}$, we have
    \begin{align} \label{eq:2DeltaspsubsetE}
    \{\Delta_{sp}(\eta^{n,(1)}(t), \eta^{n,(2)}(t)) \cap [-r,r] \neq \emptyset\}
    &=\{\exists u\in \mathcal N^n(t):Z^n_u(t)\in [-r,r]\} \notag \\
    &\subseteq \{\exists u\in \mathcal N^n(t):|Z^n_u(t)-Z^n_u(0)|\ge R-r\} \notag \\
    &\subseteq 
        \bigcup_{u \in \mathcal{U}} E^{n}_{u,\lceil (R-r)L\rceil }(t).
    \end{align}
    Let $R_{c} = e^2 mT/L + r$, so that $R > R_{c}$ implies that $\lceil (R - r)L \rceil \ge e^2mT$. By~\eqref{eq:1Nntisin} and~\eqref{eq:2DeltaspsubsetE} combined with Corollary~\ref{Paper01_corollary_almost_final_bound_spread_infection} and a union bound, for $\delta \in (0,1)$,~$t \in [0,T]$ and $R > R_{c}$, we have
    \begin{equation} \label{Paper01_almost_there_estimate_bound_spread_infection}
    \begin{aligned}
        &\mathbb{P}_{\boldsymbol{\eta}^{(1)}, \boldsymbol{\eta}^{(2)}}\Big(\Delta_{sp}(\eta^{n,(1)}(t), \eta^{n,(2)}(t)) \cap [-r,r] \neq \emptyset\Big)\\
        & \quad \leq C_{\delta,N,q_-,q_+,T} \sum_{\{u =u_0u_1\ldots u_{\ell}\in \mathcal{U}:\ell \in \mathbb N_0, u_0\le \Delta_0\}}  \delta^{|u|} e^{-(R-r)L} \\ 
        & \quad = C_{\delta,N,q_-,q_+,T} \sum_{h \in \mathbb{N}_{0}} \, \sum_{\{u =u_0\ldots u_{\ell}\in \mathcal{U}: \ell \in \mathbb N_0, u_0\le \Delta_0,\vert u \vert = h\}}  \delta^{h}e^{-(R-r)L}.
    \end{aligned}
    \end{equation}
    For $h\in \mathbb N_0$, by~\eqref{Paper01_ulam-harris-scheme} and~\eqref{Paper01_length_branch_ulam_harris_tree}, 
    and since for $k\in \mathbb N$, the number of ordered partitions of $k$ into positive integers is given by $2^{k-1}$, we have
    \begin{equation} \label{Paper01_trivial_combinatorics}
        \# \{u =u_0\ldots u_{\ell}\in \mathcal{U}: \ell \in \mathbb N_0, u_0\le \Delta_0,\vert u \vert = h\}
        =\Delta_0 (\mathds 1_{\{h=0\}}+\mathds 1_{\{h\ge 1\}}2^{h-1}). 
    \end{equation}
    Hence, the result follows by taking~$\delta \in (0,1/2)$, and then by applying~\eqref{Paper01_trivial_combinatorics} and~\eqref{eq:*Delta0defn} to~\eqref{Paper01_almost_there_estimate_bound_spread_infection}.
    \end{proof}

\appendix

\section{Appendix} \label{Paper01_subsection_preliminaries}

In this appendix we will briefly explain concepts, and state and prove results, that are used in Sections~\ref{Paper01_section_functional_analysis_feller_non_locally_compact}-\ref{Paper01_section_spread_infection}.

\subsection{Integral with respect to Poisson point processes} \label{Paper01_subsection_integral_poisson}

In this subsection, we briefly discuss some results regarding integration of càdlàg predictable processes with respect to Poisson point processes. We refer the interested reader to~\cite{bass2004stochastic} and~\cite{jacod2013limit} as some references on stochastic calculus of processes with jumps. Let $(\Omega, \mathcal{F}, \mathbb{P})$ be some probability space equipped with a filtration $\{\mathcal{F}_{t}\}_{t \geq 0}$ which satisfies the usual conditions, and let $\mathscr{R}_+$ be the Borel $\sigma$-algebra on $[0,\infty)$. We say that an integer-valued random measure $\mathcal{P}$ on $[0, \infty) \times [0, \infty)$ is a \emph{Poisson point process} with \emph{intensity measure} given by the Lebesgue measure $\lambda$ if the following properties are satisfied:
\begin{enumerate}[(i)]
    \item For any $\mathscr{A} \in \mathscr{R}_{+} \otimes \mathscr{R}_{+}$ with $\lambda(\mathscr{A}) < \infty$, we have $\mathbb{E}\left[\mathcal{P}(\mathscr{A})\right] = \lambda(\mathscr{A})$;
    \item For every $t \in [0, \infty)$ and $\mathscr{A} \in \mathscr{R}_{+} \otimes \mathscr{R}_{+}$ such that $\mathscr{A} \subset (t, \infty) \times [0,\infty)$ and $\lambda(\mathscr{A}) <  \infty$, the random variable $\mathcal{P}(\mathscr A)$ is independent of the $\sigma$-algebra $\mathcal{F}_{t}$;
    \item For each $\mathscr{A} \in \mathscr{R}_{+}$ with $\lambda(\mathscr{A}) < \infty$, the process $\left(\mathcal{P}\left([0,t] \times \mathscr{A}\right)\right)_{t \geq 0}$ is a Poisson process with rate $\lambda(\mathscr{A})$;
    \item If $\left(\mathscr{A}_{i}\right)_{i = 1}^{n}$ is an $n$-tuple of pairwise disjoint sets in $\mathscr{R}_+$ with $\lambda(\mathscr{A}_{i}) < \infty$ for each $i$, then the processes $\left(\left(\mathcal{P}([0,t] \times \mathscr{A}_{i})\right)_{t \geq 0}\right)_{i = 1}^{n}$ are independent.
\end{enumerate}

Let $H: \Omega \times [0,\infty) \times [0,\infty) \rightarrow \mathbb{R}$ be some càdlàg process adapted to the filtration $\{\mathcal{F}_{t}\}_{t \geq 0}$ such that
\begin{equation*}
    \mathbb{E}\left[\int_{0}^{\infty} \int_{0}^{\infty} \vert H(t-,y) \vert \, dt \, dy\right] < \infty \quad \textrm{and} \quad \mathbb{E}\left[\int_{0}^{\infty} \int_{0}^{\infty} H^{2}(t-,y) \, dt \, dy\right] < \infty;
\end{equation*}
then the process $\left(X(t)\right)_{t \geq 0} \defeq \left(\int_{[0,t]\times [0,\infty)} H(\tau-, y) \; \mathcal{P}(d\tau \times dy) \right)_{t \geq 0}$ is a well-defined semimartingale. If we denote the martingale and finite variation parts of $X$ by $(M(t))_{t \geq 0}$ and $(A(t))_{t \geq 0}$ respectively, then we have for all $t \geq 0$ that
\begin{equation} \label{Paper01_decomposition_semimartingale_integral_jumps}
    M(t) \defeq \int_{[0,t]\times [0,\infty)} H(\tau-, y) (\mathcal{P} - \lambda)(d\tau \times dy) \quad \textrm{and} \quad A(t) \defeq \int_{0}^{t} \int_{0}^{\infty} H(\tau-, y) \,  d\tau \, dy.
\end{equation}
We say that $M$ is the stochastic integral with respect to the \emph{compensated} Poisson process $\mathcal{P} - \lambda$. The quadratic variation $([M](t))_{t \geq 0}$ and the predictable bracket process $(\langle M \rangle (t))_{t \geq 0}$ are given by, for all $t \geq 0$,
\begin{equation} \label{Paper01_quadratic_variation_stochastic_integral_poisson_point_process}
    [M](t) = \int_{[0,t]\times [0,\infty)} H^{2}(\tau-, y)\,  \mathcal{P}(d\tau \times dy) \quad \textrm{and} \quad \langle M \rangle(t) \defeq \int_{0}^{t} \int_{0}^{\infty} H^{2}(\tau-, y) \,  d\tau \, dy.
\end{equation}
Observe that, in contrast to~$[M]$, $\langle M \rangle$ is a predictable process. Nonetheless, both $M^{2} - [M]$ and $M^{2} - \langle M \rangle$ are càdlàg martingales adapted to the filtration $\{\mathcal{F}_{t}\}_{t \geq 0}$. In particular, Itô's isometry principle yields that for all $t \geq 0$,
\begin{equation} \label{Paper01_ito_isometry_stochastic_integral_jumps}
    \mathbb{E}\Bigg[\left(\int_{[0,t]\times [0,\infty)} H(\tau-,y) \, (\mathcal{P} - \lambda)(d\tau \times dy) \right)^{2}\Bigg] = \mathbb{E}\left[\int_{0}^{t} \int_{0}^{\infty} H^{2}(t-,y) \, dt \, dy\right].
\end{equation}
We finish this subsection by stating a corollary of Itô's isometry. Since we have not found this exact estimate elsewhere, we state and prove it here for the sake of completeness.

\begin{lemma} \label{Paper01_easy_inequality_square_integral_poisson_point_process}
    Under the above conditions on $H$, 
    \begin{equation*}
    \begin{aligned}
    \mathbb{E}\Bigg[\left(\int_{[0,t]\times [0,\infty)} H(\tau-,y) \, \mathcal{P} (d\tau \times dy) \right)^{2}\Bigg] & \leq 2 \, \mathbb{E}\left[\int_{0}^{t} \int_{0}^{\infty} H^{2}(t-,y) \, dt \, dy\right] \\ & \quad \; + 2 \, \mathbb{E}\Bigg[\left(\int_{0}^{t} \int_{0}^{\infty} H(\tau-,y) \, d\tau \, dy \right)^{2}\Bigg].
\end{aligned}
\end{equation*}
\end{lemma}

\begin{proof}
    Observing that for any $a,b \in \mathbb{R}$, we have $(a + b)^{2} \leq 2(a^{2} + b^{2})$, using the decomposition of the stochastic integral into martingale and finite variation parts in~\eqref{Paper01_decomposition_semimartingale_integral_jumps} and the Itô isometry given by~\eqref{Paper01_ito_isometry_stochastic_integral_jumps}, our claim holds.
\end{proof}

\subsection{Feller semigroups and non-locally compact Polish spaces} \label{Paper01_subsection_appendix_feller_non_locally_compact_polish_space}

In this subsection, we will briefly review some aspects of the theory of Feller processes taking values in non-locally compact Polish spaces. Let $(\mathcal{S}, d_{\mathcal{S}})$ be some complete and separable metric space, and let $\mathcal{B}(\mathcal{S})$ denote the $\sigma$-algebra of Borel subsets of $\mathcal{S}$. Let $\mathscr{B}_b(\mathcal{S}, \mathbb{R})$ be the set of bounded real-valued Borel measurable functions on $\mathbb{R}$, and observe that $\mathscr{B}_b(\mathcal{S}, \mathbb{R})$ is a Banach space when equipped with the norm $\|\cdot \|_{L_\infty(\mathcal S,\mathbb R)}$, where
\begin{equation*}
     \vert \vert \phi \vert \vert_{L_{\infty}(\mathcal{S}; \mathbb{R})} \defeq \sup_{\boldsymbol{\eta} \in \mathcal{S}} \vert \phi(\boldsymbol{\eta}) \vert.
\end{equation*}
Moreover, the space of bounded real-valued continuous functions $\mathscr{C}_{b}(\mathcal{S}, \mathbb{R}) \subset \mathscr{B}_b(\mathcal{S}, \mathbb{R})$ is a complete subspace of $\mathscr{B}_b(\mathcal{S}, \mathbb{R})$ when equipped with the norm~$\vert \vert \cdot \vert \vert_{L_{\infty}(\mathcal{S}; \mathbb{R})}$, i.e.~$\mathscr{C}_{b}(\mathcal{S}, \mathbb{R})$ is also a Banach space. Following~\cite[Section~2.2, Definition~2.4]{van2011markov}, we now state the definition of Feller semigroups in this context. We highlight that similar definitions are also found in the literature under the name of generalised Feller semigroups~\cite{cuchiero2020generalized, cuchiero2023ramifications} or bi-continuous semigroups~\cite{kuhnemund2003hille}.

\begin{definition}[Feller semigroups] \label{Paper01_definition_feller_non_locally_compact}
    A family~$\{P_{t}\}_{t \geq 0}$ of bounded linear operators defined on $\mathscr{B}_b(\mathcal{S}, \mathbb{R})$ is called a Feller semigroup on $\mathscr{C}_{b}(\mathcal{S}, \mathbb{R})$ if the following properties hold:
    \begin{enumerate}[(i)]
        \item $P_{t}\phi \in \mathscr{C}_{b}(\mathcal{S}, \mathbb{R})$ for all $t \in [0, \infty)$ and $\phi \in \mathscr{C}_{b}(\mathcal{S}, \mathbb{R})$;
        \item $P_{T}P_{t} = P_{T + t}$ for all $T, t \geq 0$, and $P_{0} = \textrm{Id}$, i.e.~$P_{0}$ is the identity map;
        \item For any $t \geq 0$, $P_{t}$ is a contraction, i.e.~for any $\phi \in \mathscr{C}_{b}(\mathcal{S}, \mathbb{R})$, $\vert \vert P_{t}\phi \vert \vert_{L_{\infty}(\mathcal{S}; \mathbb{R})} \leq \vert \vert \phi \vert \vert_{L_{\infty}(\mathcal{S}; \mathbb{R})}$;
        \item $\{P_{t}\}_{t \geq 0}$ is a family of positive operators, i.e.~for $\phi \in \mathscr{C}_{b}(\mathcal{S}, \mathbb{R})$ with $\phi \geq 0$, we have $P_{t}\phi \geq 0$, for any $t \in [0, \infty)$;
        \item For all $\phi \in \mathscr{C}_{b}(\mathcal{S}, \mathbb{R})$, $T \geq 0$ and $\boldsymbol{\eta} \in \mathcal{S}$, $\lim_{t \rightarrow T} (P_{t}\phi)(\boldsymbol{\eta}) = P_{T}\phi(\boldsymbol{\eta})$;
        \item For any $\phi \in \mathscr{C}_{b}(\mathcal{S}, \mathbb{R})$, the map $(t, \boldsymbol{\eta}) \mapsto (P_{t}\phi)(\boldsymbol{\eta})$ is continuous with respect to the product topology defined on $[0, \infty) \times \mathcal{S}$.
    \end{enumerate}
\end{definition}

Observe that property~(vi) is stronger than properties~(i) and (v), and implies uniform continuity on compact time intervals and compact subsets of $\mathcal{S}$. In the context of compact or locally compact metric spaces, it is often more convenient to study the action of Feller semigroups on the space of continuous functions vanishing at infinity, which is represented by $\mathscr{C}_{0}(\mathcal{S}, \mathbb{R})$. When $\mathcal{S}$ is not locally compact, the concept of bounded continuous functions vanishing at infinity is not well defined. Another crucial difference between Feller semigroups in locally compact and non-locally compact Polish spaces is that in the former case we usually require the semigroup to be strongly continuous, i.e.~that for any $\phi \in \mathscr{C}_{0}(\mathcal{S}, \mathbb{R})$ and $T\ge 0$ we have
\begin{equation} \label{eq:Ptconv}
    \lim_{t \rightarrow T} \vert \vert P_{t}\phi - P_{T}\phi \vert \vert_{L_{\infty}(\mathcal{S}; \mathbb{R})} = 0.
\end{equation}
Observe that in contrast to~\eqref{eq:Ptconv}, property~(vi) of Definition~\ref{Paper01_definition_feller_non_locally_compact} only requires uniform continuity on compact time intervals and compact subsets of $\mathcal{S}$. Due to these differences, the theory of Feller processes in non-locally compact Polish spaces, i.e.~stochastic processes with finite dimensional distributions given by a Feller semigroup, is very different and much less studied than Feller processes in locally compact metric spaces.

For interacting particle systems with \emph{a priori} bounded particle density, it is often possible to construct the stochastic process in a compact Polish space, so that the standard theory applies. Since this article focuses on an interacting particle system without an \emph{a priori} bound on the number of particles per site and with infinitely many types of particles, we are obliged to construct the process in a non-locally compact Polish space. In the context of locally compact metric spaces, a Feller semigroup corresponds to a càdlàg strong Markov process (see~\cite[Theorem~4.2.7]{ethier2009markov}). This does not hold in full generality for Feller processes in non-locally compact spaces~\cite{beznea2024strong}, the main challenge being that one must verify the tightness of the associated probability measures in order to construct a càdlàg version of the process. We observe, however, that a semigroup $\{P_{t}\}_{t \geq 0}$ corresponding to a family of probability measures and which satisfies Definition~\ref{Paper01_definition_feller_non_locally_compact} corresponds to the finite-dimensional distributions of a Markov process. Although this result is standard (see for instance~\cite[Chapter~5]{chen2004markov}), we include the proof for the sake of completeness.

\begin{lemma} \label{Paper01_markov_property_feller_semigroup}
    Let $\{P_{t}\}_{t \geq 0}$ be a Feller semigroup satisfying the conditions of Definition~\ref{Paper01_definition_feller_non_locally_compact}, and such that for any $\boldsymbol{\eta} \in \mathcal{S}$ and $T \geq 0$, there exists a Borel probability measure $\mathbb{P}(T, \boldsymbol{\eta}, \cdot)$ on $\mathcal{S}$ such that for any $\phi \in \mathscr{C}_{b}(\mathcal{S}, \mathbb{R})$, 
    \begin{equation*}
        (P_{T}\phi)(\boldsymbol{\eta}) = \int_{\mathcal{S}} \phi(\boldsymbol{\xi}) \, \mathbb{P}(T, \boldsymbol{\eta}, d\boldsymbol{\xi}).
    \end{equation*}
    Then the map $(T, \boldsymbol{\eta}, \Gamma) \in [0, \infty) \times \mathcal{S} \times \mathscr{B}(\mathcal{S}) \mapsto \mathbb{P}(T, \boldsymbol{\eta}, \Gamma)$ is a time-homogeneous transition function, i.e.~it satisfies the following conditions:
    \begin{enumerate}[(i)]
        \item $\mathbb{P}(t, \boldsymbol{\eta}, \cdot)$ is a Borel probability measure on $\mathcal{S}$, for any $(t, \boldsymbol{\eta}) \in [0, \infty) \times \mathcal{S}$;
        \item $\mathbb{P}(0, \boldsymbol{\eta}, \cdot) = \delta_{\boldsymbol{\eta}}$, i.e.~the Dirac measure at~$\boldsymbol{\eta}$, for any $\boldsymbol{\eta} \in \mathcal{S}$;
        \item $\mathbb{P}(\cdot, \cdot, \Gamma) \in \mathscr{B}_b([0, \infty) \times \mathcal{S},\mathbb R)$, for any $\Gamma \in \mathscr{B}(\mathcal{S})$;
        \item For any $T,t \geq 0$, $\boldsymbol{\eta} \in \mathcal{S}$ and $\Gamma \in \mathscr{B}(\mathcal{S})$,
        \begin{equation*}
            \mathbb{P}(T + t, \boldsymbol{\eta}, \Gamma) = \int_{\mathcal{S}}  \mathbb{P}(T, \boldsymbol{\xi}, \Gamma) \; \mathbb{P}(t, \boldsymbol{\eta}, d \boldsymbol{\xi}).
        \end{equation*}
    \end{enumerate}
    Moreover, suppose that $(\eta(t))_{t \geq 0}$ is a stochastic process taking values in $\mathcal{S}$ whose finite-dimensional distributions are given by the expression below, which holds for any $n \in \mathbb{N}$, $\Gamma_{0}, \ldots, \Gamma_{n} \in \mathscr{B}(\mathcal{S})$ and $0 \leq t_{1} \leq \ldots \leq t_{n}$:
    \begin{equation*}
    \begin{aligned}
        & \mathbb{P}\Big(\eta(0) \in \Gamma_{0}, \, \eta(t_{1}) \in \Gamma_{1}, \ldots, \eta(t_{n}) \in \Gamma_{n} \Big) \\ & \quad = \int_{\Gamma_{0}} \, \cdots \, \int_{\Gamma_{n-1}} \mathbb{P}\Big(t_{n} - t_{n-1}, \boldsymbol{\xi}_{n-1}, \Gamma_{n}\Big) \mathbb{P}\Big(t_{n-1} - t_{n-2}, \boldsymbol{\xi}_{n-2}, d\boldsymbol{\xi}_{n-1}\Big) \cdots  \mathbb{P}\Big(t_{1}, \boldsymbol{\xi}_{0}, d\boldsymbol{\xi}_{1}\Big)\nu(d\boldsymbol{\xi}_{0}),
    \end{aligned}
    \end{equation*}
    where $\nu$ is some Borel probability measure on $\mathcal{S}$. Then $(\eta(t))_{t \geq 0}$ is a Markov process.
\end{lemma}

\begin{proof}
    By our assumption in the statement, for any $T \geq 0$ and $\boldsymbol{\eta} \in \mathcal{S}$, $\mathbb{P}(T, \boldsymbol{\eta}, \cdot)$ is a Borel probability measure on $\mathcal{S}$. Moreover, by Definition~\ref{Paper01_definition_feller_non_locally_compact}(ii), $\mathbb{P}(0, \boldsymbol{\eta}, \cdot) = \delta_{\boldsymbol{\eta}}(\cdot)$, for any $\boldsymbol{\eta} \in \mathcal{S}$, since $P_{0}$ is equal to the identity map. Furthermore, the semigroup property given by Definition~\ref{Paper01_definition_feller_non_locally_compact}(ii) implies the Chapman-Kolmogorov property for the family of transition functions, i.e.~condition~(iv) above (see for instance~\cite[Lemma~19.1]{kallenberg1997foundations}). We now verify that the map $(T, \boldsymbol{\eta}) \mapsto \mathbb{P}(T, \boldsymbol{\eta}, \Gamma)$ is Borel measurable with respect to the product topology on $[0, \infty) \times \mathcal{S}$, for any $\Gamma \in \mathscr{B}(\mathcal{S})$. By the functional monotone class theorem, it is enough to verify that the measurability holds for any closed or open subset $\Gamma \subseteq \mathcal{S}$. Observe that by Definition~\ref{Paper01_definition_feller_non_locally_compact}(vi), for any $\phi \in \mathscr{C}_{b}(\mathcal{S}, \mathbb{R})$, the map 
    \begin{equation*}
        (t, \boldsymbol{\eta}) \mapsto (P_{T}\phi)(\boldsymbol{\eta}) = \int_{\mathcal{S}} \phi(\boldsymbol{\xi}) \, \mathbb{P}(T, \boldsymbol{\eta}, d\boldsymbol{\xi})
    \end{equation*}
    is Borel measurable with respect to the product topology on $[0, \infty) \times \mathcal{S}$. Observe that since $\mathcal{S}$ is a Polish space, for any open subset $\Gamma \subseteq \mathcal{S}$, the indicator function $\boldsymbol{\eta} \mapsto \mathds{1}_{\{\boldsymbol{\eta} \in \Gamma\}}$ is the pointwise limit of a sequence of bounded continuous functions. Since the limit of a sequence of measurable functions is measurable, for any open or closed subset $\Gamma \subseteq \mathcal{S}$, we conclude that the map $(T, \boldsymbol{\eta}) \mapsto \mathbb{P}(t, \boldsymbol{\eta}, \Gamma)$ is Borel measurable with respect to the product topology on $[0, \infty) \times \mathcal{S}$, and therefore our claim holds. We refer the reader to~\cite[Theorem~4.1.1]{ethier2009markov} for the proof that a stochastic process whose finite-dimensional distributions correspond to this family of transition functions is a Markov process.
\end{proof}

We call an $\mathcal{S}$-valued stochastic process $(\eta(t))_{t\ge 0}$ whose finite-dimensional distributions correspond to a Feller semigroup a \textit{Feller process}. Lemma~\ref{Paper01_markov_property_feller_semigroup} shows that any such $(\eta(t))_{t\geq 0}$ is a Markov process. We emphasise, however, that this does not guarantee that $(\eta(t))_{t\ge 0}$ is càdlàg; see, for example, the counterexample by Beznea, C\^{\i}mpean, and R{\"o}ckner in~\cite{beznea2024strong}. We use Lemma~\ref{Paper01_markov_property_feller_semigroup} in Section~\ref{Paper01_section_functional_analysis_feller_non_locally_compact} for the construction of the spatial Muller’s ratchet.

\subsection{Topological Properties of~$\left(\mathcal{S}, d_{\mathcal{S}}\right)$} \label{Paper01_subsection_topological_properties_state_space}

In this subsection, we will characterise the space of configurations~$(\mathcal{S},d_{\mathcal{S}})$. Recall the definition of $\mathcal{S}$ and $\vert \vert \vert \cdot \vert \vert \vert_{\mathcal{S}}$ in~\eqref{Paper01_definition_state_space_formal}. We will now collect some properties of Cauchy sequences in $(\mathcal{S}, d_{\mathcal{S}})$.

\begin{lemma} \label{Paper01_cauchy_sequence_semi_metric_spaces}
    A sequence of configurations $(\boldsymbol{\eta}^{n})_{n \in \mathbb{N}}=((\eta^n_k(x))_{k\in \mathbb N_0,\, x\in L^{-1}\mathbb Z})_{n\in \mathbb N}$ is a Cauchy sequence in $(\mathcal{S}, d_{\mathcal{S}})$ if and only if the following conditions hold:
    \begin{enumerate}[(i)]
        \item  For each $x \in L^{-1}\mathbb{Z}$, $(\eta^{n}(x))_{n \in \mathbb{N}}$ is a Cauchy sequence in $\ell_{1}$;
        \item $\sup_{n\in \mathbb N} \vert \vert \vert \boldsymbol{\eta}^{n} \vert \vert \vert_{\mathcal{S}} < \infty$;
         \item For any $\varepsilon > 0$, there exists $R_{\varepsilon} > 0$ such that
        \begin{equation*}
            \sup_{n\in \mathbb N} \; \sum_{\substack{\{x \in L^{-1}\mathbb{Z}: \, \vert x \vert \geq R_{\varepsilon}\}}} \frac{\vert \vert \eta^{n}(x) \vert \vert_{\ell_{1}}}{(1 + \vert x \vert)^{2}} \leq \varepsilon;
        \end{equation*}
         \item For each $x \in L^{-1}\mathbb{Z}$, there exists $J^{(x)} \in \mathbb{N}_{0}$ such that for any $k\in \mathbb N_0$ with $k \geq J^{(x)}$, $\eta^{n}_{k}(x) = 0$ for all~$n \in \mathbb{N}$.
    \end{enumerate}
\end{lemma}

\begin{proof}
    Since the conditions are very similar to the characterisation of compact subsets of sequence spaces, we only sketch the proof and omit the details. Let ${\nu}$ be a weighted counting measure on $\mathbb N_0 \times L^{-1}\mathbb Z$ given by, for any $\mathcal{A}_{1} \times \mathcal{A}_{2} \subseteq \mathbb N_0 \times L^{-1}\mathbb Z$,
    \begin{equation*}
        {\nu}(\mathcal{A}_{1} \times \mathcal{A}_{2}) \defeq \sum_{k \in \mathbb N_0} \; \sum_{x \in L^{-1}\mathbb Z} \; \frac{\mathds{1}_{\{k \in \mathcal{A}_1, x \in \mathcal{A}_2\}}}{(1 + \vert x \vert)^2}. 
    \end{equation*}
    Then, $\mathcal{S}$ can be embedded into the space $\ell_1(\mathbb N_0 \times L^{-1}\mathbb Z; \nu)$, which is given by
    \begin{equation*}
        \ell_1(\mathbb{N}_0 \times \mathbb{Z}; {\nu}) \defeq \left\{\boldsymbol{\eta} = (\eta_k(x))_{x \in L^{-1}\mathbb{Z}, \, k \in \mathbb{N}_0} \in \mathbb{R}^{\mathbb{N}_0 \times L^{-1}\mathbb{Z}}: \sum_{x \in L^{-1}\mathbb{Z}} \; \sum_{k \in \mathbb{N}_0} \frac{\vert \eta_k(x) \vert}{(1 + \vert x \vert)^2} < \infty\right\}.
    \end{equation*}
    Then, by Vitali's convergence theorem (see for example~\cite[Exercise~6.15]{folland1999real}), and since for any $x \in L^{-1}\mathbb{Z}$ and $\boldsymbol{\eta} \in \mathcal{S}$, $\eta(x) = (\eta_{k}(x))_{k \in \mathbb{N}_{0}}$ is an integer-valued sequence, the result follows.
\end{proof}

We are now ready to state the topological properties of $(\mathcal{S}, d_{\mathcal{S}})$. 

\begin{proposition} \label{Paper01_topological_properties_state_space}
    The metric space $(\mathcal{S}, d_{\mathcal{S}})$ is complete and separable. Moreover, a subset $\mathcal{K} \subset \mathcal{S}$ is relatively compact in the topology induced by $d_{\mathcal{S}}$ if and only if all the following conditions are satisfied:
    \begin{enumerate}[(i)]
        \item $\mathcal{K}$ is $\vert \vert \vert \cdot \vert \vert \vert_{\mathcal{S}}$-bounded, i.e.~$\sup_{\boldsymbol{\eta} \in \mathcal{K}} \vert \vert \vert \boldsymbol{\eta} \vert \vert \vert_{\mathcal{S}} < \infty$;
        \item For any $\varepsilon > 0$, there exists $R_{\varepsilon} > 0$ such that
        \begin{equation} \label{eq:*etacond}
            \sup_{\boldsymbol{\eta}=(\eta(y))_{y\in L^{-1}\mathbb Z} \in \mathcal{K}} \; \sum_{\substack{\{x \in L^{-1}\mathbb{Z}: \, \vert x \vert \geq R_{\varepsilon}\}}} \frac{\vert \vert \eta(x) \vert \vert_{\ell_{1}}}{(1 + \vert x \vert)^{2}} \leq \varepsilon;
        \end{equation}
        \item For each $x \in L^{-1}\mathbb{Z}$, there exists $J^{(x)} \in \mathbb{N}_{0}$ such that for any $k\in \mathbb N_0$ with $k \geq J^{(x)}$, $\eta_{k}(x) = 0$ for all~$\boldsymbol{\eta}=(\eta_j(y))_{j\in \mathbb N_0,\, y\in L^{-1}\mathbb Z} \in \mathcal{K}$.
    \end{enumerate}
    Moreover, for any increasing sequence $(J_n)_{n\in \mathbb N}\subseteq \mathbb N$ with $\lim_{n\to \infty}J_n=\infty$, for any
    compact set $\mathcal{K} \subset \mathcal{S}$, the following limit holds:
    \begin{equation} \label{Paper01_easy_limit_bound_high_mutation_load_space}
        \lim_{n \rightarrow \infty} \sup_{\boldsymbol{\eta}=(\eta_k(y))_{k\in \mathbb N_0, \, y\in L^{-1}\mathbb Z} \in \mathcal{K}} \, \sum_{x \in L^{-1}\mathbb{Z}} \; \sum_{j = J_{n}}^{\infty} \; \frac{\eta_{j}(x)}{(1 + \vert x \vert)^{2}} = 0.
    \end{equation}
\end{proposition}

\begin{proof}
    We first prove separability. Consider the countable subset of configurations
    \begin{equation} \label{Paper01_dense_subset_state_space}
    \begin{aligned}
        \mathcal{Y} \defeq \bigcup_{i = 1}^{\infty} \; \bigcup_{j = 1}^{\infty} \left\{\boldsymbol{\xi}=(\xi_j(y))_{j\in \mathbb N_0,\, y\in L^{-1}\mathbb Z} \in \mathcal{S}: \; \xi_{k}(x) = 0 \textrm{ for } k > i \textrm{ or } \vert x \vert > jL^{-1} \right\}.
    \end{aligned}
    \end{equation}
    We claim that $\mathcal{Y}$ is dense in $(\mathcal{S}, d_{\mathcal{S}})$. Indeed, for any fixed $\boldsymbol{\eta}= (\eta_k(x))_{k\in \mathbb N_0,\, x\in L^{-1}\mathbb Z}\in \mathcal{S}$ and $\varepsilon > 0$, there exist $i,j \in \mathbb{N}$ such that 
    \begin{equation*}
        \sum_{\substack{\{x \in L^{-1}\mathbb{Z}: \, \vert x \vert > jL^{-1}\}}} \frac{\vert \vert \eta(x) \vert \vert_{\ell_{1}}}{(1 + \vert x \vert)^{2}} \leq \varepsilon,
    \end{equation*}
    and such that for $x \in L^{-1}\mathbb{Z} \, \cap \, [-jL^{-1}, jL^{-1}]$, $\eta_{k}(x) = 0$ for $k > i$. Hence, the configuration $\boldsymbol{\xi}=(\xi(x))_{x\in L^{-1}\mathbb Z} \in \mathcal{S}$ given by $\xi(x) = \eta(x)$ for $x \in L^{-1}\mathbb{Z} \cap [-jL^{-1}, jL^{-1}]$ and $\xi(x) = 0$ for $\vert x \vert > jL^{-1}$ is an element of $\mathcal{Y}$ and satisfies $\vert \vert \vert \boldsymbol{\eta} - \boldsymbol{\xi} \vert \vert \vert_{\mathcal{S}} \leq \varepsilon$. This proves separability.
    
    To prove completeness, let $(\boldsymbol{\eta}^{n})_{n \in \mathbb{N}}=((\eta^{n}_k(x)_{k\in \mathbb N_0,\,x\in L^{-1}\mathbb Z})_{n \in \mathbb{N}}$ be a Cauchy sequence in $(\mathcal{S}, d_{\mathcal{S}})$. As stated in Lemma~\ref{Paper01_cauchy_sequence_semi_metric_spaces}, for each $x \in L^{-1}\mathbb{Z}$, $(\eta^{n}(x))_{n \in \mathbb{N}}$ is a Cauchy sequence in~$\ell_{1}$. Since $\ell_{1}$ is a Banach space, for each $x \in L^{-1}\mathbb{Z}$, there exists $\eta(x)=(\eta_k(x))_{k\in \mathbb N_0} \in \ell_{1}$ such that $\eta^{n}(x) \rightarrow \eta(x)$ in~$\ell_{1}$ as~$n\rightarrow \infty$. This implies that for each $x \in L^{-1}\mathbb{Z}$ and $k \in \mathbb{N}_{0}$, we have $\eta^{n}_{k}(x) \rightarrow \eta_{k}(x)$ as $n \rightarrow \infty$ in $\mathbb{R}$. Since $(\eta_{k}^{n}(x))_{n \in \mathbb{N}}$ is an integer-valued Cauchy sequence, it must be eventually constant, and therefore $\eta(x) \in \mathbb{N}_{0}^{\mathbb{N}_{0}} \cap \ell_{1}$, for every $x \in L^{-1}\mathbb{Z}$. Define $\boldsymbol{\eta} \defeq (\eta(x))_{x \in L^{-1}\mathbb{Z}} \in (\mathbb{N}_{0}^{\mathbb{N}_{0}})^{L^{-1}\mathbb{Z}}$. We now claim that $\boldsymbol{\eta} \in \mathcal{S}$. It will suffice to establish that $\vert \vert \vert \boldsymbol{\eta} \vert \vert \vert_{\mathcal{S}} < \infty$. To see this, note that by~\eqref{Paper01_definition_state_space_formal},
    \begin{equation*}
        \vert \vert \vert \boldsymbol{\eta} \vert \vert \vert_{\mathcal{S}} = \sum_{x \in L^{-1}\mathbb{Z}} \frac{\vert \vert \eta(x) \vert \vert_{\ell_{1}}}{(1 + \vert x \vert)^{2}} 
        = \sum_{x \in L^{-1}\mathbb{Z}} \lim_{n \rightarrow \infty} \frac{\vert \vert \eta^{n}(x) \vert \vert_{\ell_{1}}}{(1 + \vert x \vert)^{2}} 
        \leq \liminf_{n \rightarrow \infty} \sum_{x \in L^{-1}\mathbb{Z}}\frac{\vert \vert \eta^{n}(x) \vert \vert_{\ell_{1}}}{(1 + \vert x \vert)^{2}}  < \infty,
    \end{equation*}
    where the first inequality follows from Fatou's lemma, and the last inequality follows from property~(ii) of Lemma~\ref{Paper01_cauchy_sequence_semi_metric_spaces}. Hence, $\boldsymbol{\eta} \in \mathcal{S}$.
    
    It then remains to verify that $\lim_{n \rightarrow \infty} \vert \vert \vert \boldsymbol{\eta}^{n} - \boldsymbol{\eta} \vert \vert \vert_{\mathcal{S}} = 0$. By property~(iii) of Lemma~\ref{Paper01_cauchy_sequence_semi_metric_spaces}, the series in the definition of $\vert \vert \vert \boldsymbol{\eta}^{n} - \boldsymbol{\eta} \vert \vert \vert_{\mathcal{S}}$ converges uniformly in $n \in \mathbb{N}$. Therefore, we get
    \begin{equation*}
    \begin{aligned}
        \lim_{n \rightarrow \infty} \; \vert \vert \vert \boldsymbol{\eta}^{n} - \boldsymbol{\eta} \vert \vert \vert_{\mathcal{S}} =  \lim_{n \rightarrow \infty} \; \sum_{\substack{x \in L^{-1}\mathbb{Z}}} \frac{\vert \vert \eta^{n}(x) - \eta(x) \vert \vert_{\ell_{1}}}{(1 + \vert x \vert)^{2}} = \sum_{\substack{x \in L^{-1}\mathbb{Z}}} \; \lim_{n \rightarrow \infty} \;  \frac{\vert \vert \eta^{n}(x) - \eta(x) \vert \vert_{\ell_{1}}}{(1 + \vert x \vert)^{2}} = 0,
    \end{aligned}
    \end{equation*}
    where for the last identity we used the fact that, by construction, for each $x\in L^{-1}\mathbb Z$, $\eta^{n}(x) \rightarrow \eta(x)$ in~$\ell_{1}$ as $n \rightarrow \infty$. Hence, $(\mathcal{S}, d_{\mathcal{S}})$ is complete, as claimed.
    
    We will now establish necessary and sufficient conditions for relative compactness. Lemma~\ref{Paper01_cauchy_sequence_semi_metric_spaces} implies that conditions~(i)-(iii) are necessary for relative compactness. The proof that they are sufficient follows along the same lines as the characterisation of relatively compact subsets of spaces of sequences~$\ell_{p}$, for $p \geq 1$ (see~\cite[Exercise~I.6]{diestel2012sequences}). Since this proof is standard, we omit the details.
    
    It remains to establish~\eqref{Paper01_easy_limit_bound_high_mutation_load_space}. For $\mathcal{K} \subset \mathcal{S}$ compact and $\varepsilon > 0$, there exists $R_{\varepsilon} > 0$ such that~\eqref{eq:*etacond} holds. For $(J_n)_{n\in \mathbb N}$ with $J_n\to \infty$ as $n \rightarrow \infty$, condition~(iii) of this proposition implies that there exists $n_{\varepsilon} \in \mathbb{N}$ such that for any $n \geq n_{\varepsilon}$ and any $\boldsymbol{\eta} =(\eta_k(x))_{k\in \mathbb N_0, \, x\in L^{-1}\mathbb Z} \in \mathcal{K}$, $\eta_{j}(x) = 0$ whenever both $j \geq J_n$ and $x \in L^{-1}\mathbb{Z} \cap [-R_\varepsilon, R_\varepsilon]$. Therefore, for $n \geq n_\varepsilon$,
    \begin{equation*}
        \sup_{\boldsymbol{\eta}=(\eta_k(y))_{k\in \mathbb N_0, \, y\in L^{-1}\mathbb Z} \in \mathcal{K}} \, \sum_{x \in L^{-1}\mathbb{Z}} \; \sum_{j = J_{n}}^{\infty} \; \frac{\eta_{j}(x)}{(1 + \vert x \vert)^{2}} \leq  \sup_{\boldsymbol{\eta}=(\eta(y))_{y\in L^{-1}\mathbb Z} \in \mathcal{K}} \; \sum_{\substack{\{x \in L^{-1}\mathbb{Z}: \, \vert x \vert \geq R_{\varepsilon}\}}} \frac{\vert \vert \eta(x) \vert \vert_{\ell_{1}}}{(1 + \vert x \vert)^{2}} \leq \varepsilon.
    \end{equation*}
    Since $\varepsilon > 0$ was arbitrary, the limit in~\eqref{Paper01_easy_limit_bound_high_mutation_load_space} holds, which completes the proof.
\end{proof}

Proposition~\ref{Paper01_topological_properties_state_space} allows us to apply the theory of Markov process in Polish spaces to stochastic processes in $(\mathcal{S}, d_{\mathcal{S}})$. Next, we state a corollary of Proposition~\ref{Paper01_topological_properties_state_space} regarding the denseness of the set of possible initial configurations $\mathcal{S}_{0}$ given by~\eqref{Paper01_set_initial_configurations} in the space of configurations $\mathcal S$.

\begin{corollary} \label{Paper01_initial_configurations_dense_subset_}
    The set $\mathcal{S}_{0}$ is a dense subset of $\mathcal{S}$ with respect to the metric $d_{\mathcal{S}}$.
\end{corollary}

\begin{proof}
    It is enough to see that the subset of configurations $\mathcal{Y}$ given by~\eqref{Paper01_dense_subset_state_space} is a subset of $\mathcal{S}_{0}$.
\end{proof}

To finish this subsection, we will briefly study the set of cylindrical bounded real-valued functions $\mathscr{C}_{b,*}^{\textrm{cyl}}(\mathcal{S}, \mathbb{R})$ defined in~\eqref{Paper01_general_definition_cylindrical_function}. Our next result shows that $\mathscr{C}_{b,*}^{\textrm{cyl}}(\mathcal{S}, \mathbb{R})$ is a suitable set of real-valued test functions on $\mathcal{S}$, since it satisfies the conditions of Assumption~\ref{Paper01_assumption_candidate_set_functions_domain_generator} (see Section~\ref{Paper01_section_functional_analysis_feller_non_locally_compact}). Recall the definition of $\mathscr{C}_*(\mathcal{S}, \mathbb{R})$ in~\eqref{Paper01_definition_Lipschitz_function}.

\begin{lemma} \label{Paper01_characterisation_lipschitz_functions_semi_metric}
     The set $\mathscr{C}_{b,*}^{\textrm{cyl}}(\mathcal{S}, \mathbb{R})$  satisfies the following conditions:
    \begin{enumerate}[(i)]
        \item The constant function $\boldsymbol{1}: \mathcal{S} \rightarrow \{1\}$ is an element of $\mathscr{C}_{b,*}^{\textrm{cyl}}(\mathcal{S}, \mathbb{R})$;
        \item $\mathscr{C}^{\textrm{cyl}}_{b,*}(\mathcal{S}, \mathbb{R})$ is an algebra over $\mathbb{R}$;
        \item $\mathscr{C}_{b,*}^{\textrm{cyl}}(\mathcal{S}, \mathbb{R})$ separates points of $\mathcal{S}$;
        \item For any $\phi \in \mathscr{C}_{b,*}^{\textrm{cyl}}(\mathcal{S}, \mathbb{R})$ and $r \geq 0$, the map $\phi_{r}: \mathcal{S} \rightarrow \mathbb{R}$ given by
        \begin{equation*}
            \phi_{r}(\boldsymbol{\eta}) \defeq \sgn( \phi
            (\boldsymbol{\eta})) \Big(\vert \phi
            (\boldsymbol{\eta})\vert \wedge r\Big)
        \end{equation*}
        is also an element of $\mathscr{C}_{b,*}^{\textrm{cyl}}(\mathcal{S}, \mathbb{R})$;
        \item $\mathscr{C}_{b,*}^{\textrm{cyl}}(\mathcal{S}, \mathbb{R}) \subseteq \mathscr{C}_*(\mathcal{S}, \mathbb{R})$.
    \end{enumerate}
\end{lemma}

\begin{proof}
    Properties~(i) and (ii) are immediate from~\eqref{Paper01_general_definition_cylindrical_function}. To see that~(iii) holds, for any $x \in L^{-1}\mathbb{Z}$, $k \in \mathbb{N}_{0}$ and $r > 0$, define the coordinate map
    \begin{equation*}
    \begin{aligned}
        \phi^{(x)}_{k,r}: \mathcal{S} & \rightarrow \mathbb{R} \\ \eta & \mapsto \phi^{(x)}_{k,r}(\eta) = \eta_{k}(x) \wedge r.
    \end{aligned}
    \end{equation*}
    Then $\phi^{(x)}_{k,r} \in \mathscr{C}^{\textrm{cyl}}_{b,*}(\mathcal{S}, \mathbb{R})$. For any $\boldsymbol{\eta}, \boldsymbol{\xi} \in \mathcal{S}$ with $\boldsymbol{\eta} \neq \boldsymbol{\xi}$, there exist $x \in L^{-1}\mathbb{Z}$, $k \in \mathbb{N}_0$ and $r > 0$ such that $\eta_{k}(x) \wedge \xi_{k}(x) < r$ and $\eta_{k}(x) \neq \xi_{k}(x)$. Hence, $\mathscr{C}^{\textrm{cyl}}_{b,*}(\mathcal{S}, \mathbb{R})$ separates points of $\mathcal{S}$. Property~(iv) is obtained by observing that for $r\ge 0$, the function $y \mapsto \sgn(y)(\vert y \vert \wedge r)$ is continuous, and that the composition of continuous functions is also continuous. To establish property~(v), recalling the definition of $\boldsymbol{e}^{(x)}_{k}$ from before~\eqref{Paper01_scaled_polynomials_carrying_capacity}, for $\phi \in \mathscr{C}^{\textrm{cyl}}_{b,*}(\mathcal{S}, \mathbb{R})$, it follows from~\eqref{Paper01_general_definition_cylindrical_function} that there exists $R = R(\phi) > 0$ such that for
    any $x \in L^{-1}\mathbb{Z}$, $k \in \mathbb{N}_{0}$ and $a \in \{1,2\}$,
     \begin{equation} \label{Paper01_bound_other_version_migration_continuous_bounded_function_cylindrical_updated_version}
     \begin{aligned}
        & \sup_{\{\boldsymbol{\zeta}=(\zeta_j(y))_{j\in \mathbb N_0,\, y\in L^{-1}\mathbb Z} \in \mathcal{S}: \zeta_{k}(x) > 0\}} \Big\vert \phi\Big(\boldsymbol{\zeta} + \boldsymbol{e}^{(x+ (-1)^{a}L^{-1})}_{k} - \boldsymbol{e}^{(x)}_{k}\Big) - \phi(\boldsymbol{\zeta}) \Big\vert  \\ 
        & \quad + \sup_{\substack{\boldsymbol{\zeta}=(\zeta_j(y))_{j\in \mathbb N_0,\, y\in L^{-1}\mathbb Z} \in \mathcal{S}}} \Big(\Big\vert \phi\Big(\boldsymbol{\zeta} + \boldsymbol{e}^{(x)}_{k}\Big) - \phi(\boldsymbol{\zeta}) \Big\vert + \mathds{1}_{\{\zeta_k(x) > 0\}}\Big\vert \phi\Big(\boldsymbol{\zeta} - \boldsymbol{e}^{(x)}_{k}\Big) - \phi(\boldsymbol{\zeta}) \Big\vert\Big) \\ 
        & \quad \leq 6 \vert \vert \phi \vert \vert_{L_{\infty}(\mathcal{S}; \mathbb{R})} \mathds{1}_{\{\vert x \vert \wedge  \vert x + (-1)^{a}L^{-1} \vert \leq R\}} \\ 
        & \quad \leq \frac{6\vert \vert \phi \vert \vert_{L_{\infty}(\mathcal{S}; \mathbb{R})}(1 + R + L^{-1})^{2(1 + \deg q_-)}}{(1 + \vert x \vert)^{2(1 + \deg q_-)}} \mathds{1}_{\{\vert x \vert \wedge  \vert x + (-1)^{a}L^{-1} \vert \leq R\}},
    \end{aligned}
    \end{equation}
    where in the last inequality we used the fact that for $x \in L^{-1}\mathbb{Z}$ such that $$\vert x \vert \wedge \vert x + L^{-1} \vert \wedge \vert x - L^{-1} \vert \leq R,$$
    we have $\vert x \vert \leq R + L^{-1}$. Then, property~(v) follows from~\eqref{Paper01_bound_other_version_migration_continuous_bounded_function_cylindrical_updated_version} and from the definition of $\mathscr{C}_*(\mathcal{S}, \mathbb{R})$ in~\eqref{Paper01_definition_Lipschitz_function}. This completes the proof.
\end{proof}

\subsection{Stochastic chain rule formula} \label{Paper01_section_integration_by_parts_formula}

In this subsection, we state a general stochastic chain rule formula. In what follows, let~$(\mathcal{S}, d_{\mathcal{S}})$ be some complete and separable metric space, let $A \subset \mathscr{C}(\mathcal{S}, \mathbb{R})$ be some vector subspace of real-valued continuous functions defined on~$\mathcal{S}$, let $\mathscr{L}\vert_{A}: A \rightarrow \mathscr{C}(\mathcal{S}, \mathbb{R})$ be a (possibly unbounded) operator and let $(\eta(t))_{t \geq 0}$ be an~$\mathcal{S}$-valued càdlàg process with initial distribution~$\nu$, which is strongly Markovian with respect to its right-continuous natural filtration~$\Big\{\mathcal{F}^{\eta}_{t+}\Big\}_{t \geq 0}$. We will assume that~$\mathscr{L}$ is the infinitesimal generator of~$(\eta(t))_{t \geq 0}$, and that~$(\eta(t))_{t \geq 0}$ is a solution of the martingale problem for~$(\mathscr{L}, A, \nu)$, i.e.~for  any~$\phi \in A$, the process~$\Big(M^{\phi}(t)\Big)_{t \geq 0}$ given by, for all $T \geq 0$,
    \begin{equation} \label{eq:Mphidefn}
        M^{\phi}(T) \defeq \phi(\eta(T)) - \phi(\eta(0)) - \int_{0}^{T} (\mathscr{L}\phi)(\eta(t-)) \, dt
    \end{equation}
    is a càdlàg martingale with respect to the filtration~$\Big\{\mathcal{F}^{\eta}_{t+}\Big\}_{t \geq 0}$.

Our next result extends a classical result concerning the action of the generator of Markov processes with transition rates that are uniformly bounded on uniformly bounded continuous real-valued functions. The classical formulation appears, for instance, in~\cite[Lemma~4.3.4]{ethier2009markov} and~\cite[Lemma~A1.5.1]{kipnis1998scaling}. A slightly different formulation of Lemma~\ref{Paper01_general_integration_by_parts_formula} is stated without proof in~\cite[Theorem~2.6.3]{demasi1991mathematical}.

\begin{lemma} \label{Paper01_general_integration_by_parts_formula}
Let~$T \geq 0$ be fixed, and let~$g \in \mathscr{C}([0,T] \times \mathcal{S}, \mathbb{R})$. Suppose the following conditions hold:
    \begin{enumerate}[(i)]
        \item $
            \sup_{t_{1},t_{2} \in [0,T]}  \mathbb{E}_{\nu}\left[\vert g(t_{1}, \eta(t_{2})) \vert \right] < \infty$;
        \item The map~$[0,T] \times \mathcal{S}\ni (t,\boldsymbol{\eta})  \mapsto \left(\frac{\partial}{\partial t} g(\cdot, \boldsymbol{\eta})\right)(t)$ is in $ \mathscr{C}([0,T] \times \mathcal{S}, \mathbb{R})$;
        \item For any $t \in [0,T]$, the map~$\mathcal{S} \ni \boldsymbol{\eta}  \mapsto g(t, \boldsymbol{\eta}) $ is in $ A$, and the map~$[0,T] \times \mathcal{S} \ni (t,\boldsymbol{\eta})  \mapsto (\mathscr{L}g)(t,\cdot)(\boldsymbol{\eta}) $ is in $ \mathscr{C}([0,T] \times \mathcal{S}, \mathbb{R})$;
        \item 
        $
    \!
    \begin{aligned}[t]
\sup_{t_{1},t_{2} \in [0,T]} \; \mathbb{E}_{\nu}\Bigg[\left(\frac{\partial}{\partial t} g(\cdot, \eta(t_{2}))\right)^{2}(t_{1})  \Bigg] + \sup_{t_{1},t_{2} \in [0,T]} \; \mathbb{E}_{\nu}\Bigg[ \Big(\mathscr{L}g(t_{1}, \cdot )\Big)^{2} (\eta(t_{2})) \Bigg] < \infty.
        \end{aligned}
    $
    \end{enumerate}
    Then the process~$\Big(M^{g}(t \wedge T)\Big)_{t \geq 0}$ given by, for all~$t \geq 0$,
    \begin{equation} \label{eq:daggerMgdefn}
    \begin{aligned}
        & M^{g}(t \wedge T) \\ & \quad \defeq g\Big( t \wedge T, \eta(t \wedge T)\Big) - g(0, \eta(0)) - \int_{0}^{t \wedge T} \left(\left(\frac{\partial}{\partial t'} g(\cdot, \eta(t'-))\right)(t') + \Big(\mathscr{L}g(t', \cdot )\Big) (\eta(t'-))\right) \, dt'
    \end{aligned}
    \end{equation}
    is a càdlàg martingale with respect to the filtration~$\Big\{\mathcal{F}^{\eta}_{t+}\Big\}_{t \geq 0}$.
\end{lemma}

The main difference between Lemma~\ref{Paper01_general_integration_by_parts_formula} and~\cite[Lemma~4.3.4]{ethier2009markov} is that the conditions of~\cite[Lemma~4.3.4]{ethier2009markov} require the map~$g$, as well as its partial derivative with respect to time, and the action of~$\mathscr{L}$ on~$g$, to be uniformly bounded. We replace these requirements by conditions~(i) and~(iv). Since we are not aware of a proof of this result in the existing literature, we prove it below. 

\begin{proof}
    The fact that~$\Big(M^{g}(t \wedge T)\Big)_{t \geq 0}$ is càdlàg follows directly from the continuity of~$g$, and the fact that the sample paths of~$(\eta(t))_{t \geq 0}$ are càdlàg. Moreover, conditions~(i) and~(iv) imply that~$M^{g}$ is an integrable process. Hence, it remains to prove that for any~$t_{1}, t_{2} \geq 0$ with $t_{2} \geq t_{1}$, we have
    \begin{equation} \label{Paper01_intergration_parts_martingale_property}
        \mathbb{E}\Big[M^{g}(t_{2} \wedge T) \Big\vert \mathcal{F}^{\eta}_{t_{1}+} \Big] = M^{g}(t_{1} \wedge T).
    \end{equation}
    Without loss of generality, we assume~$t_{2} > t_{1}$. For each $n \in \mathbb{N}$, let~$(h_{i}^{(n)})_{i = 0}^{n}$ be a partition of~$[t_{1},t_{2}]$, i.e.~suppose $t_{1} = h^{(n)}_{0} < h^{(n)}_{1} < \cdots < h^{(n)}_{n} = t_{2}$, such that
    \begin{equation*}
        \lim_{n \rightarrow \infty} \; \max_{i \in [n-1]_{0}} \; \left( h^{(n)}_{i + 1} - h^{(n)}_{i} \right)= 0.
    \end{equation*}
    By conditions~(ii) and~(iii) of this lemma, and by~\eqref{eq:Mphidefn} and the tower property of conditional expectation, we have, for every~$n \in \mathbb{N}$,
    \begin{equation} \label{Paper01_intermediate_step_integration_parts_formula}
    \begin{aligned}
         & \mathbb{E}_{\nu}\Big[g\Big(t_{2} \wedge T, \eta(t_{2} \wedge T)\Big) - g\Big(t_{1} \wedge T, \eta(t_{1} \wedge T)\Big) \Big\vert \mathcal{F}^{\eta}_{t_{1}+} \Big] \\ 
         & \quad = \mathbb{E}_{\nu}\Big[ \sum_{i=0}^{n-1} \big( g(h^{(n)}_{i+1} \wedge T,\eta(h^{(n)}_{i} \wedge T))-g(h^{(n)}_{i} \wedge T,\eta(h^{(n)}_{i} \wedge T))
         \Big\vert \mathcal{F}^{\eta}_{t_1+} \Big]\\
         &\qquad +\mathbb{E}_{\nu}\Big[ \sum_{i=0}^{n-1} \big( g(h^{(n)}_{i+1} \wedge T,\eta(h^{(n)}_{i+1} \wedge T))-g(h^{(n)}_{i+1} \wedge T,\eta(h^{(n)}_{i} \wedge T))
         \Big\vert \mathcal{F}^{\eta}_{t_1+} \Big]\\
         & \quad = \mathbb{E}_{\nu}\Bigg[\sum_{i = 0}^{n-1} \int_{h^{(n)}_{i} \wedge T}^{h^{(n)}_{i+1} \wedge T} \left(\frac{\partial}{\partial t'} g\Big(\cdot, \eta(h^{(n)}_{i} \wedge T)\Big)\right)(t') \; dt' \Bigg\vert \mathcal{F}^{\eta}_{t_{1}+}\Bigg] \\ & \quad \quad \quad +  \mathbb{E}_{\nu}\Bigg[\sum_{i = 0}^{n-1} \int_{h^{(n)}_{i} \wedge T}^{h^{(n)}_{i + 1} \wedge T} \left(\mathscr{L}g\Big(h^{(n)}_{i+1} \wedge T, \cdot\Big)\right)(\eta(t'-)) \; dt' \Bigg\vert \mathcal{F}^{\eta}_{t_{1}+}\Bigg].
    \end{aligned}
    \end{equation}
    It will suffice then to study the convergence of the terms on the right-hand side of~\eqref{Paper01_intermediate_step_integration_parts_formula} as $n \rightarrow \infty$. Observe that since $(\eta(t))_{t \geq 0}$ is a càdlàg process, for almost every $\omega \in \Omega$, there exists a compact subset $\mathscr{K}_{\omega,T} \subset S$ such that $\eta(\omega, t) \in \mathscr{K}$ $\forall t \in [0,T]$ (see~\cite[Remark~3.6.4]{ethier2009markov}). Hence, by the continuity properties of $g$, of $\frac{\partial}{\partial t}g$, and of~$\mathscr{L}g$ (see conditions~(ii) and~(iii) of this lemma), and since $(\eta(t))_{t \geq 0}$ is càdlàg, by dominated convergence the following limits hold almost surely:
    \begin{equation} \label{Paper01_first_step_limits_integration_by_parts}
    \begin{aligned}
        \lim_{n \rightarrow \infty} \; \sum_{i = 0}^{n-1} \int_{h^{(n)}_{i} \wedge T}^{h^{(n)}_{i+1} \wedge T} \left(\frac{\partial}{\partial t'} g\Big(\cdot, \eta(h^{(n)}_{i} \wedge T)\Big)\right)(t') \; dt' & = \int_{t_{1} \wedge T}^{t_{2} \wedge T} \left(\frac{\partial}{\partial t'} g(\cdot, \eta(t'-))\right)(t') \; dt', \\ \lim_{n \rightarrow \infty} \; \sum_{i = 0}^{n-1} \int_{h^{(n)}_{i} \wedge T}^{h^{(n)}_{i + 1} \wedge T} \left(\mathscr{L}g\Big(h^{(n)}_{i+1} \wedge T, \cdot\Big)\right)(\eta(t'-)) \; dt' & = \int_{t_{1} \wedge T}^{t_{2} \wedge T} \left(\mathscr{L}g(t', \cdot)\right)(\eta(t'-)) \; dt'.
    \end{aligned}
    \end{equation}
    Since uniform integrability and almost sure convergence imply convergence in~$L_{1}$, it follows from~\eqref{Paper01_first_step_limits_integration_by_parts} and~\eqref{Paper01_intermediate_step_integration_parts_formula} that to establish~\eqref{Paper01_intergration_parts_martingale_property}, it will suffice to verify that the sequences of random variables
    \begin{align}
        & \Bigg(\sum_{i = 0}^{n-1} \int_{h^{(n)}_{i} \wedge T}^{h^{(n)}_{i+1} \wedge T} \left(\frac{\partial}{\partial t'} g\Big(\cdot, \eta(h^{(n)}_{i} \wedge T)\Big)\right)(t') \; dt'\Bigg)_{n \in \mathbb{N}} \label{Paper01_sequence_derivative_time}\\  \textrm{and} \quad & \Bigg(\sum_{i = 0}^{n-1} \int_{h^{(n)}_{i} \wedge T}^{h^{(n)}_{i + 1} \wedge T} \left(\mathscr{L}g\Big(h^{(n)}_{i+1} \wedge T, \cdot\Big)\right)(\eta(t'-)) \; dt'\Bigg)_{n \in \mathbb{N}}, \label{Paper01_sequence_chain_rule_generator}
    \end{align}
    are uniformly integrable. Starting with the sequence~\eqref{Paper01_sequence_derivative_time}, observe that, for every~$n \in \mathbb{N}$,
    \begin{equation*}
    \begin{aligned}
        & \mathbb{E}_{\nu}\Bigg[\Bigg\vert \sum_{i = 0}^{n-1} \int_{h^{(n)}_{i} \wedge T}^{h^{(n)}_{i+1} \wedge T} \left(\frac{\partial}{\partial t'} g\Big(\cdot, \eta(h^{(n)}_{i} \wedge T)\Big)\right)(t') \; dt'\Bigg\vert\Bigg] \\
        & \quad \leq \sum_{i = 0}^{n-1} \int_{h^{(n)}_{i} \wedge T}^{h^{(n)}_{i+1} \wedge T} \mathbb{E}_{\nu}\Bigg[ \Bigg\vert\left(\frac{\partial}{\partial t'} g\Big(\cdot, \eta(h^{(n)}_{i} \wedge T)\Big)\right)(t') \Bigg\vert\Bigg] \; dt' \\ & \quad \leq \sum_{i = 0}^{n-1} \Bigg(\int_{h^{(n)}_{i} \wedge T}^{h^{(n)}_{i+1} \wedge T}  \; dt' \Bigg)\sup_{\tau_{1},\tau_{2} \in [0,T]} \; \mathbb{E}_{\nu}\Bigg[ \left(\frac{\partial}{\partial t} g\Big(\cdot, \eta(\tau_{2})\Big)\right)^{2}(\tau_{1}) \Bigg]^{1/2} \\ 
        & \quad \leq (t_{2} - t_{1}) \sup_{\tau_{1},\tau_{2} \in [0,T]} \; \mathbb{E}_{\nu}\Bigg[\left(\frac{\partial}{\partial t} g\Big(\cdot, \eta(\tau_{2})\Big)\right)^{2}(\tau_{1}) \Bigg]^{1/2},
    \end{aligned}
    \end{equation*}
    where for the second estimate we used Jensen's inequality. Therefore, by condition~(iv) of this lemma, we have
    \begin{equation} \label{Paper01_trivial_L1_bound_derivatibe_time}
        \sup_{n \in \mathbb{N}} \; \mathbb{E}_{\nu}\Bigg[\Bigg\vert \sum_{i = 0}^{n-1} \int_{h^{(n)}_{i} \wedge T}^{h^{(n)}_{i+1} \wedge T} \left(\frac{\partial}{\partial t'} g\Big(\cdot, \eta(h^{(n)}_{i} \wedge T)\Big)\right)(t') \; dt'\Bigg\vert\Bigg] < \infty.
    \end{equation}
    Hence, to conclude that the sequence in~\eqref{Paper01_sequence_derivative_time} is uniformly integrable, it will be enough to verify that for any~$\varepsilon > 0$, there exists~$\delta > 0$ such that for any measurable event $\mathscr{A}$ with $\mathbb{P}_{\nu}(\mathscr{A}) < \delta$, 
    \begin{equation} \label{Paper01_uniform_integrability_criterion}
        \sup_{n \in \mathbb{N}} \; \mathbb{E}_{\nu}\Bigg[\Bigg\vert \sum_{i = 0}^{n-1} \int_{h^{(n)}_{i} \wedge T}^{h^{(n)}_{i+1} \wedge T} \left(\frac{\partial}{\partial t'} g\Big(\cdot, \eta(h^{(n)}_{i} \wedge T)\Big)\right)(t') \; dt'\Bigg\vert \cdot \mathds{1}_{\mathscr{A}}\Bigg] < \varepsilon.
    \end{equation}
    Observe that by applying Fubini's theorem and then the Cauchy-Schwarz inequality, for any measurable event~$\mathscr{A}$ and any~$n \in \mathbb{N}$,
    \begin{equation*}
    \begin{aligned}
        & \mathbb{E}_{\nu}\Bigg[\Bigg\vert \sum_{i = 0}^{n-1} \int_{h^{(n)}_{i} \wedge T}^{h^{(n)}_{i+1} \wedge T} \left(\frac{\partial}{\partial t'} g\Big(\cdot, \eta(h^{(n)}_{i} \wedge T)\Big)\right)(t') \; dt'\Bigg\vert \cdot \mathds{1}_{\mathscr{A}}\Bigg] \\
        & \quad \leq \sum_{i = 0}^{n-1} \int_{h^{(n)}_{i} \wedge T}^{h^{(n)}_{i+1} \wedge T} \mathbb{E}_{\nu}\Bigg[\Bigg\vert\left(\frac{\partial}{\partial t'} g\Big(\cdot, \eta(h^{(n)}_{i} \wedge T)\Big)\right)(t')\Bigg\vert \cdot \mathds{1}_{\mathscr{A}}\Bigg] \; dt' \\
        & \quad \leq \sum_{i = 0}^{n-1} \Bigg(\int_{h^{(n)}_{i} \wedge T}^{h^{(n)}_{i+1} \wedge T}  \; dt' \Bigg)
        \sup_{\tau_{1},\tau_{2} \in [0,T]} \; \mathbb{E}_{\nu}\Bigg[ \left(\frac{\partial}{\partial t} g\Big(\cdot, \eta(\tau_{2})\Big)\right)^{2}(\tau_{1}) \Bigg]^{1/2} \mathbb{P}_{\nu}(\mathscr{A})^{1/2} \\
        & \quad \leq (t_{2} - t_{1}) \mathbb{P}_{\nu}(\mathscr{A})^{1/2} \sup_{\tau_{1},\tau_{2} \in [0,T]} \; \mathbb{E}_{\nu}\Bigg[\left(\frac{\partial}{\partial t} g\Big(\cdot, \eta(\tau_{2})\Big)\right)^{2}(\tau_{1}) \Bigg]^{1/2} ,
    \end{aligned}
    \end{equation*}
    and so by condition~(iv) of this lemma, by taking~$\delta>0$ sufficiently small,~\eqref{Paper01_uniform_integrability_criterion} holds. Therefore, by~\eqref{Paper01_trivial_L1_bound_derivatibe_time} and~\eqref{Paper01_uniform_integrability_criterion}, the sequence in~\eqref{Paper01_sequence_derivative_time} is uniformly integrable. A similar argument implies the sequence in~\eqref{Paper01_sequence_chain_rule_generator} is also uniformly integrable. Since for a sequence of uniformly integrable random variables, almost sure convergence implies convergence in~$L_{1}$, we can apply~\eqref{Paper01_first_step_limits_integration_by_parts} to take the limit as~$n \rightarrow \infty$ on the right-hand side of~\eqref{Paper01_intermediate_step_integration_parts_formula}, obtaining
    \begin{equation} \label{Paper01_almost_final_identity_integration_parts}
    \begin{aligned}
        & \mathbb{E}_{\nu}\Big[g\Big(t_{2} \wedge T, \eta(t_{2}\wedge T)\Big) - g\Big(t_{1} \wedge T, \eta(t_{1}\wedge T)\Big) \Big\vert \mathcal{F}^{\eta}_{t_1+} \Big] \\ 
        & \quad = \mathbb{E}_{\nu}\Bigg[    \int_{t_{1} \wedge T}^{t_{2} \wedge T} \left(\left(\frac{\partial}{\partial t'} g(\cdot, \eta(t'-))\right)(t') + \left(\mathscr{L}g(t', \cdot)\right)(\eta(t'-))\right) \; dt'\Bigg\vert \mathcal{F}^{\eta}_{t_{1}+}\Bigg].
    \end{aligned}
    \end{equation}
    By rearranging terms, and by the definition of the process~$M^{g}$ in~\eqref{eq:daggerMgdefn}, identity~\eqref{Paper01_almost_final_identity_integration_parts} yields~\eqref{Paper01_intergration_parts_martingale_property}, which completes the proof.
\end{proof}

\subsection{Standard large deviation results} \label{Paper01_subsection_large_deviation}

In this subsection, for ease of reference, we state and prove standard Chernoff bounds for Poisson processes and continuous-time random walks.

\begin{lemma} \label{Paper01_lemma_upper_tail_poisson}
    For $\alpha > 0$, let $Z_{\alpha}$ denote a Poisson random variable with mean $\alpha$. Then for $r \geq \alpha$, 
    \begin{equation} \label{Paper01_eq:chernoffPoisson}
         \mathbb P\left( Z_{\alpha}  \ge r\right) \le e^{r(1-\log(r/\alpha))}.
    \end{equation}
\end{lemma}

\begin{proof}
   By a Chernoff bound, we have that for $\alpha>0$ and $r\ge \alpha$,~
\begin{equation}
    \mathbb P\left( Z_\alpha \ge r\right) 
=\mathbb P\left( e^{\log(r/\alpha) Z_\alpha} \ge e^{r\log(r/\alpha)}\right)
\le e^{-r\log(r/\alpha)}e^{\alpha (\frac r \alpha -1)}
\le e^{r(1-\log(r/\alpha))},
\end{equation}
which completes the proof.
\end{proof}

We will now use Lemma~\ref{Paper01_lemma_upper_tail_poisson} to prove a large deviation estimate for
a continuous-time random walk. For $n \in \mathbb{N}$, recall the definition of $\Lambda_n$ from after~\eqref{Paper01_scaling_parameters_for_discrete_approximation}, and let $(X^n(t))_{t \geq 0}$ be a simple symmetric random walk on $\Lambda_n$ with reflecting boundaries, with total jump rate $m$.

\begin{corollary} \label{Paper01_corollary_lemma_tail_continuous_time_rw}
    Let $x \in \Lambda_n$. Then for $T > 0$ and $y \in L^{-1}\mathbb{Z}$ such that $\vert y - x \vert > e^2 mT L^{-1}$, 
    \begin{equation*}
        \sup_{n \in \mathbb{N}} \; \sup_{t \in [0,T]} \; \mathbb{P}\left(X^n(t) = y \vert \,  X^n(0) = x \right) \leq \exp\left(-L\vert y - x \vert\right).
    \end{equation*}
\end{corollary}

\begin{proof}
   Let $Z_{mT}$ be a Poisson distributed random variable with mean $mT$. By the construction of continuous-time random walks in terms of Poisson random variables, and then by Lemma~\ref{Paper01_lemma_upper_tail_poisson}, for $n\in \mathbb N$ and $t\in [0,T]$,
   \begin{equation*}
    \begin{aligned}
         \mathbb{P}\left(X^n(t) = y \vert \,  X^n(0) = x \right)  \leq \mathbb{P}(Z_{mT} \geq L\vert y - x \vert) 
         & \leq e^{L\vert y - x \vert  \left(1 - \log\left(\frac{e^2mT}{mT}\right)\right)} = e^{-L\vert y - x \vert},
    \end{aligned}
   \end{equation*}
   which completes the proof.
\end{proof}

\subsubsection*{Acknowledgements}

The authors are grateful to Matthias Birkner, Matt Roberts, and Alexandre Stauffer for helpful comments and suggestions. 
While this work was being carried out, JLdOM was supported by a scholarship from the EPSRC Centre for Doctoral Training in Statistical Applied Mathematics at Bath (SAMBa), under the project EP/S022945/1.
MO is partially supported by EPSRC research grant EP/X040089/1. 
SP is supported by a Royal Society University Research Fellowship. While part of this work was being carried out, SP was visiting SLMath as a Research Member of the Probability and Statistics of Discrete Structures program. 


\bibliographystyle{abbrv}
\bibliography{existence_muller}

\end{document}